\documentclass{amsart}
\usepackage[utf8]{inputenc}
\usepackage{amsfonts, amsmath, amsthm, amssymb}
\usepackage{graphicx}

\usepackage{a4wide}
\usepackage{rotating}
\usepackage{bbm}
\usepackage[all,cmtip]{xy}

\usepackage{dsfont}
\usepackage{color}
\usepackage{hyperref}
\usepackage{verbatim}
\usepackage{mathtools}
\usepackage{enumitem}
\usepackage{epstopdf}
\usepackage{booktabs}

\usepackage{tikz}
\usepackage[titletoc]{appendix}
\usetikzlibrary{patterns}
\usetikzlibrary{decorations.pathreplacing}
\usetikzlibrary{arrows.meta}
\usetikzlibrary{cd}
\usepackage[margin=1.0in]{geometry}
\usepackage{caption}
\usepackage{float}
\usepackage{mathtools}
\usepackage{pgfplots}
\usepackage{tkz-euclide}

\newcommand{\Z}{\mathbb{Z}}

\newcommand{\C}{\mathbb{C}}

\newcommand{\norm}[1]{\left\lVert#1\right\rVert}

\DeclareMathOperator{\id}{id}
\newcommand{\Tw}{\operatorname{Tw}}
\newcommand{\Ch}{\operatorname{Ch}}
\newcommand{\Aut}{\operatorname{Aut}}

\numberwithin{equation}{section}

\newtheorem{thm}{Theorem}[section]

\newtheorem{lemma}[thm]{Lemma}
\newtheorem{proposition}[thm]{Proposition}

\newtheorem{defn}[thm]{Definition}
\newtheorem{remark}[thm]{Remark}
\newtheorem{corollary}[thm]{Corollary}
\newtheorem{example}[thm]{Example}

\begin{document}

\title[Calabi--Yau dimensions of monotone Fukaya summands from infinitesimal gradings]{Calabi--Yau dimensions of monotone Fukaya summands from infinitesimal gradings}
\author{Yonghwan Kim}
\email{yonghkim@mit.edu}

\maketitle

\vspace{-2.1em}
\begin{abstract}
    We develop a theory of infinitesimal gradings for $\mathbb{Z}/2$-graded $A_\infty$-categories and use it to detect Calabi--Yau dimensions of monotone Fukaya summands from quantum cohomology. The dimensions we find agree with predictions from Serre functors of Kuznetsov components, such as $2$ for Gushel--Mukai and Verra fourfolds and the fractional value $\frac{n-2}{n}$ for the symplectic Grassmannian $\operatorname{IG}(2,2n)$. An infinitesimal grading induces endomorphisms on cohomological morphism spaces, and we prove a straightening theorem that realizes its weights as an actual grading compatible with the $A_\infty$-structure. We use this grading to show that Dehn twists give faithful braid group representations in the symplectic mapping class groups of Gushel--Mukai and Verra fourfolds.
\end{abstract}

\section{Introduction\label{sec:intro}}

Derived categories of Fano varieties can contain components that exhibit the geometry of a lower-dimensional Calabi--Yau variety. For example, Kuznetsov--Perry \cite{KuznetsovPerry} show that the Kuznetsov component of a Gushel--Mukai fourfold has Serre functor $[2]$, while Orlov's blowup formula \cite{Orlov} exhibits $D^b(E)$ as a semiorthogonal component of the derived category of the blowup of $\mathbb{P}^3$ along a smooth elliptic curve $E$ of degree $4$.

On the symplectic side, however, the monotone Fukaya category provides only a periodic grading \cite{Seidel00, Sheridan}, which does not determine these Calabi--Yau dimensions.

Assume that the mirror Landau--Ginzburg model $(X^\vee,W)$ is locally Morse--Bott along the corresponding critical component, which is smooth Calabi--Yau. Then this dimension is expected to appear in the Fukaya category of the corresponding Landau--Ginzburg model $(X^\vee,W)$ as the dimension of the critical component. In this case, Orlov's Knörrer periodicity shows that the same dimension is shared by the corresponding matrix factorization category $\operatorname{MF}(W-\lambda)$. Thus the other direction of mirror symmetry suggests that this dimension should be visible in the corresponding monotone Fukaya summand $\mathcal{F}(X)_\lambda$.

For example, when $X$ is the Fano threefold given by the blowup of $\mathbb{P}^3$ along a smooth elliptic quartic, Cheltsov--Przyjalkowski \cite[Family 2.25]{CP} compute the critical locus of its corresponding compactified Landau--Ginzburg model and show that it is a smooth elliptic curve with trivial normal bundle, and the superpotential is locally Morse--Bott along the normal directions, after resolving singularities in the fiber. Therefore the corresponding monotone Fukaya summand $\mathcal{F}(X)_{-1}$ is expected to carry dimension $1$ data.

We develop a theory of \emph{infinitesimal gradings} on $\mathbb{Z}/2$-graded $A_\infty$-categories and construct them on monotone Fukaya summands using the closed--open map. In the examples below, the full Fukaya summands admit homogeneous weak proper Calabi--Yau structures of degree $-2$ and $-1$, recovering the dimensions of the algebraic components above. We call the negative of this degree the \emph{infinitesimal Calabi--Yau dimension}.

We begin with the Fano threefold previously discussed. Throughout the introduction, we work over $\mathbb{C}$ and set $q=1$. All Calabi--Yau structures will be weak and proper, in the sense of \cite[Definition A.2]{Sheridan}.

\begin{thm}[Theorem \ref{thm:MM-2-25-reprise}] \label{thm:MM-2-25}
    Let $E\subset\mathbb P^3$ be a smooth elliptic curve of degree $4$, and equip $X=\mathrm{Bl}_E\mathbb P^3$ with a monotone K\"ahler form. Let $QH^*(X)_{-1}$ denote the $(-1)$-generalized eigensummand of small quantum cohomology for quantum multiplication by $c_1(X)$. With its closed-string infinitesimal grading, this summand admits an isomorphism of graded $\mathbb{C}$-algebras
    \begin{equation}
        QH^*(X)_{-1}\cong H^*(E;\mathbb{C}),
    \end{equation}
    where the right-hand side carries its ordinary cohomological grading, with weights $0,1,1,2$.
    
    The corresponding monotone Fukaya summand $\mathcal{F}(X)_{-1}$ admits an infinitesimal grading for which homogeneous Calabi--Yau structures exist precisely in infinitesimal dimension $1$.
\end{thm}

For a Gushel--Mukai fourfold, the Serre shift $[2]$ of its Kuznetsov component suggests looking for dimension $2$ on the zero Fukaya summand. Theorem \ref{thm:Gushel-Mukai} shows that infinitesimal gradings indeed detect this dimension.

This additional information also underlies our main application to symplectic topology: faithfulness of Artin group actions generated by Dehn twists. We extend arguments of Seidel--Thomas \cite{SeidelThomas} to obtain faithful Artin group actions in the monotone setting, even after the infinitesimal grading is forgotten.

\begin{thm}[Theorem \ref{thm:Gushel-Mukai-reprise}] \label{thm:Gushel-Mukai}
    Let $(X,\omega)$ be a smooth Gushel--Mukai fourfold, equipped with a monotone K\"ahler form. Then the monotone Fukaya category summand $\mathcal{F}(X)_0$ admits an infinitesimal grading for which homogeneous Calabi--Yau structures exist precisely in infinitesimal dimension $2$.

    Moreover, $X$ contains an $A_4$-configuration $L_1,\cdots,L_4$ of Lagrangian four-spheres, each of infinitesimal dimension $2$. The first three Dehn twists define an injective homomorphism
    \begin{equation}
        B_4\to\pi_0\mathrm{Symp}(X,\omega)
    \end{equation}
     to the symplectic mapping class group, sending the standard generators $\sigma_i\in B_4$ to the Dehn twists $[\tau_{L_i}]$.
\end{thm}

The same conclusions hold for smooth Verra fourfolds (Theorem \ref{thm:Verra-reprise}), which also admit Kuznetsov components of $K3$ type \cite{PerryPertusiZhao} and thus carry Serre functors $[2]$. 

Our final examples are symplectic Grassmannians $\operatorname{IG}(2,2n)$. Kuznetsov (\cite[Appendix A]{MMPS}) identifies the Kuznetsov component of $D^b(\operatorname{IG}(2,2n))$ with the derived category of representations of a type $A_{n-1}$ quiver, and thus its Serre functor satisfies $S^n\simeq [n-2]$. The following theorem recovers the resulting fractional Calabi--Yau dimension $\frac{n-2}{n}$ on the corresponding monotone Fukaya summand $\mathcal{F}(X)_0$.

\begin{thm}[Theorem \ref{thm:IG-reprise}] \label{thm:IG}
    For $n\geq 3$, let $X=\operatorname{IG}(2,2n)$ be the symplectic Grassmannian together with its monotone K\"ahler form. Then one can choose an infinitesimal grading on the monotone Fukaya summand $\mathcal{F}(X)_0$ for which the standard weak proper Calabi--Yau structure is homogeneous of dimension $\frac{n-2}{n}$. For this grading, every homogeneous weak proper Calabi--Yau structure has that dimension.
\end{thm}

These applications arise from a general construction. Suppose that quantum multiplication by $c_1(X)$ acts as scalar multiplication on the full generalized eigensummand $QH^*(X)_\lambda$:
\begin{equation}\label{eqn:intro-c1*scalar}
    c_1(X)*\vert_{QH^*(X)_\lambda}=\lambda\cdot\id.
\end{equation}
Then, we construct a geometric \emph{infinitesimal grading} on the corresponding Fukaya summand $\mathcal{F}(X)_\lambda$ using the closed--open map. The corresponding closed-string grading is computable from small quantum products and the ordinary cohomological grading; closed--open compatibility then determines the infinitesimal degrees used in our examples. After passing to graded lifts of objects and a minimal model, a straightening theorem realizes these gradings as degrees in an $A_\infty$-isomorphic model, in which each $\mu^d$ has degree $2-d$.

The remainder of this introduction presents the theory in stages building towards the proof of Theorem \ref{thm:Gushel-Mukai}.

\subsection{Infinitesimal gradings}

Our starting point is the Hochschild cocycle
\begin{equation}
    \kappa=\{(2-d)\mu^d\}
\end{equation}
associated to an uncurved $\mathbb{Z}/2$-graded $A_\infty$-category $(\mathcal{C},\mu)$. An ordinary $\mathbb{Z}$-grading determines a Hochschild cochain whose differential is $\kappa$. We define an \emph{infinitesimal grading} to be any choice of such a primitive
\begin{equation}
    \delta\eta=\kappa.
\end{equation}
The obstruction to existence is therefore the cohomology class $[\kappa]$. The cocycle $\kappa$ represents the Kodaira--Spencer class at $t=1$ of the one-parameter rescaling deformation $\mu_t^d=t^{2-d}\mu^d$ \cite[Section~3.1]{CaldararuLiTu}. This deformation enters noncommutative Hodge theory through the Getzler--Gauss--Manin connection \cite{Getzler,CaldararuLiTu}; related Euler grading operators
appear in \cite{SheridanFormulae,Hugtenburg}.

There is a second, objectwise obstruction. Given $\eta$, the length-zero component $\eta^0_L$ defines an obstruction class $[\eta^0_L]\in H^\mathrm{odd}\hom_\mathcal{C}(L,L)$ for each $L\in\mathcal{C}$. A \emph{graded lift} of $L$ is a choice of $\xi_L\in\hom^\mathrm{even}_\mathcal{C}(L,L)$ such that
\begin{equation}
    \mu^1(\xi_L)=\eta^0_L.
\end{equation}
For two graded objects, the length-one part of the infinitesimal grading, corrected by the chosen cochains $\xi_{L_0}$ and $\xi_{L_1}$, induces a \emph{grading operator}
\begin{equation}
    \Gamma:H^*\hom_\mathcal{C}(L_0,L_1)\to H^*\hom_\mathcal{C}(L_0,L_1).
\end{equation}
Its eigenvalues are the \emph{infinitesimal degrees} associated to morphism spaces. We define the \emph{graded lift} $\Tilde{\mathcal{C}}_\eta$ to be the category whose objects are pairs $(L,\xi_L)$ as above and whose morphism spaces and $A_\infty$-operations are inherited from $\mathcal{C}$. This inherits an infinitesimal grading $\Tilde{\eta}$, with trivial length-zero component. 

We show that these grading operators govern all higher products, up to $A_\infty$-isomorphism.

\begin{thm} [Theorem \ref{thm:straighten-infinitesimal-grading}] \label{thm:straighten-infinitesimal-grading-intro}
    Let $(\mathcal{C},\mu)$ be a minimal, cohomologically finite $\mathbb{Z}/2$-graded $A_\infty$-category over $\mathbb{C}$, with an infinitesimal grading $\eta$ satisfying $\eta^0\equiv0$. Then there exists an $A_\infty$-isomorphism
    \begin{equation}
        \Phi:(\mathcal{C},\{\mu\})\to(\mathcal{C},\{\mu'\})
    \end{equation}
    with $\Phi^1=\id$, such that
    \begin{equation}
        \mu'^d(\hom_{\lambda_d}\otimes\cdots\otimes\hom_{\lambda_1})\subset\hom_{\lambda_1+\cdots+\lambda_d+2-d},
    \end{equation}
    where $\hom_\lambda$ is the generalized eigenspace of $\eta^1$ with eigenvalue $\lambda$.
\end{thm}

Thus the semisimple part of $\eta^1$ defines an additional $\mathbb{C}$-valued weight grading of an $A_\infty$-isomorphic model, without any semisimplicity assumption on $\eta^1$.

Related constructions appear in Seidel--Solomon \cite{SeidelSolomon}, who use noncommutative vector fields to define equivariant objects and weight gradings in the theory of $q$-intersection numbers. Purity of these vector fields also underlies the formality criterion used to prove the formality of the arc algebra in \cite{AbouzaidSmith}.

The essential distinction is that noncommutative vector fields are Hochschild cocycles, whereas our infinitesimal gradings are generally not: they satisfy $\delta\eta=\kappa$, and are cocycles only when $\kappa$ vanishes at the chain level. The term $\kappa=\{(2-d)\mu^d\}$ encodes the degree identities for the $A_\infty$-operations, and our straightening theorem realizes these identities in an $A_\infty$-isomorphic model without imposing formality.

\subsection{Geometric infinitesimal gradings on monotone Fukaya categories}

Let $(X,\omega)$ be a closed monotone symplectic manifold of real dimension $2n$. Denote by $QH^*(X)_\lambda$ the $\lambda$-generalized eigensummand of small quantum cohomology for quantum multiplication by $c_1(X)$, and by $\mathcal{F}(X)_\lambda$ the corresponding monotone Fukaya category summand. 

The obstruction to an infinitesimal grading of $\mathcal{F}(X)_\lambda$ is determined by the closed--open image of the first Chern class.

\begin{thm}[Theorem \ref{thm:monotone-fukaya-gradeability-reprise}] \label{thm:monotone-fukaya-gradeability-intro}
    The monotone Fukaya category summand $\mathcal{F}(X)_\lambda$ admits an infinitesimal grading if and only if the following cohomology class vanishes in Hochschild cohomology $HH^*(\mathcal{F}(X)_\lambda,\mathcal{F}(X)_\lambda)$:
    \begin{equation}
        [CO_\lambda(c_1(X)-\lambda\cdot 1_X)]=0.
    \end{equation}
\end{thm}

Condition \eqref{eqn:intro-c1*scalar} means that the nilpotent part of $c_1(X)*-\lambda\id$ vanishes on the entire generalized eigensummand; it implies the cohomological vanishing above. However, existence alone does not provide the weights: different primitives of $\kappa$ can induce different grading operators. We therefore construct a particular geometric primitive $\eta_\mathrm{geom}$ and relate its induced grading operators to closed-string data.

The first comparison concerns the open--closed map. Let $\pi_\lambda$ be the idempotent projecting to $QH^*(X)_\lambda$, and let $gr_X$ denote the ordinary grading operator on $H^*(X)$: this multiplies a homogeneous element by its degree. Write $OC_\lambda:HH_*(\mathcal{F}(X)_\lambda,\mathcal{F}(X)_\lambda)\to QH^*(X)_\lambda$ for the open--closed map.

\begin{thm}[Theorem \ref{thm:OC-compatibility-reprise}] \label{thm:OC-compatibility-intro}
    Let $(X,\omega)$ be a monotone symplectic manifold that satisfies (\ref{eqn:intro-c1*scalar}). Let $\Gamma^{HH_*}_{\eta_\mathrm{geom}}$ be the induced grading operator on Hochschild homology $HH_*(\mathcal{F}(X)_\lambda,\mathcal{F}(X)_\lambda)$. Then the following holds:
    \begin{equation}
        OC_\lambda\circ\Gamma^{HH_*}_{\eta_\mathrm{geom}}=\pi_\lambda*(gr_X-n)\circ OC_\lambda.
    \end{equation}
\end{thm}

This identity makes infinitesimal weights accessible through closed-string computations. As in the usual graded setting, however, the open--closed comparison incorporates a degree shift, reflected in the term $gr_X-n$. To identify the internal grading on the closed-string summand, we turn to the closed--open map and construct a grading operator for which this map preserves weights.

\begin{thm}[Closed-string grading] \label{thm:intro-closed-string-grading}
    Assume \eqref{eqn:intro-c1*scalar}. The operator
    \begin{equation}
        \Gamma_{\mathrm{closed}}(x)
        =\pi_\lambda*gr_X(x)-gr_X(\pi_\lambda)*x,
        \quad x\in QH^*(X)_\lambda,
    \end{equation}
    is a derivation of the quantum product:
    \begin{equation}
        \Gamma_{\mathrm{closed}}(x*y)=\Gamma_{\mathrm{closed}}(x)*y+x*\Gamma_{\mathrm{closed}}(y).
    \end{equation}
    Its generalized eigenspaces define a $\mathbb{C}$-graded algebra structure on $QH^*(X)_\lambda$, with unit $\pi_\lambda$ of degree zero.
\end{thm}

This operator also admits a geometric expression through four-point Gromov--Witten invariants:
\begin{equation}
    \Gamma_\mathrm{closed}(x)=\pi_\lambda*(gr_X(x)-2GW(c_1(X),\pi_\lambda,x)^\sharp),\quad x\in QH^*(X)_\lambda,
\end{equation}
where the four-point invariant is written as a cohomology-valued operation via Poincaré duality:
\begin{equation}
    \langle GW(a,b,c)^\sharp,d\rangle=\sum_A GW_{0,4,A}(a,b,c,d).
\end{equation}

The closed--open map intertwines the closed-string grading operator with the induced operator on Hochschild cohomology.

\begin{thm}[Theorem \ref{thm:CO-compatibility-reprise}] \label{thm:CO-compatibility-intro}
    Assume \eqref{eqn:intro-c1*scalar}. Let $\Gamma^{HH^*}_{\eta_{\mathrm{geom}}}$ denote the induced grading operator on Hochschild cohomology. Then
    \begin{equation}
        CO_\lambda\circ\Gamma_{\mathrm{closed}}=\Gamma^{HH^*}_{\eta_{\mathrm{geom}}}\circ CO_\lambda.
    \end{equation}
\end{thm}

Open--closed compatibility gives a criterion for homogeneity. We use the geometric weak proper Calabi--Yau trace $\alpha\mapsto\int_X OC_\lambda(\alpha)$, together with unitality of $CO_\lambda$ and its module compatibility with $OC_\lambda$.

\begin{thm}[Proposition \ref{thm:CY-criterion-reprise}] \label{thm:intro-cy-criterion}
Assume \eqref{eqn:intro-c1*scalar}, and let $\theta\in QH^*(X)_\lambda$ be an even class invertible in $QH^*(X)_\lambda$, satisfying
\begin{equation}
    \pi_\lambda*gr_X(\theta)=(n-d)\theta.
\end{equation}
Then the functional $\Phi_\theta(\alpha)\coloneqq \int_X \theta*OC_\lambda(\alpha)$ on $HH_*(\mathcal{F}(X)_\lambda,\mathcal{F}(X)_\lambda)$ defines a homogeneous weak proper Calabi--Yau structure of infinitesimal dimension $d$ with respect to the geometric infinitesimal grading.
\end{thm}

The open--closed comparison and Poincar\'e duality give $\Phi_\theta\circ\Gamma^{HH_*}_{\eta_{\mathrm{geom}}}=d\,\Phi_\theta$.

A graded spherical object of infinitesimal dimension $d_V$ determines the dimension of any homogeneous weak proper Calabi--Yau structure on the category: nondegeneracy forces the unit to pair with the top degree class, whose degree is $d_V$. We apply this to vanishing cycles that define spherical objects in the monotone Fukaya category; their infinitesimal dimensions can be determined from the action of the closed-string grading operator on primitive cohomology.

\subsection{Spherical twists}

The existence of integer gradings in Fukaya categories plays a central role in detecting nontrivial symplectic mapping classes and establishing faithfulness of braid group actions from spherical configurations \cite{SeidelThomas}. To extend these arguments to the $\mathbb{Z}/2$-graded setting, we extend infinitesimal gradings to twisted complexes, allowing their coefficients to be decorated with weights. Twist functors along infinitesimally graded objects then admit graded lifts, with a canonical choice of objectwise grading data.

For an infinitesimally graded spherical object $V$ in the sense of Definition \ref{defn:spherical-object}, the grading operator sends the unit to zero, and therefore acts on the one-dimensional quotient
\begin{equation}
    H^*\hom_\mathcal{C}(V,V)/\langle e_V\rangle
\end{equation}
by a scalar $d_V$, which we call the \emph{infinitesimal dimension} of $V$. The canonically graded twist $T_V$ has infinitesimal weight shift $1-d_V$ on $V$, analogous to the grading shift of a Dehn twist. When $d_V\neq 1$, repeated twists produce weight shifts that the $\mathbb{Z}/2$-grading does not record.

This leads to an analogue of the classical faithfulness results for graded spherical twists \cite{SeidelThomas, BravThomas, NordskovaVolkov}. The additional issue here is to control the action after the infinitesimal grading has been forgotten. Write $D^\pi\mathcal{C}$ for the split-closed derived category, and $\mathrm{Aut}(D^\pi\mathcal{C})$ for the group of isomorphism classes of exact autoequivalences, with the infinitesimal grading forgotten.

\begin{thm}[Theorems \ref{thm:central-kernel}, \ref{thm:ADE-faithfulness-reprise}] \label{thm:ADE-faithfulness}
    Let $(\mathcal{C},\eta)$ be a cohomologically finite, infinitesimally graded $\mathbb{Z}/2$-graded $A_\infty$-category over $\mathbb{C}$. Let $\Gamma$ be a connected Dynkin diagram of type $A$, $D$, or $E$, and let $\{V_i\}$ be a $\Gamma$-configuration of infinitesimally graded spherical objects with common integer infinitesimal dimension $\neq 1$. 
    
    Then the associated Artin group action
    \begin{equation}
        \rho:B_\Gamma\to\mathrm{Aut}(D^\pi\mathcal{C}),\quad \sigma_i\mapsto T_{V_i}
    \end{equation}
    has kernel contained in the center:
    \begin{equation}
        \ker\rho\subset Z(B_\Gamma).
    \end{equation}
    Consequently, the restriction of $\rho$ to any proper parabolic subgroup $P\subset B_\Gamma$ is faithful.
\end{thm}

Thus only central elements can act trivially after the infinitesimal grading is forgotten. For a proper parabolic subgroup $P\subset B_\Gamma$ in the sense of \cite{Paris} (recalled in Definition \ref{def:braid-relations}), faithfulness follows from $P\cap Z(B_\Gamma)=\{1\}$. For a standard parabolic subgroup, we retain the action on the full ambient category; the remaining objects of the configuration serve as comparison objects for grading shifts. In particular, an $A_4$-configuration gives a faithful $B_4$-action generated by the first three twists, as used in Theorem \ref{thm:Gushel-Mukai}.

\subsection{Geometric applications}

For a Gushel--Mukai fourfold $X$, quantum multiplication by $c_1(X)$ vanishes on the full generalized zero-eigensummand $QH^*(X)_0$. We compute its closed-string grading operator and apply closed--open compatibility to determine the infinitesimal dimensions of Lagrangian vanishing spheres.

\begin{thm}
    Let $(X,\omega)$ be a smooth Gushel--Mukai fourfold equipped with a monotone K\"ahler form. Then every Lagrangian sphere has infinitesimal dimension $2$ with respect to the geometric infinitesimal grading on $\mathcal{F}(X)_0$.
\end{thm}

We obtain the $A_4$-configuration in Theorem \ref{thm:Gushel-Mukai} from an explicit degeneration with an $A_4$-singularity. Theorem \ref{thm:ADE-faithfulness} then gives a faithful $B_4$-action generated by the first three twists on the derived Fukaya category, and hence an embedding into the symplectic mapping class group. Dehn twists along four-spheres have finite order in the smooth mapping class group \cite{KeatingRandalWilliams}, so the infinitesimal degrees detect a distinction between smooth and symplectic topology.

For the elliptic example, let $X=\mathrm{Bl}_E\mathbb{P}^3$, where $E$ is a smooth elliptic curve of degree $4$. We compute that quantum multiplication by $c_1(X)$ acts as $-\id$ on the full generalized $(-1)$-eigensummand. The induced closed-string grading identifies this summand with the cohomology of the blowup center as a graded algebra:
\begin{equation}
    QH^*(X)_{-1}\cong H^*(E;\mathbb{C}).
\end{equation}
Closed--open compatibility also implies that every Lagrangian sphere $L\subset X$ with $[L]\neq 0$ in $H_3(X;\mathbb{C})$ that admits a graded lift has infinitesimal dimension $1$.

The appearance of the blowup center is also reflected in Iritani's blowup theorem \cite[Theorems 1.1, 1.2]{Iritani}. After a formal change of variables and extension of coefficients, the induced decomposition of quantum cohomology also preserves Euler vector fields. We identify the graded ring at the monotone specialization and relate its grading to the monotone Fukaya summand through the closed--open map.

The categorical Calabi--Yau statements require a further computation. In both examples, the geometric weak proper Calabi--Yau structure obtained by integrating the open--closed image is not homogeneous for the geometric infinitesimal grading, assuming surjectivity of open--closed maps onto the corresponding summands. We compute corrections that produce nondegenerate homogeneous structures on the full Fukaya summands, of degree $-2$ for the Gushel--Mukai fourfold and $-1$ for the elliptic blowup. The correcting classes are invertible elements satisfying the eigencondition of Theorem \ref{thm:intro-cy-criterion}. Thus both dimension statements arise from the same open--closed compatibility mechanism.

\subsection{Organization of the paper}

The paper is largely divided into two parts. Sections 2 and 3 develop the algebraic theory of infinitesimal gradings, while Sections 4 and 5 apply it to monotone Fukaya categories.

In Section 2, we introduce infinitesimal gradings on $A_\infty$-categories and functors and establish their basic properties. In Section 3, we extend the theory to twisted complexes and show that spherical twist functors along graded objects admit graded lifts. We then establish criteria for spherical twists to have infinite order and prove that the associated ADE Artin group actions have central kernel, yielding faithfulness on proper parabolic subgroups. In Section 4, we turn to monotone Fukaya categories and construct a geometrically defined infinitesimal grading compatible with the closed--open map. Finally, Section 5 presents concrete examples of monotone symplectic manifolds to which these results apply.

\subsection*{Acknowledgements}

First and foremost, I would like to thank my advisor Paul Seidel for his continued support and guidance throughout this project. I would also like to thank Denis Auroux, Zihong Chen, Sheel Ganatra, Sebastian Haney, Kai Hugtenburg, Samuel Johnston, Davesh Maulik, Victor Przyjalkowski, and Chris Woodward for helpful comments and discussions at different stages of this project. This work was partially supported by a Simons Dissertation Fellowship award.

\subsection*{AI disclosure}

The author fully developed the research and wrote a first draft of the manuscript. Large Language Models were then used to assist with grammatical proofreading, improving exposition and theorem statements, critically examining mathematical arguments, and checking the consistency of sign conventions before arXiv submission. They were also used to assist in combining two figures into one (Figure \ref{fig:Q_D}). The author takes full responsibility for the mathematical content and final presentation.

\section{Algebra of infinitesimal gradings} \label{sec:algebra}

In this section, we develop the theory of infinitesimal gradings for $A_\infty$-categories. We first define a Hochschild cohomology class $[\kappa]$ invariant under quasi-equivalences, which obstructs the existence of a grading. Motivated by this, we define any primitive $\eta$ of $\kappa$ to be an infinitesimal grading, and use this Hochschild cochain to define ``grading operators'' on morphism spaces and Hochschild (co)chains. These grading operators can be shown to be compatible with products and Calabi-Yau structures. After passing to the $\eta$-graded lift and a minimal model, Theorem \ref{thm:straighten-infinitesimal-grading} straightens the grading under its stated field and finiteness hypotheses. The generalized eigenspaces then define a $k$-grading compatible with the $A_\infty$-structure. If the length $1$ component is semisimple, it is the Euler operator of this grading. We conclude the section by defining corresponding notions of infinitesimal gradings for functors and natural transformations.

\subsection{Hochschild theory for $A_\infty$-algebras}\label{ssec:preliminaries}

We first review some basic notions of Hochschild theory in the setup of $A_\infty$-algebras, implicitly fixing notation and sign conventions. The same formulas carry over to $A_\infty$-categories. Our sign conventions follow \cite{Seidel08} and \cite{GanatraThesis}, up to ordering the inputs right to left.

In this section, $\mathcal{A}$ will be a $\Z/2$-graded $A_\infty$-algebra over a field $k$, unless otherwise mentioned. The notation $\lvert x\rvert$ refers to the $\mathbb{Z}/2$-degree, while the notation $\norm{x}$ refers to its reduced degree, i.e. $\norm{x}=\lvert x\rvert-1$.

\subsubsection{Hochschild cochains and the pre-Euler cocycle}

We first recall the definition of Hochschild cohomology.

\begin{defn}
    Given an $A_\infty$-algebra $\mathcal{A}$, we define its \emph{Hochschild cochain complex} to be
    \begin{equation}
        CC^*(\mathcal{A},\mathcal{A})=\prod_{k\geq0} CC^k(\mathcal{A},\mathcal{A}),\quad CC^k(\mathcal{A},\mathcal{A})=\hom(\mathcal{A}[1]^{\otimes k},\mathcal{A})
    \end{equation}
    with the differential
    \begin{equation}
    \begin{aligned}
        (\delta \alpha)^k(x_k,\cdots,x_1)&=-\sum (-1)^{\dagger_i\cdot\norm{\alpha}}\mu^{k-j+1}(x_k,\cdots,\alpha^j(x_{i+j},\cdots,x_{i+1}),x_i,\cdots,x_1)\\
        &+\sum (-1)^{\dagger_i+\norm{\alpha}}\alpha^{k-j+1}(x_k,\cdots,\mu^j(x_{i+j},\cdots,x_{i+1}),x_i,\cdots,x_1),
    \end{aligned}    
    \end{equation}
    where $\norm{\alpha}$ is the reduced degree of the Hochschild cochain $\alpha$, and $\dagger_i=\norm{x_1}+\cdots+\norm{x_i}$. 
\end{defn}

Note that the superscript $k$ stands for the length filtration: when we write either $CC^{\textrm{even}}(\mathcal{A},\mathcal{A})$ or $CC^{\textrm{odd}}(\mathcal{A},\mathcal{A})$, the superscript will stand for the $\Z/2$-grading on Hochschild cochains coming from the grading of $\mathcal{A}$.

Hochschild cohomology also admits the structure of a Gerstenhaber algebra: we spell out the bracket below.

\begin{defn}
    The \emph{Gerstenhaber bracket} on Hochschild cochains is defined by
    \begin{equation}
        [\alpha,\beta]=\alpha\circ\beta-(-1)^{(\lvert\alpha\rvert-1)(\lvert\beta\rvert-1)}\beta\circ\alpha
    \end{equation}
    where the circle product is defined by
    \begin{equation}
        (\alpha\circ\beta)^k(x_k,\cdots,x_1)=\sum_{i,j}(-1)^{\norm{\beta}\cdot\dagger_i}\alpha^{k-j+1}(x_k,\cdots,\beta^j(x_{i+j},\cdots,x_{i+1}),x_i,\cdots,x_1).
    \end{equation}
\end{defn}

The Gerstenhaber bracket satisfies the following two identities:
\begin{align}
    \delta[\alpha,\beta]&=[\delta\alpha,\beta]+(-1)^{\norm{\alpha}}[\alpha,\delta\beta], \label{eqn:Gersten-Leibniz} \\
    [\alpha,[\beta,\gamma]]&=[[\alpha,\beta],\gamma]+(-1)^{\norm{\alpha}\norm{\beta}}[\beta,[\alpha,\gamma]], \label{eqn:gersten-jacobi}
\end{align}
and therefore the shifted Hochschild cochain complex $CC^*(\mathcal{A},\mathcal{A})[1]$ is a dgLa, with graded Lie bracket induced from the Gerstenhaber bracket.

Given a $\mathbb{Z}/2$-graded $A_\infty$-algebra $\mathcal{A}$, there are two naturally defined Hochschild cocycles:

\begin{defn}
    Define the \emph{structure cocycle} to be $\mu=\{\mu^k\}\in CC^{even}(A,A)$, and the \emph{pre-Euler cocycle} to be $\kappa=\{(2-k)\mu^k\}\in CC^{even}(A,A)$. 
\end{defn}

One can indeed verify that $\mu$ and $\kappa$ are Hochschild cocycles: for example,
\begin{equation}
\begin{aligned}
    (\delta\kappa)^d&=-\sum(-1)^{\dagger_i}\mu^{d-j+1}(\cdots,\kappa^j(\cdots),\cdots)-\sum(-1)^{\dagger_i}\kappa^{d-j+1}(\cdots,\mu^j(\cdots),\cdots)\\
    &=-\sum(-1)^{\dagger_i}(2-j)\mu^{d-j+1}(\cdots,\mu^j(\cdots),\cdots)-\sum(-1)^{\dagger_i}(2-d+j-1)\mu^{d-j+1}(\cdots,\mu^j(\cdots),\cdots)\\
    &=(d-3)\sum(-1)^{\dagger_i}\mu^{d-j+1}(\cdots,\mu^j(\cdots),\cdots)=0
\end{aligned}    
\end{equation}
from the $A_\infty$-relations. One can also check that the Hochschild codifferential satisfies
\begin{equation}
    \delta=-[\mu,~].
\end{equation}

The pre-Euler cocycle $\kappa$ is related to the existence of gradings on $\mathcal{A}$ by the following proposition:

\begin{proposition} \label{prop:gradings-imply-kappa-zero}
    Suppose that the $\mathbb{Z}/2$-graded $A_\infty$-algebra $\mathcal{A}$ is defined from a $\mathbb{Z}$-graded $A_\infty$-algebra with degrees reduced modulo $2$. Then the cohomology class $[\kappa]\in HH^\mathrm{even}(\mathcal{A},\mathcal{A})$ is zero.
\end{proposition}
\begin{proof}
    Define the Euler cochain $Eu\in CC^1(\mathcal{A},\mathcal{A})$ by
    \begin{align}
        (Eu)^1(x)&=\deg(x)\cdot x,\quad x\in\mathcal{A},\\
        (Eu)^k&\equiv0,\qquad \quad\quad (k\neq1).
    \end{align}
    Then
    \begin{equation}
    \begin{aligned}
        (\delta Eu)^k(x_k,\cdots,x_1)&=-\sum_i\deg(x_i)\cdot\mu^k(x_k,\cdots,x_1)+\deg(\mu^k(x_k,\cdots,x_1))\cdot\mu^k(x_k,\cdots,x_1)\\
        &=(2-k)\cdot\mu^k(x_k,\cdots,x_1),
    \end{aligned}    
    \end{equation}
    so $\kappa=\delta (Eu)$ is a coboundary.
\end{proof}

We remark that the cocycle $\kappa$ is almost never zero at the chain level, even if the $A_\infty$-algebra $\mathcal{A}$ is $\mathbb{Z}$-graded: indeed, $\kappa$ is zero at the chain level if and only if $\mathcal{A}$ is a minimal $A_\infty$-algebra with all $\mu^k$ for $k\geq3$ trivial in characteristic 0. In characteristic $p$ the situation is more interesting, since $\mathcal{A}$ is allowed to have nontrivial operations $\mu^k$ when $k\equiv2~(\mathrm{mod}~p)$.

To provide some motivation for the theory to come, we look back at the proof of Proposition \ref{prop:gradings-imply-kappa-zero}. We notice that if there exists a direct sum decomposition
\begin{equation}
    \mathcal{A}=\bigoplus_{\lambda\in k} \mathcal{A}^\lambda
\end{equation}
such that the $A_\infty$-operations respect this decomposition, i.e.
\begin{equation}
    \mu^k(\mathcal{A}^{\lambda_k}\otimes\cdots\otimes\mathcal{A}^{\lambda_1})\subset\mathcal{A}^{\sum\lambda_i+2-k},
\end{equation}
then the cocycle $\kappa$ is a coboundary. Thus the nonvanishing of $[\kappa]$ not only obstructs $\mathbb{Z}$-gradings, but also the existence of gradings valued in the base field $k$. This leads us to the following definition.

\begin{defn} \label{defn:inf-grading-algebra}
    Let $\mathcal{A}$ be a $\mathbb{Z}/2$-graded $A_\infty$-algebra. If $[\kappa_\mathcal{A}]=0$, we define an \emph{infinitesimal grading} for $\mathcal{A}$ to be a choice of cochain $\eta\in CC^{\mathrm{odd}}(\mathcal{A},\mathcal{A})$ such that $\delta\eta=\kappa_\mathcal{A}$.
\end{defn}

For genuinely $k$-graded $A_\infty$-algebras, there is a canonical choice of an infinitesimal grading given by the Euler cochain defined above.

The remainder of this subsection introduces Hochschild homology and recalls the theory of noncommutative Cartan calculus for Hochschild invariants. Readers interested in the $A_\infty$-categorical theory of infinitesimal gradings are invited to skip ahead to subsections \ref{ssec:kappa-invariance} and \ref{ssec:infinitesimal-grading-categories}.

\subsubsection{Hochschild chains and Cartan calculus} \label{sssec:Cartan-calculus}

We will use the Lie derivative to define infinitesimal gradings for Hochschild chains, and the cap product to compare the induced gradings on Hochschild homology and cohomology. First, we recall the definition of the cup product to place these operations within the standard noncommutative Cartan calculus picture.

\begin{defn}
    For $\mathcal{A}$ as before, define the \emph{cup product} $\cup:CC^*(\mathcal{A},\mathcal{A})\otimes CC^*(\mathcal{A},\mathcal{A})\to CC^*(\mathcal{A},\mathcal{A})$ as
    \begin{equation}
        (\alpha\cup\beta)(x_k,\cdots,x_1)=\sum(-1)^{\varrho}\mu^{k-j_1-j_2+2}(x_k,\cdots,\alpha^{j_1}(x_{i_1+j_1},\cdots,x_{i_1+1}),\cdots,\beta^{j_2}(x_{i_2+j_2},\cdots,x_{i_2+1}),\cdots,x_1),
    \end{equation}
    where the sign is given by $\varrho=\norm{\alpha}(1+\norm{\beta})+\norm{\alpha}\cdot\dagger_{i_1}+\norm{\beta}\cdot\dagger_{i_2}$.
\end{defn}

The cup product satisfies the graded Leibniz rule
\begin{equation}
    \delta(\alpha\cup\beta)=\delta\alpha\cup\beta+(-1)^{\norm{\alpha}+1}\alpha\cup\delta\beta.
\end{equation}

Now we move on to the definition of Hochschild homology.

\begin{defn}
    Given an $A_\infty$-algebra $\mathcal{A}$, we define its \emph{Hochschild chain complex} by 
    \begin{equation}
        CC_*(\mathcal{A},\mathcal{A})=\bigoplus_{k\geq0} CC_k(\mathcal{A},\mathcal{A}),\quad CC_k(\mathcal{A},\mathcal{A})=\mathcal{A}\otimes \mathcal{A}[1]^{\otimes k}
    \end{equation}
    with the differential
    \begin{align}
        b(z\otimes x_k\otimes\cdots\otimes x_1)&=\sum (-1)^{\natural_i^j}\mu^{i+j+1}(x_j,\cdots,x_1,z,x_k,\cdots, 
        x_{k-i+1})\otimes x_{k-i}\otimes\cdots\otimes x_{j+1}\\
        &+\sum (-1)^{\dagger_i}(z\otimes x_k\otimes\cdots\otimes\mu^j(x_{i+j},\cdots,x_{i+1})\otimes x_{i}\otimes\cdots\otimes x_1),
    \end{align}
    where the sign is given by
    \begin{equation}
        \natural_i^j=\left(\sum_{r=j+1}^k \norm{x_r}+\norm{ z}\right)\left(\sum_{r=1}^{j} \norm{x_r}\right)+\sum_{r=j+1}^{k-i}\norm{x_r}.
    \end{equation}
\end{defn}

Now we define two operations on Hochschild chains: the cap product and the Lie derivative.

\begin{defn}
    The \emph{cap product} $\cap:CC^*(\mathcal{A},\mathcal{A})\otimes CC_*(\mathcal{A},\mathcal{A})\to CC_*(\mathcal{A},\mathcal{A})$ is defined as
    \begin{equation}
        \alpha\cap(z\otimes x_k\otimes\cdots\otimes x_1)=\sum(-1)^\vartheta\mu^{i+k-t-s+2}(x_i,\cdots,\alpha^s(x_{j+s},\cdots,x_{j+1}),\cdots,x_1,z,x_k,\cdots,x_{t+1})\otimes x_{t}\otimes\cdots\otimes x_{i+1},
    \end{equation}
    where we insert the cochain to the left of $z$: the sign is given by
    \begin{equation}
    \vartheta=1+\norm{\alpha}\dagger_j+(\dagger_i+\norm{\alpha})\left( \norm{z}+\sum_{r=i+1}^k\norm{x_r} \right)+\sum_{r=i+1}^t\norm{x_r}.
    \end{equation}
\end{defn}

The cap product equips $HH_*(A,A)$ with a module structure over $HH^*(A,A)$.

\begin{defn}
    The \emph{Lie derivative} $\mathcal{L}:CC^*(\mathcal{A},\mathcal{A})\otimes CC_*(\mathcal{A},\mathcal{A})\to CC_*(\mathcal{A},\mathcal{A})$ is defined as
    \begin{align}
        \mathcal{L}_{\alpha}(z\otimes x_k\otimes\cdots\otimes x_1)=&\sum(-1)^\varsigma\alpha^{i+j+1}(x_j,\cdots,x_1,z,x_k,\cdots,x_{k-i+1})\otimes x_{k-i}\otimes\cdots\otimes x_{j+1}\\
        +&\sum(-1)^{\norm{\alpha}\cdot\dagger_{i}} (z\otimes x_k\otimes\cdots\otimes\alpha^j(x_{i+j},\cdots,x_{i+1})\otimes x_i\otimes\cdots\otimes x_1)
    \end{align}
    with the sign
    \begin{equation}
        \varsigma=\left(\sum_{r=1}^j\norm{x_r}\right)\left(\norm{z}+\sum_{r=j+1}^k\norm{x_r}\right)+\norm{\alpha}\left(\sum_{r=j+1}^{k-i}\norm{x_r}\right).
    \end{equation}
\end{defn}

We remark that the differential for Hochschild chains is a special instance of the Lie derivative, by
\begin{equation}
    b=\mathcal{L}_\mu.
\end{equation}

The Lie derivative operation satisfies the relation
\begin{equation} \label{eqn:Lie-Leibniz}
    -\mathcal{L}_{\delta\alpha}=b\circ \mathcal{L}_\alpha-(-1)^{\norm{\alpha}}\mathcal{L}_\alpha\circ b.
\end{equation}
The Lie derivative and the cap product satisfy the following homotopy form of the noncommutative Cartan identity, with the operation $\iota^\beta_\varphi$ defined in the appendix:
\begin{equation}\label{eqn:nc-Cartan-formula}
 L_\beta(\varphi\cap\upsilon)-(-1)^{\|\beta\|\,|\varphi|}\varphi\cap L_\beta(\upsilon)={}(-1)^{\|\beta\|}[\beta,\varphi]\cap\upsilon+\iota^{\delta\beta}_\varphi(\upsilon)-(-1)^{\|\beta\|}\iota^\beta_{\delta\varphi}(\upsilon)-b\bigl(\iota^\beta_\varphi(\upsilon)\bigr)+(-1)^{\|\beta\|+\|\varphi\|}\iota^\beta_\varphi\bigl(b(\upsilon)\bigr).
\end{equation}
For closed cochains, this gives the corresponding compatibility on
Hochschild homology. At chain level, the homotopy terms must be retained.

This fits into the general picture of noncommutative Cartan calculus: indeed, in the smooth commutative setting where the HKR isomorphism applies, Hochschild cochains correspond to polyvector fields, and Hochschild chains correspond to differential forms. In this analogy, the cup product corresponds to the wedge product of polyvector fields, and the cap product corresponds to contraction of differential forms by polyvector fields. Formula \eqref{eqn:nc-Cartan-formula} expresses compatibility between the Lie derivative and contraction, up to the displayed chain homotopy.

The identities recalled in this subsection will be used in subsection \ref{sssec:HH-grading} to show that infinitesimal gradings define canonical grading operators on Hochschild (co)chains and that the cap-product pairing is compatible with these operators.

\subsection{Functors, bimodules, and invariance of $[\kappa]$ under quasi-equivalence} \label{ssec:kappa-invariance}

To show the invariance of the pre-Euler class $[\kappa]$ under quasi-equivalences, it is necessary to generalize the previous discussion to $A_\infty$-categories. All constructions of Hochschild invariants from the previous subsection carry over to $A_\infty$-categories by summing over all composable sequences of morphisms.

We first recall the notions of $A_\infty$-functors and bimodules.

\begin{defn}
    Let $\mathcal{C}$, $\mathcal{D}$ be two $\mathbb{Z}/2$-graded $A_\infty$-categories. An \emph{$A_\infty$-functor} $\mathcal{F}:\mathcal{C}\to\mathcal{D}$ is defined as a map $\mathcal{F}:Ob(\mathcal{C})\to Ob(\mathcal{D})$, together with a collection of multilinear maps
    \begin{equation}
        \mathcal{F}^d:\mathcal{C}(L_{d-1},L_d)\otimes\cdots\otimes\mathcal{C}(L_0,L_1)\to\mathcal{D}(\mathcal{F}(L_0),\mathcal{F}(L_d))
    \end{equation}
    of degree $1-d$ modulo $2$ that satisfy the $A_\infty$-functor relations:
    \begin{equation}
    \begin{aligned}
        &\sum \mu^k(\mathcal{F}^{i_k}(x_n,\cdots,x_{n-i_k+1}),\cdots,\mathcal{F}^{i_1}(x_{i_1},\cdots,x_1))\\
        =&\sum(-1)^{\dagger_j} \mathcal{F}^{n-t+1}(x_n,\cdots,\mu^t(x_{j+t},\cdots,x_{j+1}),\cdots,x_1).
    \end{aligned}    
    \end{equation}
    An $A_\infty$-functor is said to be a \emph{quasi-equivalence} if it is cohomologically fully faithful and essentially surjective.
\end{defn}

\begin{defn} \label{def:A-infty-bimodule}
    A \emph{$(\mathcal{C},\mathcal{D})$-bimodule} $\mathcal{B}$ over two $A_\infty$-categories $\mathcal{C}$, $\mathcal{D}$ is defined as an assignment of $k$-modules $\mathcal{B}(K,L)$ for each pair of objects $K\in\mathcal{C},L\in\mathcal{D}$, together with $k$-linear maps
    \begin{align} \label{eqn:A-infty-bimodule-reln}
    \begin{split}
        \mu_\mathcal{B}^{i\mid1\mid j}:~&\mathcal{D}(L_{i-1},L_i)\otimes\cdots\otimes\mathcal{D}(L_0,L_1)\\
        \otimes&~\mathcal{B}(K_0,L_0)\otimes\mathcal{C}(K_1,K_0)\otimes\cdots\otimes\mathcal{C}(K_{j},K_{j-1})\to\mathcal{B}(K_j,L_i)[1-i-j]
    \end{split}
    \end{align}
    that satisfy the $A_\infty$-bimodule relations. In particular, $(\mathcal{B}(K,L),\mu^{0\mid1\mid0})$ is a cochain complex.
\end{defn}

A complete description of the $A_\infty$-bimodule relations together with signs can be found in the discussion following \cite[Equation 2.2]{Seidel08c}, together with the definition of $A_\infty$-bimodule homomorphisms.

To each $A_\infty$-functor $\mathcal{F}:\mathcal{C}\to\mathcal{D}$, there is a naturally associated $(\mathcal{C},\mathcal{D})$-bimodule called the graph bimodule.

\begin{defn}
    Let $\mathcal{F}:\mathcal{C}\to\mathcal{D}$ be an $A_\infty$-functor. Then the \emph{graph bimodule} $\Upsilon_\mathcal{F}$ is defined as
    \begin{equation}
        \Upsilon_\mathcal{F}(K,L)\coloneqq \mathcal{D}(\mathcal{F}(K),L),
    \end{equation}
    together with the bimodule operations 
    \begin{equation}
        \mu^{i\mid1\mid j}_{\Upsilon_\mathcal{F}}:\mathcal{D}(L_{i-1},L_i)\otimes\cdots\otimes\mathcal{D}(L_0,L_1)\otimes\mathcal{D}(\mathcal{F}(K_0),L_0)\otimes\mathcal{C}(K_1,K_0)\otimes\cdots\otimes\mathcal{C}(K_{j},K_{j-1})\to\mathcal{D}(\mathcal{F}(K_j),L_i)[1-i-j]
    \end{equation}
    sending $\{y_s\in\mathcal{D}(L_{s-1},L_{s})\}$, $z\in\mathcal{D}(\mathcal{F}(K_0),L_0)$, $\{x_t\in\mathcal{C}(K_t,K_{t-1})\}$ to
    \begin{equation}
        \mu^{i\mid1\mid j}_{\Upsilon_\mathcal{F}}(y_i,\cdots,z,\cdots,x_j)=\sum\mu_\mathcal{D}^{i+r+1}(y_i,\cdots,y_1,z,\mathcal{F}^{s_r}(x_1,\cdots,x_{s_r}),\cdots,\mathcal{F}^{s_1}(x_{j-s_1+1},\cdots,x_j)).
    \end{equation}
\end{defn}

To compare a Hochschild cochain in $\mathcal{D}$ with a Hochschild cochain in $\mathcal{C}$, we must regard them as Hochschild cochains with the same bimodule coefficients. The following $(\mathcal{C},\mathcal{C})$-bimodule $\Delta_\mathcal{F}$, together with the maps $\mathcal{F}_*,\mathcal{F}^*$, will allow this by inserting a $\mathcal{C}$-cochain into the inputs of $\mathcal{F}$ and inserting the outputs of $\mathcal{F}$ into a $\mathcal{D}$-cochain.

\begin{defn} \label{defn:double-pullback-bimodule}
    Define the \emph{pullback of the diagonal bimodule}
    \begin{equation}
        \Delta_\mathcal{F}(X,Y)=\mathcal{D}(\mathcal{F}(X),\mathcal{F}(Y)),
    \end{equation}
    together with the bimodule operations
    \begin{equation}
        \mu^{i\mid1\mid j}_{\Delta_\mathcal{F}}:\mathcal{C}(Y_{i-1},Y_i)\otimes\cdots\otimes\mathcal{C}(Y_0,Y_1)\otimes\mathcal{D}(\mathcal{F}(X_0),\mathcal{F}(Y_0))\otimes\mathcal{C}(X_1,X_0)\otimes\cdots\otimes\mathcal{C}(X_{j},X_{j-1})\to\mathcal{D}(\mathcal{F}(X_j),\mathcal{F}(Y_i))[1-i-j]
    \end{equation}
    sending $\{y_s\in\mathcal{C}(Y_{s-1},Y_{s})\}$, $z\in\mathcal{D}(\mathcal{F}(X_0),\mathcal{F}(Y_0))$, $\{x_t\in\mathcal{C}(X_t,X_{t-1})\}$ to
    \begin{equation}
        \mu^{i\mid1\mid j}_{\Delta_\mathcal{F}}(y_i,\cdots,z,\cdots,x_j)=\sum\mu_\mathcal{D}^{p+q+1}(\mathcal{F}^{a_p}(y_i,\cdots),\cdots,\mathcal{F}^{a_1}(y_{a_1},\cdots,y_1),z,\mathcal{F}^{b_q}(x_1,\cdots,x_{b_q}),\cdots,\mathcal{F}^{b_1}(\cdots,x_j)).
    \end{equation}
\end{defn}

Associated to this bimodule are natural maps between Hochschild cochain complexes
\begin{align}
    \mathcal{F}_*&:CC^*(\mathcal{C},\mathcal{C})\to CC^*(\mathcal{C},\Delta_\mathcal{F})\\
    \mathcal{F}^*&:CC^*(\mathcal{D},\mathcal{D})\to CC^*(\mathcal{C},\Delta_\mathcal{F})
\end{align}
defined as
\begin{align}
    (\mathcal{F}_*\phi)^n(x_n,\cdots,x_1)&=\sum(-1)^{\norm{\phi}\cdot\dagger_i}\mathcal{F}^{n-j+1}(x_n,\cdots,\phi^j(x_{i+j},\cdots,x_{i+1}),x_i,\cdots,x_1),\label{eqn:pushfwd-functor}\\
    (\mathcal{F}^*\psi)^k(y_k,\cdots,y_1)&=\sum\psi^r(\mathcal{F}^{s_r}(y_k,\cdots,y_{k-s_r+1}),\cdots,\mathcal{F}^{s_1}(y_{s_1},\cdots,y_1)) \label{eqn:pullback-functor}.
\end{align}

Moreover, one can check that if $\mathcal{F}$ is a quasi-equivalence, the two maps $\mathcal{F}_*$ and $\mathcal{F}^*$ induce isomorphisms on the Hochschild cohomologies.

With this language, we can prove the invariance of the cohomology class of the pre-Euler cocycle $[\kappa]$ up to quasi-equivalence.

\begin{proposition} \label{prop:kappa-invariance}
    Let $\mathcal{C}$, $\mathcal{D}$ be $\mathbb{Z}/2$-graded $A_\infty$-categories, and let $\mathcal{F}:\mathcal{C}\to\mathcal{D}$ be an $A_\infty$-quasi-equivalence. Then the induced map
    \begin{equation}
        [\mathcal{F}^*]^{-1}\circ[\mathcal{F}_*]:HH^*(\mathcal{C},\mathcal{C})\to HH^*(\mathcal{D},\mathcal{D})
    \end{equation}
    matches the cohomology class $[\kappa_\mathcal{C}]$ with $[\kappa_\mathcal{D}]$.
\end{proposition}
\begin{proof}
    Define an element $H_\mathcal{F}\in CC^*(\mathcal{C},\Delta_\mathcal{F})$ as
    \begin{equation}
        H_\mathcal{F}^k=(1-k)\mathcal{F}^k.
    \end{equation}
    Then we claim that
    \begin{equation} \label{eqn:invariance-kappa}
        \mathcal{F}^*\kappa_\mathcal{D}-\mathcal{F}_* \kappa_\mathcal{C}=\delta H_\mathcal{F}
    \end{equation}
    holds at chain level. Once this is established, then
    \begin{equation}
        [\mathcal{F}^*\kappa_\mathcal{D}]=[\mathcal{F}_*\kappa_\mathcal{C}]
    \end{equation}
    in $HH^*(\mathcal{C},\Delta_\mathcal{F})$, and since $\mathcal{F}_*$ and $\mathcal{F}^*$ induce isomorphisms on cohomology we are done.

    The verification of equation (\ref{eqn:invariance-kappa}) can be done by checking
    \begin{equation}
    \begin{aligned}
        &(\mathcal{F}^*\kappa_\mathcal{D})^n(x_n,\cdots,x_1)-(\mathcal{F}_*\kappa_\mathcal{C})^n(x_n,\cdots,x_1)\\
        =&\sum(2-r)\mu_\mathcal{D}^r(\mathcal{F}^{s_r}(x_n,\cdots,x_{n-s_r+1}),\cdots,\mathcal{F}^{s_1}(x_{s_1},\cdots,x_1))\\
        -&\sum(-1)^{\dagger_i}(2-j)\mathcal{F}^{n-j+1}(x_n,\cdots,\mu^j(x_{i+j},\cdots,x_{i+1}),x_i,\cdots,x_1)
    \end{aligned}    
    \end{equation}
    and
    \begin{equation}
    \begin{aligned}
        (\delta H_\mathcal{F})^n(x_n,\cdots,x_1)=&\sum(-1)^{\dagger_i}(j-n)\mathcal{F}^{n-j+1}(x_n,\cdots,\mu^j(x_{i+j},\cdots,x_{i+1}),x_i,\cdots,x_1)\\
        -&\sum(r-n)\mu_\mathcal{D}^r(\mathcal{F}^{s_r}(x_n,\cdots,x_{n-s_r+1}),\cdots,\mathcal{F}^{s_1}(x_{s_1},\cdots,x_1)),
    \end{aligned}    
    \end{equation}
    so we conclude that
    \begin{equation}
    \begin{aligned}
        &(\mathcal{F}^*\kappa_\mathcal{D}-\mathcal{F}_*\kappa_\mathcal{C}-\delta H_\mathcal{F})^n(x_n,\cdots,x_1)\\
        =(2-n)&(\sum \mu^k(\mathcal{F}^{i_k}(x_n,\cdots,x_{n-i_k+1}),\cdots,\mathcal{F}^{i_1}(x_{i_1},\cdots,x_1))\\
        -&\sum(-1)^{\dagger_i}\mathcal{F}^{n-j+1}(x_n,\cdots,\mu^j(x_{i+j},\cdots,x_{i+1}),x_i,\cdots,x_1)),
    \end{aligned}    
    \end{equation}
    which vanishes by the $A_\infty$-functor equations.
\end{proof}

Thus, the cohomology class $[\kappa]$ provides an obstruction to finding a $\mathbb{Z}$-grading up to quasi-equivalence.

\begin{corollary}
    If the pre-Euler cohomology class $[\kappa_\mathcal{C}]\in HH^{\mathrm{even}}(\mathcal{C},\mathcal{C})$ is nonzero for a $\mathbb{Z}/2$-graded $A_\infty$-category $\mathcal{C}$, then $\mathcal{C}$ cannot be quasi-equivalent to a $\mathbb{Z}$-graded $A_\infty$-category.
\end{corollary}

We emphasize that the cocycle $\kappa$ itself is not invariant under quasi-isomorphisms: only its cohomology class is. The cochain $H_\mathcal{F}$ measures the change of $\kappa$ at chain level.

\subsection{Infinitesimal gradings for $A_\infty$-categories and graded objects} \label{ssec:infinitesimal-grading-categories}

In this subsection, we define infinitesimal gradings and explore their properties. This definition is motivated by the line of thought developed in \cite{SeidelSolomon}, \cite{Seidel13}, and \cite{Seidel14}. We note the difference in sign conventions, in that we follow the conventions from \cite{Seidel08} throughout: \cite[Section 8]{SeidelSolomon} explains how to match the signs.

\subsubsection{Infinitesimal gradings for $A_\infty$-categories}

We first generalize Definition \ref{defn:inf-grading-algebra} to $A_\infty$-categories.

\begin{defn}
    Let $\mathcal{C}$ be a $\mathbb{Z}/2$-graded $A_\infty$-category. If $[\kappa_\mathcal{C}]=0$, we define $\mathcal{C}$ to be \emph{infinitesimally gradable}, and we define an \emph{infinitesimal grading} for $\mathcal{C}$ to be a choice of cochain $\eta\in CC^{\mathrm{odd}}(\mathcal{C},\mathcal{C})$ such that $\delta\eta=\kappa_\mathcal{C}$.
\end{defn}

The length 1 part $\eta^1$ of the infinitesimal grading plays the role of the Euler operator acting on each morphism space, while the higher length parts encode the higher homotopies measuring the compatibility with the $A_\infty$-structure.

\begin{example} \label{ex:Euler-cocycle}
    Let $\mathcal{C}$ be a $\mathbb{Z}$-graded $A_\infty$-category, or more generally a $k$-graded $A_\infty$-category. Thus, for each morphism space between two objects $X,Y\in\mathcal{C}$, there exists a direct sum decomposition
    \begin{equation}
        \hom_\mathcal{C}(X,Y)=\bigoplus_{\lambda\in k} \hom_\mathcal{C}(X,Y)^\lambda
    \end{equation}
    into homogeneous subspaces. Then the Euler cochain, defined as the Hochschild cochain with length $1$ component given by
    \begin{equation}
        Eu^1:x\in\hom_\mathcal{C}(X,Y)^\lambda\mapsto \lambda\cdot x
    \end{equation}
    and trivial length $\neq1$ part, is an infinitesimal grading for $\mathcal{C}$.
\end{example}

Now we define how infinitesimal gradings define grading operators on the morphism spaces of $\mathcal{C}$. Under appropriate assumptions, this induces generalized eigenspace decompositions of the morphism spaces. Recall that since $\mathcal{C}$ is an uncurved $A_\infty$-category, $\mu^0$ is always zero, and therefore for any infinitesimal grading $\eta$,
\begin{equation}
    -\mu^1(\eta^0\vert_K)=\kappa^0\vert_K=2\mu^0\vert_K=0
\end{equation}
holds for any object $K\in\mathcal{C}$.

\begin{defn}
    Given a $\mathbb{Z}/2$-graded $A_\infty$-category $\mathcal{C}$ together with an infinitesimal grading $\eta\in CC^{\mathrm{odd}}(\mathcal{C},\mathcal{C})$, we define an object $K\in\mathcal{C}$ to be \emph{$\eta$-gradeable} if 
    \begin{equation}
        \eta^0\vert_K=\mu^1(\xi_K)
    \end{equation}
    for some $\xi_K\in\hom^{\mathrm{even}}_\mathcal{C}(K,K)$, and define the pair $(K,\xi_K)$ to be an \emph{infinitesimally graded object} with respect to $\eta$.
\end{defn}

The morphism space between two infinitesimally graded objects admits a grading operator as follows.

\begin{defn}
    Let $(\mathcal{C},\eta)$ be a $\mathbb{Z}/2$-graded $A_\infty$-category with an infinitesimal grading, and let $(K,\xi_K)$ and $(L,\xi_L)$ be two infinitesimally graded objects. Then we define the \emph{grading operator} $\Gamma_{(K,\xi_K),(L,\xi_L)}:\hom_\mathcal{C}(K,L)\to\hom_\mathcal{C}(K,L)$ as
    \begin{equation}
        \Gamma_{(K,\xi_K),(L,\xi_L)}(x)=\eta^1(x)-\mu^2(x,\xi_K)-(-1)^{\norm{x}}\mu^2(\xi_L,x).
    \end{equation}
\end{defn}

The following lemma explains how $\Gamma$ can be regarded as a ``grading operator'':

\begin{lemma} \label{lem:grading-operator-commutes-d}
    For any two infinitesimally graded objects $(K,\xi_K)$, $(L,\xi_L)$, the following identity holds:
    \begin{equation}
        \Gamma_{(K,\xi_K),(L,\xi_L)}\circ\mu^1-\mu^1\circ\Gamma_{(K,\xi_K),(L,\xi_L)}=\mu^1.
    \end{equation}
\end{lemma}
\begin{proof}
    We omit the subscript in $\Gamma$ for this proof. Since $\eta$ is an infinitesimal grading,
    \begin{equation}
        \mu^1=\kappa^1=(\delta\eta)^1=-\mu^1\circ\eta^1-\mu^2(\eta^0_L,~)-\mu^2(~,\eta^0_K)+\eta^1\circ\mu^1,
    \end{equation}
    and from the $A_\infty$-relations, we have
    \begin{align}
        &\mu^1(\mu^2(\xi_L,x))+\mu^2(\xi_L,\mu^1(x))=(-1)^{\lvert x\rvert}\mu^2(\eta^0_L,x),\\
        &\mu^1(\mu^2(x,\xi_K))-\mu^2(\mu^1(x),\xi_K)=-\mu^2(x,\eta^0_K) 
    \end{align}
    since $\mu^1(\xi_L)=\eta^0_L$ and $\mu^1(\xi_K)=\eta^0_K$. Substituting into the above, we have
    \begin{equation}
    \begin{aligned}
        \mu^1(x)=&-\mu^1\circ\eta^1(x)-(-1)^{\lvert x\rvert}\mu^1(\mu^2(\xi_L,x))-(-1)^{\lvert x\rvert}\mu^2(\xi_L,\mu^1(x))\\
        &+\mu^1(\mu^2(x,\xi_K))-\mu^2(\mu^1(x),\xi_K)+\eta^1\circ\mu^1(x)\\
        =&~\Gamma\circ\mu^1(x)-\mu^1\circ\Gamma(x).\qedhere
    \end{aligned}    
    \end{equation}
\end{proof}

In other words, the differential has weight 1 with respect to the grading operator.

We note that if $\mu^1(x)=0$, then $\mu^1\circ\Gamma(x)=0$. Therefore $\Gamma$ descends to a map on cohomology:
\begin{equation}
    \bar{\Gamma}_{(K,\xi_K),(L,\xi_L)}:H^*\hom_\mathcal{C}(K,L)\to H^*\hom_\mathcal{C}(K,L)
\end{equation}

\begin{defn}
    We define $[x]\in H^*\hom_\mathcal{C}(K,L)$ to be \emph{pure of weight $\lambda$} with respect to $(\eta,\xi_K,\xi_L)$ if 
    \begin{equation}
        \bar{\Gamma}_{(K,\xi_K),(L,\xi_L)}([x])=\lambda\cdot [x].
    \end{equation}
    By Lemma \ref{lem:grading-operator-commutes-d}, the weight of a pure cohomology class is independent of its cocycle representative.
\end{defn}

Therefore if the base field $k$ is algebraically closed and the morphism space $\hom_\mathcal{C}(K,L)$ is cohomologically finite dimensional, the eigenvalue decomposition of the grading operator induces a weight decomposition of the cohomology groups when the induced map on cohomology is semisimple.

The choice of $\xi_K$ records the choice of primitive for the object-level obstruction $\eta^0_K$. The graded lift defined below incorporates this choice into the object set. This choice makes the object-level obstruction disappear: as a consequence, the lifted infinitesimal grading has trivial length zero part, and its length 1 part exactly equals the grading operator.

\begin{defn} \label{def:lifted-category}
    Given a $\mathbb{Z}/2$-graded $A_\infty$-category $\mathcal{C}$ and an infinitesimal grading $\eta$, define the \emph{$\eta$-graded lift} of $\mathcal{C}$ to be the $A_\infty$-category $\Tilde{\mathcal{C}}_\eta$ with objects
    \begin{equation}
        Ob(\Tilde{\mathcal{C}}_\eta)=\{(K,\xi_K)\mid \mu^1(\xi_K)=\eta^0_K\},
    \end{equation}
    where $(K,\xi_K)$ is an infinitesimally graded object in $\mathcal{C}$ with respect to $\eta$, and morphisms
    \begin{equation}
        \hom_{\Tilde{\mathcal{C}}_\eta}((K,\xi_K),(L,\xi_L))=\hom_{\mathcal{C}}(K,L).
    \end{equation}
    The $A_\infty$-operations are inherited from those of $\mathcal{C}$. This naturally inherits an infinitesimal grading defined as
    \begin{equation}
        \Tilde{\eta}=\eta+\delta\xi,
    \end{equation}
    where $\xi$ is a Hochschild cochain with only a length-zero part $\xi_K$ for each object $(K,\xi_K)$. One can check
    \begin{equation}
        \delta\Tilde{\eta}=\kappa_{\Tilde{\mathcal{C}}_\eta},    
    \end{equation}
    so indeed $\Tilde{\eta}$ is an infinitesimal grading, with trivial length zero part. The grading operator is then
    \begin{equation}
         \Gamma_{(K,\xi_K),(L,\xi_L)}(x)=\Tilde{\eta}^1(x),
    \end{equation}
    which is the same as the previous definition after unpacking the definition of the Hochschild codifferential. An equivalent definition is that the objects of $\Tilde{\mathcal{C}}_\eta$ are the pairs $(K,\xi_K)$ such that $\Tilde{\eta}^0$ vanishes.
\end{defn}

We now analyze the dependence of grading operators on the infinitesimal grading. Suppose that $\eta$, $\eta'$ are two different choices of infinitesimal gradings for a $\mathbb{Z}/2$-graded $A_\infty$-category $\mathcal{C}$. For each $K\in\mathcal{C}$, $\eta^0_K-\eta'^0_K$ need not be a coboundary, so the set of $\eta$-gradeable objects need not agree with the set of $\eta'$-gradeable objects.

Thus it is reasonable to restrict to $\eta'$ that differ from $\eta$ by a coboundary $\delta\beta$, i.e. $\eta'=\eta-\delta\beta$. The underlying gradeable objects are unchanged, but their chosen primitives must be adjusted: since
\begin{equation}
    (\eta')^0_K=\eta^0_K+\mu^1(\beta^0_K),
\end{equation}
the corresponding $\eta'$-graded object is
\begin{equation}
    (K,\xi_K+\beta^0_K).
\end{equation}
With this identification, the grading operators for $\eta$ and $\eta'$ can be identified:
\begin{proposition}
    Let $(K,\xi_K),(L,\xi_L)$ be $\eta$-graded objects, and let $\eta'=\eta-\delta\beta$ be a different choice of infinitesimal grading. Then the grading operators
    \begin{align}
        \Gamma^\eta_{(K,\xi_K),(L,\xi_L)}:\hom_\mathcal{C}(K,L)\to\hom_\mathcal{C}(K,L)\\
        \Gamma^{\eta'}_{(K,\xi_K+\beta^0_K),(L,\xi_L+\beta^0_L)}:\hom_\mathcal{C}(K,L)\to\hom_\mathcal{C}(K,L)
    \end{align}
    agree on cohomology.
\end{proposition}
\begin{proof}
    Using that $\eta'=\eta-\delta\beta$, on length 1 we have
    \begin{equation}
        (\eta')^1(x)=\eta^1(x)+\mu^1\circ\beta^1(x)+\mu^2(x,\beta^0_K)+(-1)^{\norm{x}}\mu^2(\beta^0_L,x)+\beta^1\circ\mu^1(x),
    \end{equation}
    which translates to the equation
    \begin{equation}
        \Gamma^{\eta'}-\Gamma^\eta=\mu^1\circ\beta^1+\beta^1\circ\mu^1,
    \end{equation}
    completing the proof.
\end{proof}

We now demonstrate that for a general infinitesimal grading, the length 2 part $\eta^2$ together with the choices $\xi_K$ provides a chain homotopy measuring the failure of $\eta^1$ to be a chain-level derivation of $\mu^2$.

\begin{proposition}
    Let $(\mathcal{C},\eta)$ be as before, and let $(K,\xi_K)$, $(L,\xi_L)$, $(N,\xi_N)$ be infinitesimally graded objects. Then for cocycles $x\in\hom_\mathcal{C}(K,L)$ and $y\in\hom_\mathcal{C}(L,N)$,
    \begin{equation}
        [\Gamma(\mu^2(y,x))]=[\mu^2(\Gamma(y),x)]+[\mu^2(y,\Gamma(x))]
    \end{equation}
    holds in the cohomology $H^*\hom_\mathcal{C}(K,N)$. In particular, if $[x]$, $[y]$ are pure cohomology classes with weights $\lambda$, $\nu$ with respect to the grading operators, then a nonzero product $[\mu^2(y,x)]$ has weight $\lambda+\nu$.
\end{proposition}
\begin{proof}
    Since $\kappa^2=0$, the Hochschild codifferential equations for $\delta\eta=\kappa$ show
    \begin{align}
        0&=\eta^1(\mu^2(y,x))-\mu^2(\eta^1(y),x)-\mu^2(y,\eta^1(x))+(-1)^{\norm{x}}\eta^2(\mu^1(y),x)+\eta^2(y,\mu^1(x))-\mu^1(\eta^2(y,x))\\
        &-\mu^3(y,x,\eta^0_K)-\mu^3(y,\eta^0_L,x)-\mu^3(\eta^0_N,y,x)
    \end{align}
    Now using that $\mu^1(\xi_K)=\eta^0_K$, the $A_\infty$-relations show that
    \begin{equation}
    \resizebox{\textwidth}{!}{$
    \begin{aligned}
        &\mu^3(y,x,\eta^0_K)=\mu^2(\mu^2(y,x),\xi_K)-\mu^2(y,\mu^2(x,\xi_K))+\mu^3(y,\mu^1(x),\xi_K)-(-1)^{\lvert x\rvert}\mu^3(\mu^1(y),x,\xi_K)-\mu^1(\mu^3(y,x,\xi_K)),\\
        &(-1)^{\lvert x\rvert}\mu^3(y,\eta^0_L,x)=\mu^2(y,\mu^2(\xi_L,x))-(-1)^{\lvert x\rvert}\mu^2(\mu^2(y,\xi_L),x)+\mu^3(y,\xi_L,\mu^1(x))+(-1)^{\lvert x\rvert}\mu^3(\mu^1(y),\xi_L,x)+\mu^1(\mu^3(y,\xi_L,x)),\\
        &(-1)^{\lvert y\rvert+\norm{x}}\mu^3(\eta^0_N,y,x)=\mu^2(\xi_N,\mu^2(y,x))-(-1)^{\lvert x\rvert}(\mu^2(\mu^2(\xi_N,y),x)+\mu^3(\xi_N,\mu^1(y),x))+\mu^3(\xi_N,y,\mu^1(x))+\mu^1(\mu^3(\xi_N,y,x)).
    \end{aligned}
    $}
    \end{equation}
    Altogether, one has
    \begin{align}
        0&=\Gamma(\mu^2(y,x))-\mu^2(\Gamma(y),x)-\mu^2(y,\Gamma(x))+\Gamma^2(y,\mu^1(x))+(-1)^{\norm{x}}\Gamma^2(\mu^1(y),x)-\mu^1(\Gamma^2(y,x)),
    \end{align}
    where
    \begin{equation}
        \Gamma^2(y,x)=\eta^2(y,x)-\mu^3(y,x,\xi_K)-(-1)^{\norm{x}}\mu^3(y,\xi_L,x)-(-1)^{\norm{y}+\norm{x}}\mu^3(\xi_N,y,x).
    \end{equation}
    Thus, $\Gamma^2$ provides a homotopy between the two sides of the Leibniz rule for $\Gamma$. In particular, in cohomology 
    \begin{equation}
        [\Gamma(\mu^2(y,x))]=[\mu^2(\Gamma(y),x)]+[\mu^2(y,\Gamma(x))]
    \end{equation}
    holds.
\end{proof}

\begin{corollary}
    If $\mathcal{C}$ is cohomologically unital, $[\Gamma(e_L)]=0$ for the unit $e_L\in H^*\hom_\mathcal{C}(L,L)$.
\end{corollary}

The analogous statement for $\mu^3$ requires correction terms involving both $\Gamma^2$ and a new operation $\Gamma^3$ with three inputs. The following formula will show the general pattern.

For $x\in\hom_\mathcal{C}(K,L)$, $y\in\hom_\mathcal{C}(L,N)$, $z\in\hom_C(N,M)$,
\begin{align}
    -\mu^3(z,y,x)&=\Gamma(\mu^3(z,y,x))-\mu^3(\Gamma(z),y,x)-\mu^3(z,\Gamma(y),x)-\mu^3(z,y,\Gamma(x))\\
    &+\Gamma^2(z,\mu^2(y,x))+(-1)^{\norm{x}}\Gamma^2(\mu^2(z,y),x)-\mu^2(\Gamma^2(z,y),x)-\mu^2(z,\Gamma^2(y,x))\\
    &+\Gamma^3(z,y,\mu^1(x))+(-1)^{\norm{x}}\Gamma^3(z,\mu^1(y),x)+(-1)^{\norm{y}+\norm{x}}\Gamma^3(\mu^1(z),y,x)-\mu^1(\Gamma^3(z,y,x))
\end{align}
for the operation
\begin{align}
    \Gamma^3(z,y,x)&=\eta^3(z,y,x)-\mu^4(z,y,x,\xi_K)-(-1)^{\norm{x}}\mu^4(z,y,\xi_L,x)\\&-(-1)^{\norm{y}+\norm{x}}\mu^4(z,\xi_N,y,x)-(-1)^{\norm{z}+\norm{y}+\norm{x}}\mu^4(\xi_M,z,y,x).
\end{align}

Therefore as given, $\eta^1$ is not completely compatible with the $A_\infty$-structure $\mu$. In subsection \ref{ssec:straightening}, we will prove the ``straightening lemma'', which shows that after passing to the $\eta$-graded lift, one can find a quasi-equivalent $A_\infty$-category for which the $A_\infty$-operations and the induced infinitesimal grading are compatible with the generalized eigenvalue decomposition from the length $1$ part $\eta^1$. Readers who are not interested in infinitesimal gradings for Hochschild (co)chains are invited to skip to subsection \ref{ssec:straightening}.

\subsubsection{Infinitesimal gradings on Hochschild (co)chains} \label{sssec:HH-grading}

An infinitesimal grading also induces a grading operator on the Hochschild cochain and chain complexes. To do so, we recall the definition of the length operator for an $A_\infty$-algebra $\mathcal{A}$; the same definition works verbatim for $A_\infty$-categories.

\begin{defn}
    The \emph{length operator} $\Lambda$ on Hochschild cochains and chains sends $\alpha^k\in \hom(\mathcal{A}[1]^{\otimes k},\mathcal{A})$ and $z\otimes x_k\otimes\cdots\otimes x_1\in \mathcal{A}\otimes\mathcal{A}[1]^{\otimes k}$ to
    \begin{align}
        \Lambda(\alpha^k)&= k\cdot\alpha^k,\\
        \Lambda(z\otimes x_k\otimes\cdots\otimes x_1)&=k\cdot (z\otimes x_k\otimes\cdots\otimes x_1).
    \end{align}
\end{defn}

Now we define the Hochschild cochain and chain grading operators.

\begin{defn} \label{defn:HH-grading-operators}
    Given an infinitesimally graded $A_\infty$-category $(\mathcal{C},\eta)$, we define the \emph{Hochschild cochain grading operator} $\Gamma^{CC}$ as
    \begin{equation}
        \Gamma^{CC}(\alpha)=[\eta,\alpha]+\Lambda(\alpha),
    \end{equation}
    and the \emph{Hochschild chain grading operator} $\Gamma_{CC}$ as
    \begin{equation}
        \Gamma_{CC}(\beta)=\mathcal{L}_\eta(\beta)-\Lambda(\beta).
    \end{equation}
\end{defn}

We verify that their commutator with the differential is the differential.

\begin{proposition}
    The operators $\Gamma^{CC}$, $\Gamma_{CC}$ defined above are grading operators; i.e., they satisfy
    \begin{align}
        \Gamma^{CC}\circ\delta-\delta\circ\Gamma^{CC}&=\delta,\\
        \Gamma_{CC}\circ b-b\circ\Gamma_{CC}&=b.
    \end{align}
\end{proposition}
\begin{proof}
    Recall from Equation (\ref{eqn:gersten-jacobi}) that
    \begin{equation}
        [\mu,[\eta,\alpha]]=[[\mu,\eta],\alpha]+[\eta,[\mu,\alpha]],
    \end{equation}
    and since $\delta=-[\mu,~]$, we obtain
    \begin{equation} \label{eqn:eta-jacobi}
        [\kappa,\alpha]=\delta[\eta,\alpha]-[\eta,\delta\alpha].
    \end{equation}
    Now for a Hochschild cochain $\alpha$, we find that
    \begin{align}
        (\Lambda\circ\delta \alpha)^k(\cdots)&=-\sum (-1)^{\dagger_i\cdot\norm{\alpha}}k\mu^{k-j+1}(\cdots,\alpha^j(\cdots),x_i,\cdots)\\
        &+\sum (-1)^{\dagger_i+\norm{\alpha}}k\alpha^{k-j+1}(\cdots,\mu^j(\cdots),x_i,\cdots),
    \end{align}
    and
    \begin{align}
        (\delta\circ\Lambda \alpha)^k(\cdots)&=-\sum (-1)^{\dagger_i\cdot\norm{\alpha}}j\mu^{k-j+1}(\cdots,\alpha^j(\cdots),x_i,\cdots)\\
        &+\sum (-1)^{\dagger_i+\norm{\alpha}}(k-j+1)\alpha^{k-j+1}(\cdots,\mu^j(\cdots),x_i,\cdots),
    \end{align}
    and finally,
    \begin{align}
        -[\kappa,\alpha]^k(\cdots)&=-\sum (-1)^{\dagger_i\cdot\norm{\alpha}}(2-k+j-1)\mu^{k-j+1}(\cdots,\alpha^j(\cdots),x_i,\cdots)\\
        &+\sum (-1)^{\dagger_i+\norm{\alpha}}(2-j)\alpha^{k-j+1}(\cdots,\mu^j(\cdots),x_i,\cdots).
    \end{align}
    Thus, we conclude that
    \begin{equation}
        \Lambda\circ\delta-\delta\circ\Lambda-[\kappa,~]=\delta,
    \end{equation}
    which together with (\ref{eqn:eta-jacobi}) proves the desired result. For Hochschild chains, recall from (\ref{eqn:Lie-Leibniz})
    \begin{equation} \label{eqn:Lie-eta}
        -\mathcal{L}_{\kappa}=b\circ \mathcal{L}_\eta-\mathcal{L}_\eta\circ b,
    \end{equation}
    and by direct computation
    \begin{align}
        \mathcal{L}_{\kappa}(z\otimes x_k\otimes\cdots\otimes x_1)=&\sum(-1)^{\natural_i^j}(1-i-j)\mu^{i+j+1}(x_j,\cdots,x_1,z,x_k,\cdots,x_{k-i+1})\otimes x_{k-i}\otimes\cdots\otimes x_{j+1}\\
        +&\sum (-1)^{\dagger_i}(2-j)(z\otimes x_k\otimes\cdots\otimes\mu^j(x_{i+j},\cdots,x_{i+1})\otimes x_i\otimes\cdots\otimes x_1),
    \end{align}
    and
    \begin{align}
        \Lambda\circ b(z\otimes x_k\otimes\cdots\otimes x_1)&=\sum (-1)^{\natural_i^j}(k-i-j)(\mu^{i+j+1}(x_j,\cdots,x_1,z,x_k,\cdots, 
        x_{k-i+1})\otimes x_{k-i}\otimes\cdots\otimes x_{j+1})\\
        &+\sum (-1)^{\dagger_i}(k-j+1)(z\otimes x_k\otimes\cdots\otimes\mu^j(x_{i+j},\cdots,x_{i+1})\otimes x_{i}\otimes\cdots\otimes x_1),
    \end{align}
    and finally
    \begin{align}
        b\circ\Lambda(z\otimes x_k\otimes\cdots\otimes x_1)&=\sum (-1)^{\natural_i^j}k\cdot(\mu^{i+j+1}(x_j,\cdots,x_1,z,x_k,\cdots, 
        x_{k-i+1})\otimes x_{k-i}\otimes\cdots\otimes x_{j+1})\\
        &+\sum (-1)^{\dagger_i}k\cdot(z\otimes x_k\otimes\cdots\otimes\mu^j(x_{i+j},\cdots,x_{i+1})\otimes x_{i}\otimes\cdots\otimes x_1).
    \end{align}
    From this we conclude that
    \begin{equation}
        b\circ\Lambda-\Lambda\circ b+\mathcal{L}_\kappa=b,
    \end{equation}
    which combined with (\ref{eqn:Lie-eta}) concludes our proof.
\end{proof}

\begin{example}
    As a sanity check, we compute the above operators for a $\mathbb{Z}$-graded $A_\infty$-category $\mathcal{C}$ with its canonical infinitesimal grading given by the Euler cochain $Eu$. Then for a Hochschild cochain $\alpha$,
    \begin{align}
        [Eu,\alpha]^k(x_k,\cdots,x_1)&=\deg(\alpha^k)\cdot\alpha^k(x_k,\cdots,x_1)\\
        \Lambda(\alpha)^k&=k\cdot\alpha^k,    
    \end{align}
    and thus we compute that
    \begin{equation}
        \Gamma^{CC}(\alpha)^k=(\deg(\alpha^k)+k)\cdot\alpha^k,
    \end{equation} 
    which is the associated grading for the component $\alpha^k$ in the Hochschild cochain complex. Similarly for a Hochschild chain $z\otimes x_k\otimes\cdots\otimes x_1$, we find
    \begin{align}
        \mathcal{L}_{Eu}(z\otimes x_k\otimes\cdots\otimes x_1)&=(\deg(z)+\sum\deg(x_i))\cdot (z\otimes x_k\otimes\cdots\otimes x_1),\\
        \Lambda(z\otimes x_k\otimes\cdots\otimes x_1)&=k\cdot(z\otimes x_k\otimes\cdots\otimes x_1),
    \end{align}
    and therefore
    \begin{equation}
        \Gamma_{CC}(z\otimes x_k\otimes\cdots\otimes x_1)= \left(\deg{ z}+\sum(\deg{x_i}-1)\right)\cdot (z\otimes x_k\otimes\cdots\otimes x_1)
    \end{equation}
    which is indeed its associated grading in the Hochschild chain complex.
\end{example}

The induced grading operators can be shown to be invariant under changing the infinitesimal grading $\eta$ by a coboundary.

\begin{proposition}
    Let $(\mathcal{C},\eta)$ be as before, and let $\eta'=\eta+\delta\beta$ be another choice of infinitesimal grading for $\mathcal{C}$. Then the grading operators $\Gamma^{CC}$, $\Gamma_{CC}$ for $\eta$ and $\eta'$ induce the same map on (co)homology.
\end{proposition}
\begin{proof}
    From (\ref{eqn:Gersten-Leibniz}),
    \begin{equation}
        [\delta\beta,\phi]=\delta[\beta,\phi]+[\beta,\delta\phi]
    \end{equation}
    shows that $[\beta,~]$ is an odd-degree homotopy between the two grading operators $\Gamma^{CC}$. Similarly,
    \begin{equation}
        -\mathcal{L}_{\delta\beta}=b\circ\mathcal{L}_\beta+\mathcal{L}_\beta\circ b    
    \end{equation}
    shows that $-\mathcal{L}_\beta$ is an odd-degree homotopy between the two grading operators $\Gamma_{CC}$. In particular, the induced maps on (co)homology agree.
\end{proof}

Therefore the induced grading operators on Hochschild homology and cohomology only depend on the equivalence class of the infinitesimal grading modulo Hochschild coboundaries.

From Definition \ref{defn:HH-grading-operators}, we can define Hochschild cohomology and homology classes of pure weight with respect to an infinitesimal grading. This notion of pure weight is independent of the chain-level representative of the (co)homology class. The induced actions of $\Gamma^{CC}$ on cohomology and $\Gamma_{CC}$ on homology are related by the cap product:

\begin{proposition} \label{prop:cap-product-matches-grading}
    Let $(\mathcal{C},\eta)$ be as before. Then for a cocycle $\phi\in CC^*(\mathcal{C},\mathcal{C})$ and a cycle $\upsilon\in CC_*(\mathcal{C},\mathcal{C})$,
    \begin{equation} \label{eqn:cap-product-matches-grading}
        [\Gamma^{CC}(\phi)\cap\upsilon]+[\phi\cap\Gamma_{CC}(\upsilon)]=[\Gamma_{CC}(\phi\cap\upsilon)]
    \end{equation}
    in $HH_*(\mathcal{C},\mathcal{C})$.
\end{proposition}
\begin{proof}
    From the noncommutative Cartan formula (\ref{eqn:nc-Cartan-formula}) specialized to cycle and cocycle inputs, one has
    \begin{equation}
        \mathcal{L}_\eta(\phi\cap\upsilon)-\phi\cap\mathcal{L}_\eta\upsilon=[\eta,\phi]\cap\upsilon+\iota^\kappa_\phi(\upsilon)-b(\iota^\eta_\phi(\upsilon)).
    \end{equation}
    The failure of (\ref{eqn:cap-product-matches-grading}) to hold at chain level for a cycle input $\upsilon=z\otimes x_k\otimes\cdots\otimes x_1$ can be measured by $\kappa=\{(2-k)\mu^k\}$ as follows: the sign $\vartheta_\phi$ can be found in \eqref{eqn:Appendix-sign};
    \begin{equation}
        \Lambda(\phi)\cap\upsilon-\phi\cap\Lambda(\upsilon)+\Lambda(\phi\cap\upsilon)=\sum(-1)^{\vartheta_\phi}\kappa(x_i,\cdots,\phi^s(x_{j+s},\cdots,x_{j+1}),\cdots,z,\cdots,x_{t+1})\otimes\cdots\otimes x_{i+1}.
    \end{equation}
    Let us define the right-hand side of the above equation as an operation $\iota_\phi^\kappa(\upsilon)$. Then one can check that
    \begin{equation} \label{eqn:App-1}
        b\circ\iota^\eta_\phi-(-1)^{\norm{\phi}+\norm{\eta}}\iota^\eta_\phi\circ b=\iota_\phi^{\delta\eta}-(-1)^{\norm{\eta}}\iota^\eta_{\delta\phi}-\mathcal{L}_{\eta}(\phi\cap~)+(-1)^{\norm{\eta}\cdot\lvert\phi\rvert}\phi\cap\mathcal{L}_\eta+(-1)^{\norm{\eta}}[\eta,\phi]\cap~,
    \end{equation}
    proved in the appendix as Lemma \ref{lem:App-1}. But now since $\kappa=\delta\eta$ and $\phi$ is a cocycle, we conclude that
    \begin{equation}
        \Gamma^{CC}(\phi)\cap\upsilon+\phi\cap\Gamma_{CC}(\upsilon)-\Gamma^{CC}(\phi\cap\upsilon)=b\circ\iota^\eta_\phi(\upsilon)-\iota^\eta_\phi\circ b(\upsilon),
    \end{equation}
    which is exact since $\upsilon$ is a cycle. Therefore (\ref{eqn:cap-product-matches-grading}) holds on the level of homology.
\end{proof}

Now we recall the definition of Calabi-Yau structures from \cite[Definition A.2]{Sheridan}.

\begin{defn}
    Let $\mathcal{C}$ be a $\mathbb{Z}/2$-graded $A_\infty$-category, whose morphism spaces are cohomologically finite-dimensional. A \emph{weak proper Calabi-Yau structure} $\sigma$ for $\mathcal{C}$ of parity $d$ is a chain map
    \begin{equation}
        \sigma:CC_*(\mathcal{C},\mathcal{C})\to k[-d]
    \end{equation}
    such that for every pair of objects $K,L$, the composition
    \begin{equation}
        [\sigma]\circ\iota\circ[\mu^2]:H^*\hom_\mathcal{C}(K,L)\otimes H^*\hom_\mathcal{C}(L,K)\to H^*\hom_\mathcal{C}(L,L)\to HH_*(\mathcal{C},\mathcal{C})\to k
    \end{equation}
    is a perfect pairing. The map induced by the inclusion of length-zero Hochschild chains $\iota$ sends an endomorphism cocycle to its corresponding length 0 chain.
\end{defn}

Given an infinitesimal grading, one can define the infinitesimal dimension of a Calabi-Yau structure.

\begin{defn}
    Let $(\mathcal{C},\eta)$ be a $\mathbb{Z}/2$-graded $A_\infty$-category equipped with an infinitesimal grading, and let $\sigma$ be a weak proper Calabi-Yau structure for $\mathcal{C}$. We define $\sigma$ to be of \emph{infinitesimal dimension $n$} with respect to $\eta$ if
    \begin{equation}
        [\sigma]\circ(\overline{\Gamma}_{CC}-n\cdot\id)=0.
    \end{equation}
    Under semisimplicity assumptions, $\sigma$ is cohomologically supported on the eigenvalue $n$ summand of the grading operator $\Gamma_{CC}$.
\end{defn}

Recall from \cite[Lemma A.2]{Sheridan} that a weak proper Calabi-Yau structure $\sigma$ defines a nondegenerate pairing between Hochschild homology and cohomology by the map
    \begin{equation} \label{eqn:HH^*HH_*-pairing}
        \sigma\circ\cap:CC^*(\mathcal{C},\mathcal{C})\otimes CC_*(\mathcal{C},\mathcal{C})\to CC_*(\mathcal{C},\mathcal{C})\to k.
    \end{equation}

The following proposition is the infinitesimal analogue of Poincaré duality.

\begin{proposition}\label{prop:inf-grading-cap}
    Let $(\mathcal{C},\eta)$ be as before, and let $\sigma$ be a weak proper Calabi-Yau structure of $\eta$-infinitesimal dimension $n$. Then for any cocycle $\phi\in CC^*(\mathcal{C},\mathcal{C})$ and cycle $\upsilon\in CC_*(\mathcal{C},\mathcal{C})$ of pure weights that pair nontrivially by (\ref{eqn:HH^*HH_*-pairing}), their weights must add up to $n$.
\end{proposition}
\begin{proof}
    Since we assumed that $\phi$ and $\upsilon$ are pure, there exist $\lambda$, $\nu$ such that
    \begin{equation}
        [\Gamma^{CC}(\phi)]=\lambda\cdot[\phi],~[\Gamma_{CC}(\upsilon)]=\nu\cdot[\upsilon].
    \end{equation}
    From Proposition \ref{prop:cap-product-matches-grading} we have
    \begin{equation}
        [\Gamma^{CC}(\phi)\cap\upsilon]+[\phi\cap\Gamma_{CC}(\upsilon)]=[\Gamma_{CC}(\phi\cap\upsilon)],
    \end{equation}
    but since the Calabi-Yau structure $\sigma$ satisfies
    \begin{equation}
        [\sigma]\circ(\Gamma_{CC}-n\cdot\id)=0,
    \end{equation}
    we conclude
    \begin{align}
        (\lambda+\nu)[\sigma(\phi\cap\upsilon)] &= [\sigma(\Gamma_{CC}(\phi)\cap\upsilon+\phi\cap\Gamma_{CC}(\upsilon))]\\
        &= [\sigma(\Gamma_{CC}(\phi\cap\upsilon))]\\
        &= n[\sigma(\phi\cap\upsilon)].
    \end{align}
    Since the pairing $\sigma(\phi\cap\upsilon)$ is assumed to be nonzero, the weights $\lambda$ and $\nu$ add up to $n$.
\end{proof}

The existence of a nonzero gradeable object also forces the infinitesimal dimension of a homogeneous weak proper Calabi--Yau structure to be unique.

\begin{proposition}\label{prop:CY-dim-unique}
    Let $(\mathcal{C},\eta)$ be a c-finite infinitesimally graded $A_\infty$-category over a field of characteristic zero. Suppose there exists a $\eta$-gradeable nonzero object $L$. Then all weak proper Calabi--Yau structures on $\mathcal{C}$ homogeneous with respect to $\eta$ have the same infinitesimal dimension, given by the formula
    \begin{equation}
        d=\frac{2\mathrm{Tr}(\Gamma_L)}{\dim_k H^*\hom_\mathcal{C}(L,L)}
    \end{equation}
    for the grading operator $\Gamma_L:\hom_\mathcal{C}(L,L)\to\hom_\mathcal{C}(L,L)$ for a choice of graded lift $(L,\xi_L)$.
\end{proposition}
\begin{proof}
    A weak proper Calabi--Yau structure induces a pairing
    \begin{equation}
        \Phi_L:H^*\hom_\mathcal{C}(L,L)\otimes H^*\hom_\mathcal{C}(L,L)\to k
    \end{equation}
    and homogeneity of infinitesimal dimension $d$ shows that
    \begin{equation}
        \Phi_L(\Gamma_Lx,y)+\Phi_L(x,\Gamma_Ly)=d\Phi_L(x,y).
    \end{equation}
    Therefore setting the matrix for the linear transformation $\Gamma_L$ as $A$ and the bilinear form $\Phi_L$ as $B$, we have
    \begin{equation}
        A^tB+BA=d\cdot B,
    \end{equation}
    and since $B$ is invertible we multiply by $B^{-1}$ on both sides and take the trace.
\end{proof}

\subsection{Straightening of infinitesimal gradings} \label{ssec:straightening}

We prove the ``straightening lemma'' in this subsection. Throughout this subsection, the underlying base field $k$ will be assumed to be algebraically closed.

We first introduce the definition of a straightened infinitesimally graded category.

\begin{defn}
    Let $(\mathcal{C},\eta)$ be an infinitesimally graded minimal $A_\infty$-category  such that $\eta^0\equiv 0$. For each $K,L\in\mathcal{C}$ and $\lambda\in k$, define $\hom_\mathcal{C}(K,L)_\lambda$ to be the $\lambda$-generalized eigenspace of the operator
    \begin{equation}
        \eta^1:\hom_\mathcal{C}(K,L)\to \hom_\mathcal{C}(K,L).
    \end{equation}
    We say that $(\mathcal{C},\eta)$ is \emph{$\eta^1$-straightened} if for every $d\geq 2$, collection of objects $X_0,\cdots,X_{d+1}\in\mathcal{C}$, and collection of generalized eigenvalues $(\lambda_1,\cdots,\lambda_{d+1})\in k^{d+1}$, 
    \begin{align}
        \mu^{d+1}\bigl( \hom_\mathcal{C}(X_{d},X_{d+1})_{\lambda_{d+1}}\otimes\cdots\otimes\hom_\mathcal{C}(X_0,X_1)_{\lambda_1} \bigr)&\subset\hom_\mathcal{C}(X_0,X_{d+1})_{\lambda_1+\cdots+\lambda_{d+1}+1-d},\\
        \eta^{d}\bigl( \hom_\mathcal{C}(X_{d-1},X_d)_{\lambda_d}\otimes\cdots\otimes\hom_\mathcal{C}(X_0,X_1)_{\lambda_1} \bigr)&\subset\hom_\mathcal{C}(X_0,X_d)_{\lambda_1+\cdots+\lambda_{d}+1-d}.
    \end{align}
    In words, the operations $\mu^{d+1}$ and $\eta^d$ preserve the generalized eigenvalue decomposition for $\eta^1$ up to a shift by $1-d$. We remark that since $\mathcal{C}$ is minimal, this condition is automatic for $d=1$.
\end{defn}

The assumptions for the above definition can be achieved by passing to a minimal model.

\begin{lemma}
    For an infinitesimally graded $A_\infty$-category, let $\Tilde{\mathcal{C}_\eta}$ be the lifted category with objects $(K,\xi_K)$. Then there exists a minimal infinitesimally graded $A_\infty$-category $(\mathcal{C}',\eta')$ quasi-equivalent to $(\Tilde{\mathcal{C}_\eta},\Tilde{\eta})$, such that $(\eta')^0\equiv0$.
\end{lemma}
\begin{proof}
    By the homological perturbation lemma, there exists a quasi-equivalence
    \begin{equation}
        \mathcal{F}:\Tilde{\mathcal{C}_\eta}\to\mathcal{C}'
    \end{equation}
    identity on objects, where $\mathcal{C}'$ is minimal. Now pick a quasi-inverse $\mathcal{G}_0$ to $\mathcal{F}$ again identity on objects. Since $\mathcal{C}'$ is minimal, $\mathcal{F}\circ\mathcal{G}_0$ has invertible linear part, and thus is an $A_\infty$-diffeomorphism. We define
    \begin{equation}
        \mathcal{G}=\mathcal{G}_0\circ(\mathcal{F}\circ\mathcal{G}_0)^{-1},
    \end{equation}
    which is a strict right inverse to $\mathcal{F}$. Then the formula for transporting an infinitesimal grading from (\ref{eqn:invariance-kappa}) shows that for the map $\mathcal{G}_*^\Delta$ defined in (\ref{eqn:pushfwd-G}),
    \begin{equation}
        \eta'=\mathcal{G}_*^\Delta\circ\mathcal{F}_*(\Tilde{\eta})+\mathcal{G}_*^\Delta(H_\mathcal{F})
    \end{equation}
    is an infinitesimal grading for $\mathcal{C}'$. Since all functors under consideration have trivial length zero part and $\Tilde{\eta}^0\equiv0$, we conclude that the transported infinitesimal grading $\eta'$ satisfies $(\eta')^0\equiv0$.
\end{proof}

Therefore the minimality and triviality of $\eta^0$ hypotheses can be achieved, after passing to the lifted category $\Tilde{\mathcal{C}}_\eta$. Under these hypotheses, we state the straightening lemma.

\begin{thm} \label{thm:straighten-infinitesimal-grading}
    Let $(\mathcal{C},\eta)$ be a $\mathbb{Z}/2$-graded $A_\infty$-category with an infinitesimal grading, and assume that the base field $k$ is algebraically closed. We further assume that $\eta^0\equiv0$, and that $\mathcal{C}$ is minimal with finite-dimensional morphism spaces. Then there exists an $A_\infty$-diffeomorphism $\Phi:(\mathcal{C},\{\mu\})\to(\mathcal{C},\{\mu'\})$ such that the pair consisting of $\mu'$ and the induced infinitesimal grading $\eta'$ is $(\eta')^1$-straightened. Moreover, the length 1 parts $\eta^1$ and $(\eta')^1$ agree.
\end{thm}

We have chosen to present the theorem in this way to emphasize that the choice of gradeable objects $K$ and corresponding choices of primitives $\xi_K$ are fixed under the $A_\infty$-diffeomorphism chosen in the theorem. 

In particular, straightening produces a $k$-grading without requiring semisimplicity.

\begin{corollary}
Under the hypotheses of Theorem \ref{thm:straighten-infinitesimal-grading}, the straightened $A_\infty$-category $(\mathcal C,\mu')$ admits a $k$-grading given by the generalized eigenspaces of $\eta^1$.
\end{corollary}
\begin{proof}
Since the morphism spaces are finite-dimensional and $k$ is algebraically closed, their generalized eigenspaces give direct-sum decompositions. The straightening condition says precisely that $\mu'^d$ has weight $2-d$ with respect to these decompositions.
\end{proof}

If $\eta^1$ is semisimple, it agrees with the Euler operator of this grading. In general, their difference is nilpotent on each morphism space.

We now begin the proof of Theorem \ref{thm:straighten-infinitesimal-grading}. The key idea is isolated in the following lemma:

\begin{lemma}
    Let $(\mathcal{C},\eta)$ be an infinitesimally graded minimal $A_\infty$-category such that $\eta^0\equiv0$, and $\eta^1$ acts semisimply on all morphism spaces. Suppose there exists some $d$ such that for any $2\leq l\leq d-1$, the pair $(\mu^{l+1},\eta^l)$ is $\eta^1$-straightened. Then there exists an $A_\infty$-diffeomorphism $\Phi_d:\mathcal{C}\to\mathcal{C}$ such that the pair consisting of the induced $A_\infty$-structure $\mu'$ and infinitesimal grading $\eta'$ is $\eta^1$-straightened for all lengths $2\leq l\leq d$.
\end{lemma}
\begin{proof}
    We start with the case $d=2$, and show how to straighten the pair $(\mu^3,\eta^2)$.

    By the semisimplicity assumption on $\eta^1$, we may pick a weight decomposition for each morphism space $\hom_\mathcal{C}(K,L)$. Thus, we may decompose $\mu^3$ and $\eta^2$ according to the $\eta^1$-weights of the morphism spaces of $\mathcal{C}$, and write
    \begin{align}
        \mu^3=\sum \mu^3_{\lambda},&~\mu^3_{\lambda}\vert_{L_0,L_1,L_2,L_3}:\hom_\mathcal{C}(L_2,L_3)\otimes\hom_\mathcal{C}(L_1,L_2)\otimes\hom_{\mathcal{C}}(L_0,L_1)\to\hom_\mathcal{C}(L_0,L_3)\langle-1+\lambda\rangle, \label{eqn:mu-weight-decomp} \\
        \eta^2=\sum \eta^2_{\lambda},&~\eta^2_{\lambda}\vert_{K_0,K_1,K_2}:\hom_\mathcal{C}(K_1,K_2)\otimes\hom_\mathcal{C}(K_0,K_1)\to\hom_\mathcal{C}(K_0,K_2)\langle-1+\lambda\rangle, \label{eqn:eta-weight-decompt}
    \end{align}
    where $\langle\lambda\rangle$ means the shift in the $\eta^1$-eigenvalues. Since $\eta$ is an infinitesimal grading for $\mathcal{C}$, we have
    \begin{equation}\label{eqn:Z-Z/2-obstruction}
        -[\eta^1,\mu^{3}]-\mu^{3}=\sum(-1)^{\dagger_j}\eta^2(\cdots,\mu^2(~,~),\cdots)-\mu^2(~,\eta^2(\cdots))-\mu^2(\eta^2(\cdots),~),
    \end{equation}
    and decomposing by the $\eta^1$-weights of the morphism spaces of $\mathcal{C}$, we obtain that
    \begin{equation}\label{eqn:Z-Z/2-length-decomp}
        -\lambda\mu^{3}_\lambda(z,y,x)=\eta^2_\lambda(z,\mu^2(y,x))+(-1)^{\norm{x}}\eta^2_\lambda(\mu^2(z,y),x)-\mu^2(z,\eta^2_\lambda(y,x))-\mu^2(\eta^2_\lambda(z,y),x),
    \end{equation}
    which we abbreviate as
    \begin{equation}
        -\lambda\mu^3_\lambda=\sum(-1)^{\dagger_j}\eta^2_{\lambda}(\cdots,\mu^2(~,~),\cdots)-\mu^2(~,\eta^2_{\lambda}(\cdots))-\mu^2(\eta^2_{\lambda}(\cdots),~).
    \end{equation}

     Now define an $A_\infty$-functor $\mathcal{F}$ as
    \begin{equation}\label{eqn:grading-straight-functor}
        \mathcal{F}^1=\id,~\mathcal{F}^2=\sum_{\lambda\neq0}\frac{1}{\lambda}\eta^2_{\lambda}.
    \end{equation}
    Then this functor induces a new $A_\infty$-structure $\mu_{new}$ on the image
    \begin{equation}
        \mathcal{F}:(\mathcal{C},\mu)\to(\mathcal{C},\mu_{new}),
    \end{equation}
    and we claim that $\mu^3_{new}$ preserves the $\eta^1$-weights.

    First, since $\mu^1$ is zero and $\mathcal{F}^1=\id$, we obtain that
    \begin{equation}
        \mu^1_{new}=\mu^1,\quad \mu^2_{new}=\mu^2.
    \end{equation}
    For $\mu^3_{new}$, the $A_\infty$-functor equations show that
    \begin{equation}
        \mu^3_{new}=\sum(-1)^{\dagger_j}\mathcal{F}^2(\cdots,\mu^2(~,~),\cdots)+\mu^3-\mu^2_{new}(~,\mathcal{F}^2(~,))-\mu^2_{new}(\mathcal{F}^2(~,~),~)
    \end{equation}
    and from (\ref{eqn:Z-Z/2-length-decomp}) we conclude that
    \begin{equation}
        \mu^3_{new}=\mu^3_0,
    \end{equation}
    which is the summand of $\mu^3$ that preserves the $\eta^1$-weights of the morphism spaces of $\mathcal{C}$.

    We now examine the induced infinitesimal grading $\eta_{new}$. Recall the maps from (\ref{eqn:pushfwd-functor}), (\ref{eqn:pullback-functor}):
    \begin{equation}
    \begin{tikzcd}
        CC^*(\mathcal{C},\mathcal{C}) \arrow[r,"\mathcal{F}_*"] &CC^*(\mathcal{C},\Delta_\mathcal{F}) &CC^*(\mathcal{C}_\mathcal{F},\mathcal{C}_\mathcal{F}) \arrow[l,"\mathcal{F}^*"'] 
    \end{tikzcd}
    \end{equation}
    for $\mathcal{C}_\mathcal{F}=(\mathcal{C},\mu_{new})$. Choose a chain-level inverse functor $\mathcal{G}$ for $\mathcal{F}$, which must satisfy
    \begin{equation}
        \mathcal{G}^1=\mathrm{Id},\quad \mathcal{G}^2=-\mathcal{F}^2.
    \end{equation}
    Then there exists a chain map
    \begin{equation} \label{eqn:pushfwd-G}
        \mathcal{G}_*^{\Delta_{\mathcal{F}}}:CC^*(\mathcal{C},\Delta_\mathcal{F})\to CC^*(\mathcal{C}_\mathcal{F},\mathcal{C}_\mathcal{F})
    \end{equation}
    defined as
    \begin{equation}
        \mathcal{G}_*^{\Delta_\mathcal{F}}(\alpha)^k(y_k,\cdots,y_1)=\sum\alpha^r\Bigl(\mathcal{G}^{s_r}(y_k,\ldots,y_{k-s_r+1}),\ldots,\mathcal{G}^{s_1}(y_{s_1},\ldots,y_1)\Bigr),
    \end{equation}
    which is a chain-level inverse to $\mathcal{F}^*$. Then the induced infinitesimal grading is defined as
    \begin{equation}
        \eta_{new}=\mathcal{G}_*^{\Delta_\mathcal{F}}\circ(\mathcal{F}_*\eta+H_\mathcal{F})
    \end{equation}
    for $H_\mathcal{F}=\{(1-d)\mathcal{F}^d\}$. For length 1, we compute that
    \begin{equation}
        \eta^1_{new}=\eta^1(\mathcal{G}^1\circ\mathcal{F}^1(~))=\eta^1,
    \end{equation}
    so indeed the length $1$ part is unchanged. The formula for the length 2 part reads
    \begin{equation}
        \eta^2_{new}=(\mathcal{F}_*\eta)^1\circ \mathcal{G}^2(~,~)+(\mathcal{F}_*\eta)^2(\mathcal{G}^1(~),\mathcal{G}^1(~))-\mathcal{F}^2(\mathcal{G}^1(~),\mathcal{G}^1(~)),
    \end{equation}
    and using the formula for $\mathcal{F}_*$ we can further compute that
    \begin{align}
        \eta^2_{new}&= (\mathcal{F}_*\eta)^1\circ \mathcal{G}^2(~,~)+(\mathcal{F}_*\eta)^2(~,~)-\mathcal{F}^2(~,~)\\
        &=\mathcal{F}^1(\eta^1\circ\mathcal{G}^2(~,~))+\mathcal{F}^2(\eta^1(~),~)+\mathcal{F}^2(~,\eta^1(~))+\mathcal{F}^1(\eta^2(~,~))-\mathcal{F}^2(~,~)
    \end{align}
    and using that $\mathcal{F}^1=\id$, $\mathcal{G}^2=-\mathcal{F}^2$,
    \begin{align}
        \eta^2_{new}&=-\eta^1(\mathcal{F}^2(~,~))+\mathcal{F}^2(\eta^1(~),~)+\mathcal{F}^2(~,\eta^1(~))+\eta^2(~,~)-\mathcal{F}^2(~,~)\\
        &=\eta^2+\sum_{\lambda\neq0}\frac{1}{\lambda}\left(\eta^2_\lambda(\eta^1(~),~)+\eta^2_\lambda(~,\eta^1(~))-\eta^1(\eta^2_\lambda(~,~))-\eta^2_\lambda\right)\\
        &= \sum_\lambda \eta^2_\lambda-\sum_{\lambda\neq0}\eta^2_\lambda= \eta^2_0
    \end{align}
    from (\ref{eqn:eta-weight-decompt}) and the definition of the functor $\mathcal{F}$. Therefore, we have checked that the choice (\ref{eqn:grading-straight-functor}) of $\mathcal{F}$ maps $(\mu^3,\eta^2)$ to $(\mu^3_0,\eta^2_0)$, the summands that preserve the $\eta^1$-weights of $\mathcal{C}$.

    We now consider the case for $d>2$. Again take a decomposition for $\lambda\in k$
    \begin{align}
        \mu^{d+1}=\sum \mu^{d+1}_{\lambda},&~\mu^{d+1}_{\lambda}\vert_{L_0,\cdots,L_{d+1}}:\hom_\mathcal{C}(L_d,L_{d+1})\otimes\cdots\otimes\hom_{\mathcal{C}}(L_0,L_1)\to\hom_\mathcal{C}(L_0,L_{d+1})\langle-d+1+\lambda\rangle, \\
        \eta^d=\sum \eta^d_{\lambda},&~\eta^d_{\lambda}\vert_{K_0,\cdots,K_d}:\hom_\mathcal{C}(K_{d-1},K_d)\otimes\cdots\otimes\hom_\mathcal{C}(K_0,K_1)\to\hom_\mathcal{C}(K_0,K_d)\langle-d+1+\lambda\rangle.
    \end{align}
    We now look at the $\eta^1$-weight decomposition for the equation for $\delta\eta=\kappa$ at length $d+1$:
    \begin{equation}
        -[\eta^1,\mu^{d+1}]-(d-1)\mu^{d+1}=\sum_{2\leq j\leq d}(-1)^{\dagger_i}\eta^j(\cdots,\mu^{d-j+2}(\cdots),x_i,\cdots)-\sum_{2\leq j\leq d}\mu^j(\cdots,\eta^{d-j+2}(\cdots),\cdots).
    \end{equation}
    Since the operations $(\mu^j,\eta^{j-1})$ are already compatible with the $\eta^1$-weight decomposition for all $j\leq d$, the $\lambda\neq0$ weights entirely come from the $\mu^{d+1}$ and $\eta^d$ contributions. Therefore,
    \begin{equation}
        -\lambda\mu^{d+1}_\lambda=\sum(-1)^{\dagger_j}\eta^d_{\lambda}(\cdots,\mu^2(~,~),\cdots)-\mu^2(~,\eta^d_{\lambda}(\cdots))-\mu^2(\eta^d_{\lambda}(\cdots),~)
    \end{equation}
    still holds. Therefore by a similar computation, the $A_\infty$-functor
    \begin{equation}\label{eqn:grading-straight-functor-2}
        \mathcal{F}^1=\id,~\mathcal{F}^d=\sum_{\lambda\neq0}\frac{1}{\lambda}\eta^d_{\lambda}
    \end{equation}
    induces an $A_\infty$-structure on $\mathcal{C}$ such that
    \begin{equation}
        \mu^j_{new}=\mu^j
    \end{equation}
    for $j\leq d$, and
    \begin{equation}
        \mu^{d+1}_{new}=\mu^{d+1}_0
    \end{equation}
    the $\eta^1$-weight compatible part. Again pick a chain-level inverse $\mathcal{G}$, which satisfies $\mathcal{G}^1=\id$, $\mathcal{G}^j\equiv0$ for $2\leq j\leq d-1$, and $\mathcal{G}^d=-\mathcal{F}^d$. Then the induced infinitesimal grading
    \begin{equation}
        \eta_{new}=\mathcal{G}_*^{\Delta_\mathcal{F}}\circ(\mathcal{F}_*\eta+H_\mathcal{F})
    \end{equation}
    satisfies
    \begin{equation}
        \eta^j_{new}=\eta^j
    \end{equation}
    for all $j\leq d-1$, and at length $d$
    \begin{equation}
        \eta^d_{new}=(\mathcal{F}_*\eta)^1\circ \mathcal{G}^d(~,~)+(\mathcal{F}_*\eta)^d(\mathcal{G}^1(~),\cdots,\mathcal{G}^1(~))+(1-d)\mathcal{F}^d(\mathcal{G}^1(~),\cdots,\mathcal{G}^1(~)).
    \end{equation}
    Then following through the computations, we see that
    \begin{align}
        \eta^d_{new}&=-\eta^1(\mathcal{F}^d(~,~))+\sum\mathcal{F}^d(\cdots,\eta^1(~),\cdots)+\eta^d+(1-d)\mathcal{F}^d\\
        &=\eta^d+\sum_{\lambda\neq0}\frac{1}{\lambda}\left(\sum\eta^d_\lambda(\cdots,\eta^1(~),\cdots)-\eta^1(\eta^d_\lambda(~,~))+(1-d)\eta^d_\lambda\right)\\
        &= \sum_\lambda \eta^d_\lambda-\sum_{\lambda\neq0}\eta^d_\lambda= \eta^d_0
    \end{align}
    which concludes the proof.
\end{proof}

The non-semisimple case is a slight extension of the above argument: the additional feature is to invert the nonzero-eigenvalue part, including its nilpotent contribution.

\begin{lemma} \label{lem:straight-inductive-step}
    Let $(\mathcal{C},\eta)$ be an infinitesimally graded minimal cohomologically finite $A_\infty$-category such that $\eta^0\equiv0$. Suppose there exists some $d$ such that for any $2\leq l\leq d-1$, the pair $(\mu^{l+1},\eta^l)$ is $\eta^1$-straightened. Then there exists an $A_\infty$-diffeomorphism $\Phi_d:\mathcal{C}\to\mathcal{C}$ such that the pair consisting of the induced $A_\infty$-structure $\mu'$ and infinitesimal grading $\eta'$ is $\eta^1$-straightened for all lengths $2\leq l\leq d$.
\end{lemma}
\begin{proof}
    First define the space of multilinear maps 
    \begin{equation} \label{eqn:V_d}
         V^{d}_{X_0,\cdots,X_d}=\hom_k(\bigotimes_{i=0}^{d-1}\hom(X_i,X_{i+1}),\hom(X_0,X_{d})),\quad V^d=\prod_{X_0,\cdots,X_{d}}V^{d}_{X_0,\cdots,X_d}
    \end{equation}
    and define the defect operator $D^d:V^d\to V^d$ as
    \begin{equation}
        D^d(\alpha)=\eta^1\circ\alpha-\sum\alpha(\cdots,\eta^1(~),\cdots).
    \end{equation}
    Because of the grading shift difference for $\mu$ and $\eta$, we define
    \begin{align}
        D^d_\mu&=D^d+(d-2)\id,\\
        D^d_\eta&=D^d+(d-1)\id.
    \end{align}
    In particular, the defect operator restricted to each component
    \begin{equation}
        D^d_\mu\vert_{X_0,\cdots,X_d}:V^{d}_{X_0,\cdots,X_d}\to V^{d}_{X_0,\cdots,X_d}
    \end{equation}
    admits a generalized eigenspace decomposition over an algebraically closed base field: recall that under our minimality and c-finiteness assumptions, morphism spaces are finite dimensional. Therefore we may write
    \begin{equation}
        \mu^{d+1}=\sum\mu^{d+1}_\lambda,
    \end{equation}
    where, for each $(X_0,\cdots,X_{d+1})$, $\mu^{d+1}_\lambda$ is the generalized $\lambda$-eigenspace component for $D^{d+1}_\mu$. Define $\eta^d_\lambda$ componentwise in the same way, using $D^d_\eta$.

    By the assumptions that $(\mu^{l+1},\eta^l)$ are straightened for $l\leq d-1$, we again have for $\lambda\neq0$,
    \begin{equation} \label{eqn:lambda-defect}
        -D^{d+1}_\mu(\mu_\lambda^{d+1})=\sum(-1)^{\dagger_j}\eta^d_{\lambda}(\cdots,\mu^2(~,~),\cdots)-\mu^2(~,\eta^d_{\lambda}(\cdots))-\mu^2(\eta^d_{\lambda}(\cdots),~),
    \end{equation}
    but since we do not have semisimplicity anymore, the left-hand side cannot be written as $-\lambda\mu^{d+1}_\lambda$. Instead, write
    \begin{equation}
        D^{d+1}_\mu(\mu^{d+1}_\lambda)=\lambda\mu_\lambda^{d+1}+N_\lambda^\mu(\mu^{d+1}_\lambda).
    \end{equation}
    Here all generalized eigenspaces are understood componentwise for each object tuple. The operators $N_\lambda^\mu$ and $N_\lambda^\eta$ are the restrictions of $D_\mu^{d+1}-\lambda\mathrm{id}$ and $D_\eta^d-\lambda\mathrm{id}$ to the corresponding generalized $\lambda$-eigenspaces. These restrictions are nilpotent on each component, since the component spaces in \eqref{eqn:V_d} are finite-dimensional. We retain the product notation $V^d$ and $V^{d+1}$ introduced above. Similarly, we can write
    \begin{equation}
        D^d_\eta(\eta^d_\lambda)=\lambda\eta_\lambda^d+N_\lambda^\eta(\eta^d_\lambda).
    \end{equation}
    Then from (\ref{eqn:lambda-defect}), we obtain
    \begin{equation}
        -(\lambda\cdot\id+N_\lambda^\mu)\circ\mu^{d+1}_\lambda=\partial_{\mu^2}(\eta^d_\lambda)
    \end{equation}
    for the map $\partial_{\mu^2}:V^d\to V^{d+1}$ defined as
    \begin{equation}
        \partial_{\mu^2}\phi=\sum(-1)^{\dagger_{i-1}}\phi(\cdots,\mu^2(x_{i+1},x_i),\cdots)-\mu^2(~,\phi(\cdots))-\mu^2(\phi(\cdots),~).
    \end{equation}
    Independently, from the equation $(\delta\eta)^2=\kappa^2=0$, we obtain
    \begin{equation}
        \eta^1(\mu^2(~,~))-\mu^2(\eta^1(~),~)-\mu^2(~,\eta^1(~))=0.
    \end{equation}
    A direct calculation using the above equation gives
    \begin{equation} \label{eqn:App-2}
        D^{d+1}_\mu\circ\partial_{\mu^2}=\partial_{\mu^2}\circ D^d_\eta,
    \end{equation}
    which also shows after projecting to generalized eigenspaces
    \begin{equation}
        (\lambda\cdot\id+N^\mu_\lambda)\circ\partial_{\mu^2}=\partial_{\mu^2}\circ(\lambda\cdot\id+N^\eta_\lambda),
    \end{equation}
    since $\partial_{\mu^2}$ preserves generalized eigenspaces. Now define the functor
    \begin{equation}
        \mathcal{F}^1=\id,\quad \mathcal{F}^d=\sum_{\lambda\neq0} (\lambda\cdot\id+N^\eta_\lambda)^{-1}\eta^d_\lambda.
    \end{equation}
    All inverses and sums in this definition are taken componentwise on the corresponding generalized eigenspaces. Each component is finite-dimensional, so the nonzero-eigenvalue blocks are invertible and only finitely many eigenvalues contribute. Thus $\mathcal{F}^d$ is well-defined. We first compute
    \begin{align}
        \mu^{d+1}_{new}&=\partial_{\mu^2}\mathcal{F}^d+\mu^{d+1}\\
        &=\mu^{d+1}+\sum_{\lambda\neq0}\partial_{\mu^2}\circ(\lambda\cdot\id+N^\eta_\lambda)^{-1}\eta^d_\lambda\\
        &=\mu^{d+1}+\sum_{\lambda\neq0}(\lambda\cdot\id+N^\mu_\lambda)^{-1}\partial_{\mu^2}\eta^d_\lambda\\
        &=\mu^{d+1}-\sum_{\lambda\neq0}\mu^{d+1}_\lambda=\mu^{d+1}_0.
    \end{align}

    For the induced infinitesimal grading, the key computation is at length $d$, where
    \begin{align}
        \eta^d_{new}&=-\eta^1(\mathcal{F}^d(\cdots))+\sum\mathcal{F}^d(\cdots,\eta^1(~),\cdots)+\eta^d+(1-d)\mathcal{F}^d\\
        &=\eta^d-D^d_\eta(\mathcal{F}^d)\\
        &=\eta^d-\sum_{\lambda\neq0}(\lambda\cdot\id+N^\eta_\lambda)\mathcal{F}^d_\lambda\\
        &=\eta^d-\sum_{\lambda\neq0}\eta^d_\lambda=\eta^d_0. \qedhere
    \end{align}
\end{proof}

Using the above lemma, we will now prove Theorem \ref{thm:straighten-infinitesimal-grading}.
\begin{proof}
    The proof will be by induction. For each $d$, by Lemma \ref{lem:straight-inductive-step}, we may find an $A_\infty$-functor
    \begin{equation}
        \Phi_d:\mathcal{C}\to\mathcal{C}
    \end{equation}
    such that the induced $A_\infty$-structure $\{\mu_d^{j+1}\}$ and induced infinitesimal grading $\{\eta_d^j\}$ at each step $d$ are $\eta^1$-straightened for all lengths $j\leq d$. One also notes that the $A_\infty$-diffeomorphism $\Phi_d$ preserves $\mu^1,\cdots,\mu^d$ and $\eta^1,\cdots,\eta^{d-1}$. Therefore, the composition
    \begin{equation}
        \Phi=\cdots\circ\Phi_3\circ\Phi_2
    \end{equation}
    is well-defined, and stabilizes the operations $\mu^{d+1}$ and $\eta^d$. Therefore $\Phi$ is a well-defined $A_\infty$-diffeomorphism, which satisfies the conclusion of the theorem.
\end{proof}

The assumption that $k$ is algebraically closed is only used to obtain the eigenspace decompositions over the base field. For example, if all the eigenvalues of $\eta^1$ lie in $\mathbb{Z}$, as in the $\mathbb{Z}$-graded case, the proof works over $\mathbb{Q}$.

\subsection{Infinitesimally graded functors}

We now develop the infinitesimal grading theory for functors. To motivate the theory, recall that for a degree 0 functor $\mathcal{F}:\mathcal{C}\to\mathcal{D}$ between $\mathbb{Z}$-graded $A_\infty$-categories, we have
\begin{equation}
    (\mathcal{F}^*\mathrm{Eu}_\mathcal{D}-\mathcal{F}_*\mathrm{Eu}_\mathcal{C})^d=(1-d)\mathcal{F}^d,
\end{equation}
since the length $d$ part $\mathcal{F}^d$ has degree $1-d$. This motivates the following definition:

\begin{defn}
    For two infinitesimally graded $A_\infty$-categories $(\mathcal{C},\eta_\mathcal{C})$ and $(\mathcal{D},\eta_\mathcal{D})$, we define an $A_\infty$-functor $\mathcal{F}:\mathcal{C}\to\mathcal{D}$ to be \emph{infinitesimally gradeable} if
    \begin{equation} \label{eqn:inf-graded-functor}
        \mathcal{F}_*\eta_\mathcal{C}-\mathcal{F}^*\eta_\mathcal{D}+H_\mathcal{F}=\delta\gamma_\mathcal{F}
    \end{equation}
    for some $\gamma_\mathcal{F}\in CC^{even}(\mathcal{C},\Delta_\mathcal{F})$, and call the pair $(\mathcal{F},\gamma_\mathcal{F})$ an \emph{infinitesimally graded functor}.
\end{defn}

The cochain $H_\mathcal{F}$ is defined as $\{(1-d)\mathcal{F}^d\}$: recall from (\ref{eqn:invariance-kappa}) that
\begin{equation}
    \delta(\mathcal{F}_*\eta_\mathcal{C}-\mathcal{F}^*\eta_\mathcal{D}+H_\mathcal{F})=0,
\end{equation}
so the above definition indeed makes sense.

For example, the identity functor can always be equipped with a canonical infinitesimal grading given by the zero cocycle. However, we can also equip it with any other even degree cocycle as its infinitesimal grading.

We now assess how the infinitesimal grading of functors acts on the grading operator on morphism spaces. From the length zero part of (\ref{eqn:inf-graded-functor}), 
\begin{equation}
    \eta^0_\mathcal{F(L)}-\mathcal{F}^1(\eta^0_L)=\mu^1(\gamma_\mathcal{F}^0\vert_{\mathcal{F}(L)}).
\end{equation}
\begin{defn}
    For an infinitesimally graded $A_\infty$-functor $(\mathcal{F},\gamma_\mathcal{F})$, define the \emph{induced infinitesimal grading} on an object $(\mathcal{F}(L),\xi_{\mathcal{F}(L)})$ as
    \begin{equation}
        \xi^\mathcal{F}_{\mathcal{F}(L)} =\mathcal{F}^1(\xi_L)+\gamma_\mathcal{F}^0\vert_{\mathcal{F}(L)}.
    \end{equation}
\end{defn}

\begin{proposition}
    Given an infinitesimally graded $A_\infty$-functor $(\mathcal{F},\gamma_\mathcal{F})$ and two infinitesimally graded objects $(K,\xi_K)$, $(L,\xi_L)$, the following diagram commutes:
    \begin{equation}
    \begin{tikzcd}
        H^*\hom_\mathcal{C}((K,\xi_K),(L,\xi_L)) \arrow[rr,"\mathcal{F}^1"] \arrow[d,"\Gamma"'] &&H^*\hom_\mathcal{D}((\mathcal{F}(K),\xi^\mathcal{F}_{\mathcal{F}(K)}),(\mathcal{F}(L),\xi^\mathcal{F}_{\mathcal{F}(L)})) \arrow[d,"\Gamma^\mathcal{F}"]\\
        H^*\hom_\mathcal{C}((K,\xi_K),(L,\xi_L)) \arrow[rr,"\mathcal{F}^1"'] &&H^*\hom_\mathcal{D}((\mathcal{F}(K),\xi^\mathcal{F}_{\mathcal{F}(K)}),(\mathcal{F}(L),\xi^\mathcal{F}_{\mathcal{F}(L)}))
    \end{tikzcd}.    
    \end{equation}
    In words, $\mathcal{F}^1$ matches the induced grading operator on $H^*\hom_\mathcal{C}((K,\xi_K),(L,\xi_L))$ with that on the morphism space for the objects $\mathcal{F}(K)$ and $\mathcal{F}(L)$ with shifted infinitesimal gradings $\xi^\mathcal{F}$ involving $\gamma^0_\mathcal{F}$ terms.
\end{proposition}
\begin{proof}
    From the length $1$ part of (\ref{eqn:inf-graded-functor}), for $x\in\hom_\mathcal{C}(K,L)$,
    \begin{align}
        &(\mathcal{F}^*\eta)^1(x)-(\mathcal{F}_*\eta)^1(x)=-\mathcal{F}^1(\eta^1(x))-\mathcal{F}^2(\eta^0_L,x)-\mathcal{F}^2(x,\eta^0_K)+\eta^1(\mathcal{F}^1(x))\\
        =-&(\delta\gamma_\mathcal{F})^1(x)=\mu^1(\gamma_\mathcal{F}^1(x))+\mu^2(\mathcal{F}^1(x),\gamma^0_\mathcal{F}\vert_K)+(-1)^{\norm{x}}\mu^2(\gamma^0_\mathcal{F}\vert_L,\mathcal{F}^1(x))+\gamma^1_\mathcal{F}(\mu^1(x)).
    \end{align}
    
    Now from the $A_\infty$-functor equations, we have
    \begin{align}
        \mathcal{F}^2(x,\eta^0_K)&=\mu^2(\mathcal{F}^1(x),\mathcal{F}^1(\xi_K))+\mu^1(\mathcal{F}^2(x,\xi_K))-\mathcal{F}^1(\mu^2(x,\xi_K))+\mathcal{F}^2(\mu^1(x),\xi_K),\\(-1)^{\norm{x}}\mathcal{F}^2(\eta^0_L,x)&=\mu^2(\mathcal{F}^1(\xi_L),\mathcal{F}^1(x))+\mu^1(\mathcal{F}^2(\xi_L,x))-\mathcal{F}^1(\mu^2(\xi_L,x))-\mathcal{F}^2(\xi_L,\mu^1(x))
    \end{align}

    From these equations, we conclude that
    \begin{equation} \label{eqn:functor-commutes-grading}
        \mathcal{F}^1(\Gamma(x))-\Gamma^\mathcal{F}(\mathcal{F}^1(x))
        =\mu^1(\Xi(x))+\Xi(\mu^1(x))
    \end{equation}
    for the two formulae
    \begin{align}
        \Gamma^\mathcal{F}(y)&=\eta^1(y)-\mu^2(y,\mathcal{F}^1(\xi_K)+\gamma^0_\mathcal{F}\vert_K)-(-1)^{\norm{y}}\mu^2(\mathcal{F}^1(\xi_L)+\gamma^0_\mathcal{F}\vert_L,y),\\
        \Xi(x)&=-\gamma^1_\mathcal{F}(x)-\mathcal{F}^2(x,\xi_K)-(-1)^{\norm{x}}\mathcal{F}^2(\xi_L,x).
    \end{align}
    The right hand side of (\ref{eqn:functor-commutes-grading}) vanishes in cohomology.
\end{proof}

The $A_\infty$-functors from $\mathcal{C}$ to $\mathcal{D}$ form an $A_\infty$-category $fun(\mathcal{C},\mathcal{D})$, as described in \cite[Section 1(d)]{Seidel08}. The morphisms are given by $A_\infty$-pre-natural transformations, whose definition we now introduce.

\begin{defn}
    For two $A_\infty$-functors $\mathcal{F},\mathcal{G}:\mathcal{C}\to\mathcal{D}$, a \emph{natural transformation of $\mathbb{Z}/2$-degree g} from $\mathcal{F}$ to $\mathcal{G}$ is a family of maps
    \begin{equation}
        T^d:\hom_\mathcal{C}(X_{d-1},X_d)\otimes\cdots\otimes\hom_\mathcal{C}(X_0,X_1)\to\hom_\mathcal{D}(\mathcal{F}X_0,\mathcal{G}X_d)[g-d]
    \end{equation}
    that satisfy the $A_\infty$-natural transformation equations
    \begin{align}
        &\sum (-1)^\aleph\mu^k(\mathcal{G}^{n_k}(x_d,\cdots,x_{d-n_k+1}),\cdots,T^{n_i}(x_{n_1+\cdots+n_i},\cdots,x_{n_1+\cdots+n_{i-1}+1}),\cdots,\mathcal{F}^{n_1}(x_{n_1},\cdots,x_1))\\
        &+\sum(-1)^{\dagger_i+g}~T^{d-k+1}(x_d,\cdots,\mu^{k}_{\mathcal{C}}(x_{i+k},\cdots,x_{i+1}),x_i,\cdots,x_1)=0
    \end{align}
    for the sign
    \begin{equation}
        \aleph=(g-1)(\norm{x_{n_1+\cdots+n_{i-1}}}+\cdots+\norm{x_1}).
    \end{equation}
\end{defn}

The morphism space in $fun(\mathcal{C},\mathcal{D})$ can be interpreted as Hochschild cochains valued in a bimodule
\begin{equation}
    \hom_{fun(\mathcal{C},\mathcal{D})}(\mathcal{F},\mathcal{G})=CC^*(\mathcal{C},\prescript{}{\mathcal{F}}{\Delta}_\mathcal{G}),
\end{equation}
with the coefficient differential $\delta=-\mu^1_{fun}$ for the $A_\infty$-bimodule
\begin{equation}
    \prescript{}{\mathcal{F}}{\Delta}_\mathcal{G}(X,Y)=\hom_{\mathcal{D}}(\mathcal{F}X,\mathcal{G}Y) 
\end{equation}
defined similarly to Definition \ref{defn:double-pullback-bimodule}.

We wish to understand the pair $(\mathcal{F},\gamma_\mathcal{F})$ as an infinitesimally graded object of $fun(\mathcal{C},\mathcal{D})$. As such, the next task would be to construct a grading operator on the morphism spaces.

Our motivation comes again from the $\mathbb{Z}$-graded case. Suppose that $\mathcal{F},\mathcal{G}:\mathcal{C}\to\mathcal{D}$ are degree zero functors, and the natural transformation $T$ has degree $g$. To recover this degree, we compute that
\begin{align}
    \mathrm{Eu}^1_\mathcal{D}(T^d(x_d,\cdots,x_1))=(\sum\deg{(x_i)}+g-d) T^d(x_d,\cdots,x_1),\\
    \sum T^d(x_d,\cdots,\mathrm{Eu}^1_\mathcal{C}(x_i),\cdots,x_1)=(\sum\deg{(x_i)}) T^d(x_d,\cdots,x_1)
\end{align}
for the integer degree $g$. Therefore, the difference of the two left-hand sides plus the length term
\begin{equation}
    \mathrm{Eu}^1_\mathcal{D}(T^d)-\sum T^d(\cdots,\mathrm{Eu}^1_\mathcal{C},\cdots)+dT^d=gT^d
\end{equation}
recovers the degree of the natural transformation $T$.

This leads us to define two new operations $\iota_T,\rho_T$ associated to any pre-natural transformation $T$:
\begin{equation}
\begin{tikzcd}
    CC^*(\mathcal{C},\mathcal{C}) \arrow[r,"\rho_T"] &CC^*(\mathcal{C},\prescript{~}{\mathcal{F}}{\Delta}_{\mathcal{G}})  & CC^*(\mathcal{D},\mathcal{D}) \arrow[l,"\iota_T"'].
\end{tikzcd}    
\end{equation}
The operation $\iota_T(\eta_\mathcal{D})$ will apply the infinitesimal gradings of the target to the outputs of $T$, while the operation $\rho_T(\eta_\mathcal{C})$ will apply the infinitesimal gradings of the domain to the inputs of $T$.

\begin{defn}
    Given two $A_\infty$-functors $\mathcal{F},\mathcal{G}:\mathcal{C}\to\mathcal{D}$, and a pre-natural transformation $T\in CC^*(\mathcal{C},\prescript{~}{\mathcal{F}}{\Delta}_{\mathcal{G}})$, define the two operations $\iota_T:CC^*(\mathcal{D},\mathcal{D})\to CC^*(\mathcal{C},\prescript{~}{\mathcal{F}}{\Delta}_{\mathcal{G}})$ and $\rho_T:CC^*(\mathcal{C},\mathcal{C})\to CC^*(\mathcal{C},\prescript{~}{\mathcal{F}}{\Delta}_{\mathcal{G}})$ by
    \begin{equation}
        (\iota_T\phi)^k(x_k,\cdots,x_1)=\sum(-1)^{\norm{T}\cdot\dagger}\phi^{s+t+1}(\mathcal{G}^{j_s}(x_k,\cdots,x_{k-j_s+1}),\cdots,T^l(x_q,\cdots,x_{i_1+\cdots+i_t+1}),\cdots,\mathcal{F}^{i_1}(x_{i_1},\cdots,x_1))
    \end{equation}
    for the sign
    \begin{equation}
        \dagger=\norm{x_{i_1+\cdots+i_t}}+\cdots+\norm{x_1},
    \end{equation}
    and
    \begin{equation}
        (\rho_T\phi)^k(x_k,\cdots,x_1)=\sum(-1)^{\norm{\phi}\cdot\dagger_i}T^{k-j+1}(x_k,\cdots,\phi^j(x_{i+j},\cdots,x_{i+1}),x_i,\cdots,x_1).
    \end{equation}
\end{defn}

For example, the length 0 and 1 parts of the above operations are
\begin{align}
    (\iota_T\phi)^0_L&=\phi^1(T^0_L),\\
    (\iota_T\phi)^1(x)&=\phi^1(T^1(x))+\phi^2(\mathcal{G}^1(x),T^0_K)+(-1)^{\norm{T}\cdot\norm{x}}\phi^2(T^0_L,\mathcal{F}^1(x)),\\
    (\rho_T\phi)^0_L&=T^1(\phi^0_L),\\
    (\rho_T\phi)^1(x)&= T^1(\phi^1(x))+T^2(x,\phi^0_K)+(-1)^{\norm{\phi}\cdot\norm{x}}T^2(\phi^0_L,x)
\end{align}
for $x\in\hom_\mathcal{C}(K,L)$.

With this definition, we can define the grading operator.
\begin{defn}
     Let $(\mathcal{C},\eta_\mathcal{C})$, $(\mathcal{D},\eta_\mathcal{D})$ be infinitesimally graded $A_\infty$-categories, and let $(\mathcal{F},\gamma_\mathcal{F})$, $(\mathcal{G},\gamma_\mathcal{G})$ be infinitesimally graded functors from $(\mathcal{C},\eta_\mathcal{C})$ to $(\mathcal{D},\eta_\mathcal{D})$. Then for any $A_\infty$-pre-natural transformation $T\in CC^*(\mathcal{C},\prescript{~}{\mathcal{F}}{\Delta}_{\mathcal{G}})$, we define the \emph{grading operator}
     \begin{equation}
         \Gamma(T)=\iota_T(\eta_\mathcal{D})-\rho_T(\eta_\mathcal{C})+\Lambda(T)-\mu^2_{fun}(T,\gamma_\mathcal{F})-(-1)^{\norm{T}}\mu^2_{fun}(\gamma_\mathcal{G},T),
     \end{equation}
     for the length operator $\Lambda$ defined such that $\Lambda(T^d)=d\cdot T^d$.
\end{defn}

We now prove that the above grading operator commutes with the differential in the expected way.

\begin{proposition}
    For the objects $\mathcal{C},\mathcal{D},\mathcal{F},\mathcal{G},T$ defined above,
    \begin{equation}
        \Gamma\circ\mu^1_{fun}-\mu^1_{fun}\circ\Gamma=\mu^1_{fun}.
    \end{equation}
\end{proposition}
\begin{proof}
    The key identity to prove is
    \begin{align}
        &\iota_{\mu^1T}(\eta_\mathcal{D})-\rho_{\mu^1T}(\eta_\mathcal{C})+\Lambda(\mu^1T)-\mu^1_{fun}(\iota_T\eta_\mathcal{D}-\rho_T\eta_\mathcal{C}+\Lambda(T))\\
        =&\mu^1_{fun}T-\mu^2_{fun}(\mathcal{G}_*\eta_\mathcal{C}-\mathcal{G}^*\eta_\mathcal{D}+H_\mathcal{G},T)-\mu^2_{fun}(T,\mathcal{F}_*\eta_\mathcal{C}-\mathcal{F}^*\eta_\mathcal{D}+H_\mathcal{F}).
    \end{align}
    We show three independent identities, which add up to the above identity:
    \begin{align}
        \iota_{\mu^1T}(\eta_\mathcal{D})-\mu^1_{fun}(\iota_T\eta_\mathcal{D})&=\iota_T(\kappa_\mathcal{D})+\mu^2_{fun}(\mathcal{G}^*\eta_\mathcal{D},T)+\mu^2_{fun}(T,\mathcal{F}^*\eta_\mathcal{D}), \label{eqn:App-3}\\ 
        \rho_{\mu^1T}(\eta_\mathcal{C})-\mu^1_{fun}(\rho_T\eta_\mathcal{C})&=(-1)^{\norm{T}}\rho_T(\kappa_\mathcal{C})+\mu^2_{fun}(\mathcal{G}_*\eta_\mathcal{C},T)+\mu^2_{fun}(T,\mathcal{F}_*\eta_\mathcal{C}), \label{eqn:App-4}\\
        \Lambda(\mu^1T)-\mu^1_{fun}(\Lambda T)&=\mu^1T-\iota_T(\kappa_\mathcal{D})+(-1)^{\norm{T}}\rho_T(\kappa_\mathcal{C})-\mu^2_{fun}(H_\mathcal{G},T)-\mu^2_{fun}(T,H_\mathcal{F}). \label{eqn:App-5}
    \end{align}
    The proofs of the above equalities appear in the appendix as lemmas \ref{lem:App-3}, \ref{lem:App-4}, and \ref{lem:App-5}.
\end{proof}

We are ready to define graded isomorphisms between infinitesimally graded functors.

\begin{defn}
    Let $(\mathcal{C},\eta_\mathcal{C})$, $(\mathcal{D},\eta_\mathcal{D})$ be infinitesimally graded $A_\infty$-categories, and let $(\mathcal{F},\gamma_\mathcal{F})$, $(\mathcal{G},\gamma_\mathcal{G})$ be infinitesimally graded functors from $(\mathcal{C},\eta_\mathcal{C})$ to $(\mathcal{D},\eta_\mathcal{D})$. Then $\mathcal{F}$, $\mathcal{G}$ are \emph{graded isomorphic} if there exists an invertible degree zero natural transformation $T:\mathcal{F}\to \mathcal{G}$ such that
    \begin{equation}
        \iota_T(\eta_\mathcal{D})-\rho_T(\eta_\mathcal{C})+\Lambda(T)-\mu^2_{fun}(T,\gamma_\mathcal{F})+\mu^2_{fun}(\gamma_\mathcal{G},T)
    \end{equation}
    is a Hochschild coboundary. We say that the natural transformation $T$ has \emph{infinitesimal degree $\lambda$} if
    \begin{equation}
        [\iota_T(\eta_\mathcal{D})-\rho_T(\eta_\mathcal{C})+\Lambda(T)-\mu^2_{fun}(T,\gamma_\mathcal{F})+\mu^2_{fun}(\gamma_\mathcal{G},T)]=\lambda\cdot[T]
    \end{equation}
    holds for some $\lambda\in k$.
\end{defn}

We unpack the above definition at length zero. Evaluating at an object $L\in\mathcal{C}$,
\begin{equation} \label{eqn:grading-agree}
    \eta^1(T^0_L)-T^1(\eta^0_L)-\mu^2(T^0_L,\gamma^0_\mathcal{F}\vert_{\mathcal{F}(L)})+\mu^2(\gamma^0_\mathcal{G}\vert_{\mathcal{G}(L)},T^0_L)=\mu^1(H_L)
\end{equation}
for some element $H_L$.

The following proposition shows that graded isomorphic functors induce the same grading operators on cohomology-level morphism spaces.

\begin{proposition} \label{prop:isomorphic-functors-same-grading}
    Let $(\mathcal{F},\gamma_\mathcal{F})$ and $(\mathcal{G},\gamma_\mathcal{G})$ be infinitesimally graded endofunctors that are graded isomorphic. Then for infinitesimally graded objects $(K,\xi_K)$, $(L,\xi_L)$, the following diagram commutes:
    \begin{equation}
    \begin{tikzcd}
        H^*\hom_\mathcal{C}((K,\xi_K),(\mathcal{F}(L),\xi^{\mathcal{F}}_{\mathcal{F}(L)})) \arrow[rr, "{\mu^2(T^0_L,{-})}"] \arrow[d,"\Gamma"'] &&H^*\hom_\mathcal{C}((K,\xi_K),(\mathcal{G}(L),\xi^{\mathcal{G}}_{\mathcal{G}(L)})) \arrow[d,"\Gamma"]\\
        H^*\hom_\mathcal{C}((K,\xi_K),(\mathcal{F}(L),\xi^{\mathcal{F}}_{\mathcal{F}(L)})) \arrow[rr, "{\mu^2(T^0_L,{-})}"] &&H^*\hom_\mathcal{C}((K,\xi_K),(\mathcal{G}(L),\xi^{\mathcal{G}}_{\mathcal{G}(L)})) 
    \end{tikzcd}.    
    \end{equation}
\end{proposition}
\begin{proof}
    We need to compare the two compositions for a cocycle $x\in\hom_\mathcal{C}(K,\mathcal{F}(L))$
    \begin{equation}  \label{eqn:grading-agree-0}
        \mu^2(T^0_L,\eta^1(x)-\mu^2(x,\xi_K)-(-1)^{\norm{x}}\mu^2(\xi^{\mathcal{F}}_{\mathcal{F}(L)},x))-\eta^1(\mu^2(T^0_L,x))+\mu^2(\mu^2(T^0_L,x),\xi_K)+(-1)^{\norm{x}}\mu^2(\xi^{\mathcal{G}}_{\mathcal{G}(L)},\mu^2(T^0_L,x))
    \end{equation}
     at chain level and show they agree up to a coboundary.

    The $A_\infty$-relation
    \begin{equation}
        \mu^2(T^0_L,\mu^2(x,\xi_K))-\mu^2(\mu^2(T^0_L,x),\xi_K)+\mu^3(T^0_L,x,\eta^0_K)-\mu^3(T^0_L,\mu^1(x),\xi_K)+\mu^1(\mu^3(T^0_L,x,\xi_K))=0
    \end{equation}
    shows that up to coboundary
    \begin{equation} \label{eqn:grading-agree-1}
         -\mu^2(T^0_L,\mu^2(x,\xi_K))+\mu^2(\mu^2(T^0_L,x),\xi_K)\simeq\mu^3(T^0_L,x,\eta^0_K).
    \end{equation}
    From $(\delta\eta)^2=\kappa^2=0$, we find that 
    \begin{align}
        0=(\delta\eta)^2(T^0_L,x)=&-\mu^1(\eta^2(T^0_L,x))-\mu^2(T^0_L,\eta^1(x))-\mu^2(\eta^1(T^0_L),x)-\mu^3(T^0_L,x,\eta^0_K)\\&-
        \mu^3(T^0_L,\eta^0_{\mathcal{F}(L)},x)-\mu^3(\eta^0_{\mathcal{G}(L)},T^0_L,x)+\eta^1(\mu^2(T^0_L,x))+\eta^2(T^0_L,\mu^1(x))
    \end{align}
    since $T^0_L$ is a cocycle. Thus up to coboundary,
    \begin{equation} \label{eqn:grading-agree-2}
        \mu^2(T^0_L,\eta^1(x))-\eta^1(\mu^2(T^0_L,x))\simeq -\mu^2(\eta^1(T^0_L),x)-\mu^3(T^0_L,x,\eta^0_K)\\-
        \mu^3(T^0_L,\eta^0_{\mathcal{F}(L)},x)-\mu^3(\eta^0_{\mathcal{G}(L)},T^0_L,x).
    \end{equation}
    Finally, the $A_\infty$-relations show
    \begin{align}
        \mu^2(T^0_L,\mu^2(\xi^{\mathcal{F}}_{\mathcal{F}(L)},x))+(-1)^{\norm{x}}\mu^2(\mu^2(T^0_L,\xi^{\mathcal{F}}_{\mathcal{F}(L)}),x)+\mu^1(\mu^3(T^0_L,\xi^{\mathcal{F}}_{\mathcal{F}(L)},x))&\\+\mu^3(T^0_L,\xi^{\mathcal{F}}_{\mathcal{F}(L)},\mu^1(x))+(-1)^{\norm{x}}\mu^3(T^0_L,\eta^0_{\mathcal{F}(L)},x)=0
    \end{align}
    and similarly for $\mathcal{G}$. Thus up to coboundary,
    \begin{align} \label{eqn:grading-agree-3}
        -(-1)^{\norm{x}}\mu^2(T^0_L,\mu^2(\xi^{\mathcal{F}}_{\mathcal{F}(L)},x))+(-1)^{\norm{x}}\mu^2(\xi^{\mathcal{G}}_{\mathcal{G}(L)},\mu^2(T^0_L,x))&\simeq \mu^2(\mu^2(T^0_L,\xi^{\mathcal{F}}_{\mathcal{F}(L)}),x)+\mu^3(T^0_L,\eta^0_{\mathcal{F}(L)},x)\\
        -&\mu^2(\mu^2(\xi^{\mathcal{G}}_{\mathcal{G}(L)},T^0_L),x)+\mu^3(\eta^0_{\mathcal{G}(L)},T^0_L,x).
    \end{align}
    
    Now adding up (\ref{eqn:grading-agree-1}), (\ref{eqn:grading-agree-2}), (\ref{eqn:grading-agree-3}), we conclude that (\ref{eqn:grading-agree-0}) is equal up to coboundary to
    \begin{equation}
        -\mu^2(\eta^1(T^0_L)-\mu^2(T^0_L,\xi^{\mathcal{F}}_{\mathcal{F}(L)})+\mu^2(\xi^{\mathcal{G}}_{\mathcal{G}(L)},T^0_L),x).
    \end{equation}

    Now the $A_\infty$-natural transformation relations show that
    \begin{equation}
        \mu^1(T^1(\xi_L))+T^1(\mu^1(\xi_L))+\mu^2(\mathcal{G}^1(\xi_L),T^0_L)-\mu^2(T^0_L,\mathcal{F}^1(\xi_L))=0,
    \end{equation}
    and (\ref{eqn:grading-agree}) shows that
    \begin{equation}
        \eta^1(T^0_L)-T^1(\eta^0_L)-\mu^2(T^0_L,\gamma^0_\mathcal{F}\vert_{\mathcal{F}(L)})+\mu^2(\gamma^0_\mathcal{G}\vert_{\mathcal{G}(L)},T^0_L)
    \end{equation}
    is a coboundary, so adding these two equations shows that
    \begin{equation}
        \eta^1(T^0_L)-\mu^2(T^0_L,\xi^{\mathcal{F}}_{\mathcal{F}(L)})+\mu^2(\xi^{\mathcal{G}}_{\mathcal{G}(L)},T^0_L)\simeq T^1(\eta^0_L)-\mu^2(T^0_L,\mathcal{F}^1(\xi_L))+\mu^2(\mathcal{G}^1(\xi_L),T^0_L)\simeq-\mu^1(T^1(\xi_L))
    \end{equation}
    is a coboundary, which concludes the proof.
\end{proof}

Consequently, $\mathcal{F}$ and $\mathcal{G}$ cannot be graded isomorphic if, for two infinitesimally graded objects $K,L$, the grading operators on $H^*\hom(K,\mathcal{F}(L))$ and $H^*\hom(K,\mathcal{G}(L))$ admit no intertwining linear isomorphism.

We end this section by showing that compositions of infinitesimally graded functors are again infinitesimally graded.

\begin{lemma}\label{lem:composition-of-functors}
    Let $(\mathcal{F},\gamma_\mathcal{F}):(\mathcal{C},\eta_\mathcal{C})\to(\mathcal{D},\eta_\mathcal{D})$ and $(\mathcal{G},\gamma_\mathcal{G}):(\mathcal{D},\eta_\mathcal{D})\to(\mathcal{E},\eta_\mathcal{E})$ be infinitesimally graded functors. Then the composition $\mathcal{G}\circ\mathcal{F}$ inherits an infinitesimal grading, given by
    \begin{equation}
        \gamma_{\mathcal{G}\circ\mathcal{F}}=\mathcal{G}_\sharp(\gamma_\mathcal{F})+\mathcal{F}^\sharp(\gamma_\mathcal{G})
    \end{equation}
    for the chain maps $\mathcal{G}_\sharp:CC^*(\mathcal{C},\Delta_{\mathcal{F}})\to CC^*(\mathcal{C},\Delta_{\mathcal{G}\circ\mathcal{F}})$ given by
    \begin{equation}
        (\mathcal{G}_\sharp\alpha)^d(x_d,\cdots,x_1)=\sum(-1)^{\norm{\alpha}\cdot\dagger_m}\mathcal{G}^{t+1}(\mathcal{F}^{i_t}(\cdots),\cdots,\alpha^l(\cdots),\mathcal{F}^{i_s}(x_m,\cdots),\cdots,\mathcal{F}^{i_1}(\cdots)),
    \end{equation}
    and $\mathcal{F}^\sharp:CC^*(\mathcal{D},\Delta_{\mathcal{G}})\to CC^*(\mathcal{C},\Delta_{\mathcal{G}\circ\mathcal{F}})$ given by
    \begin{equation}
        (\mathcal{F}^\sharp\beta)^d(y_d,\cdots,y_1)=\sum\beta^r(\mathcal{F}^{j_r}(y_d,\cdots,y_{d-j_r+1}),\cdots,\mathcal{F}^{j_1}(y_{j_1},\cdots,y_1)).
    \end{equation}
    In particular, the length zero part for an object $L\in\mathcal{C}$ can be computed as
    \begin{equation}
        \gamma^0_{\mathcal{G}\circ\mathcal{F}}\vert_L=\mathcal{G}^1(\gamma^0_\mathcal{F}\vert_L)+\gamma^0_\mathcal{G}\vert_{\mathcal{F}(L)}.
    \end{equation}
\end{lemma}
\begin{proof}
    To show that $\gamma_{\mathcal{G}\circ\mathcal{F}}$ is an infinitesimal grading, we check
    \begin{equation} \label{eqn:inf-grading-functor-composition-identity}
        (\mathcal{G}\circ\mathcal{F})_*\eta_{\mathcal{C}}-(\mathcal{G}\circ\mathcal{F})^*\eta_{\mathcal{E}}+H_{\mathcal{G}\circ\mathcal{F}}=\delta\gamma_{\mathcal{G}\circ\mathcal{F}}.
    \end{equation}
    The left-hand side is equal to
    \begin{equation}
        \mathcal{G}_\sharp(\mathcal{F}_*\eta_\mathcal{C})-\mathcal{F}^\sharp(\mathcal{G}^*\eta_\mathcal{E})+H_{\mathcal{G}\circ\mathcal{F}}.
    \end{equation}
    We claim the following two identities:
    \begin{align}
        &H_{\mathcal{G}\circ\mathcal{F}}=\mathcal{G}_\sharp(H_\mathcal{F})+\mathcal{F}^\sharp(H_\mathcal{G}), \label{eqn:inf-grading-functor-composition-1} \\
        &\mathcal{G}_\sharp(\mathcal{F}^*\eta_\mathcal{D})=\mathcal{F}^\sharp(\mathcal{G}_*\eta_\mathcal{D}). \label{eqn:inf-grading-functor-composition-2}
    \end{align}
    Once these identities are established, (\ref{eqn:inf-grading-functor-composition-identity}) follows from
    \begin{align}
        \mathcal{G}_\sharp(\mathcal{F}_*\eta_\mathcal{C})-\mathcal{F}^\sharp(\mathcal{G}^*\eta_\mathcal{E})+H_{\mathcal{G}\circ\mathcal{F}}&=  \mathcal{G}_\sharp(\mathcal{F}_*\eta_\mathcal{C})-\mathcal{F}^\sharp(\mathcal{G}^*\eta_\mathcal{E})+\mathcal{G}_\sharp(H_\mathcal{F})+\mathcal{F}^\sharp(H_\mathcal{G})\\
        &=\mathcal{G}_\sharp(\mathcal{F}_*\eta_\mathcal{C}-\mathcal{F}^*\eta_\mathcal{D}+H_\mathcal{F})+\mathcal{F}^\sharp(\mathcal{G}_*\eta_\mathcal{D}-\mathcal{G}^*\eta_\mathcal{E}+H_\mathcal{G})\\
        &=\mathcal{G}_\sharp(\delta\gamma_\mathcal{F})+\mathcal{F}^\sharp(\delta\gamma_\mathcal{G}) = \delta(\gamma_{\mathcal{G}\circ\mathcal{F}}).
    \end{align}
    We first check the second identity (\ref{eqn:inf-grading-functor-composition-2}). By definition,
    \begin{align}
        \mathcal{F}^\sharp(\mathcal{G}_*\eta_\mathcal{D})^k(x_k,\cdots,x_1) &=\sum(\mathcal{G}_*\eta_\mathcal{D})^r(\mathcal{F}^{i_r}(x_k,\cdots,x_{k-i_r+1}),\cdots,\mathcal{F}^{i_1}(x_{i_1},\cdots,x_1))\\
        &=\sum\mathcal{G}^{r-s+1}(\mathcal{F}^{i_r}(\cdots),\cdots,\eta_\mathcal{D}^s(\mathcal{F}^{i_t}(\cdots),\cdots,\mathcal{F}^{i_{t-s+1}}(\cdots)),\cdots,\mathcal{F}^{i_1}(\cdots))\\
        &=\mathcal{G}_\sharp(\mathcal{F}^*\eta_\mathcal{D})^k(x_k,\cdots,x_1).
    \end{align}
    For the first identity (\ref{eqn:inf-grading-functor-composition-1}), recall that
    \begin{equation}
        H_{\mathcal{G}\circ\mathcal{F}}^d=(1-d)(\mathcal{G}\circ\mathcal{F})^d=(1-d)\sum\mathcal{G}^r(\mathcal{F}^{i_r}(\cdots),\cdots,\mathcal{F}^{i_1}(\cdots)).
    \end{equation}
    The corresponding term in $\mathcal{G}_\sharp(H_\mathcal{F})$ contributes
    \begin{equation}
        \sum(\sum(1-i_t)\cdot\mathcal{G}^r(\mathcal{F}^{i_r}(\cdots),\cdots,\mathcal{F}^{i_1}(\cdots))),
    \end{equation}
    and $\mathcal{F}^\sharp(H_\mathcal{G})$ contributes
    \begin{equation}
        \sum(1-r)\cdot\mathcal{G}^r(\mathcal{F}^{i_r}(\cdots),\cdots,\mathcal{F}^{i_1}(\cdots)),
    \end{equation}
    and since $\sum i_t=d$ the contributions add up correctly.
\end{proof}

\section{Infinitesimal gradings for twisted complexes}

\subsection{Twisted complexes}

We now extend our notion of infinitesimal gradings to twisted complexes. We first recall the definition of a twisted complex.

\begin{defn} \label{def:twisted-cpx}
    A \emph{twisted complex} $\mathcal{L}=(I,\{L_i\},\{V_i\},\delta_\mathcal{L})$ in an $A_\infty$-category $\mathcal{C}$ consists of a finite ordered index set $I$, a collection of objects $L_i\in\mathcal{C}$ indexed by $i\in I$, $\mathbb{Z}/2$-graded $k$-vector spaces $V_i$ of finite rank indexed by $i\in I$, and a collection of morphisms $\delta_\mathcal{L}=\{\delta_{ij}\in\hom^1(V_i\otimes L_i,V_j\otimes L_j)\}$ such that $\delta_{ij}=0$ if $i\geq j$ and $\delta_\mathcal{L}$ satisfies the following Maurer-Cartan equation:
    \begin{equation}
        \sum_{d\geq1} \mu^d(\delta_{\mathcal{L}},\cdots,\delta_{\mathcal{L}},\delta_{\mathcal{L}})=0,
    \end{equation}
    which is shorthand for the equation
    \begin{equation}
        \sum_{d\geq1}\sum_{p=i_0<\cdots<i_d=q}\mu^d(\delta_{i_{d-1}i_d},\delta_{i_{d-2}i_{d-1}},\cdots,\delta_{i_0i_1})=0
    \end{equation}
    for all $(p,q)$. Because of the lower triangular condition for $\delta_\mathcal{L}$, the above sum is finite.
\end{defn}

Twisted complexes form an $A_\infty$-category $\Tw\mathcal{C}$ with $A_\infty$-operations $\mu^d_{Tw}$, and $\Tw\mathcal{C}$ is cohomologically unital if $\mathcal{C}$ is. Moreover, we can explicitly construct a map that extends Hochschild cochains in $\mathcal{C}$ to $\Tw\mathcal{C}$:

\begin{defn} \label{def:cc-extends-Tw}
    For a Hochschild cochain $\beta$, we define the Hochschild cochain $\Tw\beta\in CC^*(\Tw\mathcal{C},\Tw\mathcal{C})$ as
    \begin{equation}
       (\Tw\beta)^k(x_k,\cdots,x_1)=\sum_{j_0,\cdots,j_k}\beta^{k+j_0+\cdots+j_k}(\underbrace{\delta,\cdots,\delta}_{j_k},x_k,\cdots,x_1,\underbrace{\delta,\cdots,\delta}_{j_0}).
    \end{equation}
    Then $\Tw:CC^*(\mathcal{C},\mathcal{C})\to CC^*(\Tw\mathcal{C},\Tw\mathcal{C})$ is a chain map, and is a quasi-isomorphism for c-unital $\mathcal{C}$.
\end{defn}

This definition is somewhat misleading, in that the $x_k$ in the above definition are actually morphisms between twisted complexes. If we decompose the $x_k$ as direct sums of elements $\phi_k\otimes z_k$, for $\phi_k\in\hom(V_{k-1},V_k)$ and $z_k\in\hom_\mathcal{C}(L_{k-1},L_k)$, the formula for the additive extension is
\begin{equation}
    (\Sigma\beta)^k(\phi_k\otimes z_k,\cdots,\phi_1\otimes z_1)=(-1)^{\bowtie} \phi_k\circ\cdots\circ\phi_1\otimes\beta^k(z_k,\cdots,z_1)
\end{equation}
for the Koszul sign $\bowtie=\sum_{i<j}\lvert\phi_i\rvert\norm{z_j}$, and $\Tw\beta$ is obtained by taking a sum over the above formula, where Maurer-Cartan elements can be inserted in all possible places.

For this section, we will use the symbol $\boldsymbol{\delta}$ for the Hochschild cochain, and $-[\mu_{Tw},~]$ for the codifferential of the Hochschild cochain complex.

\begin{defn} \label{defn-bold-delta}
    For $\Tw\mathcal{C}$, define the Hochschild cochain $\boldsymbol{\delta}$ to have only a nontrivial length-$0$ part defined as
    \begin{equation}
        \boldsymbol{\delta}^0\vert_\mathcal{L}=\delta_\mathcal{L}\in\hom_{\Tw\mathcal{C}}(\mathcal{L},\mathcal{L}).
    \end{equation}
    Then in $CC^*(\Tw\mathcal{C},\Tw\mathcal{C})$, we compute
    \begin{equation}
        -[\mu_{Tw},\boldsymbol{\delta}]^k(x_k,\cdots,x_1)=\sum_{j_0,\cdots,j_k}-(j_0+\cdots+j_k)\mu^{k+j_0+\cdots+j_k}(\underbrace{\delta,\cdots,\delta}_{j_k},x_k,\cdots,x_1,\underbrace{\delta,\cdots,\delta}_{j_0}),
    \end{equation}
    since the Hochschild cochain $\delta$ can be chosen to occur in $(j_0+\cdots+j_k)$ places, where the other inputs of $\boldsymbol{\delta}$ are coming from the differential of twisted complexes.
\end{defn}

The following proposition shows that an infinitesimal grading for $\mathcal{C}$ extends to an infinitesimal grading for $\Tw\mathcal{C}$.

\begin{proposition} \label{prop:infinitesimal-grading-Tw}
    For an infinitesimally graded $A_\infty$-category $(\mathcal{C},\eta)$, the $A_\infty$-category of twisted complexes $\Tw\mathcal{C}$ inherits a canonical infinitesimal grading $\eta_{Tw}$ defined as
    \begin{equation}
        \eta_{Tw}=\Tw\eta-\boldsymbol{\delta}.
    \end{equation}
\end{proposition}
\begin{proof}
    We first compute that
    \begin{equation}
        \kappa_{Tw}^k(x_k,\cdots,x_1)=\sum_{j_0,\cdots,j_k}(2-k)\mu^{k+j_0+\cdots+j_k}(\underbrace{\delta,\cdots,\delta}_{j_k},x_k,\cdots,x_1,\underbrace{\delta,\cdots,\delta}_{j_0}),
    \end{equation}
    while
    \begin{equation}
        -[\mu_{Tw},\Tw\eta]=\Tw(\kappa)=\sum_{j_0,\cdots,j_k}(2-k-j_0-\cdots-j_k)\mu^{k+j_0+\cdots+j_k}(\underbrace{\boldsymbol{\delta},\cdots,\boldsymbol{\delta}}_{j_k},x_k,\cdots,x_1,\underbrace{\boldsymbol{\delta},\cdots,\boldsymbol{\delta}}_{j_0}).
    \end{equation}
    We conclude that
    \begin{equation}
        \kappa_{\Tw\mathcal{C}}=\Tw\kappa+[\mu_{Tw},\boldsymbol{\delta}],
    \end{equation}
    which shows that our choice of $\eta_{Tw}$ defines an infinitesimal grading.
\end{proof}

We now define infinitesimally graded twisted complexes.

\begin{defn}
    An \emph{infinitesimally graded twisted complex} $(\mathcal{L},\xi_{\mathcal{L}})$ for the infinitesimal grading $\eta_{Tw}$ is a pair of a twisted complex $\mathcal{L}\in\Tw\mathcal{C}$, together with an element $\xi_\mathcal{L}\in\hom^\mathrm{even}_{\Tw\mathcal{C}}(\mathcal{L},\mathcal{L})$ such that
    \begin{equation}
        \mu^1_{Tw}(\xi_{\mathcal{L}})=\eta^0_{Tw}\vert_{\mathcal{L}}.
    \end{equation}
\end{defn}

For example, infinitesimally graded objects $(L,\xi_L)$ of $(\mathcal{C},\eta)$ are infinitesimally graded objects in $(\Tw\mathcal{C},\eta_{Tw})$.

To understand the above definition better, we look at shifts $K[1]$ for objects $K\in\mathcal{C}$, which are defined as $V\otimes K\in\Tw\mathcal{C}$ for a rank 1 $k$-vector space $V$ of degree $1$ with chosen generator $v$. 

\begin{defn} \label{def:mapping-cone}
    For an even cocycle $\alpha\in\hom_\mathcal{C}(K,L)$, the \emph{mapping cone of $\alpha$} is defined to be the twisted complex
    \begin{equation}
        Cone(\alpha)=(V\otimes K)\oplus L,\quad \delta_{Cone(\alpha)}= \begin{pmatrix}
            0 &0\\ s\otimes\alpha &0
        \end{pmatrix}
    \end{equation}
    for the degree $1$ vector space $V$ of rank $1$, and the dual $s=v^\vee$.
\end{defn}

We claim that mapping cones of pure weight morphisms define infinitesimally graded objects. For a $\mathbb{Z}$-graded $A_\infty$-category, the mapping cone is genuinely graded after shifting the source summand so that the mapping cone differential has degree 1. The following proposition is the infinitesimal analogue of this construction.

\begin{proposition}
    Let $(K,\xi_{K})$, $(L,\xi_L)$ be infinitesimally graded objects in $(\mathcal{C},\eta)$, and assume that there exists a homotopy unit $e_K\in\hom_{\mathcal{C}}(K,K)$. Suppose that a cocycle $\alpha\in\hom^{even}_{\mathcal{C}}(K,L)$ satisfies
    \begin{equation}
        \Gamma_{K,L}(\alpha)=\lambda\alpha+\mu^1(\beta)
    \end{equation}
    for some $\beta\in\hom_{\mathcal{C}}(K,L)$, and let the homotopy unit $e_K$ satisfy
    \begin{equation}
        \mu^2(\alpha,e_K)-\alpha=\mu^1(h)
    \end{equation}
    for some $h\in\hom_\mathcal{C}(K,L)$. Then the mapping cone $Cone(\alpha)$, defined as the twisted complex $K[1]\oplus L$ with internal differential $\alpha$, admits an infinitesimal grading defined as
    \begin{equation}
        \xi_{Cone(\alpha)}=\begin{pmatrix}
            \id_V\otimes(\xi_K+(\lambda-1)e_K) &0\\
            s\otimes(\beta-(\lambda-1)h) &\xi_L
        \end{pmatrix}.
    \end{equation}
    Since the differential $K[1]\to L$ must have weight $1$, a morphism of weight $\lambda$ is corrected by shifting the grading of the $K[1]$-summand by $\lambda-1$.
\end{proposition}
\begin{proof}
    We need to show that $\mu^1_{Tw}(\xi_{Cone(\alpha)})=\eta^0_{Tw}\vert_{Cone(\alpha)}$. First, we compute for $\eta^0_{Tw}=(\Tw\eta)^0-\boldsymbol{\delta}$
    \begin{equation} \label{eqn:mapping-cone-1}
        (\Tw\eta)^0=\begin{pmatrix}
            \id_V\otimes\eta^0_{K} &0\\
            s\otimes\eta^1(\alpha) &\eta^0_L
        \end{pmatrix},\quad \delta_{Cone(\alpha)}=\begin{pmatrix}
            0 &0 \\ s\otimes\alpha &0
        \end{pmatrix},
    \end{equation}
    and 
    \begin{equation} \label{eqn:mapping-cone-2}
        \mu^1_{Tw}(\xi_{Cone(\alpha)})=\begin{pmatrix}
            \id_V\otimes\mu^1(\xi_{K}) &0\\
            s\otimes\left(\mu^1(\beta-(\lambda-1)h)-\mu^2(\xi_L,\alpha)+\mu^2(\alpha,\xi_{K}+(\lambda-1)e_K)\right) &\mu^1(\xi_L)
        \end{pmatrix}
    \end{equation}
    since the signs work out as
    \begin{align}
        \mu_{Tw}^2(s\otimes\alpha,\id_V\otimes\xi_K)&=s\otimes\mu^2(\alpha,\xi_K)\\
        \mu^2_{Tw}(\id\otimes\xi_L,s\otimes\alpha)&=(-1)^{\lvert s\rvert\cdot\norm{\xi_L}}s\otimes\mu^2(\xi_L,\alpha)=-s\otimes\mu^2(\xi_L,\alpha).
    \end{align}

    Now since
    \begin{equation}
        \lambda\alpha+\mu^1(\beta)=\eta^1(\alpha)-\mu^2(\alpha,\xi_K)+\mu^2(\xi_L,\alpha),
    \end{equation}
    we conclude that the difference of the two terms in (\ref{eqn:mapping-cone-1}) agrees with (\ref{eqn:mapping-cone-2}).
\end{proof}

\subsection{Twisted complexes with cochain complex coefficients}

The previous discussion of infinitesimal gradings for twisted complexes is incomplete in that it does not incorporate the grading of the vector spaces $V_i$ appearing in each twisted complex. To remedy this, we will first introduce twisted complexes with cochain complex coefficients, as in \cite[Definition 3.6]{Kim}.

\subsubsection{The category $\Tw_{\Ch}\mathcal{C}$ and the object $T_V\mathcal{L}$}

We first define the additive enlargement $\Sigma_{\Ch}\mathcal{C}$.

\begin{defn}
    Given an $A_\infty$-category $\mathcal{C}$, the \emph{additive enlargement} $\Sigma_{\Ch}\mathcal{C}$ is defined to have objects
    \begin{equation}
        Ob(\Sigma_{\Ch}\mathcal{C})=\{\bigoplus_{i\in I}C^*_i\otimes L_i\}
    \end{equation}
    where each $L_i\in\mathcal{C}$, $C_i$ is a $\mathbb{Z}/2$-graded cochain complex of finitely generated $k$-modules, and $I$ is a finite index set. Morphisms are defined by
    \begin{equation}
        \hom_{\Sigma_{\Ch} \mathcal{C}} (\bigoplus_{i\in I}C^*_i\otimes L_i,\bigoplus_{j\in J}D^*_j\otimes K_j)=\bigoplus_{i,j}\hom_k(C_i^*,D_j^*)\otimes\hom_\mathcal{C}(L_i,K_j),
    \end{equation}
    where $\hom_k(C_i^*,D_j^*)$ denotes the morphism spaces between these two cochain complexes in the dg category of cochain complexes over $k$. Then $\Sigma_{\Ch}\mathcal{C}$ admits an $A_\infty$-structure induced from $\mathcal{C}$. Explicitly,
    \begin{align}
        \mu^1_{\Sigma_{\Ch}\mathcal{C}}(\phi\otimes x)&= (-1)^{\lvert x\rvert+\lvert\phi\rvert}\partial_{Ch}(\phi)\otimes x + \phi\otimes \mu^1_{\mathcal{C}}(x),\\
        \mu^k_{\Sigma_{\Ch}\mathcal{C}}(\phi_k\otimes x_k,\cdots,\phi_1\otimes x_1)&= (-1)^{\bowtie} (\phi_k\circ\cdots\circ\phi_1) \otimes \mu^k_\mathcal{C}(x_k,\cdots,x_1),\quad (k\geq2)
    \end{align}
    for the sign $\bowtie=\sum_{i<j}\lvert\phi_i\rvert\cdot\norm{x_j}$. 
\end{defn}

A primitive morphism $\phi\otimes x: C^*\otimes K\to D^*\otimes L$ appearing in the above definition consists of a linear map $\phi\in\hom_k(C^*,D^*)$ and an element $x\in\hom_\mathcal{C}(K,L)$; the tensor $\phi\otimes x$ has total degree $\lvert\phi\rvert+\lvert x\rvert$. The sign $\bowtie$ appears when one commutes $\phi_i$ across $x_j$: the $x_j$ are in reduced degrees. For example,
\begin{equation}
    \mu^2(\phi_2\otimes x_2,\phi_1\otimes x_1)=(-1)^{\lvert\phi_1\rvert\norm{x_2}} \phi_2\circ\phi_1\otimes\mu^2(x_2,x_1).
\end{equation}

Now we introduce internal differentials. Define the category $\Xi\mathcal{C}$ to be the category with objects
\begin{equation}
    Ob(\Xi\mathcal{C})=\{(\mathcal{L}=\bigoplus_{i\in I}C^*_i\otimes L_i,\delta_\mathcal{L})\},
\end{equation}
where $\mathcal{L}\in\Sigma_{Ch}\mathcal{C}$, and $\delta\in\hom^1_{\Sigma_{\Ch}\mathcal{C}}(\mathcal{L},\mathcal{L})$, such that each object carries a decreasing filtration of finite length
\begin{equation}
    \mathcal{L}=F^0\mathcal{L}\supset F^1\mathcal{L} \supset\cdots\supset F^n\mathcal{L}=0
\end{equation}
for which $\delta_\mathcal{L}$ strictly increases the filtration: $\delta_\mathcal{L}(\mathcal{F}^i\mathcal{L})\subset\mathcal{F}^{i+1}\mathcal{L}$. The morphism spaces are inherited from $\Sigma_{\Ch}\mathcal{C}$:
\begin{equation}
    \hom_{\Xi\mathcal{C}}(\mathcal{L},\mathcal{K})=\hom_{\Sigma_{\Ch}\mathcal{C}}(\mathcal{L},\mathcal{K}).
\end{equation}

Then $\Xi\mathcal{C}$ is a curved $A_\infty$-category, with structure maps defined by inserting $\delta$ as for twisted complexes:
\begin{equation} \label{eqn:Tw_Ch-Ainfty}
    \mu^k_{\Xi\mathcal{C}}(x_k,\cdots,x_1)=\sum_{i_0,\cdots,i_k}\mu_{\Sigma_{\Ch}\mathcal{C}}^{k+i_0+\cdots+i_k}\big(~\overbrace{\delta_{\mathcal{L}_k},\cdots,\delta_{\mathcal{L}_k}}^{i_k},x_k,\overbrace{\delta_{\mathcal{L}_{k-1}},\delta_{\mathcal{L}_{k-1}},\cdots,\delta_{\mathcal{L}_{k-1}}}^{i_{k-1}},\cdots,x_1,\overbrace{\delta_{\mathcal{L}_0},\cdots,\delta_{\mathcal{L}_0}}^{i_0}~\big).
\end{equation}

\begin{defn}
    We define the \emph{$A_\infty$-category of twisted complexes with cochain complex coefficients} $\Tw_{\Ch}\mathcal{C}$ to be the full subcategory of $\Xi\mathcal{C}$ with objects whose curvature term $\mu^0$ is zero. The $A_\infty$-structure maps are inherited from those of $\Xi\mathcal{C}$, defined in (\ref{eqn:Tw_Ch-Ainfty}).
\end{defn}

When $\mathcal{C}$ is cohomologically unital, $\Tw_{\Ch}\mathcal{C}$ is also cohomologically unital (see \cite[Proposition 3.7]{Kim}). Similarly when $\mathcal{C}$ is strictly unital, so is $\Tw_{\Ch}\mathcal{C}$:

\begin{lemma}
    When $\mathcal{C}$ admits a strict unit $e_L$ for each object $L\in\mathcal{C}$, any object $\mathcal{L}=(\bigoplus C^*_i\otimes L_i,\delta_{\mathcal{L}})\in\Tw_{\Ch}\mathcal{C}$ also admits a strict unit given by
    \begin{equation}
        e_{\mathcal{L}}=\bigoplus \id_{C_i^*}\otimes e_{L_i}.
    \end{equation}
\end{lemma}
\begin{proof}
    Since $e_{L_i}$ are strict units, any $A_\infty$-operation of length $\geq3$ involving $e_{L_i}$ vanishes. Therefore it remains to check two chain-level identities
    \begin{equation}
        \mu^1_{\Tw_{\Ch}}(e_\mathcal{L})=0,\quad
        \mu^2(x,e_\mathcal{K})=(-1)^{\lvert x\rvert}\mu^2_{\Tw_{\Ch}}(e_\mathcal{L},x)=x.
    \end{equation}
    The first identity can be checked as
    \begin{equation}
        \mu^1_{\Tw_{\Ch}}(e_\mathcal{L})=\mu^1_{\Sigma_{\Ch}}(e_\mathcal{L})+\mu^2_{\Sigma_{\Ch}}(e_\mathcal{L},\delta_\mathcal{L})+\mu^2_{\Sigma_{\Ch}}(\delta_\mathcal{L},e_\mathcal{L}),
    \end{equation}
    and for each primitive morphism $\phi_{ij}\otimes x_{ij}\in\hom(C^*_i\otimes L_i, C^*_j\otimes L_j)$, the above equation evaluates as
    \begin{align}
        \mu^2(\phi_{ij}\otimes x_{ij},\id_{C^*_i}\otimes e_{L_i})+ \mu^2(\id_{C^*_j}\otimes e_{L_j},\phi_{ij}\otimes x_{ij})&=\phi_{ij}\otimes\mu^2(x_{ij},e_{L_i})+(-1)^{\lvert\phi_{ij}\rvert}\phi_{ij}\otimes\mu^2(e_{L_j},x_{ij})\\
        &=\phi_{ij}\otimes x_{ij}+(-1)^{\lvert\phi_{ij}\rvert+\lvert x_{ij}\rvert}\phi_{ij}\otimes x_{ij}=0,
    \end{align}
    since $\lvert\delta_\mathcal{L}\rvert=\lvert\phi_{ij}\rvert+\lvert x_{ij}\rvert$ is odd. A similar computation on primitive morphisms confirms the $\mu^2$ identities.
\end{proof}

From this point, we require that every morphism complex of $\mathcal{C}$ is finite-dimensional as a vector space: this allows arbitrary twisted complexes with coefficients in the morphism spaces of $\mathcal{C}$ to define objects of $\Tw_{\Ch}\mathcal{C}$.

There is a canonical morphism called the evaluation map, which is central to the definition of twist functors.

\begin{defn}  \label{def:ev}
    For two objects $V,L\in\mathcal{C}$ and a $\mathbb{Z}/2$-graded cochain complex $C^*$ over $k$, the \emph{evaluation map}
    \begin{equation}
        ev_{V,C\otimes L}\in\hom((C^*\otimes \hom_\mathcal{C}(V,L))\otimes V, C^*\otimes L)
    \end{equation}
    is the canonical map $\id_{C^*}\otimes ev_{V,L}$, where the identification
    \begin{equation}
        \hom_\mathcal{C}(V,L)^\vee\otimes\hom_\mathcal{C}(V,L)\cong\hom_k(\hom_\mathcal{C}(V,L),\hom_\mathcal{C}(V,L))
    \end{equation}
    matches the identity on the right-hand side with $ev_{V,L}$ on the left-hand side. Similarly, for any
    \begin{equation}
        \mathcal{L}=(\bigoplus_{i\in I}C^*_i\otimes L_i,\delta_\mathcal{L})\in\Tw_{\Ch}\mathcal{C},
    \end{equation}
    the \emph{evaluation map}
    \begin{equation}
        ev_{V,\mathcal{L}}\in\hom_{\Tw_{\Ch}\mathcal{C}}(\hom_{\Tw_{\Ch}\mathcal{C}}(V,\mathcal{L})\otimes V, \mathcal{L)}
    \end{equation}
    is given as the direct sum of the individual evaluation maps $ev_{V,C^*_i\otimes L_i}$ defined above. In particular, the
    \begin{equation}
        \hom_{\Tw_{\Ch}\mathcal{C}}((C^*_i\otimes\hom_\mathcal{C}(V,L_i))\otimes V, C_j^*\otimes L_j)
    \end{equation}
    component of the evaluation map is zero when $i\neq j$.
\end{defn}

We show that the evaluation map is a cocycle.

\begin{lemma}
    For any $V\in\mathcal{C}$ and $\mathcal{L}\in\Tw_{\Ch}\mathcal{C}$, the evaluation map
    \begin{equation}
        ev_{V,\mathcal{L}}\in\hom_{\Tw_{\Ch}\mathcal C}\bigl(\hom_{\Tw_{\Ch}\mathcal C}(V,\mathcal L)\otimes V,\mathcal{L}\bigr)^0
    \end{equation}
    is a cocycle.
\end{lemma}
\begin{proof}
    The canonical identification of vector spaces
    \begin{equation}
        \hom_{\Tw_{\Ch}\mathcal C}\bigl(\hom_{\Tw_{\Ch}\mathcal C}(V,\mathcal L)\otimes V,\mathcal{L}\bigr)\cong \hom_k\bigl(\hom_{\Tw_{\Ch}\mathcal C}(V,\mathcal L),\hom_{\Tw_{\Ch}\mathcal C}(V,\mathcal L)\bigr)
    \end{equation}
    is actually an identification of cochain complexes, since it comes from the tensor-hom adjunction in the category $\Ch_k$. Under this identification, $ev_{V,\mathcal{L}}$ corresponds to $\id$ on the right, which is a cocycle. Therefore $ev_{V,\mathcal{L}}$ is also a cocycle.
\end{proof}

We remark that whenever $\hom_{\Tw_{\Ch}\mathcal C}(V,\mathcal L)$ is used as a coefficient cochain complex, its differential is given by
\begin{equation} \label{eqn:T_VL-cochain-differential}
    \partial(z)=(-1)^{\lvert z\rvert}\mu^1_{\Tw_{\Ch}\mathcal{C}}(z),
\end{equation}
which is an essential point in the above proof.

Using the shift of the evaluation morphism, we will define the twist functor on the level of objects. We first see the toy case when $\mathcal{L}$ consists of a single object $L\in\mathcal{C}$.

\begin{defn} \label{def:T_VL}
    Given two objects $V$, $L$, define the \emph{shifted evaluation map}
    \begin{equation}
        ev_{V,L}[1]:\hom_\mathcal{C}(V,L)[1]\otimes V\to L
    \end{equation}
    by the identification with $ev_{V,L}$ under
    \begin{equation}
        \hom_k(\hom_\mathcal{C}(V,L),\hom_\mathcal{C}(V,L))^0\cong\hom_k(\hom_\mathcal{C}(V,L)[1],\hom_\mathcal{C}(V,L))^1.
    \end{equation}
    The twisted complex $T_VL\in\Tw_{\Ch}\mathcal{C}$ is defined as
    \begin{equation}
        T_VL=(\hom_\mathcal{C}(V,L)[1]\otimes V)\oplus L,
    \end{equation}
    with the twisting differential given by
    \begin{equation}
        \delta_{T_VL}=\begin{pmatrix}
            0 &0\\ ev_{V,L}[1] &0
        \end{pmatrix}.
    \end{equation}
\end{defn}

The above definition can be identified as the mapping cone construction of Definition \ref{def:mapping-cone} applied to two objects in $\Sigma_{\Ch}\mathcal{C}$. With this interpretation, we generalize the above definition to arbitrary $\mathcal{L}\in\Tw_{\Ch}\mathcal{C}$.

\begin{defn}
    Given an object $V\in\mathcal{C}$ and a twisted complex with cochain complex coefficients
    \begin{equation}
        \mathcal{L}=\bigoplus_{i\in I} C^*_i\otimes L_i\in\Tw_{\Ch}\mathcal{C},
    \end{equation}
    define the \emph{algebraic twist} $T_V\mathcal{L}\in\Tw_{\Ch}\mathcal{C}$ given by the mapping cone of $ev_{V,\mathcal{L}}$.
\end{defn}

Since $ev_{V,\mathcal{L}}$ is a cocycle, our definition indeed defines an object of $\Tw_{\Ch}\mathcal{C}$.

\subsubsection{The algebraic twist functor $T_V$}

We now assume that $\mathcal{C}$ is strictly unital.

The functor $T_V$ will be defined as the mapping cone of a natural transformation $ev:S_V\to\id$ in the $A_\infty$-functor category. We first define the functor $S_V$.

\begin{defn} \label{defn:S_V}
    Given an object $V\in\mathcal{C}$, the functor $S_V:\Tw_{\Ch}\mathcal{C}\to\Tw_{\Ch}\mathcal{C}$ is defined as
    \begin{equation}
        S_V(\mathcal{L})=\hom_{\Tw_{\Ch}\mathcal{C}}(V,\mathcal{L})\otimes V,
    \end{equation}
    with trivial twisting differential $\delta$, and cochain complex differential as in (\ref{eqn:T_VL-cochain-differential}). For $x_i\in\hom_{\Tw_{\Ch}\mathcal{C}}(\mathcal{L}_{i-1},\mathcal{L}_i)$, the length $d$ component
    \begin{equation}
        S_V^d:\hom_{\Tw_{\Ch}\mathcal{C}}(V,\mathcal{L}_0)\otimes V\to\hom_{\Tw_{\Ch}\mathcal{C}}(V,\mathcal{L}_d)\otimes V
    \end{equation}
    is defined by 
    \begin{equation}
        S_V^d(x_d,\cdots,x_1)=Y^d(x_d,\cdots,x_1)\otimes e_V,
    \end{equation}
    where $Y^d$ denotes the unshifted map
    \begin{equation}
        Y^d(x_d,\cdots,x_1)(z)=(-1)^{\lvert z\rvert}\mu^{d+1}_{\Tw_{\Ch}\mathcal{C}}(x_d,\cdots,x_1,z).
    \end{equation}
\end{defn}

In the lemma below, we verify that $S_V$ defines an $A_\infty$-functor.

\begin{lemma}
    The $A_\infty$-functor relations hold for $S_V$ defined above.
\end{lemma}
\begin{proof}
    The $A_\infty$-functor relation for $S_V$ with $d$ inputs is
    \begin{equation}
        \sum\mu^r(S_V^{s_r}(x_d,\cdots),\cdots,S_V^{s_1}(\cdots,x_1))=\sum(-1)^{\dagger_i}S_V^{d-m+1}(x_d,\cdots,\mu^m(x_{i+m},\cdots,x_{i+1}),x_i,\cdots,x_1).
    \end{equation}
    Since $S_V$ involves the strict identity $e_V$, any $\mu^r$ operation on the left hand side for $r\geq3$ vanishes. Thus the left-hand side consists of the sum of the $\mu^1$-term
    \begin{equation}
        \mu^1_{\Sigma_{\Ch}\mathcal{C}}(Y^d(x_d,\cdots,x_1)\otimes e_V)=(-1)^{\lvert Y(x_d,\cdots,x_1)\rvert}\partial_{\Ch}(Y^d(x_d,\cdots,x_1))\otimes e_V
    \end{equation}
    and the sum coming from $r=2$: incorporating the Koszul sign for the tensor product, this is
    \begin{equation}
        (\sum(-1)^{\lvert Y^j(x_j,\cdots,x_1)\rvert}\cdot Y^{d-j}(x_d,\cdots,x_{j+1})\circ Y^{j}(x_j,\cdots,x_1) )\otimes e_V
    \end{equation}
    and the right-hand side is
    \begin{equation}
        \sum_{m=1}^d\sum_{i=0}^{d-m}
        (-1)^{\dagger_i} Y^{d-m+1}(x_d,\cdots,\mu^m(x_{i+m},\cdots,x_{i+1}),x_i,\cdots,x_1)\otimes e_V.
    \end{equation}
    Thus we only need to verify that the cochain coefficients agree: their evaluations at $z$ are
    \begin{align}
        &\sum\mu^{d-j+1}(x_d,\cdots,x_{j+1},\mu^{j+1}(x_j,\cdots,x_1,z))\\
        =&\sum(-1)^{\norm{z}+1+\dagger_i}\mu^{d-m+2}(x_d,\cdots,\mu^m(x_{i+m},\cdots,x_{i+1}),x_i,\cdots,x_1,z)
    \end{align}
    which is the decomposition of the $A_\infty$-relations for $(x_d,\cdots,x_1,z)$ decomposed into the terms where the inner $\mu$-operation incorporates $z$, and the terms that do not.
\end{proof}

Our next step is defining the natural transformation $ev$.

\begin{defn}
    For an object $V\in\mathcal{C}$, the \emph{evaluation natural transformation} $ev$ is defined on objects $\mathcal{L}\in\Tw_{\Ch}\mathcal{C}$ as
    \begin{equation}
        ev^0_\mathcal{L}=ev_{V,\mathcal{L}}:S_V(\mathcal{L)}\to\mathcal{L},
    \end{equation}
    and all higher components $ev^d\equiv 0$ for $d\geq1$.
\end{defn}
\begin{proposition}
    The $ev$ defined above satisfies the $A_\infty$-natural transformation equations.
\end{proposition}
\begin{proof}
    Recall that the length $d$ part of the natural transformation equation reads
    \begin{equation}
        (\mu^1_{fun}(ev))^d(x_d,\cdots,x_1)=\sum(-1)^{\dagger_i}\mu^{d-i+t+1}_{\Tw_{\Ch}\mathcal{C}}(x_d,\cdots,x_{i+1},ev^0_{\mathcal{L}_i},S_V^{j_t}(x_i,\cdots,x_{i-j_t+1}),\cdots,S_V^{j_1}(x_{j_1},\cdots,x_1)),
    \end{equation}
    for $x_i\in\hom_{\Tw_{\Ch}\mathcal{C}}(\mathcal{L}_{i-1},\mathcal{L}_i)$. Since the $S_V^d$ all involve the strict unit $e_V$, any $\mu^k$ for $k\geq3$ in the above sum vanishes. Therefore the above sum reduces to
    \begin{equation} \label{eqn:ev-nat-trans}
        \mu^{d+1}_{\Tw_{\Ch}\mathcal{C}}(x_d,\cdots,x_1,ev^0_{\mathcal{L}_0})+(-1)^{\dagger_d}\mu^2_{\Tw_{\Ch}\mathcal{C}}(ev^0_{\mathcal{L}_d},S_V^d(x_d,\cdots,x_1)).
    \end{equation}
    Under the identification
    \begin{equation}
        \hom_{\Tw_{\Ch}\mathcal{C}}(\hom_{\Tw_{\Ch}\mathcal{C}}(V,\mathcal{L}_0)\otimes V,\mathcal{L}_d)\cong \hom_k(\hom_{\Tw_{\Ch}\mathcal{C}}(V,\mathcal{L}_0),\hom_{\Tw_{\Ch}\mathcal{C}}(V,\mathcal{L}_d)),
    \end{equation}
    the first term in (\ref{eqn:ev-nat-trans}) corresponds to $Y^d$, while the second term corresponds to $-Y^d$, which cancel each other.
\end{proof}

Therefore the mapping cone of the natural transformation defines an $A_\infty$-functor, which we will call $T_V$.

\begin{defn}
    Given an object $V\in\mathcal{C}$, the \emph{twist functor along $V$}, $T_V:\Tw_{\Ch}\mathcal{C}\to\Tw_{\Ch}\mathcal{C}$ is defined as the mapping cone of the even degree natural transformation $ev:S_V\to\id$. On the level of objects,
    \begin{equation} \label{eqn:twist-functor}
        T_V\mathcal{L}=\hom_{\Tw_{\Ch}\mathcal{C}}(V,\mathcal{L})[1]\otimes V\oplus \mathcal{L},
    \end{equation}
    with the internal differential
    \begin{equation}
        \delta_{T_V\mathcal{L}}=\begin{pmatrix}
            0 &0\\ ev_{V,\mathcal{L}}[1] &\delta_\mathcal{L}
        \end{pmatrix}.
    \end{equation}
    For a morphism $x\in\hom_{\Tw_{\Ch}\mathcal{C}}(\mathcal{L}_0,\mathcal{L}_1)$, the length 1 part can be computed as
    \begin{equation}
        T_V^1(x)=\begin{pmatrix}
            \Phi^1(x)\otimes e_V &0\\
            0 &x
        \end{pmatrix},
    \end{equation}
    and for $x_i\in\hom_{\Tw_{\Ch}\mathcal{C}}(\mathcal{L}_{i-1},\mathcal{L}_i)$, the length $d$ part for $d\geq 2$ is computed as
    \begin{equation}
        T_V^d(x_d,\cdots,x_1)=\begin{pmatrix}
            \Phi^d(x_d,\cdots,x_1)\otimes e_V &0\\
            0 &0
        \end{pmatrix}
    \end{equation}
    for the shifted coefficient summand $\Phi^d=Y^d[1]$.
\end{defn}   

One advantage of defining the twist functor as a mapping cone, instead of defining it directly as above, is that the comparison with the twist functor $\mathcal{T}_V$ defined for $A_\infty$-modules in \cite[Section 3a]{Seidel08} is immediate.

\begin{proposition}
    Via the Yoneda embedding functor $\mathcal{Y}:\Tw_{\Ch}\mathcal{C}\to$ mod-$\mathcal{C}$, the twist functor $T_V$ is naturally quasi-isomorphic to the algebraic twist $\mathcal{T}_{\mathcal{Y}_V}$ defined for $A_\infty$-modules.
\end{proposition}
\begin{proof}
    Recall that the twist functor for $A_\infty$-modules is given for objects as
    \begin{equation}
        \mathcal{T}_{\mathcal{Y}_V}(\mathcal{M})=\hom_{mod-\mathcal{C}}(\mathcal{Y}_V,\mathcal{M})[1]\otimes\mathcal{Y}_V\oplus \mathcal{M},
    \end{equation}
    which is the mapping cone of the evaluation map for $A_\infty$-modules
    \begin{equation}
        ev_{mod}:\hom_{\mathrm{mod}-\mathcal C}(\mathcal Y_V,\mathcal{M})\otimes\mathcal Y_V\longrightarrow\mathcal{M}.
    \end{equation}
    Now from the Yoneda lemma, we have a quasi-isomorphism
    \begin{equation}
        y_\mathcal{L}:\hom_{\Tw_{\Ch}\mathcal{C}}(V,\mathcal{L})\simeq\hom_{mod-\mathcal{C}}(\mathcal{Y}_V,\mathcal{Y}_\mathcal{L}), 
    \end{equation}
    and this induces a map
    \begin{equation}
        y_\mathcal{L}\otimes\id:\mathcal{Y}_{S_V(\mathcal{L})}=\hom_{\Tw_{\Ch}\mathcal{C}}(V,\mathcal{L})\otimes\mathcal{Y}_V\to \hom_{mod-\mathcal{C}}(\mathcal{Y}_V,\mathcal{Y}_\mathcal{L})\otimes\mathcal{Y}_V.
    \end{equation}
    The key claim is that the diagram below naturally commutes:
    \begin{equation}
    \begin{tikzcd}
        \mathcal{Y}_{S_V(\mathcal{L})} \arrow[rr,"\mathcal{Y}(ev_{V,\mathcal{L}})"] \arrow[d,"y_\mathcal{L}\otimes\id"'] &&\mathcal{Y_\mathcal{L}} \arrow[d,"="]\\
        \hom_{mod-\mathcal{C}}(\mathcal{Y}_V,\mathcal{Y}_\mathcal{L})\otimes \mathcal{Y}_V \arrow[rr,"ev_{mod}"'] &&\mathcal{Y_\mathcal{L}}
    \end{tikzcd}    
    \end{equation}
    Under the identification
    \begin{equation}
        \hom_{mod-\mathcal{C}}(\hom_{\Tw_{\Ch}\mathcal{C}}(V,\mathcal{L})\otimes \mathcal{Y}_V,\mathcal{Y_\mathcal{L}})\cong\hom_k(\hom_{\Tw_{\Ch}\mathcal{C}}(V,\mathcal{L}),\hom_{mod-\mathcal{C}}(\mathcal{Y}_V,\mathcal{Y}_\mathcal{L})),
    \end{equation}
    $\mathcal{Y}(ev_{V,\mathcal{L}})$ corresponds to $y_\mathcal{L}$ since $ev_{V,\mathcal{L}}$ corresponds to the identity before Yoneda. Also, the composition
    \begin{equation}
        ev_{mod}\circ y_\mathcal{L}\otimes\id
    \end{equation}
    corresponds to $y_\mathcal{L}$ by the definition of the evaluation map for modules.
    
    Since the Yoneda functor preserves shifts and direct sums, taking mapping cones shows that
    \begin{equation}
        \mathcal{T}_{\mathcal{Y}_V}\circ\mathcal{Y}\simeq\mathcal{Y}\circ T_V,
    \end{equation}
    which is moreover a natural quasi-isomorphism.
\end{proof}

\subsection{Infinitesimal gradings with cochain complex coefficients} \label{ssec:Tw_Ch}

The goal of this subsection is to develop a framework in which the algebraic twist functor $T_V$ along an infinitesimally graded object $(V,\xi_V)$ is infinitesimally graded, and moreover carries an infinitesimal grading inducing the specified object-level grading. Our exposition for this subsection will assume that $\mathcal{C}$ is strictly unital, which simplifies computation.

\subsubsection{The total infinitesimal grading of $\Tw_{\Ch}\mathcal{C}$} 

Recall that our goal is to take into account the grading that each coefficient cochain complex carries. As a first step, we define an infinitesimal grading on the category of cochain complexes. The $A_\infty$ structure on the category of $\mathbb{Z}/2$-graded cochain complexes is defined as
\begin{align}
    \mu^1_{\Ch}(\phi)&=(-1)^{\lvert\phi\rvert}\partial_{\Ch}(\phi),\\
    \mu^2_{\Ch}(\phi_2,\phi_1)&=(-1)^{\lvert\phi_1\rvert}\phi_2\circ\phi_1.
\end{align}

\begin{proposition} \label{prop:eta-ch}
    The $A_\infty$-category of cochain complexes admits an infinitesimal grading $\eta_{\Ch}$, which is trivial in lengths $\geq1$, with length-zero component given by minus the differential of the cochain complex
    \begin{equation}
        \eta^0\vert_{C^*}=-\partial_C.
    \end{equation}
\end{proposition}
\begin{proof}
    For length $1$, we compute for $\phi\in\hom_{\Ch}(C^*,D^*)$,
    \begin{align}
        (\delta\eta)^1(\phi)&=-\mu^2_{\Ch}(\phi,-\partial_C)-\mu^2_{\Ch}(-\partial_D,\phi)\\
        &= (-1)^{\deg(\partial_C)}\phi\circ\partial_C+(-1)^{\deg(\phi)}\partial_D\circ\phi = \mu^1_{\Ch}(\phi).
    \end{align}
    And for length $k\geq2$, $(\delta\eta)^k\equiv0$ since all the $\mu^{k+1}$ are trivial.
\end{proof}

This infinitesimal grading will be assembled with the infinitesimal grading of $\mathcal{C}$ to produce an infinitesimal grading of $\Tw_{\Ch}\mathcal{C}$. We first explain how both $\eta_{\Ch}$ and $\eta_\mathcal{C}$ extend to the additive extension $\Sigma_{\Ch}\mathcal{C}$.

\begin{defn}
    Define the \emph{additive extensions} $\Sigma\eta_{\Ch}$, $\Sigma\eta_\mathcal{C}$ as follows: on a twisted complex with cochain coefficients
    \begin{equation}
        \mathcal{L}=\bigoplus_{i\in I} C^*_i\otimes L_i,
    \end{equation}
    the length zero parts are given as
    \begin{align}
        \Sigma\eta_{\Ch}^0\vert_{\mathcal{L}}&=\bigoplus_{i\in I}(-\partial_{C^*_i}\otimes e_{L_i}),\\
        \Sigma\eta_\mathcal{C}^0\vert_{\mathcal{L}}&=\bigoplus_{i\in I}(\id_{C^*_i}\otimes\eta^0\vert_{L_i}), 
    \end{align}
    and the positive length parts are given for primitive inputs as
    \begin{align}
        &\Sigma\eta_{\Ch}^d\equiv0,\\
        &\Sigma\eta^d_{\mathcal{C}}(\phi_d\otimes x_d,\cdots,\phi_1\otimes x_1)=(-1)^{\bowtie} \phi_d\circ\cdots\circ\phi_1\otimes\eta^d(x_d,\cdots,x_1),
    \end{align}
    for the Koszul sign $\bowtie=\sum_{i<j}\lvert\phi_i\rvert\cdot\norm{x_j}$.
\end{defn}

To extend these to Hochschild cochains of $\Tw_{\Ch}\mathcal{C}$, we now add the contributions of the twisting differentials. First, we give a general definition for the extension of any Hochschild cochain of $\mathcal{C}$, similar to Definition \ref{def:cc-extends-Tw}.

\begin{defn}
    For a Hochschild cochain $\beta\in CC^*(\mathcal{C},\mathcal{C})$, define the Hochschild cochain $\overline{\Tw}\beta\in CC^*(\Tw_{\Ch}\mathcal{C},\Tw_{\Ch}\mathcal{C})$ as
    \begin{equation}
        \overline{\Tw}\beta^k(x_k,\cdots,x_1)=\sum_{j_0,\cdots,j_k}\Sigma\beta^{k+j_0+\cdots+j_k}(\underbrace{\delta_{\mathcal{L}_k},\cdots,\delta_{\mathcal{L}_k}}_{j_k},x_k,\cdots,x_1,\underbrace{\delta_{\mathcal{L}_0},\cdots,\delta_{\mathcal{L}_0}}_{j_0}).
    \end{equation}
\end{defn}

Since $\eta_{\Ch}$ is supported only on length zero, this extends as
\begin{equation}
    \Tw\eta_{\Ch}^0\vert_{\mathcal{L}}=\bigoplus_{i\in I}(-\partial_{C^*_i}\otimes e_{L_i}),\quad \Tw\eta_{\Ch}^d\equiv 0.
\end{equation}
Here $\Tw\eta_{\Ch}$ accounts for the coefficient differentials. The coefficient degree operators enter through the object gradings introduced below.

We claim that these assemble together to define an infinitesimal grading of $\Tw_{\Ch}\mathcal{C}$.

\begin{proposition}
    Let $(\mathcal{C},\eta_\mathcal{C})$ be an infinitesimally graded category. Then the $A_\infty$-category of twisted complexes with cochain complex coefficients $\mathcal{T}=\Tw_{\Ch}\mathcal{C}$ inherits a canonical infinitesimal grading
    \begin{equation}
        \eta_\mathcal{T}=\overline{\Tw}\eta+\Tw\eta_{\Ch}-\boldsymbol{\delta},
    \end{equation}
    for the Hochschild cochain $\boldsymbol{\delta}$ defined in Definition \ref{defn-bold-delta}.
\end{proposition}
\begin{proof}
    We first compute the pre-Euler class of $\Tw_{\Ch}\mathcal{C}$. There are three contributions: the internal differential of $\Ch_k$, the $A_\infty$-operations of $\mathcal{C}$, and the twisting differentials. The proof of Proposition \ref{prop:infinitesimal-grading-Tw} computes the contribution of the latter two, so we conclude that
    \begin{equation}
        \kappa_{\Tw_{\Ch}\mathcal{C}}= \overline{\Tw}(\kappa_\mathcal{C})+\Phi+[\mu_{\Tw_{\Ch}\mathcal{C}},\boldsymbol{\delta}],
    \end{equation}
    where $\Phi$ is the contribution of $\Ch_k$ given as a Hochschild cochain with only length 0 and 1 components: it evaluates on a primitive input as
    \begin{equation}
        \Phi^1(\phi\otimes z)=(-1)^{\lvert\phi\rvert+\lvert z\rvert}\partial_{\Ch}(\phi)\otimes z,
    \end{equation}
    and the length zero part is given as 
    \begin{equation}
        \Phi^0_\mathcal{L}=\Phi^1(\delta_\mathcal{L}).
    \end{equation}
    It remains for us to prove that
    \begin{equation}
        -[\mu_{\Tw},\Tw\eta_{\Ch}]=\Phi.
    \end{equation}
    Since $\eta_{\Ch}$ only involves strict units, it cannot compose nontrivially with more than two elements in morphism spaces of $\mathcal{C}$. Therefore it is enough to restrict to when
    \begin{equation}
        \mathcal{L}=C^*_1\otimes L_1\oplus C^*_2\otimes L_2,\quad \delta_\mathcal{L}=\sum \phi_j\otimes x_j.
    \end{equation}
    For length zero, signs coming from the Leibniz rule show that
    \begin{align}
        -\mu^1(\Tw\eta_{\Ch}^0\vert_{\mathcal{L}})&=\mu^2(\delta_\mathcal{L},\partial_{C_1}\otimes e_{L_1})+\mu^2(\partial_{C_2}\otimes e_{L_2},\delta_\mathcal{L})\\
        &= \sum(-1)^{\norm{x_j}}\phi_j\circ\partial_{C_1}\otimes \mu^2(x_j,e_{L_1})+\sum(-1)^{\lvert\phi_j\rvert}\partial_{C_2}\circ \phi_j\otimes\mu^2(e_{L_2},x_j)\\
        &=\sum(-1)^{\lvert x_j\rvert+\lvert\phi_j\rvert}(\partial_{\Ch}\phi_j)\otimes x_j.
    \end{align}
    For length one, similarly, pick a primitive input $\phi\otimes z\in\hom(C^*\otimes K, D^*\otimes L)$. Then
    \begin{align}
        -&\mu^2_{\Tw_{\Ch}}(\phi\otimes z,-\partial_C\otimes e_K)-\mu^2_{\Tw_{\Ch}}(-\partial_D\otimes e_L,\phi\otimes z)\\
        &=(-1)^{\norm{z}}\phi\circ\partial_C\otimes\mu^2(z,e_K)+(-1)^{\lvert\phi\rvert}\partial_D\circ\phi\otimes\mu^2(e_L,z)\\
        &=(-1)^{\lvert\phi\rvert+\lvert z\rvert}(\partial_{\Ch}\phi)\otimes z,
    \end{align}
    and for length $k\geq2$,
    \begin{equation}
        -[\mu_{\Tw_{\Ch}\mathcal{C}},\Tw\eta_{\Ch}]^k(\phi_k\otimes z_k,\cdots,\phi_1\otimes z_1)=-\sum_i \mu^{k+1}_{\Tw_{\Ch}}(\phi_k\otimes z_k,\cdots,-\partial_{C_i}\otimes e_{L_i},\cdots,\phi_1\otimes z_1)=0,    
    \end{equation}
    since $e_{L_i}$ are strict units.
\end{proof}

\subsubsection{Infinitesimally graded objects of $\Tw_{\Ch}\mathcal{C}$} We now analyze which objects of $(\mathcal{T}=\Tw_{\Ch}\mathcal{C},\eta_\mathcal{T})$ are infinitesimally graded. Instead of providing a complete answer, we describe a large class of objects that admit an infinitesimal grading.

As such, we first define infinitesimally graded cochain complexes.

\begin{defn} \label{defn:inf-graded-cochain-cpx}
    An \emph{infinitesimally graded cochain complex} is a $\mathbb{Z}/2$-graded cochain complex $(C^*,\partial_C)$ together with a linear map $\psi_C\in\hom_{\Ch}(C^*,C^*)$ of even degree such that
    \begin{equation}
        \psi_C\circ\partial_C-\partial_C\circ\psi_C=\partial_C.
    \end{equation}
    The \emph{$A_\infty$-category of infinitesimally graded cochain complexes} has objects $(C^*,\psi_C)$ and the same morphisms and $A_\infty$ structure as $\Ch_k$. We also define a length zero Hochschild cochain $\xi^0$ as
    \begin{equation}
        \xi^0_{(C^*,\psi_C)}=\psi_C\in\hom_{\Ch}(C^*,C^*).
    \end{equation}
\end{defn}

We remark that $\Ch_k$ admits a canonical infinitesimal grading $\eta_{\Ch}$, while $\psi_C$ is an extra datum which records how an object $C^*$ is infinitesimally gradeable. For example, a $\mathbb{Z}$-graded cochain complex admits a canonical infinitesimal grading, multiplication by degree.

\begin{proposition} \label{prop:eta-ch}
    The pair $(C^*,\psi_C)$ defines an infinitesimally graded object with respect to the infinitesimal grading $\eta_{\Ch}$, and the induced grading operator on morphism spaces is given on $\hom_{\Ch}(C^*,D^*)$ by
    \begin{equation}
        \Gamma(\phi)=\psi_D\circ\phi-\phi\circ\psi_C.
    \end{equation}
\end{proposition}
\begin{proof}
    We compute
    \begin{equation}
        \mu^1_{\Ch}(\psi_C)=\partial_C\circ\psi_C-\psi_C\circ\partial_C=-\partial_C,
    \end{equation}
    and the induced grading operator can be computed as
    \begin{align}
        \Gamma(\phi)&= \eta^1(\phi)-\mu^2_{\Ch}(\phi,\psi_C)-(-1)^{\norm{\phi}}\mu^2(\psi_D,\phi)\\
        &=-(-1)^{\deg(\psi_C)}\phi\circ\psi_C-(-1)^{\norm{\phi}}\cdot(-1)^{\deg(\phi)}\psi_D\circ\phi\\
        &=\psi_D\circ\phi-\phi\circ\psi_C.  \qedhere
    \end{align}
\end{proof}

We restrict our attention to objects of $\Tw_{\Ch}\mathcal{C}$ whose coefficient cochain complexes are infinitesimally graded.

\begin{defn}
    Define the full subcategory of \emph{twisted complexes with infinitesimally graded cochain complex coefficients} to consist of the objects
    \begin{equation}
        \mathcal{L}=\bigoplus C^*_i\otimes L_i\in\Tw_{\Ch}\mathcal{C}
    \end{equation}
    for which each cochain complex $C_i$ admits an infinitesimal grading as in Definition \ref{defn:inf-graded-cochain-cpx}.
\end{defn}

We claim that all twisted complexes as above with trivial twisting differential are graded with respect to the infinitesimal grading $\eta_\mathcal{T}$.

\begin{proposition}
    Suppose that each object $L_i$ and each coefficient cochain complex $C^*_i$ for the object
    \begin{equation}
        (\mathcal{L}=\bigoplus_{i\in I}C^*_i\otimes L_i, \delta_\mathcal{L}\equiv0)
    \end{equation}
    is infinitesimally graded, i.e. admits an element $\xi_{L_i}$ as before and a map $\psi_i:C^*_i\to C^*_i$ such that
    \begin{equation}
        \psi_i\circ\partial_{C_i}-\partial_{C_i}\circ\psi_i=\partial_{C_i}.
    \end{equation}
    Then $\mathcal{L}$ with the above choices canonically defines an infinitesimally graded object of $(\Tw_{\Ch}\mathcal{C},\eta_\mathcal{T})$ given by
    \begin{equation}
        \xi_\mathcal{L}=\bigoplus_{i\in I}(\psi_i\otimes e_{L_i}+\id\otimes\xi_{L_i}).
    \end{equation}
\end{proposition}
\begin{proof}
    The claim to verify is
    \begin{equation}
        \mu^1_{\Sigma_{\Ch}\mathcal{C}}(\xi_\mathcal{L})=\eta_\mathcal{T}^0\vert_\mathcal{L}.
    \end{equation}
    First, the left-hand side evaluates on the component $C^*_i\otimes L_i$ as
    \begin{equation}
        \partial_{\Ch}\psi_i\otimes e_{L_i}+\id\otimes\mu^1(\xi_{L_i})=-\partial_i\otimes e_{L_i}+\id\otimes\eta^0_{L_i}.
    \end{equation}
    Since we assume that $\mathcal{L}$ has $\delta_\mathcal{L}=0$, the right-hand side computes on the same component as
    \begin{equation}
        -\partial_i\otimes e_{L_i}+\id\otimes\eta^0_{L_i}.
    \end{equation}
    Thus they agree.
\end{proof}

However, when the twisting differential $\delta_\mathcal{L}$ is nonzero, the twisted complex $(\mathcal{L},\delta_\mathcal{L})$ need not be infinitesimally graded even if each $C^*_i$ and $L_i$ are individually graded. 

We now introduce the category of graded twisted complexes used below.

\begin{defn}
    Set $\mathcal{T}=\Tw_{\Ch}\mathcal{C}$. We use the following decorated restriction of its graded lift. An object of $\Tilde{\mathcal{T}}$  consists of a twisted complex
    \begin{equation}
        \mathcal{L}=\big(\bigoplus_i (C^*_i,\psi_i)\otimes (L_i,\xi_i),\delta_\mathcal{L}\big)
    \end{equation}
    with chosen infinitesimal gradings $\psi_i$ for $C_i^*$ and $\xi_i$ of $L_i$, and an additional choice $\xi_\mathcal{L}\in\hom_\mathcal{T}^0(\mathcal{L},\mathcal{L})$ such that
    \begin{equation}
        \mu^1_\mathcal{T}(\xi_\mathcal{L})=\eta^0_\mathcal{T}\vert_\mathcal{L}.
    \end{equation}
    The total grading $\xi_L$ is part of the data: it is not determined by the given gradings. Morphism spaces and $A_\infty$-operations are inherited from $\mathcal{T}$. 

    The lifted infinitesimal grading on $\Tilde{\mathcal{T}}$ is given by
    \begin{equation}
        \Tilde{\eta}=\eta_\mathcal{T}+\delta\xi
    \end{equation}
    for the length zero Hochschild cochain $\xi$ associated to each object.
\end{defn}

We explicitly record the grading of a shifted morphism complex.

\begin{lemma}
    When $V,(\mathcal{L},\delta_\mathcal{L})$ are infinitesimally graded objects of $(\Tw_{Ch}\mathcal{C})$ with its induced infinitesimal grading, then the morphism complex
    \begin{equation}
        \hom(V,\mathcal{L)}[1]
    \end{equation}
    defines a graded cochain complex with
    \begin{equation}
        \psi_{\hom(V,\mathcal{L})[1]}=\Gamma_{V,\mathcal{L}}-\id,
    \end{equation}
    for the grading operator $\Gamma_{V,\mathcal{L}}$ of $\Tw_{\Ch}\mathcal{C}$.
\end{lemma}

In the following subsection, we will show that twist functors along infinitesimally graded objects inherit a natural infinitesimal grading: this will imply that the twist of an infinitesimally graded object will again be infinitesimally graded.

\subsubsection{The infinitesimal grading of the twist functor}

Now fix an infinitesimally graded object $(V,\xi_V)$ of $(\mathcal{C},\eta_\mathcal{C})$. Our goal for this subsection is to prove that the algebraic twist $T_V$ admits an infinitesimal grading such that the length zero part assigned to each object is canonical. To be precise, we will show that the induced functor $T_{(V,\xi_V)}$ on the graded lift $\Tilde{T}$ can be infinitesimally graded.

\begin{defn}
    We define the \emph{twist functor along $(V,\xi_V)\in(\Tilde{C},\eta$)} to be the endofunctor $T_{(V,\xi_V)}:\Tilde{T}\to\Tilde{T}$ defined as
    \begin{equation}
        T_{(V,\xi_V)}\mathcal{L}=\hom_{\mathcal{T}}(V,\mathcal{L})[1]\otimes V\oplus\mathcal{L}.
    \end{equation}
    The new coefficient cochain complex has differential and grading
    \begin{align}
        \partial(sx)&=(-1)^{\norm{x}}s\mu^1_\mathcal{T}(x),\\
        \psi_{\hom_\mathcal{T}(V,\mathcal{L})[1]}(sx)&=s(\Tilde{\eta}^1(x)-x).
    \end{align}
    It carries an internal twisting differential and infinitesimal grading
    \begin{equation}\label{eqn:xi_TVL}
        \delta_{T_V\mathcal{L}}=\begin{pmatrix}
            0 &0\\
            ev_{\mathcal{L}} &\delta_{\mathcal{L}}
        \end{pmatrix},\qquad \xi_{T_V\mathcal{L}}=\begin{pmatrix}
        \id\otimes\xi_V+\psi_{\hom(V,\mathcal{L})[1]}\otimes e_V &0\\
       0 & \xi_\mathcal{L}
    \end{pmatrix}.
    \end{equation}
    for the shifted evaluation $ev_\mathcal{L}$.
\end{defn}

We first verify that this indeed defines an object in $\Tilde{T}$.

\begin{lemma}
    \begin{equation}
        \mu_\mathcal{T}^1(\xi_{T_V\mathcal{L}})={\eta^0_\mathcal{T}}\vert_{T_V\mathcal{L}}
    \end{equation}
\end{lemma}
\begin{proof}
    We first compute that
    \begin{equation}
        \eta^0_{\mathcal{T}}\vert_{T_V\mathcal{L}}=\begin{pmatrix}
            \id\otimes\eta^0_V-\partial_{\hom(V,\mathcal{L})[1]}\otimes e_V &\\
            \mu_\mathcal{T}^2(ev_\mathcal{L},\eta^1_\mathcal{T})-ev_\mathcal{L} & \eta^0_\mathcal{L}
        \end{pmatrix}.
    \end{equation}
    By gradeability of $\mathcal{L}$, the bottom right matrix entry agrees. The upper left agrees since
    \begin{equation}
        \mu^1_{\mathcal{T}}(\id\otimes\xi_V)=\id\otimes\eta^0_V,\quad \mu^1_\mathcal{T}(\psi\otimes e_V)=-\partial\otimes e_V,    
    \end{equation}
    so the key computation is for the lower left term. 
    
    We first introduce the following notation. For a linear map $\Upsilon: \hom_\mathcal{T}(V,\mathcal{L}_0)\to\hom_\mathcal{T}(V,\mathcal{L}_1)$, we define the \emph{evaluation morphism associated to $\Upsilon$} by
    \begin{equation}
        \Upsilon^{ev}\coloneqq \mu^2_\mathcal{T}(ev_{\mathcal{L}_1},[sx\mapsto s\Upsilon(x)]\otimes e_V).
    \end{equation}
    Equivalently, we apply the morphism $\Upsilon$ to the Hom coefficient and evaluate the resulting morphism. The resulting definition incorporates the Koszul and suspension signs, and simplifies the notation.
    
    With this notation, we compute the lower left term as
    \begin{align}
        \mu^1_\mathcal{T}(\xi_{T_V\mathcal{L}})_{21}&=\mu^2_\mathcal{T}(ev_\mathcal{L},\id\otimes \xi_V)+\mu^2_\mathcal{T}(ev_\mathcal{L},\psi_{\hom_\mathcal{T}(V,\mathcal{L})[1]}\otimes e_V)+\mu^2_\mathcal{T}(\xi_\mathcal{L},ev_\mathcal{L})\\
        &=(x\mapsto \mu^2_\mathcal{T}(x,\xi_V)+\Tilde{\eta}_{\mathcal{T}}^1(x)-x+(-1)^{\norm{x}}\mu^2_\mathcal{T}(\xi_\mathcal{L},x))^{ev}
    \end{align}
    for a homogeneous input $x\in\hom_\mathcal{T}(V,\mathcal{L})$. Since the lifted grading is defined as
    \begin{equation}
        \Tilde{\eta}_{\mathcal{T}}^1(x)=\eta_\mathcal{T}^1(x)-\mu^2_\mathcal{T}(x,\xi)-(-1)^{\norm{x}}\mu_\mathcal{T}^2(\xi,x),
    \end{equation}
    we indeed find that this agrees with the bottom left term of $\eta^0_{\mathcal{T}}\vert_{T_V\mathcal{L}}$.
\end{proof}

Now we are ready to prove that $T_{(V,\xi_V)}$ is a gradeable functor.

\begin{thm}\label{thm:T_V-grading}
    Given an infinitesimally graded object $(V,\xi_V)\in(\Tilde{C},\eta)$, the induced twist functor $T_{(V,\xi_V)}:\Tilde{T}\to\Tilde{T}$ admits an infinitesimal grading with trivial length zero part.
\end{thm}
\begin{proof}
    The obstruction to grading is given by
    \begin{equation}
        {T_{(V,\xi_V)}}_*\Tilde{\eta}_{\mathcal{T}}-T_{(V,\xi_V)}^*\Tilde{\eta}_\mathcal{T}+H_{T_{(V,\xi_V)}}.
    \end{equation}
    We compute the obstruction in the two potentially nonzero matrix entries: the upper-left entry measures compatibility with the grading on the $\hom$ coefficient, while the lower left entry measures compatibility with evaluation. Recall that the twist functor is given by
    \begin{equation}
        T_V^1(x_1)={\begin{pmatrix}
            \Phi^1(x_1)\otimes e_V &0\\
            0 & x_1
        \end{pmatrix}},\qquad T_V^d(x_d,\cdots,x_1)=\begin{pmatrix}
            \Phi^d(x_d,\cdots,x_1)\otimes e_V &0\\
             0 &0
        \end{pmatrix}
    \end{equation}
    for the equation given by $\Phi^d=Y^d[1]$:
    \begin{equation}
        \Phi^d(x_d,\cdots,x_1)(sz)=(-1)^{\dagger_d+\norm{z}}s\mu^{d+1}_\mathcal{T}(x_d,\cdots,x_1,z).
    \end{equation}
    By definition, the pushforward of $\Tilde{\eta}$ is computed by inserting $\eta$ into $T_V$, and therefore
    \begin{equation}
    (T_*\Tilde{\eta}_{\mathcal{T}})^d(x_d,\cdots,x_1)=\begin{pmatrix}
        \sum_j \Phi^{d-j+1}(x_d,\cdots,\Tilde{\eta}^j_{\mathcal{T}}(x_{i+j},\cdots,x_{i+1}),\cdots,x_1)\otimes e_V &0\\
        0 &\Tilde{\eta}_\mathcal{T}^d(x_d,\cdots,x_1)
    \end{pmatrix}.
\end{equation}
    We omit the $(V,\xi_V)$ subscript for the twist functor for the remainder of this proof, for readability. For the pullback, we apply the components of $T_V$ to consecutive blocks of inputs then apply $\Tilde{\eta}$, giving
    \begin{equation}
        (T_{(V,\xi_V)}^*\Tilde{\eta}_\mathcal{T})^d(x_d,\cdots,x_1)=\sum\Tilde{\eta}_\mathcal{T}^r(T^{s_r}(x_d,\cdots),\cdots,T^{s_1}(\cdots,x_1)).
    \end{equation}
    Since the output will be a $2\times2$ matrix, we interpret the terms one by one. The two easy terms to check are the right entries
    \begin{equation}
        (T_V^*\Tilde{\eta}_\mathcal{T})^d_{12}=0,\qquad (T_V^*\Tilde{\eta}_\mathcal{T})^d_{22}=\Tilde{\eta}^d_\mathcal{T}.
    \end{equation}
    For the upper left entry, we compute that
    \begin{align}
        (T_V^*\Tilde{\eta}_\mathcal{T})^d_{11}(x_d,\cdots,x_1)&=(\psi_d\circ\Phi^d(x_d,\cdots,x_1)-\Phi^d(x_d,\cdots,x_1)\circ\psi_0)\otimes e_V\\
        &+\sum(-1)^k\Phi^{s_r}(x_d,\cdots)\cdots\Phi^{s_1}(\cdots,x_1)\otimes\eta_\mathcal{C}^r(e_V,\cdots,e_V)
    \end{align}
    for the sign
    \begin{equation}
        k=\sum_{i=1}^r(r-i)(\dagger^{s_1+\cdots+s_i}_{s_1+\cdots+s_{i-1}+1}+1)
    \end{equation}
    coming from moving $\Phi^{s_i}$ past the $r-i$ inputs of $e_V$. The first term comes from the $\delta\xi$ contribution of the graded lift
    \begin{equation}
        \Tilde{\eta}_\mathcal{T}=\eta_\mathcal{T}+\delta\xi,
    \end{equation}
    which vanishes when $r\geq2$, since it contributes a $\mu^{r+1}$ operation with the $r$ inputs and one strict unit. Therefore the only contribution comes from $r=1$, when
    \begin{align}
        ((\delta\xi)^1(T_V^d))_{11}&=-\mu^2_\mathcal{T}(\Phi^d\otimes e_V,\id\otimes\xi_V+\psi_0\otimes e_V)-(-1)^{\dagger_d}\mu^2_\mathcal{T}(\id\otimes\xi_V+\psi_d\otimes e_V,\Phi^d\otimes e_V)\\
        &=-\Phi^d\otimes\xi_V-(\Phi^d\circ\psi_0)\otimes e_V-(-1)^{\dagger_d}(-1)^{\dagger_d+1}(\Phi^d\otimes\xi_V+\psi_d\circ\Phi^d\otimes e_V)\\
        &=(\psi_d\circ\Phi^d-\Phi^d\circ\psi_0)\otimes e_V,
    \end{align}
    where the sign correction is from the degree of $\Phi^d$
    \begin{equation}
        \lvert\Phi^d\rvert=1+\dagger_d.
    \end{equation}
    For the lower left entry, the evaluation morphism is inserted exactly once. As such, we compute
    \begin{align}
        (T_V^*\Tilde{\eta}_\mathcal{T})^d_{21}(x_d,\cdots,x_1)&= (x\mapsto (-1)^{\dagger_d}\Tilde{\eta}^{d+1}_\mathcal{T}(x_d,\cdots,x_1,x))^{ev} \label{eqn:cancelation-T_V}\\
        &+\sum_j\Tilde{\eta}_\mathcal{T}^{d-j+r+1}(x_d,\cdots,x_{j+1},ev_{\mathcal{L}_j},\Phi^{s_r}(x_j,\cdots)\otimes e_V,\cdots,\Phi^{s_1}(\cdots,x_1)\otimes e_V).\label{eqn:cancellation-T_V}
    \end{align}
    This concludes the computation of the obstruction term. We summarize the computations:
    \begin{align}
        (T_*\Tilde{\eta}_\mathcal{T}-T^*\Tilde{\eta}_\mathcal{T}+H_{T_V})^d_{11}(x_d,\cdots,x_1)&=\sum \Phi^{d-j+1}(x_d,\cdots,\Tilde{\eta}_\mathcal{T}^j(x_{i+j},\cdots,x_{i+1}),\cdots,x_1)\otimes e_V \label{eqn:3.144}\\
        &+(-\psi_d\circ\Phi^d(x_d,\cdots,x_1)+\Phi^d(x_d,\cdots,x_1)\circ\psi_0+(1-d)\Phi^d)\otimes e_V\\
        &-\sum(-1)^k\Phi^{s_r}(x_d,\cdots)\cdots\Phi^{s_1}(\cdots,x_1)\otimes\eta_\mathcal{C}^r(e_V,\cdots,e_V),\label{eqn:3.146}
    \end{align}
    and
    \begin{align}
         (T_*\Tilde{\eta}_\mathcal{T}-T^*\Tilde{\eta}_\mathcal{T}+H_{T_V})^d_{21}(x_d,\cdots,x_1)=&-(x\mapsto (-1)^{\dagger_d}\Tilde{\eta}^{d+1}_\mathcal{T}(x_d,\cdots,x_1,x))^{ev}\\
         &- \sum_j\Tilde{\eta}_\mathcal{T}^{d-j+r+1}(x_d,\cdots,x_{j+1},ev_{\mathcal{L}_j},\Phi^{s_r}(x_j,\cdots)\otimes e_V,\cdots,\Phi^{s_1}(\cdots,x_1)\otimes e_V),
    \end{align}
    and the other two entries vanish.

    We now claim that the above obstruction admits a primitive. To do so, we first define the projection operator
    \begin{equation}
        p_\mathcal{L}=\begin{pmatrix}
            \id_{\hom_\mathcal{T}(V,\mathcal{L})[1]}\otimes e_V &0\\0&0
        \end{pmatrix}
    \end{equation}
    projecting to the shifted summand of $T_V\mathcal{L}$, and check that
    \begin{equation}\label{eqn:coeff}
        \mu^1_\mathcal{T}(p_\mathcal{L})=\begin{pmatrix}
            0 &0\\ ev_\mathcal{L} &0
        \end{pmatrix}.
    \end{equation}

    We define the primitive $\xi_{T_V}$ to have trivial length zero part, and 
    \begin{align}
        \xi_{T_V}^d(x_d,\cdots,x_1)&=+\begin{pmatrix}
            [sx\mapsto (-1)^{\dagger_d}s\Tilde{\eta}_\mathcal{T}^{d+1}(x_d,\cdots,x_1,x)]\otimes e_V &0 \\0 &0
        \end{pmatrix}\\
        &-\sum_{s_1+\cdots+s_r=d}\Tilde{\eta}^{r+1}_\mathcal{T}(T^{s_r}(x_d,\cdots),\cdots,T^{s_2},\begin{pmatrix}
            \Phi^{s_1}\otimes e_V &0\\0&0
        \end{pmatrix},p_{\mathcal{L}_0}).
    \end{align}
    The first term comes from Yoneda, while the second term compensates for the difference between its induced action on the coefficient summand and the grading inherited from the category of twisted complexes. For this, we use the projection $p_\mathcal{L}$ onto the coefficient summand: its differential is the evaluation map, as seen in \eqref{eqn:coeff}.

    Expanding out the matrix notation, we have
    \begin{align}
         (\xi_{T_V}^d(x_d,\cdots,x_1))_{11}&= +[sx\mapsto (-1)^{\dagger_d}s\Tilde{\eta}_\mathcal{T}^{d+1}(x_d,\cdots,x_1,x)]\otimes e_V \label{eqn:first-summand} \\
         &-\sum(-1)^k(\Phi^{s_r}(x_d,\cdots)\circ\cdots\circ\Phi^{s_1}(\cdots,x_1))\otimes\Tilde{\eta}_\mathcal{T}^{r+1}(e_V,\cdots,e_V)
    \end{align}
    and
    \begin{equation}
        (\xi_{T_V}^d(x_d,\cdots,x_1))_{21}=-\sum_j \Tilde{\eta}_\mathcal{T}^{d-j+r+2}(x_d,\cdots,x_{j+1},ev_{\mathcal{L}_j},\Phi^{s_r}(x_j,\cdots)\otimes e_V,\cdots,\Phi^{s_1}(\cdots,x_1)\otimes e_V,p_{\mathcal{L}_0}).
    \end{equation}
    Now we claim that
    \begin{equation}\label{eqn:hom-identity-T_V}
        \delta(\xi_{T_V})=-\mu^1_{fun}(\xi_{T_V})= {T_V}_*\Tilde{\eta}_{\mathcal{T}}-T_{V}^*\Tilde{\eta}_\mathcal{T}+H_{T_{(V,\xi_V)}}.
    \end{equation}
    We first verify the upper-left entry. By the definition of the functor differential,
    \begin{align}
        (\mu^1_{fun}(\xi_{T_V}))^d_{11}&= \sum_k(-1)^{\dagger_{s_1+\cdots+s_{k-1}}}\mu^r_\mathcal{T}(\Phi^{s_r}\otimes e_V,\cdots,\Phi^{s_{k+1}}\otimes e_V,(\xi^{s_k}_{T_V})_{11},\Phi^{s_{k-1}}\otimes e_V,\cdots,\Phi^{s_1}\otimes e_V)\\&+\sum_{j=1}^{d}\sum_{i=0}^{d-j}(-1)^{\dagger_i}\Bigl(\xi_{T_V}^{d-j+1}(x_d,\ldots,\mu_{\mathcal T}^{j}(x_{i+j},\ldots,x_{i+1}),\ldots,x_1)\Bigr)_{11}.
    \end{align}
    We first simplify the first term: when $r=1$ we have the contribution
    \begin{equation}
        \mu^1_\mathcal{T}((\xi_{T_V}^d)_{11}),
    \end{equation}
    and when $r=2$ we have the contribution 
    \begin{equation}
        \sum\mu^2_\mathcal{T}(\Phi^{d-j}\otimes e_V,(\xi^j_{T_V})_{11})+\sum(-1)^{\dagger_j}\mu^2_\mathcal{T}((\xi^{d-j}_{T_V})_{11},\Phi^j\otimes e_V).
    \end{equation}
    For $r\geq 3$, the contribution involves $\mu^r$ with the strict unit $e_V$ involved, so vanishes. 

    Now we substitute in our definition of $\xi_{T_V}$ to compute the formula. The first term
    \begin{equation}
        [sx\mapsto (-1)^{\dagger_d}s\Tilde{\eta}_\mathcal{T}^{d+1}(x_d,\cdots,x_1,x)]\otimes e_V
    \end{equation}
    contributes $\Theta\otimes e_V$, for the map $\Theta$ sending $sx$ to
    \begin{align}
        (-1)^{\dagger_d+\norm{x}}s\bigl[ &\mu^1_\mathcal{T}(\Tilde{\eta}_\mathcal{T}^{d+1}(x_d,\cdots,x_1,x))-\Tilde{\eta}^{d+1}_\mathcal{T}(x_d,\cdots,x_1,\mu^1_\mathcal{T}(x))\\
        +&\sum\mu^{d-r+1}_\mathcal{T}(x_d,\cdots,x_{r+1},\Tilde{\eta}_\mathcal{T}^{r+1}(x_r,\cdots,x_1,x))\\
        -&\sum\Tilde{\eta}_\mathcal{T}^{d-r+1}(x_d,\cdots,x_{r+1},\mu^{r+1}(x_r,\cdots,x_1,x))\\
        -& \sum(-1)^{\dagger_i+\norm{x}} \Tilde{\eta}_\mathcal{T}^{d-j+2} (x_d,\cdots,\mu^j_\mathcal{T}(x_{i+j},\cdots,x_{i+1}),\cdots,x_1,x)
        \bigr]
    \end{align}
    and the partition term contributes
    \begin{align}
        &\sum_r(-1)^{k-1}\{[\sum_j(-1)^{r-j+\dagger_{s_1+\cdots+s_{j-1}}}\Phi^{s_r}\circ\cdots\Phi^{s_{j+1}}\circ((-1)^{\dagger_{s_1+\cdots+s_{j-1}+1}^{s_1+\cdots+s_j}}\partial_{Ch}\Phi^{s_j} \label{eqn:T_V-grading-1}\\
        +& (-1)^{\dagger_{s_1+\cdots+s_{j-1}+1}^{s_1+\cdots+s_{j-1}+i}} \Phi^{s_j-m+1}(\cdots,\mu^m_\mathcal{T}(\cdots),\cdots)\circ\Phi^{s_{j-1}}\cdots\circ\Phi^{s_1}]\otimes \Tilde{\eta}_\mathcal{T}^{r+1}(e_V,\cdots,e_V)) \label{eqn:T_V-grading-2} \\
        +& \Phi^{s_r}\circ\cdots\circ\Phi^{s_1}\otimes \mu^1_\mathcal{T}(\Tilde{\eta}_\mathcal{T}^{r+1}(e_V,\cdots,e_V))\}\\
        +&\sum(-1)^k (1+(-1)^r)\Phi^{s_r}\circ\cdots\circ\Phi^{s_1}\otimes \Tilde{\eta}_\mathcal{T}^r(e_V,\cdots,e_V)
    \end{align}
    where $k$ stands for the sign convention from before. The $1+(-1)^r$ terms come from composition from the left and right, i.e. the $r=2$ terms we described.

    From the $A_\infty$-functor equations for $\Phi^d$, we have
    \begin{align}
        &(-1)^{\dagger_d}\partial_{Ch}\Phi^d(x_d,\cdots,x_1)+\sum(-1)^{\dagger_i}\Phi^{d-j+1}(x_d,\cdots,\mu^j_\mathcal{T}(x_{i+j},\cdots,x_{i+1}),x_i,\cdots,x_1)\\
        =&-\sum(-1)^{\dagger_i}\Phi^{d-i}(x_d,\cdots,x_{i+1})\circ\Phi^i(x_i,\cdots,x_1).
    \end{align}
    Applying this to (\ref{eqn:T_V-grading-1}) and (\ref{eqn:T_V-grading-2}), we simplify their sum to
    \begin{equation}
        \sum(-1)^k(\sum_{i=1}^{r-1}(-1)^{i})\Phi^{s_r}\circ\cdots\circ\Phi^{s_1}\otimes \Tilde{\eta}_\mathcal{T}^r(e_V,\cdots,e_V),
    \end{equation}
    and therefore the total contribution from the partition terms is
    \begin{align}
        &-\sum(-1)^k\Phi^{s_r}\circ\cdots\circ\Phi^{s_1}\otimes\mu^1_\mathcal{C}(\Tilde{\eta}^{r+1}(e_V,\cdots,e_V)) \label{eqn:T_V-grading-3}\\
        &+\sum(-1)^k (1+(-1)^r+\sum_{i=1}^{r-1}(-1)^i)\Phi^{s_r}\circ\cdots\circ\Phi^{s_1}\otimes\Tilde{\eta}^{r}(e_V,\cdots,e_V).
    \end{align}

    Now from the Hochschild equation
    \begin{equation}
        \delta\Tilde{\eta}_\mathcal{T}=\kappa_\mathcal{T},
    \end{equation}
    if we input strict units, then the right-hand side vanishes. Thus
    \begin{equation}
        0=\mu^1_\mathcal{T}(\Tilde{\eta}_\mathcal{T}^{r+1}(e_V,\cdots,e_V))+(1+(-1)^{r+1}-\sum_{i=0}^{r-1}(-1)^i)\Tilde{\eta}^r_\mathcal{T}(e_V,\cdots,e_V)
    \end{equation}
    where the right term is from adjacent multiplications. And therefore
    \begin{equation}
        -\mu^1_\mathcal{T}(\Tilde{\eta}_\mathcal{T}^{r+1}(e_V,\cdots,e_V))=(\sum_{i=0}^{r-1}(-1)^i)\Tilde{\eta}^r_\mathcal{T}(e_V,\cdots,e_V),
    \end{equation}
    and plugging this into (\ref{eqn:T_V-grading-3}), we conclude that the contribution from the partition terms is
    \begin{equation}\label{eqn:3.178}
        \sum(-1)^k\Phi^{s_r}\circ\cdots\circ\Phi^{s_1}\otimes\Tilde{\eta}^{r}(e_V,\cdots,e_V).
    \end{equation}
    This is the negative of the unit-input contribution \eqref{eqn:3.146}, as required. It remains to compute the contribution from the first summand \eqref{eqn:first-summand}. Since this summand applies $\Tilde{\eta}$ to the external inputs together with the $\hom$ coefficient $x$, we apply the Hochschild identity to $(x_d,\cdots,x_1,x)$. Since $ (\delta\Tilde{\eta}_\mathcal{T})^{d+1}=(1-d)\mu^{d+1}_\mathcal{T}$,
    \begin{align}
        (d-1)\mu^{d+1}_\mathcal{T}&=\mu^1_\mathcal{T}(\Tilde{\eta}_\mathcal{T}^{d+1})-\Tilde{\eta}^{d+1}_\mathcal{T}(x_d,\cdots,x_1,\mu^1(x))\\
        &+\sum\mu^{d-j+1}_\mathcal{T}(x_d,\cdots,x_{j+1},\Tilde{\eta}^{j+1}_\mathcal{T}(x_j,\cdots,x_1,x))\\
        &-\sum\Tilde{\eta}^{d-j+1}_\mathcal{T}(x_d,\cdots,x_{j+1},\mu^{j+1}_\mathcal{T}(x_j,\cdots,x_1,x))\\
        &-\sum(-1)^{\dagger_i+\norm{x}}\Tilde{\eta}^{d-j+2}_\mathcal{T}(x_d,\cdots,\mu^j_\mathcal{T}(x_{i+j},\cdots,x_{i+1}),x_i,\cdots,x_1,x)\\
        &+\sum\mu^{d-j+2}_\mathcal{T}(x_d,\cdots,\Tilde{\eta}^j_\mathcal{T}(x_{i+j},\cdots,x_{i+1}),\cdots,x_1,x)\\
        &+\mu^{d+1}_\mathcal{T}(x_d,\cdots,x_1,\Tilde{\eta}^1_\mathcal{T}(x))-\Tilde{\eta}^1_\mathcal{T}(\mu^{d+1}_\mathcal{T}(x_d,\cdots,x_1,x)).
    \end{align}
    Thus the $\Theta$-contribution sends $sx$ to
    \begin{align}
        (-1)^{\dagger_d+\norm{x}}s\bigl(& \Tilde{\eta}^1_\mathcal{T}(\mu^{d+1}_\mathcal{T}(x_d,\cdots,x_1,x))-\mu^{d+1}_\mathcal{T}(x_d,\cdots,x_1,\Tilde{\eta}^1_\mathcal{T}(x))\\
        &-\sum\mu^{d-j+2}_\mathcal{T}(x_d,\cdots,\Tilde{\eta}^j_\mathcal{T}(x_{i+j},\cdots,x_{i+1}),\cdots,x_1,x)\\
        &+(d-1)\mu^{d+1}_\mathcal{T}(x_d,\cdots,x_1,x)\bigr).
    \end{align}
    Now rewriting this in the $\Phi$ notation, we have
    \begin{align}
        (\psi_d\circ\Phi^d)(sx)&=(-1)^{\dagger_d+\norm{x}}s(\Tilde{\eta}^1_\mathcal{T}(\mu^{d+1}_\mathcal{T}(x_d,\cdots,x_1,x))-\mu^{d+1}_\mathcal{T}(x_d,\cdots,x_1,x)),\\
        (\Phi^d\circ\psi_0)(sx)&=(-1)^{\dagger_d+\norm{x}}s(\mu^{d+1}_\mathcal{T}(x_d,\cdots,x_1,\Tilde{\eta}^1_\mathcal{T}(x))-\mu^{d+1}_\mathcal{T}(x_d,\cdots,x_1,x)),
    \end{align}
    and
    \begin{align}
        &\Phi^{d-j+1}(x_d,\cdots,\Tilde{\eta}^j_\mathcal{T}(x_{i+j},\cdots,x_{i+1}),\cdots,x_1)(sx)\\
        =&(-1)^{\dagger_d+\norm{x}}s\mu^{d-j+2}_\mathcal{T}(x_d,\cdots,\Tilde{\eta}^j_\mathcal{T}(x_{i+j},\cdots,x_{i+1}),\cdots,x_1,x).
    \end{align}
    Together with \eqref{eqn:3.178}, these identities give the negative of \eqref{eqn:3.144}--\eqref{eqn:3.146}, completing the upper-left component of \eqref{eqn:hom-identity-T_V}. 

    We now check compatibility with evaluation, encoded by the lower-left entry. The pushforward and $H_{T_{(V,\xi_V)}}$ have vanishing lower-left entries, so the required identity reduces to
    \begin{align} \label{eqn:T_V-grading-21term}
        (\mu^1_{fun}(\xi_{T_V}))^d_{21}=(T_V^*\Tilde{\eta}_\mathcal{T})^d_{21}.
    \end{align}
    
    By definition, the left-hand side is
    \begin{align}
        (\mu^1_{fun}\xi_{T_V})^d_{21}&=\sum\mu^{d-j+1}_\mathcal{T}(x_d,\cdots,x_{j+1},(\xi_{T_V}^j(x_j,\cdots,x_1))_{21})\\
        &+\sum\mu_\mathcal{T}^{d-j+2}(x_d,\cdots,x_{j+1},ev_{\mathcal{L}_j},(\xi_{T_V}^j(x_j,\cdots,x_1))_{11})\\
        &+\sum(-1)^{\dagger_j}\mu^2_\mathcal{T}((\xi^{d-j}_{T_V}(x_d,\cdots,x_{j+1}))_{21},\Phi^j(x_j,\cdots,x_1)\otimes e_V)\\
        &+\sum(-1)^{\dagger_i}\xi_{T_V}^{d-j+1}(x_d,\cdots,\mu^j_\mathcal{T}(x_{i+j},\cdots,x_{i+1}),x_i,\cdots,x_1)_{21}.
    \end{align}
    After using the definition of $\xi_{T_V}$, this expands as
    \begin{align}
        &-\sum\mu^{d-j+1}_\mathcal{T}(x_d,\cdots,\Tilde{\eta}_\mathcal{T}^{j-i+r+2}(x_j,\cdots,ev_{\mathcal{L}_i},\Phi^{s_r}\otimes e_V,\cdots,\Phi^{s_1}\otimes e_V,p_{L_0})) \label{eqn:T_V-grading-6}\\
        &+\sum\mu_\mathcal{T}^{d-j+2}(x_d,\cdots,ev_{\mathcal{L}_j},[sx\mapsto(-1)^{\dagger_j}s\Tilde{\eta}_\mathcal{T}^{j+1}(x_j,\cdots,x_1,x)]\otimes e_V) \label{eqn:T_V-grading-5}\\
        &-\sum(-1)^k \mu_\mathcal{T}^{d-j+2}(x_d,\cdots,ev_{\mathcal{L}_j},\Phi^{s_r}\circ\cdots\circ\Phi^{s_1}\otimes\Tilde{\eta}^{r+1}_\mathcal{T}(e_V,\cdots,e_V)) \label{eqn:T_V-grading-7}\\
        &-\sum(-1)^{\dagger_i}\mu^2_\mathcal{T}(\Tilde{\eta}_\mathcal{T}^{d-j+r+2}(x_d,\cdots,ev_{\mathcal{L}_j},\Phi^{s_r}\otimes e_V,\cdots,\Phi^{s_1}\otimes e_V,p_{\mathcal{L}_i}),\Phi^i\otimes e_V)\label{eqn:T_V-grading-4}\\
        &-\sum (-1)^{\dagger_i}\Tilde{\eta}_\mathcal{T}^{d-q-j+r+3}(x_d,\cdots,\mu^j_\mathcal{T}(x_{i+j},\cdots,x_{i+1}),\cdots,ev_{\mathcal{L}_q},\Phi^{s_r}\otimes e_V,\cdots,p_{\mathcal{L}_0})\label{eqn:T_V-grading-8}\\
        &-\sum(-1)^{\dagger_i}\Tilde{\eta}_\mathcal{T}^{d-j+r+2}(x_d,\cdots,ev_{\mathcal{L}_j},\Phi^{s_r}\otimes e_V,\cdots,\Phi^{s_a}(\cdots,\mu^m_\mathcal{T}(\cdots),x_i,\cdots),\cdots,p_{\mathcal{L}_0}). \label{eqn:T_V-grading-9}
    \end{align}
    First we simplify (\ref{eqn:T_V-grading-5}) by noting that the strict unit vanishes for $\mu^{\geq3}$, and thus the sum simplifies to
    \begin{equation}\label{eqn:T_V-grading-204}
        \mu^2_\mathcal{T}(ev_{\mathcal{L}_d},[sx\mapsto(-1)^{\dagger_d}s\Tilde{\eta}^{d+1}_\mathcal{T}(x_d,\cdots,x_1,x)]\otimes e_V).
    \end{equation}
    The next place to attack is (\ref{eqn:T_V-grading-4}), by using that
    \begin{equation}
        \mu^2_\mathcal{T}(\Tilde{\eta}^d_\mathcal{T}(,\cdots,p_{\mathcal{L}_j}),\Phi^j\otimes e_V)=(-1)^{\lvert\Phi^j\rvert}\Tilde{\eta}^d_\mathcal{T}(,\cdots,\Phi^j\otimes e_V).
    \end{equation}
    This simplifies the fourth term (\ref{eqn:T_V-grading-4}) as
    \begin{equation} \label{eqn:T_V-grading-206}
        \sum\Tilde{\eta}_\mathcal{T}^{d-j+r+1}(x_d,\cdots,ev_{\mathcal{L}_j},\Phi^{s_r}\otimes e_V,\cdots,\Phi^{s_1}\otimes e_V),\quad r\geq2.
    \end{equation}
    Equation \eqref{eqn:T_V-grading-204} gives \eqref{eqn:cancelation-T_V}, while \eqref{eqn:T_V-grading-206} gives the terms of \eqref{eqn:cancellation-T_V} with at least two $\Phi$ blocks. It remains to show that the four remaining terms add up to supply the terms of \eqref{eqn:cancellation-T_V} with one $\Phi$ block. We again use
    \begin{equation}
        \sum(\delta\Tilde{\eta}_\mathcal{T})^{r+1}(T^{s_r},\cdots,T^{s_2},\Phi^{s_1}\otimes e_V,p_{\mathcal{L}_0})_{21}=0
    \end{equation}
    from strict unitality. First, the $\mu\circ\Tilde{\eta}$ terms are
    \begin{align}
        &\sum\mu_\mathcal{T}^{d-j+1}(x_d,\cdots,\Tilde{\eta}_\mathcal{T}^{j-q+r+2}(x_j,\cdots,ev_{\mathcal{L}_q},\Phi^{s_r}\otimes e_V,\cdots,\Phi^{s_1}\otimes e_V,p_{\mathcal{L}_0}))\\
        +&\sum\mu_\mathcal{T}^{d-j+2}(x_d,\cdots,ev_{\mathcal{L}_j},\Tilde{\eta}_\mathcal{T}^{r+1}(\Phi^{s_r}\otimes e_V,\cdots,\Phi^{s_1}\otimes e_V,p_{\mathcal{L}_0})),
    \end{align}
    and the second term further simplifies to
    \begin{equation}
        \sum(-1)^k\mu^{d-j+2}(x_d,\cdots,x_{j+1},ev_{L_j},(\Phi^{s_r}\circ\cdots\circ\Phi^{s_1})\otimes\Tilde{\eta}_\mathcal{T}^{r+1}(e_V,\cdots,e_V)),
    \end{equation}
    which thus accounts for (\ref{eqn:T_V-grading-7}). 

    Now for the $-\Tilde{\eta}\circ\mu$ terms, we first analyze the case when the $\mu$ is inserted to the left of evaluation. The terms to the right of $\mu$ are
    \begin{equation}
        x_i,\cdots,ev_{\mathcal{L}_q},\Phi^{s_r}\otimes e_V,\cdots,\Phi^{s_1}\otimes e_V,p_{\mathcal{L}_0},
    \end{equation}
    and therefore the sign contribution is
    \begin{equation}
        \dagger^i_{q+1}+\sum(\lvert\Phi^{s_t}\rvert-1)+\norm{p_{\mathcal{L}_0}}=\dagger_i+1.
    \end{equation}
    Therefore, the contribution is
    \begin{equation}
        \sum(-1)^{\dagger_i}\Tilde{\eta}_\mathcal{T}^{d-q-j+r+3}(x_d,\cdots,\mu^j_\mathcal{T}(x_{i+j},\cdots),x_i,\cdots,ev_{\mathcal{L}_q},\Phi^{s_r}\otimes e_V,\cdots,\Phi^{s_1}\otimes e_V,p_{\mathcal{L}_0}),
    \end{equation}
    which accounts for the fifth term (\ref{eqn:T_V-grading-8}).

    Finally, the remaining $\Tilde{\eta}\circ\mu$ terms insert $\mu$ among the $\Phi\otimes e_V$ terms, but by strict unitality only $\mu^1$ and $\mu^2$ survive by the formulae
    \begin{align}
        \mu^1_\mathcal{T}\begin{pmatrix}
            \Phi^s\otimes e_V &0 \\ 0 &0
        \end{pmatrix}
        &=\begin{pmatrix}
            (-1)^{\lvert \Phi^s\rvert}\partial_{Ch}\Phi^s\otimes e_V &0\\
            \mu^2_\mathcal{T}(ev_{L_s},\Phi^s\otimes e_V) &0
        \end{pmatrix}, \label{eqn:cancellation-2}\\
        \mu^2_\mathcal{T}(\Phi^{s_2}\otimes e_V,\Phi^{s_1}\otimes e_V)&=(-1)^{\lvert\Phi^{s_1}\rvert}\Phi^{s_2}\circ\Phi^{s_1}\otimes e_V.
    \end{align}
    These terms cancel the final term (\ref{eqn:T_V-grading-9}) by the $A_\infty$-functor equation for $\Phi$:
    \begin{align}
        &\sum(-1)^{\dagger_{s_1+\cdots+s_{j-1}}}\Tilde{\eta}_\mathcal{T}^{d-q+r+2}(x_d,\cdots,ev_{\mathcal{L}_q},\Phi^{s_r}\otimes e_V,\cdots,\\&((-1)^{\lvert\Phi^{s_j}\rvert}\partial_{Ch}\Phi^{s_j})+\sum_{a=1}^{s_j-1} (-1)^{\lvert\Phi^{a}\rvert}\Phi^{s_j-a}\circ\Phi^a)\otimes e_V,\cdots,p_{\mathcal{L}_0})\\
        =&\sum(-1)^{\dagger_i}\Tilde{\eta}^{d-q+r+2}_\mathcal{T}(x_d,\cdots,ev_{\mathcal{L}_q},\Phi^{s_r}\otimes e_V,\cdots,\\&\Phi^{s_j-a+1}(\cdots,\mu^a_\mathcal{T}(x_{i+a},\cdots),x_i,\cdots)\otimes e_V,\Phi^{s_{j-1}}\otimes e_V,\cdots,\Phi^{s_1}\otimes e_V,p_{\mathcal{L}_0}). 
    \end{align}
    The lower left entry of \eqref{eqn:cancellation-2} supplies the remaining one $\Phi$ block contribution to \eqref{eqn:cancellation-T_V}, with positive sign. This completes the lower left component of the desired identity.
    
    Therefore we have shown that our choice of $\xi_{T_V}$ provides an infinitesimal grading for the functor $T_{(V,\xi_V)}$.
    \end{proof}
In conclusion, the functor $T_V$ admits some infinitesimal grading, but the construction given above is not canonical. Nevertheless, the object-level grading is the specified grading \eqref{eqn:xi_TVL}, independent of the higher components of the functor grading. We summarize the conclusion below.

\begin{defn} \label{defn:T_V-canonical-grading-objects}
    For each infinitesimally graded object $(V,\xi_V)$ and $k\geq1$, $T_V^k$ admits an infinitesimal grading. For each infinitesimally graded object $(\mathcal{L},\xi_\mathcal{L})$, the induced object-level grading 
    \begin{equation}
        \xi^{can}_{T_V^k\mathcal{L}}
    \end{equation}
    is canonical: we call this the \emph{canonical object-level grading} associated to an algebraic twist. Here we equip iterated twists with the infinitesimal grading obtained by successive applications of the construction above.
\end{defn}

\subsection{Spherical objects and spherical algebraic twists}

The discussion in the previous sections used strict units to give explicit formulas for algebraic twists and their infinitesimal gradings. We now apply these constructions to spherical objects in a proper, c-unital, infinitesimally graded $A_\infty$-category. For the computations below, we fix a strictly unital model.

First pass to a minimal model of $\mathcal{C}$. By \cite[Lemma 2.1]{Seidel08}, as corrected in the erratum \cite{SeidelErratum}, there exists a formal diffeomorphism $\mathcal{C}\to\mathcal{C}'$ such that $\mathcal{C}'$ is strictly unital. Using the invariance results for infinitesimal gradings established above, equip the resulting model $(\mathcal C_{\min},\eta_{\min})$ with the transported infinitesimal grading and object gradings. Let
\begin{equation}
    \Phi:(\mathcal C_{\min},\eta_{\min})\to(\mathcal C,\eta)
\end{equation}
be the chosen quasi-equivalence preserving infinitesimal gradings. 

The proposition below records the grading choices for this model and the comparison with the twist functors for the original category.

\begin{proposition} \label{prop:transport-grading-unital}
    Let $(V,\xi_V)$ be an infinitesimally graded object of $\mathcal C_{\min}$. Then:
    \begin{enumerate}
        \item For every $k\geq 1$, the twist $T_V^k$ admits an infinitesimal grading whose induced object gradings are obtained by iterating the twist-grading construction of Section \ref{ssec:Tw_Ch}. These gradings can be chosen by fixing a grading as in Theorem \ref{thm:T_V-grading} and taking its composites as in Lemma \ref{lem:composition-of-functors}.
        \item The equivalence $\widehat\Phi$ of module categories induced by $\Phi$ admits a natural equivalence of underlying functors
        \begin{equation}
            \widehat\Phi\circ T_{\mathcal{Y}_V}^{k}\simeq T_{\mathcal{ Y}_{\Phi(V)}}^{k}\circ\widehat\Phi.
        \end{equation}
        As a consequence, the relations of the corresponding twist autoequivalences of each category are preserved under the equivalence between the two models.
    \end{enumerate}
\end{proposition}

Throughout this subsection, we work in the fixed minimal strictly unital model, denoted by $\mathcal{C}$. For each twist functor we fix the grading given by the proposition, and for its iterates we use the gradings obtained by composition. The cohomology level computations below only involve the resulting object gradings and induced grading operators on cohomology. In particular, no higher components of the chosen functor gradings are needed.

\subsubsection{Spherical objects and infinitesimal dimension}

\begin{defn} \label{defn:spherical-object}
    An object $V$ in a proper, cohomologically unital $\mathbb{Z}/2$-graded $A_\infty$-category $\mathcal{C}$ is \emph{spherical of $\mathbb{Z}/2$-dimension n} if the following conditions hold. All operations below are understood on cohomology.
    \begin{enumerate}
        \item There exists a homogeneous class $[V]$ of degree $n$ such that
        \begin{equation}
            H^*\hom_\mathcal{C}(V,V)=k\langle e_V\rangle\oplus k\langle [V]\rangle,\qquad \mu^2_\mathcal{C}([V],[V])=0
        \end{equation} 
        for the cohomological unit $e_V$.
        \item The linear functional
        \begin{equation}
            \int_V:H^*\hom_{\mathcal C}(V,V)\to k, \qquad \int_V e_V=0,\quad \int_V[V]=1,
        \end{equation}
        has degree $-n$.
        \item For every object $L\in\mathcal{C}$ and $i\in\mathbb{Z}/2$, the pairing
        \begin{equation}
            H^{n-i}\hom_\mathcal{C}(L,V)\otimes H^i\hom_\mathcal{C}(V,L)\xrightarrow{\ \mu^2\ }H^n\operatorname{hom}_{\mathcal C}(V,V)\xrightarrow{\ \int_V\ } k
        \end{equation}
        is nondegenerate.
    \end{enumerate}
\end{defn}

Recall that the algebraic twist along a spherical object admits an \emph{adjoint twist}, as defined in \cite[Section 5g]{Seidel08}. In particular, when $V$ is spherical, $T_V$ is invertible.

We define the infinitesimal dimension of a spherical object.

\begin{defn}
    An infinitesimally graded spherical object $(V,\xi_V)$ has \emph{infinitesimal dimension $n$} if
    \begin{equation}
        [\Gamma_{V,V}([V])]=n\cdot[V]
    \end{equation}
    holds in $H^*\hom_\mathcal{C}(V,V)$.
\end{defn}

For example, spherical objects of odd $\mathbb{Z}/2$-dimension always have well-defined infinitesimal dimension for degree reasons. For spherical objects of even $\mathbb{Z}/2$-dimension the following argument shows that their infinitesimal dimension is well-defined. Write
\begin{equation}
    \Gamma_{V,V}([V])=a\cdot e_V+\lambda_V\cdot[V].
\end{equation}
Then, by the derivation property of $\Gamma_{V,V}$, we have
\begin{equation}
    2a[V]=\mu^2(\Gamma_{V,V}([V]),[V])+\mu^2([V],\Gamma_{V,V}([V]))=\Gamma_{V,V}(\mu^2([V],[V]))=0
\end{equation}
on cohomology, where the last equality is from Definition \ref{defn:spherical-object}. Therefore over a base field $k$ of characteristic not equal to $2$, we conclude that $a=0$, and therefore the infinitesimal dimension of a spherical object is always well-defined.

We will see that for monotone Fukaya categories, the open-closed map together with the algebra structure of the corresponding summand of quantum cohomology can be used to compute the infinitesimal dimensions of Lagrangian spheres.

\begin{lemma} \label{lem:spherical-dim-shift}
    For an infinitesimally graded spherical object $(V,\xi_V)$ of infinitesimal dimension $\lambda_V$, the following diagram commutes:
    \begin{equation} \label{eqn:comm-diag}
    \begin{tikzcd}
        H^*\hom(V,\mathcal{L}) \arrow[d,"\iota_*"] \arrow[rr,"\Gamma_{V,\mathcal{L}}+(\lambda_V-1)\cdot\id"] &&H^*\hom(V,\mathcal{L}) \arrow[d,"\iota_*"]\\
        H^*\hom(V,T_V\mathcal{L}) \arrow[rr,"\Gamma_{V,T_V\mathcal{L}}"']  && H^*\hom(V,T_V\mathcal{L})
    \end{tikzcd},    
    \end{equation}
    where the map $\iota:\hom(V,\mathcal{L})\to\hom(V,T_V\mathcal{L})$ with parity $n-1$ sends
    \begin{equation}
        \iota(x)=sx\otimes [V]+(-1)^{n-1}s\mu^2(x,[V])\otimes e_V.
    \end{equation}
    In words, under $\iota_*$, the weights on $H^*\hom(V,T_VL)$ are those on $H^*\hom(V,L)$ increased by $\lambda_V-1$.
\end{lemma}
\begin{proof}
    We first verify that for a cocycle $x$, $\iota(x)$ is a cocycle. We compute
    \begin{equation}
        \mu^1_{\Tw_{\Ch}}(\iota(x))=(0,(-1)^{\norm{x}}\mu^2(x,[V])+(-1)^{n-1}(-1)^{\lvert x\rvert+n-1}\mu^2(x,[V])),
    \end{equation}
    where the first term is zero because every term in $\iota(x)$ is a cocycle.

    Now using Proposition \ref{prop:transport-grading-unital}, we move to a minimal strictly unital model for $\mathcal{C}$. Then there exists a filtration
    \begin{equation}
        \hom(V,T_VL)=(\hom(V,L)[1]\otimes e_V)\oplus \hom(V,L)\oplus\hom(V,L)[1]\otimes[V],
    \end{equation}
    where the first two summands form an acyclic pair. Thus, the induced complex on the $\hom(V,L)[1]\otimes[V]$ summand can be identified to be quasi-isomorphic to $\hom(V,T_VL)$, and this identifies $\iota(x)$ in $\hom(V,T_VL)$ with $sx\otimes[V]$.

    Therefore we may compute the grading operator on this chain model. Then for $x\in\hom(V,\mathcal{L})$, using that
    \begin{equation}
    \xi_{T_V\mathcal{L}=}\begin{pmatrix}
        \id\otimes\xi_V+\psi_{\hom(V,\mathcal{L})[1]}\otimes e_V &0\\
        0 &\xi_\mathcal{L}
    \end{pmatrix},    
    \end{equation}
    the summand $E_L=(\hom(V,L)[1]\otimes e_V)\oplus \hom(V,L)$ is acyclic and preserved by the grading operator. Thus on the quotient $\hom(V,T_VL)/E_L$, we compute for the induced operator on the quotient
    \begin{equation}
        \Bar{\Gamma}(sx\otimes[V])=sx\otimes \Gamma_{V,V}([V])+s(\Gamma_{V,L}(x)-x)\otimes[V].
    \end{equation}
    Since we assume that $V$ has infinitesimal dimension $\lambda_V$, we conclude that on cohomology
    \begin{equation}
        [\overline{\Gamma}_{V,T_V\mathcal{L}}(sx\otimes [V])]=(\lambda_V-1)[sx\otimes[V]]+s[\Gamma_{V,\mathcal{L}}(x)\otimes[V]].
    \end{equation}
    Here all classes are taken in $H^*(\hom(V,T_VL)/E_L)$. Under the quotient quasi-isomorphism, $[\iota(x)]$ maps to $[sx\otimes[V]]$, so this identity proves the commutativity of \eqref{eqn:comm-diag}.
\end{proof}

When $K\in\mathcal{C}$ pairs cohomologically trivially to $V$, the algebraic twist preserves the grading operator.

\begin{lemma} \label{lem:acyclic-twist-preserves-grading}
    Let $(V,\xi_V)$ be an infinitesimally graded object in $(\mathcal{C},\eta)$, and let $(K,\xi_K)$ be an infinitesimally graded object such that $H^*\hom(K,V)=0$. Then for any infinitesimally graded twisted complex $\mathcal{L}$, the isomorphism
    \begin{equation}
        \iota_*:H^*\hom(K,\mathcal{L})\xhookrightarrow{}H^{*}\hom(K,T_V\mathcal{L})
    \end{equation}
    preserves the grading operators on both sides.
\end{lemma}
\begin{proof}
    Since we assume that $H^*\hom(K,V)=0$, the inclusion
    \begin{equation}
        \iota:\hom(K,\mathcal{L})\xhookrightarrow{}\hom(V,\mathcal{L})[1]\otimes\hom(K,V)\oplus\hom(K,\mathcal{L)}\cong \hom(K,T_V\mathcal{L})
    \end{equation}
    is an isomorphism on cohomology. Since the cohomology-level grading operator is independent of the chain representative of the cohomology class, it is enough to compute it for the classes in the image of $\iota$. Then
    \begin{align}
        [\Gamma_{K,T_V\mathcal{L}}(0,x)]&=[\eta_\mathcal{T}^1(0,x)-\mu^2_{\Tw_{\Ch}}((0,x),\xi_K)-(-1)^{\norm{x}}\mu^2_{\Tw_{\Ch}}(\xi_{T_V\mathcal{L}},(0,x))]\\
        &=[\eta_\mathcal{T}^1(x)-\mu^2_{\Tw_{\Ch}}(x,\xi_K)-(-1)^{\norm{x}}\mu^2(\xi_{\mathcal{L}},x)]\\
        &=[\Gamma_{K,\mathcal{L}}(x)]. \qedhere
    \end{align}
\end{proof}

Therefore, using Lemma \ref{lem:composition-of-functors} we arrive at the following conclusion.

\begin{proposition} \label{prop:inf-grading-T_V^k}
    Let $(V,\xi_V)$ be an infinitesimally graded spherical object of infinitesimal dimension $\lambda_V$, and $(\mathcal{L},\xi_\mathcal{L}),(K,\xi_K)$ be infinitesimally graded objects such that $H^*\hom(K,V)=0$. Then for every $k\geq1$, under the quasi-isomorphisms
    \begin{equation}
        \iota^k_{V,T_V^kL}:H^*\hom(V,\mathcal{L})\cong H^{*+k(n-1)}\hom(V,T_V^k\mathcal{L}),
    \end{equation}
    \begin{equation}
        \Gamma_{V,T_V^k\mathcal{L}}\circ\iota^k=\iota^k\circ(\Gamma_{V,\mathcal{L}}+k\cdot(\lambda_V-1)\id),
    \end{equation}
    and
    \begin{equation}
        \Gamma_{K,{T_V^k}_{\mathcal{L}}}=\Gamma_{K,\mathcal{L}}.
    \end{equation}
\end{proposition}

This will be our tool to show the infinitude of spherical twists in the next subsection.

\subsubsection{Orders of spherical twists and faithful braid group actions}

This is the final subsection of our algebraic foundations of infinitesimal gradings, where we use the previously developed machinery to show that spherical twists have infinite order, and moreover generate faithful braid group actions.

We first need a small lemma.

\begin{lemma} \label{lem:spherical-mult-trivial}
    For a spherical object $L$ and an object $K$ such that $H^*\hom(K,L)$ is 1-dimensional, the map
    \begin{equation}
        [\mu^2([L],~)]:H^*\hom(K,L)\to H^*\hom(K,L)
    \end{equation}
    is trivial.
\end{lemma}
\begin{proof}
    For homogeneous $x\in H^*\hom_\mathcal{C}(K,L)$, we have in cohomology
    \begin{equation}
        \mu^2([L],\mu^2([L],x))=(-1)^{\lvert x\rvert}\mu^2(\mu^2([L],[L]),x)=0,
    \end{equation}
    where the second equality follows from Definition \ref{defn:spherical-object}. Since $H^*\hom_\mathcal{C}(K,L)$ is one-dimensional over $k$, we have that $[\mu^2([L],~)]$ is multiplication by a scalar $c$, and by the above equality we conclude $c=0$, equivalently
    \begin{equation}
        \mu^2([L],x)=0
    \end{equation}
    in cohomology.
\end{proof}

We finally prove our main theorem that spherical twists along infinitesimally graded spherical objects have infinite order, even in a $\mathbb{Z}/2$-graded $A_\infty$-category.

\begin{thm} \label{thm:dehn-twist-reprise}
    Let $K,L,V$ be an $A_3$-configuration (in this order) of infinitesimally graded spherical objects in $\mathcal{C}$ defined over a characteristic zero field $k$, and suppose that $V$ has well-defined infinitesimal dimension $\lambda_V\neq1$. Then 
    \begin{equation}
        T_V^m\not\simeq\id
    \end{equation}
    for every $m\neq0$.
\end{thm}
\begin{proof}
    Assume the existence of such $K,L,V$ and $m$ as above. Since $T_V$ admits an inverse, we may assume that $m>0$. Let $\gamma^{(m)}$ be the grading on $T_V^m$ obtained by composing the fixed grading on $T_V$, as in Proposition \ref{prop:transport-grading-unital}. Then, by transporting $\gamma^{(m)}$ across the assumed natural transformation
    \begin{equation}
        T:T_V^m\simeq\id,
    \end{equation}
    we obtain an infinitesimal grading $\gamma$ of the identity functor such that $T$ is a graded isomorphism of functors. Although the infinitesimal grading $\gamma^{(m)}$ is not canonical, its induced infinitesimal grading $\xi^{can}_{T_V^m\mathcal{L}}$ on objects $T_V^m\mathcal{L}$ is canonical. Thus, the computations of the cohomological grading operators below will be independent of the particular choice of $\gamma^{(m)}$ chosen above.

    Now, such a grading $\gamma$ of the identity functor assigns cochains $\gamma^0_K$, $\gamma^0_L$, and $\gamma^0_V$ to each corresponding spherical object: we will show that such a choice is impossible.

    First, note that $\gamma^0_V$ and $\gamma^0_L$ are cocycles in the endomorphism algebra of a spherical object, and thus we may write
    \begin{equation}
        [\gamma^0_L]=a_L\cdot [e_L]+b_L\cdot[L].
    \end{equation}
    Since $H^*\hom(K,L)$ and $H^*\hom(V,L)$ are one-dimensional, by Lemma \ref{lem:spherical-mult-trivial}, multiplication by the $b_L[L]$-part acts trivially on both spaces. 

    Recall from Proposition \ref{prop:isomorphic-functors-same-grading} that the induced gradings on the cohomological-level morphism spaces
    \begin{equation}
        H^*\hom_{\Tw_{\Ch}\mathcal{C}}(X,T_V^mY),\quad H^*\hom_\mathcal{C}(X,Y)
    \end{equation}
    each induced from $(T_V^m,\gamma^{(m)})$ and $(\id,\gamma)$ will agree, up to a multiplication by the natural transformation element $T^0_Y$. In particular, on cohomologically one-dimensional morphism spaces, this will be a multiplication by a nonzero scalar, and thus cancels in the comparison of grading operators.

    We first look at the morphism space $\hom(K,T_V^mL)$. Since $\hom_\mathcal{C}(K,V)$ is acyclic, Proposition \ref{prop:inf-grading-T_V^k} shows that the induced grading operator on cohomology equals
    \begin{equation}
        [\Gamma_{K,L}(x)]=[\eta^1(x)-\mu^2(x,\xi_K)-(-1)^{\norm{x}}\mu^2(\xi_L+\gamma^0_L,x)]
    \end{equation}
    through the isomorphism $H^*\hom(K,T_V^mL)\cong H^*\hom(K,L)$. Therefore we conclude that the shift by $\mu^2([\gamma^0_L],~)$ on the 1-dimensional morphism space $H^*\hom(K,L)$ is trivial. By Lemma \ref{lem:spherical-mult-trivial}, multiplication by the $b_L[L]$ part acts trivially, and thus we conclude that $a_L=0$.

    Now look at the morphism space
    \begin{equation}
        \hom_{\Tw_{\Ch}\mathcal{C}}(V,T_V^mL).
    \end{equation}
    By Proposition \ref{prop:inf-grading-T_V^k}, we know that the induced grading operator on cohomology equals
    \begin{equation}
        [\Gamma_{V,L}]+m\cdot(\lambda_V-1)\id
    \end{equation}
    for the infinitesimal dimension $\lambda_V$. On the other hand, since $a_L=0$, we conclude that the induced grading operator on the morphism space $H^*\hom(V,L)$ given by the graded functor $(\id,\gamma)$ is
    \begin{equation}
        [\Gamma_{V,L}].
    \end{equation}
    Since multiplication by $[T^0_L]$ intertwines the grading operators of
    \begin{equation}
        H^*\hom(V,T_V^mL)\cong H^*\hom(V,L),
    \end{equation}
    we conclude that
    \begin{equation}
        [\Gamma_{V,L}]+m\cdot(\lambda_V-1)\id=[\Gamma_{V,L}],
    \end{equation}
    and therefore
    \begin{equation}
        m\cdot(\lambda_V-1)=0.
    \end{equation}
    But since we assume that $\lambda_V\neq 1$ and the characteristic of the base field is zero, this provides a contradiction.
\end{proof}

We now move on to the second application of our infinitesimal gradings, constructing braid group actions. First we set up notation.

\begin{defn}\label{def:braid-relations}
    The \emph{braid group on $m+1$ strands} $B_{m+1}$ is generated by the elements $\sigma_1,\cdots,\sigma_m$ that satisfy the braid relations:
    \begin{align}
        \sigma_i\sigma_j&=\sigma_j\sigma_i\qquad\quad\quad(|i-j|>1),\\
        \sigma_i\sigma_{i+1}\sigma_i&=\sigma_{i+1}\sigma_i\sigma_{i+1}\quad(1\leq i\leq m-1).
    \end{align}
    In general for a finite simple graph $\Gamma$, the associated \emph{Artin group} $B_\Gamma$ is generated by $\{\sigma_v,v\in V(\Gamma)\}$ with the relations
    \begin{align}
        \sigma_v\sigma_w&=\sigma_w\sigma_v &&(vw\not\in E(\Gamma)),\\
    \sigma_v\sigma_w\sigma_v&=\sigma_w\sigma_v\sigma_w &&(vw\in E(\Gamma)),
    \end{align}
    where $E(\Gamma)$ is the edge set. In particular, $B(A_m)$ is the braid group $B_{m+1}$. 
    
    For $I\subset V(\Gamma)$,  we call the subgroup
    \begin{equation}
        B_I=\langle \sigma_v\mid v\in I\rangle\subseteq B_\Gamma
    \end{equation}
    a \emph{standard parabolic subgroup}, and its conjugates \emph{parabolic subgroups}.
\end{defn}

This Artin group will act as $A_\infty$-autoequivalences on the $A_\infty$-category of right $\mathcal{C}$-modules.

\begin{defn}
    The \emph{derived autoequivalence group} $\pi_0\Aut(\mathcal{C})$ for an $A_\infty$-category $\mathcal{C}$ is the group of isomorphism classes of $A_\infty$-autoequivalences of $\mathcal{C}$, up to equivalence by invertible degree zero natural transformations.
\end{defn}

Recall the graded algebra associated to a $\Gamma$-configuration.

\begin{defn} 
    Let $\Gamma$ be a connected ADE Dynkin diagram and $n\in\mathbb{Z}$. Let $\overline{\Gamma}$ be the doubled quiver, with one arrow in each direction for each edge. Assign opposite arrows integer degrees $d,n-d$ that sum to $n$, and the idempotents associated to each vertex degree zero.

    The \emph{graded zigzag algebra} associated to $\Gamma$ is given by
    \begin{equation}
        A_{\Gamma,n}=k\overline{\Gamma}/J_\Gamma
    \end{equation}
    for the ideal $J_\Gamma$ generated by
    \begin{align}
        &(v\vert w\vert u)\quad\text{for }v\neq u,\\
        &(v\vert w\vert v)-(v\vert u\vert v)\quad \text{for }u,w\text{ adjacent to }v,\\
        &\text{all paths of length }3.
    \end{align}
    For $\Gamma=A_1$, we define $A_{A_1,n}=k[x]/(x^2)$ with $\lvert x\rvert=n$.
\end{defn}

When $\Gamma=A_m$, this recovers the quiver $\Gamma_{m,n}$ and algebra $A_{m,n}$ of \cite[Section~4.3]{SeidelThomas}. The choice of arrow degrees is suppressed from the notation, since the properties used below are independent of this choice. Intrinsic formality in type $A$ for $n\geq2$ was proved by Seidel--Thomas \cite[Lemma~4.21]{SeidelThomas}. For $n=2$, Etg\"u--Lekili show the type $D$ formality statement in \cite[Corollary 45]{EtguLekili}. The remaining cases $E_6,E_7,E_8$ in characteristic zero follow from the computations of Lekili--Ueda \cite[Sections 5.4--5.6]{LekiliUeda}.

The following lemma extends the statement to every integer $n\neq1$. We remark that when $n=1$, \cite[Section 4.5]{SeidelThomas} notes that $A_{3,1}$ is not intrinsically formal.

\begin{lemma} \label{lem:A_m,n-intrinsic-formal}
    Let $\Gamma$ and $A_{\Gamma,n}$ be as above. Over a field of characteristic zero, the algebra $A_{\Gamma,n}$ is intrinsically formal whenever $n\in\mathbb{Z}\setminus\{1\}$. Moreover, the formality equivalence can be chosen to preserve the original $\mathbb{Z}/2$-grading.
\end{lemma}
\begin{proof}
    For $\Gamma=A_1$, a degree argument excludes any nontrivial higher operations $\mu^r,r\geq 3$. Excluding this case, we will imitate the proof of \cite[Lemma 4.21]{SeidelThomas} and compute Hochschild cohomology. We will write $A=A_{\Gamma,n}$, and our base ring $R$ will be given as
    \begin{equation}
        R=\oplus_{v} k e_v,\quad v\in V(\Gamma).
    \end{equation}
    Let $A^+$ be the ideal of $A_{\Gamma,n}$ spanned by paths of positive length. The $r$th term of the Hochschild cochain complex is given by
    \begin{equation}
        CC^r(A,A[2-r])=\hom^0_{R-R}((A^+)^{\otimes_R r},A[2-r]).
    \end{equation}

    Because $\Gamma$ is a tree, we may regrade the degrees such that each arrow has degree $n/2$ when $n$ is even. When $n$ is odd, we assign each arrow and its opposite the degrees
    \begin{equation}
        \lfloor\frac{n}{2}\rfloor,\quad\lceil\frac{n}{2}\rceil
    \end{equation}
    such that the degrees of two composable arrows add up to $n$. Then a nonzero cochain in $CC^r(A,A[2-r])$ with total input length $l$ and output length $j$ must satisfy
    \begin{equation}
        \frac{n}{2}(l-j)=r-2.
    \end{equation}
    Moreover, by positive length we know that $l\geq r$, and since all length $\geq3$ paths are trivial in $A$ we know that $j\leq 2$. Therefore, for degree reasons, the relevant group $CC^r(A,A[2-r])$ vanishes when $r>2$ and $n\neq1,2$, and intrinsic formality follows from the same argument as \cite[Theorem 4.7]{SeidelThomas}.

    For $n=2$, intrinsic formality follows from the three references cited before the statement of the lemma.

    We now prove the stronger assertion that the formality equivalence can be chosen to preserve the original $\mathbb{Z}/2$ gradings. When $n\neq 2$, the preceding weight argument excludes all higher operations, so the assertion is immediate.

    For $n=2$, regrade such that each arrow has weight $1$. The only possibility of a nontrivial $\mu^d$ operation comes from one that takes $d$ composable arrows to a spherical class. Thus each $\mu^d$ has the correct parity. Now when choosing the sequence of $A_\infty$-diffeomorphisms eliminating $\mu^d$ for each $d$ inductively, we may replace each eliminating Hochschild cochain by its odd degree component, since the weight-graded obstruction equation respects the original parity decomposition. The resulting $A_\infty$-diffeomorphism will preserve the $\mathbb{Z}/2$-grading, giving the desired equivalence.
\end{proof}

With this, we can prove that the kernel of the Artin group action on $\mathcal{C}\text{-}\mathrm{mod}$ is central.

\begin{thm} \label{thm:central-kernel}
    Let $\Gamma$ be a connected Dynkin diagram of type $A$, $D$, or $E$, and let $(\mathcal{C},\eta)$ be a proper infinitesimally graded $A_\infty$-category over an algebraically closed field $k$ of characteristic zero. Suppose that $\mathcal{C}$ contains a spherical $\Gamma$-configuration $\{V_v\}_{v\in V(\Gamma)}$ of infinitesimally graded objects, all of infinitesimal dimension $n\in\mathbb{Z}\setminus\{1\}$.

    Then the spherical twists define a homomorphism
    \begin{equation}
        \rho_\Gamma:B_\Gamma\to\pi_0\bigl(\mathrm{Aut}(\mathcal{C}\text{-}\mathrm{mod})\bigr)
    \end{equation}
    sending the standard generator $\sigma_v$ to $[T_{V_v}]$. Its kernel is contained in the center:
    \begin{equation}
        \ker\rho_\Gamma\subseteq Z(B_\Gamma).
    \end{equation}
\end{thm}
\begin{proof}
    First, define $\rho:B_\Gamma\to\pi_0\bigl(\mathrm{Aut}(\mathcal{C}\text{-}\mathrm{mod})\bigr)$ to send the generator $\sigma_i$ to the twist functor $T_{\mathcal{Y}_{V_i}}$, associated to the Yoneda image of $V_i$. To show that this is well-defined, one needs to check that the twist functors $T_{\mathcal{Y}_{V_i}}$ satisfy the relations in Definition \ref{def:braid-relations}: this follows from the arguments of \cite[Propositions 2.12, 2.13]{SeidelThomas}.

    Let 
    \begin{equation}
        \mathcal{D}=\langle \mathcal{Y}_{V_1},\cdots, \mathcal{Y}_{V_m}\rangle\subset\mathcal{C}\text{-}\mathrm{mod}
    \end{equation}
    be the split-closed triangulated subcategory generated by these objects. The twists and their inverses preserve $\mathcal{D}$, giving a restricted action
    \begin{equation}
        \overline{\rho}:B_\Gamma\to \pi_0\Aut(\mathcal{D}).
    \end{equation}
    It suffices to show that $\ker\overline{\rho}$ is contained in the center $Z(B_\Gamma)$.
    
    Define the $A_\infty$-algebra
    \begin{equation}
        \mathcal{A}=\bigoplus_{i,j\in V(\Gamma)}\hom_{\mathcal{C}}(V_i,V_j).
    \end{equation}
    Morita invariance gives an equivalence of $\mathbb{Z}/2$-graded categories
    \begin{equation}
        \mathcal{D}\cong\mathrm{Perf}^{\mathbb{Z}/2}(\mathcal{A}),
    \end{equation}
    where the right-hand side denotes the category of perfect supergraded right $A_\infty$-modules over $\mathcal{A}$. This sends each $\mathcal{Y}_{V_j}$ to $e_j\mathcal{A}=\bigoplus_i \hom_\mathcal{C}(V_i,V_j)$, or $\mathcal{P}_j$ in the notation of \cite{SeidelThomas}. Moreover, under this identification, the twist functor $T_{\mathcal{Y}_{V_i}}$ is identified with the twist functor $t_i$ on $\mathrm{Perf}^{\mathbb{Z}/2}(\mathcal{A})$. Therefore, we obtain the homomorphism
    \begin{equation}
        \Tilde{\rho}:B_\Gamma\to \pi_0\Aut(\mathrm{Perf}^{\mathbb{Z}/2}\mathcal{A}),
    \end{equation}
    and it is enough for us to show that its kernel is contained in $Z(B_\Gamma)$.

    To do so we pass to the graded lift $\Tilde{\mathcal{C}}_\eta$, and consider its full subcategory consisting of all lifts of the objects $V_v$. Since $\Gamma$ is a tree and each $V_v$ has integer infinitesimal dimension, we may choose a lift $(V_v,\xi_{V_v})$ for each vertex $v$ such that the induced grading operators on all cohomological morphism spaces have integer eigenvalues. We now restrict to the full subcategory of these chosen lifts, pass to a minimal model, and apply Theorem \ref{thm:straighten-infinitesimal-grading}.

    With these choices, the grading operator on $\mathcal{A}$ has integer eigenvalues. Since $k$ has characteristic zero, the compatibility with these weights gives an integer weight grading for which $\mu^d$ has degree $2-d$. Thus the cohomology algebra $H^*(\mathcal{A})$, equipped with its integer weights, is isomorphic to the zigzag algebra $A_{\Gamma,n}$ with corresponding degrees assigned to each path.

    Therefore by Lemma \ref{lem:A_m,n-intrinsic-formal}, we conclude that $\mathcal{A}$ is formal as a $\mathbb{Z}/2$-graded algebra. Thus,
    \begin{equation}\label{eqn:morita}
        \mathrm{Perf}^{\mathbb{Z}/2}(\mathcal{A})\cong \mathrm{Perf}^{\mathbb{Z}/2}(A_{\Gamma,n})
    \end{equation}
    for the chosen degrees assigned to paths and the $\mathbb{Z}/2$-degree inherited from $\mathcal{C}$. Under this equivalence, each $\mathcal{P}_i$ corresponds to $e_i A_{\Gamma,n}$, and $t_i$ corresponds to the twist functor along that projective.

    We clarify that $A_{\Gamma,n}$ need not have superdegrees induced from its integer degrees. For the following argument, we choose an auxiliary integer lift of its inherited $\mathbb{Z}/2$-grading from the equivalence \eqref{eqn:morita}, identifying it with some $A_{\Gamma,N}$ for some integer $N$ equal to either $2$ or $3$ with the same parity as the spherical classes $[V_i]$.
    
    Now take an element $\beta\in\ker\Tilde{\rho}$, and let $W_\beta$ be the corresponding word generated by the graded twists. We remark that $\beta$ is a relation between the twist functors before the graded lift: therefore the reduction of $W_\beta$ mod $2$ is isomorphic to the identity functor. 
    
    Its action on each projective shows that
    \begin{equation}
        W_\beta(P_i)\cong P_i[2r_i]
    \end{equation}
    for some $r_i\in\mathbb{Z}$. In particular, conjugation by $W_\beta$ preserves each standard graded twist. By faithfulness of Artin group actions generated by graded spherical configurations \cite[Theorem 1]{NordskovaVolkov}, it follows that $\beta$ commutes with every standard generator $\sigma_i$ of $B_\Gamma$, and therefore $\beta\in Z(B_\Gamma)$.
\end{proof}

\begin{thm}\label{thm:ADE-faithfulness-reprise}
    Under the hypotheses of Theorem \ref{thm:central-kernel}, the restriction of $\rho_\Gamma$ to every proper parabolic subgroup of $B_\Gamma$ is faithful. In particular, if $I\subsetneq V(\Gamma)$ and $\Gamma'$ is the full subdiagram on $I$, then the twists $\{T_{V_v}\mid v\in I\}$ define a faithful action
    \begin{equation}
        B_{\Gamma'}\to \pi_0\operatorname{Aut}(\mathcal C\text{-}\mathrm{mod}).
    \end{equation}
\end{thm}
\begin{proof}
    Identify $B_{\Gamma'}$ with the standard parabolic subgroup generated by $\{\sigma_v\mid v\in I\}$. Theorem \ref{thm:central-kernel} shows that
    \begin{equation}
        \ker\left(\rho_\Gamma\vert_{B_\Gamma'}\right)\subseteq B_{\Gamma'}\cap Z(B_\Gamma).
    \end{equation}
    Since a proper parabolic subgroup of an irreducible Artin group of finite type has trivial intersection with the center \cite[Corollary 2.6]{Paris}, the restricted action is faithful. The same conclusion holds for conjugate parabolic subgroups, since $\ker\rho_\Gamma$ is normal.
\end{proof}

This theorem establishes faithfulness for Artin group actions induced by spherical configurations of infinitesimally graded objects that can be extended to a larger configuration satisfying the same hypotheses. In particular, the above theorem does not cover faithfulness for the full $E_8$ action.

We will see in the next section that for summands of monotone Fukaya categories, the algebra structure of the corresponding summand in quantum cohomology, together with the homology classes of the Lagrangian spheres, can determine the infinitesimal dimensions of these spheres.

\section{Geometry}
\subsection{Monotone Fukaya categories}

We first recall the definition of the monotone Fukaya category with its open--closed and closed--open string maps, following \cite[Section 2]{Sheridan} and \cite[Section 5]{RitterSmith}. We work over a base field $k$ throughout.

Let $(X,\omega)$ be a closed symplectic manifold of dimension $2n$. We assume that $X$ is monotone, normalizing the symplectic form as
\begin{equation}
    [\omega]=c_1(X)\in H^2(X;\mathbb{R}),
\end{equation}
and that $X$ carries a smooth anticanonical divisor $D$. The Lagrangians $L$ we consider will be assumed to be closed, connected, oriented Lagrangians equipped with spin structures, such that $\pi_1(L)\to\pi_1(X)$ is trivial. We additionally require that $L$ is monotone, i.e. satisfies
\begin{equation}
    [\omega]=\frac{1}{2}\mu_L\in H^2(X,L;\mathbb{R})
\end{equation}
for the Maslov class $\mu_L:H_2(X,L)\to\mathbb{Z}$. To this Lagrangian, we associate the datum of an almost complex structure $J_L$ such that the moduli space $\mathcal{M}(L)$ of $J_L$-holomorphic disks of Maslov index 2 with boundary on $L$ and one boundary marked point is regular. We say that $W(L)=\lambda\in k$ if
\begin{equation}\label{eqn:J_L}
    ev_*[\mathcal{M}(L)]=\lambda\cdot[L]\in H_n(L)
\end{equation}
for some $\lambda\in k$. In general, we equip $L$ with the additional datum of a rank $1$ $k$-local system $\mathcal{E}$, and take the evaluation map twisted by the local system, to define a monotone Lagrangian pair $(L,\mathcal{E})$.

Having defined the objects of the monotone Fukaya category, we explain how to define the morphisms. For each pair of Lagrangians $(L_0,L_1)$ satisfying $W(L_0)=W(L_1)=\lambda$, we fix a regular Floer datum $(H_t,J_t)_{t\in[0,1]}$ as in \cite[Section 2.3]{Sheridan} such that $J_0=J_{L_0}$ and $J_1=J_{L_1}$. The Floer cochain complex $CF^*(L_0,L_1)$ is defined as a $\mathbb{Z}/2$-graded complex generated by time 1 Hamiltonian chords of $H_t$, and its differential counts solutions to the Floer equation for $(H_t,J_t)$ with Lagrangian boundary conditions given by $L_0,L_1$. Similarly, one can count polygons with coherently chosen regular perturbation data, agreeing with the fixed Floer data at the strip-like ends and extending the prescribed data near boundary strata, to define the $A_\infty$-operations.

\begin{defn}
    The \emph{monotone Fukaya category summand} $\mathcal{F}(X)_\lambda$ is a $\mathbb{Z}/2$-graded uncurved $A_\infty$-category with objects given by Lagrangians $L$ with $W(L)=\lambda$, and $A_\infty$-structure defined by counting holomorphic polygons with appropriate Lagrangian boundary conditions. This uncurved $A_\infty$-category $\mathcal{F}(X)_\lambda$ is cohomologically unital, chain-level proper, and has base field $k$.
\end{defn}

We use a Morse function $h_L:L\to\mathbb{R}$ to define $CF^*(L,L)$, whose generators have indices between $0$ and $n$. With this choice, a chain-level representative $e_L$ for the geometric unit can be defined from the moduli of discs with one output marked point with boundary on $L$.

To construct the open--closed and closed--open maps, we first define the chain model for $QH^*(X)$. Fix a Morse-Smale pair $(f,g)$ on $X$ with a single minimum, and define its Morse cochain complex by $CM^*(X)$. The differential counts gradient flowlines, and the quantum product is defined by holomorphic spheres with three marked points, two of which are connected to half-lines $(-\infty,0]$ and one to $[0,\infty)$. The moduli space of coupled Morse and Floer solutions for this configuration can be chosen to be regular, and these counts induce the quantum product on cohomology.

The unit $1_X\in CM^*(X)$ is represented by the minimum: we claim that one can choose small perturbation data such that
\begin{equation}\label{eqn:strict-unit-chain}
    1_X*x=x
\end{equation}
holds at chain level for any $x\in CM^*(X)$. First, assume that the domain-dependent perturbation for the almost complex structure on the sphere is small. When the holomorphic sphere carries a constant almost complex structure, since the intersection condition with the stable manifold of the minimum is an open condition, that intersection condition may be forgotten. As such, the resulting family admits a $\mathbb{C}^*$-rescaling action, and thus cannot be rigid if it represents a nontrivial homology class. This chain level triviality continues to hold for domain-dependent almost complex structures whose perturbation is small.

Thus we may assume that only constant spheres contribute to the product $1_X*~$. Assume that a rigid $Y$-shaped configuration has left input $1_X$, right input $x$ and output $y$. When only the left gradient flowline equation is perturbed, and the other two legs satisfy the gradient flowline equation for the chosen Morse--Smale pair $(f,g)$, then it admits a rescaling $\mathbb{R}$-action moving the intersection point with the remaining leg. Thus the configuration cannot be rigid unless it is constant: thus there are no rigid configurations when $x\neq y$, and exactly one rigid configuration when $x=y$. This persists for small perturbations. Hence we may assume that \eqref{eqn:strict-unit-chain} holds at chain level for our chosen datum. A similar argument shows that one can choose coherent perturbation data for the \emph{quantum Massey products} $m_3$ such that
\begin{equation}
    m_3(1_X,1_X,x)=0,\quad m_3(1_X,D,x)=0,\quad \forall x\in CM^*(X).
\end{equation}

For $CO_\lambda^d$, the domain is defined as the disc with $d$ input boundary marked points, $1$ output boundary marked point, and $1$ interior marked point joined to a half-line $(-\infty,0]$. Similarly for $OC_\lambda$, the domain has input boundary marked points, one of which is distinguished, and an interior marked point attached to a half-line $[0,\infty)$. One chooses coherent Floer and Morse perturbation data to ensure transversality of the relevant moduli spaces and coherence with the Hochschild (co)chain models. We refer to \cite[Sections 5.4, 5.5]{RitterSmith} for the precise construction. For the proofs of these properties, see \cite[Propositions 2.1, 2.2 and Lemmas 2.3, 2.4]{Sheridan}. We summarize the resulting structure below:

\begin{proposition}
    There exist two chain maps
    \begin{align}
        OC_\lambda&:CC_*(\mathcal{F}(X)_\lambda,\mathcal{F}(X)_\lambda)\to CM^{*+n}(X),\\
        CO_\lambda&:CM^*(X) \to CC^*(\mathcal{F}(X)_\lambda,\mathcal{F}(X)_\lambda)  
    \end{align}
    such that on cohomology $CO_\lambda$ is a unital algebra homomorphism, and $OC_\lambda$ is a $QH^*(X)$-module homomorphism, i.e. for $\upsilon\in HH_*(\mathcal{F}(X)_\lambda,\mathcal{F}(X)_\lambda)$ and $x\in QH^*(X)$,
    \begin{equation}
        OC_\lambda(CO_\lambda(x)\cap\upsilon)=x*OC_\lambda(\upsilon).
    \end{equation}
    Moreover, the open--closed map equips the monotone Fukaya summand $\mathcal{F}(X)_\lambda$ with a weak proper Calabi--Yau structure given by
    \begin{equation}
        \Phi(\upsilon)=\int_X OC_\lambda(\upsilon).
    \end{equation}
\end{proposition}

We next identify the Hochschild obstruction to an infinitesimal grading on a monotone Fukaya category summand.

\begin{thm} \label{thm:monotone-fukaya-gradeability-reprise}
    Assume that the base field $k$ contains the real numbers $\mathbb{R}$. The monotone Fukaya category summand $\mathcal{F}(X)_\lambda$ admits an infinitesimal grading if and only if the following cohomology class vanishes in Hochschild cohomology $HH^*(\mathcal{F}(X)_\lambda,\mathcal{F}(X)_\lambda)$:
    \begin{equation}
        [CO_\lambda(c_1(X)-\lambda\cdot 1_X)]=0.
    \end{equation}
\end{thm}
\begin{proof}
    Let $\kappa$ be the pre-Euler cochain of the monotone Fukaya category $\mathcal{F}(X)_\lambda$. We first fix an object $L_0$ and a point $p_0\in L_0$, and for each $L$ choose a path $\gamma_L$ from $p_0$ to a point of $L$, together with a lift to the Lagrangian Grassmannian. The reference path from $K$ to $L$ is the concatenation of $\gamma_K^{-1}$ with $\gamma_L$. To each generator $x\in CF^*(K,L)$, we assign a reference capping $w$ by the reference path chosen above, which determines the value
    \begin{equation}
        g(x)=\deg(x,w)-2\langle\omega,w\rangle
    \end{equation}
    for the Maslov degree $\deg(x,w)$. Since choosing a different capping changes the Maslov degree and twice the symplectic area by the same amount, this is independent of $w$. By monotonicity, we see that for a rigid polygon $u\in\mathcal{M}(y;x_d,\cdots,x_1)$,
    \begin{equation}
        0=\operatorname{ind}(u)=2\langle\omega,u\rangle+g(y)-\sum g(x_i)+d-2,
    \end{equation}
    and therefore we conclude that
    \begin{equation}
        \kappa=\delta G+2E
    \end{equation}
    for the Hochschild cochain $G$ supported at length one, as $G^1(x)=g(x)\cdot x$, and the Hochschild cochain $E$ weighting each polygon count for $\mu$ by its area, and $E^0\equiv0$. By monotonicity we see that
    \begin{equation}
        [E]=[CO_\lambda(c_1(X))]-\lambda\cdot 1_{HH},
    \end{equation}
    the second term correcting the contribution of Maslov two discs. Therefore it follows that
    \begin{equation}
        [\kappa]=2[CO_\lambda(c_1(X)-\lambda\cdot 1_X)],
    \end{equation}
    which is the obstruction to an infinitesimal grading.
\end{proof}  
If all Lagrangians $L$ are disjoint from a representative $D$ of $c_1(X)$ and admit lifts to the relative Fukaya category $\mathcal{F}(X,D)$, this follows from \cite[Proposition 3.8]{Sheridan} without the assumption on the base field.

We also record a lemma, which will be used to determine the open--closed image of Lagrangian spheres.

\begin{lemma}[{\cite[Lemma 3.2]{Smith}}]\label{lem:smith-lemma}
    For $L\in\mathcal{F}(X)_\lambda$, the composition of the length zero parts
    \begin{equation}
        OC_L^0\circ CO_L^0:QH^*(X)\to HF^*(L,L)\to QH^{*+n}(X)
    \end{equation}
    is equal to the quantum product with $PD_X[L]$.
\end{lemma}

The remainder of this subsection prepares the cyclic open--closed map and homotopies used in Section \ref{ssec:OC-compatibility}: $H$ for the open--closed module property, $K$ for closed--open multiplicativity, and $S$ for their compatibility. Readers familiar with their construction are invited to skip to Section \ref{ssec:quantum-cohomology} and especially Definition \ref{def:geometric-inf-grading}.

We first recall the construction of \emph{cyclic open-closed maps}, following \cite[Section 4]{Ganatra23} and \cite[Section 2]{GanatraSheridan}, beginning by adjoining strict units to the monotone Fukaya category.

\begin{defn}\label{defn:strict-unit}
    The \emph{unit enlarged monotone Fukaya category summand} $\mathcal{F}(X)^+_\lambda$ is defined to have the same objects as $\mathcal{F}(X)_\lambda$, and for two distinct objects $K,L$ the same morphism spaces
    \begin{equation}
        \hom_{\mathcal{F}(X)_\lambda^+}(K,L)=\hom_{\mathcal{F}(X)_\lambda}(K,L).
    \end{equation}
    For each object $L$, we adjoin a formal generator to its endomorphism space
    \begin{equation}
        \hom_{\mathcal{F}(X)_\lambda^+}(L,L)=CF^*(L,L)\oplus k\langle e_L^+\rangle
    \end{equation}
    with $\lvert e_L^+\rvert=0$, and set $\mu^1(e_L^+)=0$. The $A_\infty$-structure is extended from $\mathcal{F}(X)_\lambda$ such that $e_L^+$ satisfies the axioms for the strict unit
    \begin{align}
        (-1)^{\lvert x\rvert}\mu^2(e_L^+,x)=\mu^2(x,e_K^+)=x,\\
        \mu^k(\cdots,e_L^+,\cdots)=0,\quad k\geq3.
    \end{align}
    This definition works for any $A_\infty$-category $\mathcal{C}$.
\end{defn}

The Hochschild chain complex of $\mathcal{F}(X)_\lambda^+$ contains a subcomplex of \emph{degenerate elements} defined as
\begin{equation}
    \mathcal{D}=\oplus k\langle e_L^+\rangle \oplus\bigoplus k\langle z\otimes x_d\otimes\cdots\otimes e^+\otimes\cdots\otimes x_1\rangle.
\end{equation}

\begin{defn}
    For an $A_\infty$-category $\mathcal{C}$, its \emph{non-unital Hochschild complex} $CC_*^{nu}(\mathcal{C},\mathcal{C})$ is defined as the quotient
    \begin{equation}
        CC_*^{nu}(\mathcal{C},\mathcal{C})=CC_*(\mathcal{C}^+,\mathcal{C}^+)/\mathcal{D}.
    \end{equation}
    As a graded vector space, this is naturally isomorphic to two copies $CC_*(\mathcal{C},\mathcal{C})\oplus CC_*(\mathcal{C},\mathcal{C})[1]$, where the second copy is represented by chains with $e^+$ in the distinguished input. This is equipped with the differential
    \begin{equation}
        b^{nu}(\upsilon,\omega)=(b(\upsilon)+D(\omega),b'(\omega))
    \end{equation}
    for the maps sending $\omega=z\otimes x_d\otimes\cdots\otimes x_1$ to
    \begin{align}
        b'(\omega)&= \sum(-1)^{\dagger_s} z\otimes x_d\otimes\cdots \otimes x_{s+j+1}\otimes\mu^j(x_{s+j},\cdots,x_{s+1})\otimes x_s\otimes \cdots\otimes x_1 \\
        &+\sum(-1)^{\dagger_{d-j}}\mu^{j+1}(z,x_d,\cdots,x_{d-j+1})\otimes x_{d-j}\otimes\cdots\otimes x_1,\\
        D(\omega)&=(-1)^{(\dagger_2^d+\lvert z\rvert)\cdot\lvert x_1\rvert+1} x_1\otimes z\otimes\cdots\otimes x_2+(-1)^{\dagger_d+\lvert z\rvert} z\otimes x_d\otimes\cdots\otimes x_1.
    \end{align}
\end{defn}

When $\mathcal{C}$ is c-unital, the inclusion of $CC_*(\mathcal{C},\mathcal{C})$ into the non-unital complex is a quasi-isomorphism (\cite[Lemma 3.2]{Ganatra23}). We also recall the definition of the Connes B operator.

\begin{defn}
    The \emph{Connes B operator} $B:CC_*^{nu}(\mathcal{C},\mathcal{C})\to CC_*^{nu}(\mathcal{C},\mathcal{C})$ is defined as
    \begin{equation}
        B(z\otimes x_n\otimes\cdots\otimes x_1)=\sum (-1)^{(\lvert z\rvert+\dagger_{j+1}^n)(1+\dagger_j)}e^+\otimes x_j\otimes\cdots\otimes x_1\otimes z\otimes x_n\otimes\cdots\otimes x_{j+1},
    \end{equation}
    and $B$ vanishes on the second summand. Henceforth we write $b$ for $b^{nu}$. This operator satisfies the mixed complex identities
    \begin{equation}
        \lvert B\rvert=-1,\quad B^2=0,\quad bB+Bb=0.
    \end{equation}
    For cochains of the original category (extended by zero on inputs containing a strict unit), the cap product is evaluated in the canonical representatives in each summand. This satisfies the cap product identity
    \begin{equation}\label{eqn:Connes-B}
        \iota_{\phi\cup\psi}B=\iota_\phi\iota_\psi B.
    \end{equation}
\end{defn}

We now introduce the cyclic open--closed map. We do not need its precise construction, but emphasize that a genuine distinguished Hochschild input is represented by a boundary puncture with a strip-like end. An external distinguished unit is represented instead by an unpunctured boundary marked point. This marked point remains part of the domain data.

\begin{defn}
    The \emph{cyclic open--closed map} $OC^-:CC_*^{nu}(\mathcal{F}(X)_\lambda,\mathcal{F}(X)_\lambda)[[u]]\to CM^{*+n}(X)[[u]]$, for $\lvert u\rvert=2$, is defined by
    \begin{equation}
        OC^-=OC+u\cdot OC^1+\cdots
    \end{equation}
    and satisfies
    \begin{equation}
        \partial_X\circ OC^k=OC^k\circ b+OC^{k-1}\circ B,\quad k\geq 1.
    \end{equation}
\end{defn}

To define a universal consistent choice of Floer data with the adjoined strict unit, we first construct the moduli spaces of domains defining these $A_\infty$-operations, following \cite[Section 10]{GanatraThesis}.

Let $\mathcal{R}^d$ be the moduli space of holomorphic discs with $d$ input boundary marked points and one output boundary marked point for $d\geq2$, and pick a subset $F\subset\{1,\cdots,d\}$.

\begin{defn}
    Let $\mathcal{R}^{d,F}$ denote the moduli space of discs with $d$ input boundary marked points and one output boundary marked point, with the inputs indexed by $F$ designated as forgotten. Such discs are called \emph{$f$-stable} if $d-\lvert F\rvert\geq2$, and \emph{f-semistable} if $d-\lvert F\rvert=1$. For $f$-stable domains, there exists a \emph{forgetful map} $\mathcal{R}^{d,F}\to\mathcal{R}^{d-\lvert F\rvert}$, which extends to the compactification $\overline{\mathcal{R}^{d,F}}$. We call the image of the forgetful map for an $f$-stable disc $S$ its \emph{associated reduced surface} $\widehat{S}$. 
\end{defn}

For an $f$-stable disc, we transport the Floer datum assigned to its associated reduced surface by the forgetful map. For a domain whose stabilization contracts components to a boundary point $p$, we extend the solution on the reduced surface constantly over the contracted components, and assign a constant zero Hamiltonian and a constant almost complex structure $J_p$ as the perturbation data. If the reduced surface is a strip, as for an $f$-semistable disc, we pull back the fixed translation-invariant datum on that surface via the forgetful map. The moduli spaces of solutions to the Floer equations for this datum are smooth of the expected dimension. For $d>2$, with at least one $e_L^+$ input, the Floer data pulled back via forgetful maps prevent rigid contributions. We make an exception for the case when all inputs are $e_L^+$, for which the underlying disc is not (semi)stable after forgetting marked points, and declare the operation $\mu^d(e_L^+,\cdots,e_L^+)$ to be zero for $d\geq 3$. This justifies that the Floer datum extends to incorporate the strict unit $e_L^+$.

\begin{figure}[h]
\centering
\begin{tikzpicture}[line width=0.85pt, line cap=round]

\draw (0,0) circle (2);

\draw (0,-0.9) circle (1.1);

\draw[line width=1.1pt] (0.963,1.667) -- (1.038,1.797);
\draw[line width=1.1pt] (-1.667,0.963) -- (-1.797,1.038);
\draw[line width=1.1pt] (1.859,-0.498) -- (2.004,-0.537);

\draw[line width=1.1pt] (0,0.125) -- (0,0.275);
\draw[line width=1.1pt] (-0.929,-0.467) -- (-1.065,-0.403);

\draw[line width=1.1pt] (0,-2.075) -- (0,-1.925);
\node[below=5pt] at (0,-2) {$p_\infty$};
\node at (0,2.45) {$\mathbb{D}$};

\draw[->] (2.55,0) -- (3.85,0)
  node[midway,above] {$\sim$};

\node at (6.65,2.45) {$\mathbb{H}$};
\draw[-] (4.35,-0.85) -- (8.95,-0.85);
\draw[-] (4.35,0.75) -- (8.95,0.75);
\node[below=5pt] at (6.65,-0.85) {$\mathbb{R}$};
\node[above=5pt] at (6.65,0.75) {$\operatorname{Im}w=1$};

\draw[line width=1.1pt] (5.057,-0.925) -- (5.057,-0.775);
\draw[line width=1.1pt] (6.323,-0.925) -- (6.323,-0.775);
\draw[line width=1.1pt] (7.356,-0.925) -- (7.356,-0.775);

\draw[line width=1.1pt] (6.650,0.675) -- (6.650,0.825);
\draw[line width=1.1pt] (8.066,0.675) -- (8.066,0.825);

\node at (6.65,-2.27) {$p_\infty\longmapsto\infty$};

\end{tikzpicture}
\caption{The domain of a horocycle moduli space, and its biholomorphic image in the upper half plane. The tangency point $p_\infty$ is sent to $\infty$ in $\mathbb{H}$. Morse flowlines are attached to interior marked points on the horocycle, but are not drawn in this figure.}
\label{fig:horocycle}
\end{figure}
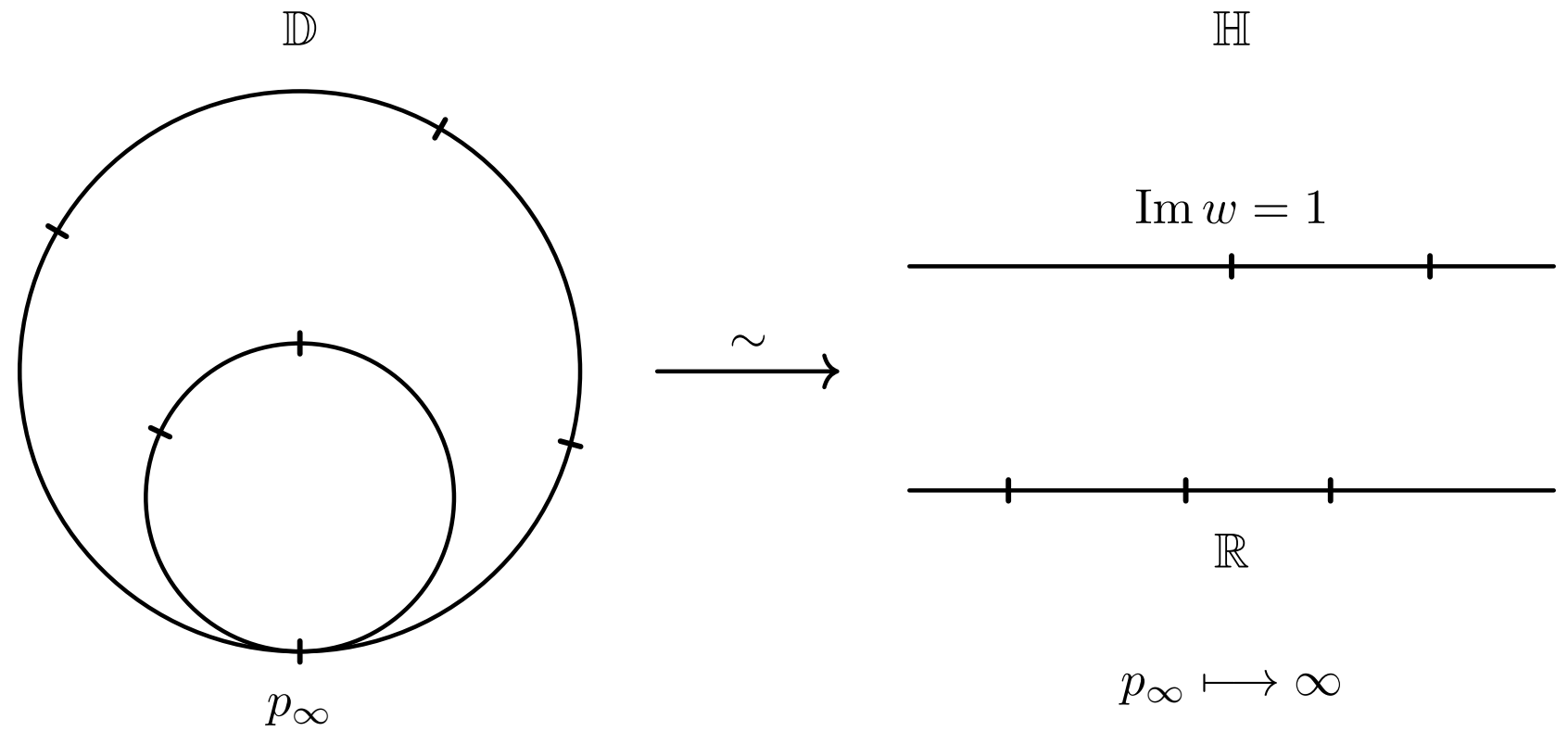

We also introduce horocycle moduli spaces, following \cite{SolomonTukachinsky,Hugtenburg}. The domains parametrized by these moduli spaces are holomorphic discs with one boundary marked point specified as $p_\infty$, decorated by a horocycle tangent to the boundary at $p_\infty$. Each such disc is biholomorphic to the upper half-plane via a map sending the horocycle to a horizontal line and $p_\infty$ to $\infty$. See Figure \ref{fig:horocycle}. We constrain $k$ interior marked points to lie on the horocycle in a prescribed order, at which we impose the incidence conditions.

We adapt the inductive construction of a universal consistent choice of Floer data from \cite[Proposition 10.8]{GanatraThesis}, and construct such data for open--closed, closed--open, and horocycle moduli spaces with strict units. The relevant families are ordered by the number of interior markings (for horocycle moduli), and boundary inputs, with the Floer data for the unstable domains fixed in advance. For a disc component with one boundary node, no ordinary inputs, and no interior marking other than a direct divisor incidence, we use the $J_L$-holomorphic disc datum from \eqref{eqn:J_L}. On its adjacent component, we pull back the prescribed data after forgetting the attaching point and stabilizing. On each boundary stratum, we prescribe the previously chosen perturbation data for each component and extend to a neighborhood of the stratum by gluing. We extend these data from a neighborhood of the boundary over the interior, and choose the remaining perturbations generically. We apply this construction simultaneously for the open--closed, closed--open, horocycle, and connection homotopy families that appear, and henceforth assume that such moduli spaces are regular and satisfy the algebraic equations from codimension $1$ boundary analysis.

\begin{figure}[h]
\centering
\resizebox{\linewidth}{!}{%
\begin{tikzpicture}[
  line width=0.8pt, line cap=round, line join=round,
  morse/.style={densely dashed, line width=0.7pt, ->, >=stealth}
]

\draw (0,1.15) circle (1.65);
\draw (0,0.325) circle (0.825);

\draw (-1.368,1.940) -- (-1.490,2.010);
\draw (1.368,1.940) -- (1.490,2.010);
\draw (0,-0.575) -- (0,-0.425);
\node[below=5pt] at (0,-0.5) {$p_\infty$};

\draw (-0.055,1.095) -- (0.055,1.205);
\draw (-0.055,1.205) -- (0.055,1.095);
\draw[morse,shorten <=2pt] (0,1.15) -- (0,3.25);
\node[above=3pt] at (0,3.25) {$\mathrm{out}$};

\draw (0.659,0.683) -- (0.769,0.793);
\draw (0.659,0.793) -- (0.769,0.683);
\draw[morse,shorten >=2pt] (2.35,0.738) -- (0.714,0.738);
\node[right=3pt] at (2.35,0.738) {$A$};

\draw[->] (-1.10,-0.36) -- (-4.33,-2.31);
\draw[->] (0,-1.3) -- (0,-2.08);
\draw[->] (1.10,-0.36) -- (4.42,-2.12);

\begin{scope}[shift={(-1.4,-0.65)}]
\draw (-3.9,-3.5) circle (1.65);
\draw (-3.9,-4.325) circle (0.825);

\draw (-5.268,-2.710) -- (-5.390,-2.640);
\draw (-2.532,-2.710) -- (-2.410,-2.640);
\draw (-3.9,-5.225) -- (-3.9,-5.075);
\node[below=5pt] at (-3.9,-5.15) {$p_\infty$};

\draw (-3.955,-3.555) -- (-3.845,-3.445);
\draw (-3.955,-3.445) -- (-3.845,-3.555);
\draw[morse,shorten <=2pt] (-3.9,-3.5) -- (-3.9,-1.28);
\node[above=3pt] at (-3.9,-1.28) {$\mathrm{out}$};

\draw (-3.241,-3.967) -- (-3.131,-3.857);
\draw (-3.241,-3.857) -- (-3.131,-3.967);
\draw[morse,shorten >=2pt] (-1.94,-3.912) -- (-3.186,-3.912);
\node[below=4pt] at (-2.16,-3.912) {$\partial A$};
\end{scope}

\begin{scope}[yshift=-0.65cm]
\draw (0,-3.5) circle (1.65);
\draw (-1.368,-2.710) -- (-1.490,-2.640);
\draw (1.368,-2.710) -- (1.490,-2.640);
\draw (0,-5.225) -- (0,-5.075);
\node[below=5pt] at (0,-5.15) {$p_\infty$};

\filldraw[fill=white,draw=black]
  (0,-3.60)
  .. controls (-0.25,-3.37) and (-0.83,-2.96) .. (-0.87,-2.43)
  .. controls (-0.93,-1.94) and (-0.57,-1.64) .. (-0.05,-1.65)
  .. controls (0.48,-1.64) and (0.90,-1.91) .. (0.87,-2.43)
  .. controls (0.82,-2.96) and (0.26,-3.37) .. (0,-3.60)
  -- cycle;
\fill (0,-3.60) circle (0.035);

\draw (-0.335,-2.355) -- (-0.225,-2.245);
\draw (-0.335,-2.245) -- (-0.225,-2.355);
\draw[morse,shorten <=2pt] (-0.28,-2.30) -- (-1.36,-2.04);
\node[above left=2pt] at (-1.36,-2.04) {$\mathrm{out}$};

\draw (0.225,-2.355) -- (0.335,-2.245);
\draw (0.225,-2.245) -- (0.335,-2.355);
\draw[morse,shorten >=2pt] (1.36,-2.04) -- (0.28,-2.30);
\node[above right=2pt] at (1.36,-2.04) {$A$};
\end{scope}

\begin{scope}[shift={(1.4,-0.65)}]
\draw (3.65,-2.7) circle (0.9);
\draw (3.65,-4.15) circle (0.55);
\draw (4.85,-4.15) circle (0.65);
\fill (3.65,-3.6) circle (0.035);
\fill (4.20,-4.15) circle (0.035);

\draw (2.923,-2.280) -- (2.819,-2.220);
\draw (4.490,-2.700) -- (4.610,-2.700);
\draw (3.110,-3.343) -- (3.033,-3.435);
\draw (3.595,-2.755) -- (3.705,-2.645);
\draw (3.595,-2.645) -- (3.705,-2.755);
\draw[morse,shorten <=2pt] (3.65,-2.7) -- (3.65,-1.28);
\node[above=3pt] at (3.65,-1.28) {$\mathrm{out}$};

\draw (3.372,-4.547) -- (3.297,-4.654);
\node[below left=2pt] at (3.297,-4.654) {$p_\infty$};

\draw (4.85,-4.88) -- (4.85,-4.72);
\draw (4.795,-4.205) -- (4.905,-4.095);
\draw (4.795,-4.095) -- (4.905,-4.205);
\draw[morse,shorten >=2pt] (6.18,-4.15) -- (4.85,-4.15);
\node[right=3pt] at (6.18,-4.15) {$A$};
\end{scope}

\end{tikzpicture}
}
\caption{The OC module homotopy and three boundary configurations. Dashed arrows indicate Morse flowlines. The middle picture shows the nodal configuration when the two horocycle marked points collide. The collision configurations are drawn at zero internal Morse length; the attached Morse-length families are suppressed.}
\label{fig:oc-module-homotopy}
\end{figure}
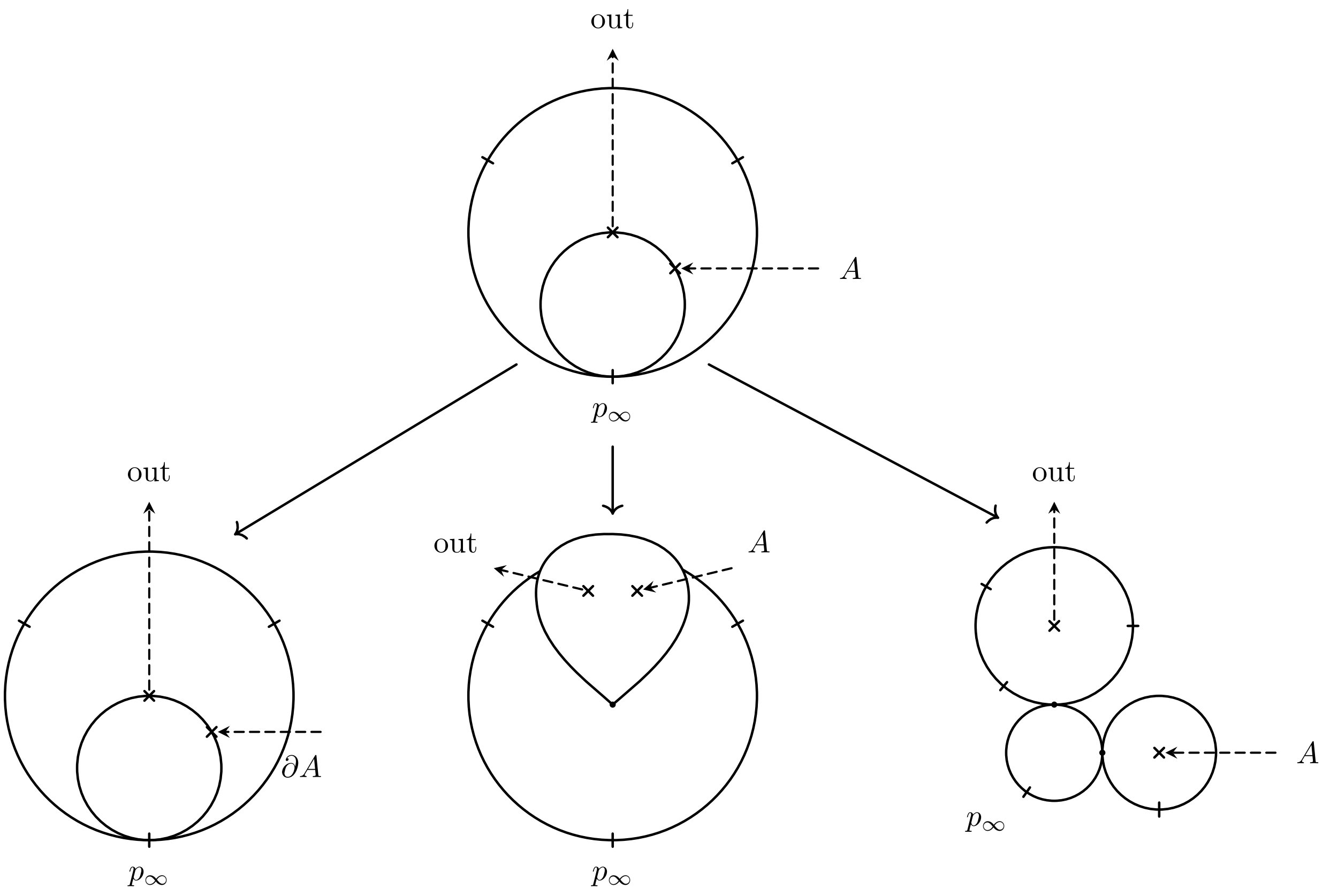

For a generator of the Hochschild chain complex $x_0\otimes\cdots\otimes x_d$, consider a disc with $d+1$ input boundary marked points, and a horocycle tangent to the puncture corresponding to $x_0$. Two interior marked points $p,p_{out}$ are placed on this horocycle, ordered as in Figure \ref{fig:oc-module-homotopy}. Attach an incoming half-line at $p$ and an outgoing half-line at $p_{out}$.

For generators $a,y$ of the Morse complex $CM^*(X)$, let
\begin{equation}
    \mathcal{M}_H(y;x_0,\cdots,x_d;a)
\end{equation}
be the moduli space of solutions to the coupled Morse and Floer equations on this domain, with boundary asymptotics $x_0,\cdots,x_d$, and Morse asymptotics to $a$ and $y$. We use coherent perturbation data to ensure that this moduli space is transversely cut out.

Throughout, the horocycle operations include their Morse-length completions. At a collision face where a sphere attaches to a disc, we attach the family obtained by replacing the interior node with a gradient trajectory of length between zero and infinity. The zero-length configuration agrees with the original nodal face, while the infinite-length boundary gives the corresponding composition in the Morse model. For $H_A$ and $K_{A,A'}$, this adds a one-parameter family; for $S_{A,A'}$, we attach the corresponding two-parameter families along its collision faces, using the same component data on their intersections. We orient these attachments so that the common faces cancel. The moduli spaces and counts below refer to these completed families.

We define the operations $H_A, S_{A,A'}$ below on both summands of the non-unital Hochschild complex. On the second summand, the distinguished input is represented by a moving unpunctured boundary marked point, retained as the horocycle tangency point. We choose the perturbation data compatibly with its endpoint identifications, whose contributions give the connecting operator $D$. Thus the module-homotopy identity for $H_A$ holds on the full non-unital complex.

\begin{defn}
    For $A\in CM^*(X)$, define the algebraic operation of degree $n+\lvert A\rvert-1$
    \begin{equation}
        H_A:CC_*^{nu}(\mathcal{F}(X)_\lambda,\mathcal{F}(X)_\lambda)\to CM^{*+n}(X)
    \end{equation}
    as
    \begin{equation}
        H_A(x_0\otimes\cdots\otimes x_d)=\sum_y\lvert\mathcal{M}^0_H(y;x_0,\cdots,x_d;a)\rvert y
    \end{equation}
    summing over the rigid moduli spaces $\mathcal{M}^0_H(y;x_0,\cdots,x_d;a)$ first for a critical point $A=a$, then extending linearly to arbitrary Morse cochains.
\end{defn}

The codimension $1$ boundary of the compactified $1$-dimensional moduli space $\overline{\mathcal{M}_H}(y;x_0,\cdots,x_d;a)$ is drawn in Figure \ref{fig:oc-module-homotopy}, which shows that
\begin{equation}\label{eqn:H_A}
    \partial_X\circ H_A+(-1)^{\lvert A\rvert}H_A\circ b+H_{\partial_X A}=OC\circ\iota_{CO(A)}-A*OC.
\end{equation}
For a cocycle $A$, this exhibits $H_A$ as the chain homotopy between the two module actions, proving the $QH^*(X)$-module property of the open--closed map on homology.

By exchanging the roles of $p_{out}$ and $p_\infty$, we may similarly define an operation which measures the homotopy for the algebra property of the closed--open map. For generators $a,a'$ of the Morse complex $CM^*(X)$, let
\begin{equation}
    \mathcal{M}_K(z;x_1,\cdots,x_d;a,a')
\end{equation}
be the moduli space of solutions to the coupled Morse and Floer equations on this domain, with boundary asymptotics $z,x_1,\cdots,x_d$, and Morse asymptotics to $a$ and $a'$.

\begin{defn}
    For $A,A'\in CM^*(X)$, define the Hochschild cochain
    \begin{equation}
        K_{A,A'}\in CC^*(\mathcal{F}(X)_\lambda,\mathcal{F}(X)_\lambda)
    \end{equation}
    as
    \begin{equation}
        K_{A,A'}(x_1\otimes\cdots\otimes x_d)=\sum_z\lvert\mathcal{M}^0_K(z;x_1,\cdots,x_d;a,a')\rvert z
    \end{equation}
    summing over the rigid moduli spaces $\mathcal{M}^0_K(z;x_1,\cdots,x_d;a,a')$ first for Morse critical points $a,a'$ representing $A,A'$, then extending bilinearly. For cocycles $A,A'$, this satisfies
    \begin{equation}
        \delta K_{A,A'}=CO(A*A')-CO(A)\cup CO(A').
    \end{equation}
\end{defn}

\begin{figure}[h]
\centering
\resizebox{\linewidth}{!}{%
\begin{tikzpicture}[
  line width=0.8pt, line cap=round, line join=round,
  morse/.style={densely dashed,line width=0.7pt,->,>=stealth},
  leglabel/.style={font=\small,inner sep=2pt}
]

\draw[line width=1pt]
  (-5,0) -- (-2.5,3.7) -- (2.5,3.7) -- (5,0)
  -- (2.5,-3.7) -- (-2.5,-3.7) -- cycle;

\draw (0,0) circle (2.3);
\fill (-1.38,1.84) circle (0.055);
\fill (1.38,1.84) circle (0.055);
\draw (0,-1.15) circle (1.15);
\draw (-0.055,-0.055) -- (0.055,0.055);
\draw (-0.055,0.055) -- (0.055,-0.055);
\draw (0.865,-0.515) -- (0.975,-0.405);
\draw (0.865,-0.405) -- (0.975,-0.515);
\draw (1.095,-1.205) -- (1.205,-1.095);
\draw (1.095,-1.095) -- (1.205,-1.205);
\draw (0,-2.375) -- (0,-2.225);
\node[below=5pt] at (0,-2.3) {$p_\infty$};

\draw[morse,shorten <=2pt] (0,0) -- (0,2.95);
\node[above=3pt] at (0,2.95) {$\mathrm{out}$};
\draw[morse,shorten >=2pt] (3.10,-0.46) -- (0.92,-0.46);
\node[right=3pt] at (3.10,-0.46) {$A$};
\draw[morse,shorten >=2pt] (3.05,-1.65) -- (1.15,-1.15);
\node[right=3pt] at (3.05,-1.65) {$A'$};

\draw (-0.35,6.3) circle (1);
\draw (-0.35,5.8) circle (0.5);
\draw (-0.35,4.7) circle (0.6);
\draw (0.8,4.7) circle (0.55);
\fill (-0.95,4.70) circle (0.045);
\fill (-0.35,4.10) circle (0.045);
\fill (-0.35,5.3) circle (0.028);
\fill (0.25,4.7) circle (0.028);
\draw (-0.395,6.255) -- (-0.305,6.345);
\draw (-0.395,6.345) -- (-0.305,6.255);
\draw (0.005,6.055) -- (0.095,6.145);
\draw (0.005,6.145) -- (0.095,6.055);
\draw (0.755,4.655) -- (0.845,4.745);
\draw (0.755,4.745) -- (0.845,4.655);
\draw (-0.883,4.303) -- (-0.778,4.378);
\node[font=\small,inner sep=1pt,below left=2pt]
  at (-0.83,4.34) {$p_\infty$};

\draw[morse,shorten <=2pt] (-0.35,6.3) -- (-0.35,7.95);
\node[leglabel,above=3pt] at (-0.35,7.95) {$\mathrm{out}$};
\draw[morse,shorten >=2pt] (1.45,6.55) -- (0.05,6.1);
\node[leglabel,right=3pt] at (1.45,6.55) {$A$};
\draw[morse,shorten >=2pt] (2.15,4.7) -- (0.8,4.7);
\node[leglabel,right=3pt] at (2.15,4.7) {$A'$};

\draw (6.25,3.35) circle (0.8);
\draw (5.25,4.1) circle (0.45);
\draw (7.25,4.1) circle (0.45);
\draw (6.25,2.12) circle (0.43);
\fill (5.45,3.35) circle (0.045);
\fill (7.05,3.35) circle (0.045);
\fill (5.61,3.83) circle (0.028);
\fill (6.89,3.83) circle (0.028);
\fill (6.25,2.55) circle (0.028);
\draw (5.205,4.055) -- (5.295,4.145);
\draw (5.205,4.145) -- (5.295,4.055);
\draw (7.205,4.055) -- (7.295,4.145);
\draw (7.205,4.145) -- (7.295,4.055);
\draw (6.205,2.075) -- (6.295,2.165);
\draw (6.205,2.165) -- (6.295,2.075);
\draw (5.738,2.838) -- (5.630,2.730);
\node[font=\small,inner sep=1pt,below left=2pt]
  at (5.684,2.784) {$p_\infty$};

\draw[morse,shorten <=2pt] (5.25,4.1) -- (5.25,5.35);
\node[leglabel,above=3pt] at (5.25,5.35) {$\mathrm{out}$};
\draw[morse,shorten >=2pt] (8.50,4.55) -- (7.25,4.1);
\node[leglabel,right=3pt] at (8.50,4.55) {$A$};
\draw[morse,shorten >=2pt] (7.60,1.65) -- (6.25,2.12);
\node[leglabel,right=3pt] at (7.60,1.65) {$A'$};

\draw (6.1,-1.85) circle (0.8);
\draw (6.1,-3.55) circle (0.9);
\fill (5.30,-1.85) circle (0.045);
\fill (6.90,-1.85) circle (0.045);
\draw (6.1,-3.1) circle (0.45);
\fill (6.1,-2.65) circle (0.028);
\draw (6.1,-1.12) -- (6.1,-0.98);
\draw (6.055,-1.895) -- (6.145,-1.805);
\draw (6.055,-1.805) -- (6.145,-1.895);
\draw (6.505,-3.145) -- (6.595,-3.055);
\draw (6.505,-3.055) -- (6.595,-3.145);
\draw (5.605,-3.145) -- (5.695,-3.055);
\draw (5.605,-3.055) -- (5.695,-3.145);

\draw[morse,shorten <=2pt] (6.1,-1.85) -- (7.25,-0.85);
\node[leglabel,above right=2pt] at (7.25,-0.85) {$\mathrm{out}$};
\draw[morse,shorten >=2pt] (7.75,-3.1) -- (6.55,-3.1);
\node[leglabel,right=3pt] at (7.75,-3.1) {$A$};
\draw[morse,shorten >=2pt] (4.50,-3.1) -- (5.65,-3.1);
\node[leglabel,left=3pt] at (4.50,-3.1) {$A'$};

\draw (0,-5.85) circle (1.05);
\fill (-0.84,-5.22) circle (0.045);
\fill (0.63,-6.69) circle (0.045);
\draw (0,-6.375) circle (0.525);
\filldraw[fill=white,draw=black]
  (0.42,-6.06)
  .. controls (0.12,-5.72) and (0.11,-5.21) .. (0.30,-4.91)
  .. controls (0.52,-4.48) and (1.06,-4.40) .. (1.24,-4.85)
  .. controls (1.48,-5.41) and (1.06,-5.94) .. (0.42,-6.06)
  -- cycle;
\fill (0.42,-6.06) circle (0.028);
\draw (-0.045,-5.895) -- (0.045,-5.805);
\draw (-0.045,-5.805) -- (0.045,-5.895);
\draw (0.525,-5.145) -- (0.615,-5.055);
\draw (0.525,-5.055) -- (0.615,-5.145);
\draw (0.835,-5.395) -- (0.925,-5.305);
\draw (0.835,-5.305) -- (0.925,-5.395);

\draw[morse,shorten <=2pt] (0,-5.85) -- (-1.65,-5.60);
\node[leglabel,left=3pt] at (-1.65,-5.60) {$\mathrm{out}$};
\draw[morse,shorten >=2pt] (1.85,-4.35) -- (0.57,-5.1);
\node[leglabel,above right=2pt] at (1.85,-4.35) {$A$};
\draw[morse,shorten >=2pt] (2.30,-5.55) -- (0.88,-5.35);
\node[leglabel,right=3pt] at (2.30,-5.55) {$A'$};

\draw (-6.2,-3.25) circle (1.1);
\fill (-7.08,-3.91) circle (0.045);
\fill (-5.32,-3.91) circle (0.045);
\filldraw[fill=white,draw=black]
  (-6.20,-3.50)
  .. controls (-6.83,-3.11) and (-7.15,-2.56) .. (-6.99,-1.88)
  .. controls (-6.88,-1.29) and (-6.21,-1.08) .. (-5.74,-1.49)
  .. controls (-5.11,-2.04) and (-5.35,-3.00) .. (-6.20,-3.50)
  -- cycle;
\fill (-6.20,-3.50) circle (0.028);
\draw[dashed,line width=0.55pt]
  (-6.2,-2.04) ellipse [x radius=0.74,y radius=0.19];
\draw (-6.245,-1.895) -- (-6.155,-1.805);
\draw (-6.245,-1.805) -- (-6.155,-1.895);
\draw (-6.645,-2.245) -- (-6.555,-2.155);
\draw (-6.645,-2.155) -- (-6.555,-2.245);
\draw (-5.935,-2.257) -- (-5.845,-2.167);
\draw (-5.935,-2.167) -- (-5.845,-2.257);

\draw[morse,shorten <=2pt] (-6.2,-1.85) -- (-6.2,-0.55);
\node[leglabel,above=3pt] at (-6.2,-0.55) {$\mathrm{out}$};
\draw[morse,shorten >=2pt] (-7.95,-2.65) -- (-6.6,-2.2);
\node[leglabel,left=3pt] at (-7.95,-2.65) {$A$};
\draw[morse,shorten >=2pt] (-4.65,-2.65) -- (-5.89,-2.212);
\node[leglabel,right=3pt] at (-4.65,-2.65) {$A'$};

\draw (-6.25,2.35) circle (0.95);
\fill (-7.01,1.78) circle (0.045);
\fill (-5.49,1.78) circle (0.045);
\draw (-6.25,1.875) circle (0.475);
\filldraw[fill=white,draw=black]
  (-6.25,2.35)
  .. controls (-6.72,2.60) and (-6.98,3.19) .. (-6.84,3.58)
  .. controls (-6.73,4.02) and (-6.13,4.18) .. (-5.85,3.83)
  .. controls (-5.40,3.29) and (-5.72,2.59) .. (-6.25,2.35)
  -- cycle;
\fill (-6.25,2.35) circle (0.028);
\draw (-6.575,3.265) -- (-6.485,3.355);
\draw (-6.575,3.355) -- (-6.485,3.265);
\draw (-6.015,3.265) -- (-5.925,3.355);
\draw (-6.015,3.355) -- (-5.925,3.265);
\draw (-5.915,2.115) -- (-5.825,2.205);
\draw (-5.915,2.205) -- (-5.825,2.115);

\draw[morse,shorten <=2pt] (-6.53,3.31) -- (-7.80,4.10);
\node[leglabel,above left=2pt] at (-7.80,4.10) {$\mathrm{out}$};
\draw[morse,shorten >=2pt] (-4.65,3.65) -- (-5.97,3.31);
\node[leglabel,right=3pt] at (-4.65,3.65) {$A$};
\draw[morse,shorten >=2pt] (-4.65,2.10) -- (-5.87,2.16);
\node[leglabel,right=3pt] at (-4.65,2.10) {$A'$};

\node[above=5pt] at (-6.2,0.5) {$A*H_{A'}$};
\node[above=5pt] at (-6.1,-5.5) {$m_3(A,A',~)\circ OC$};
\node[above=5pt] at (6.2,-5.5) {$OC\circ \iota_{K_{A,A'}}$};
\node[above=5pt] at (0.1,-7.7) {$H_{A*A'}$};
\node[above=5pt] at (-2.65,5.5) {$H_A\circ\iota_{CO(A')}$};

\end{tikzpicture}%
}
\caption{The three-point horocycle family and its six depicted boundary
configurations. Dashed arrows denote incoming Morse legs at $A,A'$ and
outgoing Morse legs at $\mathrm{out}$. Again the attached Morse-length families are suppressed.}
\label{fig:three-point-horocycle-hexagon}
\end{figure}
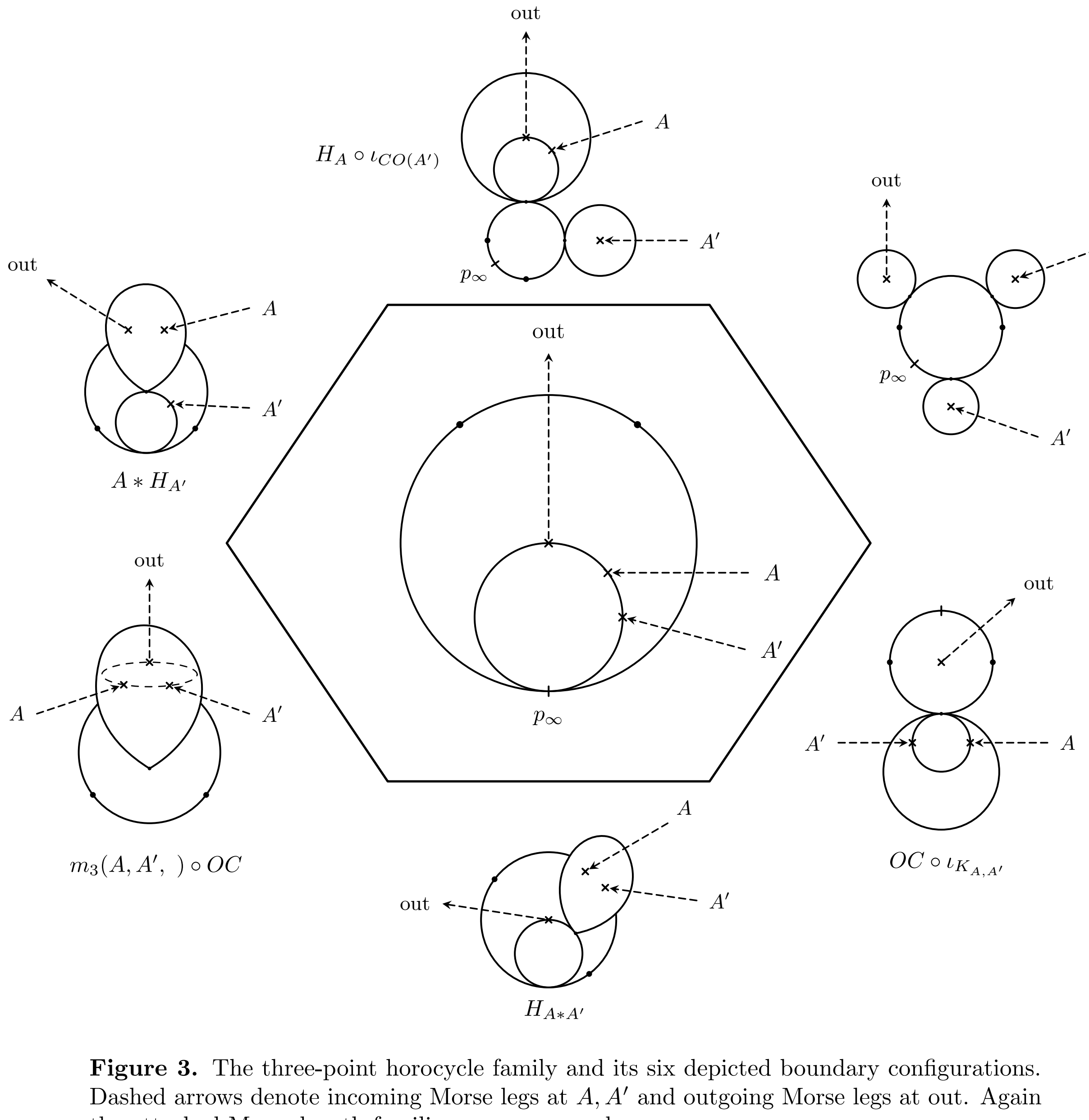

We can similarly consider the domain with the same $d+1$ boundary marked points, a horocycle tangent at the distinguished puncture, and three interior marked points $p_1,p_2,p_{out}$ on the horocycle ordered as in Figure \ref{fig:three-point-horocycle-hexagon}. Denote the corresponding moduli space of coupled Morse and Floer equations by
\begin{equation}
    \mathcal{M}_S(y;x_0,\cdots,x_d;a,a').
\end{equation}

\begin{defn}
    For $A,A'\in CM^*(X)$, define the algebraic operation
    \begin{equation}
        S_{A,A'}:CC_*^{nu}(\mathcal{F}(X)_\lambda,\mathcal{F}(X)_\lambda)\to CM^{*+n}(X)
    \end{equation}
    as
    \begin{equation}
         S_{A,A'}(x_0\otimes\cdots\otimes x_d)=\sum_y\lvert\mathcal{M}^0_S(y;x_0,\cdots,x_d;a,a')\rvert y
    \end{equation}
    again, first for Morse generators then extending bilinearly to arbitrary cochains $A,A'$.
\end{defn}

The codimension $1$ boundaries of the compactified $1$-dimensional moduli space $\overline{\mathcal{M}_S}(y;x_0,\cdots,x_d;a,a')$ are depicted in Figure \ref{fig:three-point-horocycle-hexagon}. We additionally precompose this algebraic identity with the Connes $B$ operator, which inserts a strict unit in the distinguished input. In the upper right configuration consisting of four discs, the strict unit is an input into the $A_\infty$-operation $\mu^k$ for $k\geq 3$, so its contribution is zero. Therefore we conclude that for even cocycles $A,A'$,
\begin{align}\label{eqn:S_A,A'}
    &(A*H_{A'}-H_{A*A'}+m_3(A,A',~)\circ OC+H_{A}\circ\iota_{CO(A')}+OC\circ\iota_{K_{A,A'}})\circ B\\
    &= \partial_X\circ S_{A,A'}\circ B-(-1)^{\lvert S_{A,A'}\rvert}S_{A,A'}\circ b\circ B\\&=\partial_X\circ S_{A,A'}\circ B+(-1)^{\lvert S_{A,A'}\rvert}S_{A,A'}\circ B\circ b=\partial_\mathrm{hom}(S_{A,A'}\circ B),
\end{align}
where we used that $b\circ B+B\circ b=0$.

\subsection{Closed-string grading operator}\label{ssec:quantum-cohomology}

We define a closed--string grading operator on each generalized eigensummand of quantum cohomology. Under the hypothesis that quantum multiplication with $c_1$ acts as a scalar on the whole generalized eigensummand, this operator is a derivation. We then compare it with the Hochschild grading operators induced by geometric infinitesimal gradings, and use this comparison to define notions of infinitesimal dimensions.

\begin{defn}\label{def:closed-string-grading-operator}
    Let $QH^*(X)_\lambda$ be the $\lambda$-generalized eigensummand of quantum cohomology with respect to quantum multiplication with $c_1$, and let $\pi_\lambda$ be the idempotent corresponding to that summand. Define the four-point Gromov--Witten invariant $GW(a,b,c)^\sharp$ as a cohomology-valued operation via Poincaré duality:
    \begin{equation}
        \langle GW(a,b,c)^\sharp,d\rangle=\sum_A GW_{0,4,A}(a,b,c,d).
    \end{equation}
    We define the \emph{closed-string grading operator} $\Gamma_\lambda:QH^*(X)_\lambda\to QH^*(X)_\lambda$ as
    \begin{equation}
        \Gamma_\lambda(x)=\pi_\lambda*(\mathrm{gr}_X(x)-2GW(c_1,\pi_\lambda,x)^\sharp),
    \end{equation}
    where the operation $\mathrm{gr}_X:QH^*(X)\to QH^*(X)$ acts on the underlying cohomology $H^*(X)$ by $gr_X(x)=r\cdot x$ for $x\in H^r(X)$, and is extended linearly. For example, $\mathrm{gr}_X(1_X)=0$, and $\mathrm{gr}_X(c_1)=2c_1$.
\end{defn}

In general, an additional insertion in a Gromov--Witten invariant measures the variation of the quantum product in the corresponding direction in big quantum cohomology. In particular,
\begin{equation}
    GW(c_1,\pi_\lambda,x)^\sharp=\partial_{c_1}(\pi_\lambda*_tx)\vert_{t=0},
\end{equation}
where the inputs are fixed, i.e. do not change in the $t$ direction. The divisor axiom expresses this differentiation as a derivative in the Novikov direction:
\begin{equation}\label{eqn:GW_4-qdq}
    GW(c_1,\pi_\lambda,x)^\sharp=q\frac{d}{dq}(\pi_\lambda*_q x)\vert_{q=1}.
\end{equation}

\begin{example} \label{ex:cp2-infinitesimal-grading}
We first compute the closed-string grading for the quantum cohomology of $\mathbb{CP}^2$. Using anticanonical Novikov coefficients, its small quantum cohomology ring is
\begin{equation}
    QH^*(\mathbb{CP}^2;k)=\frac{k[H,q]}{H^{3}-q^3}.
\end{equation}
Assume that $\mathrm{char}(k)\neq3$, and that $k$ contains a primitive third root of unity $\zeta$. At $q=1$, the idempotent corresponding to the eigenvalue $3$ of quantum multiplication by $c_1=3H$ is
\begin{equation}
    \pi_3=\frac{1+H+H^2}{3}.
\end{equation}
The divisor axiom gives
\begin{align}
    2q\frac{d}{dq}(\pi_3*_q\pi_3)\vert_{q=1}&=\frac{2}{9}q\frac{d}{dq}\left( 1+2H+3H^2+2q^3+q^3H \right)\vert_{q=1}\\
    &=\frac{2}{3}(2+H),
\end{align}
and therefore the closed-string grading of the idempotent $\pi_3$ is
\begin{equation}
    \Gamma_3(\pi_3)=\pi_3*\left(\frac{1}{3}(2H+4H^2)-\frac{2}{3}(2+H)\right)=0,
\end{equation}
where the last equality follows from $\pi_3*H=\pi_3$. The same computation applies to the other two idempotents for the eigenvalues $3\zeta,3\zeta^2$. Thus the closed-string grading operator is zero on each one-dimensional summand.
\end{example}

The above phenomenon is explained by the following proposition:

\begin{proposition}\label{prop:idempotent-zero}
    The closed-string grading of an idempotent for the eigensummands of $c_1*:QH^*(X)\to QH^*(X)$ is zero. In particular, one-dimensional summands of the $c_1*$ operator on quantum cohomology carry closed-string grading zero.
\end{proposition}
\begin{proof}
    The key identity we prove is
    \begin{equation} \label{eqn:grading-defect}
        2q\frac{d}{dq}(\pi_\lambda*_q x)\vert_{q=1}=\mathrm{gr}_X(\pi_\lambda)*x+\pi_\lambda*\mathrm{gr}_X(x)-\mathrm{gr}_X(\pi_\lambda*x).
    \end{equation}
    To show this, it is enough to show that for homogeneous elements $\alpha$, $\beta$,
    \begin{equation}
        2q\frac{d}{dq}(\alpha*_q\beta)\vert_{q=1}=gr(\alpha)*\beta+\alpha*gr(\beta)-gr(\alpha*\beta),
    \end{equation}
    which is true because both sides measure the failure of quantum multiplication (with Novikov variable $q$ set to 1) to preserve cohomological degree. From this we compute
    \begin{align}
        \Gamma_\lambda(\pi_\lambda)&=\pi_\lambda*(gr(\pi_\lambda)-2q\frac{d}{dq}(\pi_\lambda*_q \pi_\lambda)\vert_{q=1})\\
        &=\pi_\lambda*(gr(\pi_\lambda)-2gr(\pi_\lambda)*\pi_\lambda+gr(\pi_\lambda*\pi_\lambda))\\
        &=2\pi_\lambda*(gr(\pi_\lambda)-gr(\pi_\lambda)*\pi_\lambda)=0,
    \end{align}
    since $\pi_\lambda*\pi_\lambda=\pi_\lambda$.
\end{proof}

The grading-defect identity \eqref{eqn:grading-defect} allows us to express this four-point invariant entirely in terms of the cohomological degree operator and the small quantum product.

\begin{proposition}
    Let $QH^*(X)_\lambda$ be a generalized eigensummand as before, and let $\pi_\lambda$ be its idempotent. Then for $x\in QH^*(X)_\lambda$,
    \begin{equation}
        \Gamma_\lambda(x)=\pi_\lambda*\mathrm{gr}_X(x)-\mathrm{gr}_X(\pi_\lambda)*x.
    \end{equation}
\end{proposition}
\begin{proof}
    Since $x\in QH^*(X)_\lambda$, we have that
    \begin{equation}
        \pi_\lambda*x=x.
    \end{equation}
    Now from (\ref{eqn:grading-defect}),
    \begin{align}
        2q\frac{d}{dq}(\pi_\lambda*_qx)\vert_{q=1}&=gr(\pi_\lambda)*x+\pi_\lambda*gr(x)-gr(\pi_\lambda*x)\\
        &=gr(\pi_\lambda)*x+\pi_\lambda*gr(x)-gr(x),
    \end{align}
    and therefore
    \begin{align}
        \Gamma_\lambda(x)&=\pi_\lambda*(gr(x)-2q\frac{d}{dq}(\pi_\lambda*_q x)\vert_{q=1})\\
        &= \pi_\lambda*(gr(x)-gr(\pi_\lambda)*x-\pi_\lambda*gr(x)+gr(x))\\
        &=\pi_\lambda*gr(x)-gr(\pi_\lambda)*x
    \end{align}
    since $\pi_\lambda*\pi_\lambda=\pi_\lambda$ and $\pi_\lambda*x=x$.
\end{proof}

Assume that quantum multiplication by $c_1$ acts as $\lambda\cdot\id$ on the whole generalized eigensummand $QH^*(X)_\lambda$. Equivalently, the corresponding idempotent satisfies
\begin{equation}\label{eqn:c1-acts-scalarly}
    c_1(X)*\pi_\lambda=\lambda\pi_\lambda.
\end{equation}

Under this stronger assumption, we can prove that the closed-string grading operator is a derivation.

\begin{proposition}\label{prop:derivation}
    Under the hypotheses of \eqref{eqn:c1-acts-scalarly}, the closed-string grading operator is a derivation, i.e. for $x,y\in QH^*(X)_\lambda$,
    \begin{equation}
        \Gamma_\lambda(x*y)=\Gamma_\lambda(x)*y+x*\Gamma_\lambda(y).
    \end{equation}
\end{proposition}
\begin{proof}
    For the proof we omit the subscript for the grading operator $\Gamma$. Recall from \eqref{eqn:grading-defect} that
    \begin{equation}
        \mathrm{gr}_X(x*y)-\mathrm{gr}_X(x)*y-x*\mathrm{gr}_X(y)=-2GW(c_1,x,y)^\sharp.
    \end{equation}
    Since
    \begin{align}
        &\Gamma(x*y)-\Gamma(x)*y-x*\Gamma(y)
        \\=&\pi_\lambda*(\mathrm{gr}_X(x*y)-\mathrm{gr}_X(x)*y-x*\mathrm{gr}_X(y))+\mathrm{gr}_X(\pi_\lambda)*x*y\\
        =&-2\pi_\lambda*GW(c_1,x,y)^\sharp+\mathrm{gr}_X(\pi_\lambda)*x*y,
    \end{align}
    it is enough to show that
    \begin{equation}
        2\pi_\lambda*GW(c_1,x,y)^\sharp=\mathrm{gr}_X(\pi_\lambda)*x*y.
    \end{equation}
    We obtain this identity by differentiating the associativity equation in big quantum cohomology and using the scalar action of $c_1$ on the generalized eigensummand. Take the equation in big quantum cohomology
    \begin{equation}
        (c_1-\lambda\cdot 1_X)*_t(\pi_\lambda*_ty)=((c_1-\lambda\cdot1_X)*_t\pi_\lambda)*_ty,
    \end{equation}
    differentiate in the $x$ direction, and evaluate at $t=0$:
    \begin{align}
        &GW(c_1,x,\pi_\lambda*y)^\sharp+(c_1-\lambda\cdot1_X)*GW(\pi_\lambda,x,y)^\sharp\\=&GW(c_1-\lambda\cdot1_X,\pi_\lambda,x)^\sharp*y+GW((c_1-\lambda\cdot1_X)*\pi_\lambda,x,y)^\sharp.
    \end{align}
    Multiplying both sides by $\pi_\lambda$ and using $c_1*\pi_\lambda-\lambda\pi_\lambda=0$ gives
    \begin{equation}
        \pi_\lambda*GW(c_1,x,y)^\sharp=GW(c_1,\pi_\lambda,x)^\sharp*y,
    \end{equation}
    since the $1_X$ contribution to a four-point Gromov--Witten invariant vanishes due to the fundamental class axiom. Now using \eqref{eqn:grading-defect}:
    \begin{equation}
        \mathrm{gr}_X(\pi_\lambda*x)-\mathrm{gr}_X(\pi_\lambda)*x-\pi_\lambda*\mathrm{gr}_X(x)=-2GW(c_1,\pi_\lambda,x)^\sharp
    \end{equation}
    and quantum multiplying both sides by $y$, we obtain
    \begin{equation}
        2GW(c_1,\pi_\lambda,x)^\sharp*y=\mathrm{gr}_X(\pi_\lambda)*x*y
    \end{equation}
    which concludes the proof. 
\end{proof}

Thus quantum multiplication adds the generalized weights of $\Gamma_\lambda$ acting on $QH^*(X)_\lambda$.

We now introduce geometric infinitesimal gradings on the monotone Fukaya summand $\mathcal{F}(X)_\lambda$. For the construction and open--closed compatibility proof below, we restrict to objects that admit a graded lift to the relative Fukaya category $\mathcal{F}(X,D)$, and fix such lifts. We use the integer degrees determined by such lifts, with $\lvert q\rvert=2$, and let $gr$ denote the Euler operator for the $\mathbb{Z}$-grading of $\mathcal{F}(X,D)$. We use the uncurved monotone summand, retaining the $q$-weights of its operations; the strict units $e_L^+$ are from Definition \ref{defn:strict-unit}.

We will show that the closed--open map intertwines their induced grading operators on Hochschild cohomology with $\Gamma_\lambda$. The construction requires a Hochschild primitive of $q\partial_q\mu$, which we obtain using closed--open operations and horocycle homotopies.

To begin, choose cocycles $D,\Pi_\lambda$ representing the cohomology classes $c_1(X),\pi_\lambda\in H^*(X)$. Then we may choose a bounding cochain $C_\lambda$ such that
\begin{equation}
    \partial_XC_\lambda=\Pi_\lambda*D-\lambda\Pi_\lambda.
\end{equation} 
Now since both $CO(\Pi_\lambda)$ and $CO(1_X)$ represent the unit in Hochschild cohomology before adjoining the strict unit, there exists a Hochschild cochain $\beta_\Pi$ such that
\begin{equation}
    \delta\beta_\Pi=CO(\Pi_\lambda)-CO(1_X).
\end{equation}
Similarly since the closed--open map preserves products, there exists a Hochschild cochain $K_{\Pi,D}$ such that
\begin{equation}
    \delta K_{\Pi,D}=CO(\Pi_\lambda*D)-CO(\Pi_\lambda)\cup CO(D),
\end{equation}
and moreover this can be defined geometrically as a count of holomorphic discs with interior incidence conditions on a horocycle. We can similarly define $K_{1,D}$ and $K_{1,1}$ as a count of holomorphic discs, for the previous choice of chain level unit $1_X$. 

Finally we may define a Hochschild cochain $K_{D,1}^\mathrm{inc}$ from the same horocycle family as $K_{D,1}$, but we replace the incoming Morse trajectory at the $D$-marked point by direct incidence $u(z_D)\in D$. This cochain satisfies
\begin{equation}
    \delta K_{D,1}^\mathrm{inc}=CO(D)-q\partial_q\mu\cup CO(1_X)-\lambda CO(1_X),
\end{equation}
instead of the cup product homotopy. Since $[CO(D)]=\lambda[CO(1_X)]$ holds in $HH^*(\mathcal{F}(X)_\lambda,\mathcal{F}(X)_\lambda)$, the above equation implies that
\begin{equation}
    [q\partial_q\mu]\cup[CO(1_X)]=0.
\end{equation}
Since $[CO(1_X)]$ is the unit, we may choose a Hochschild cochain $\beta_0$ such that
\begin{equation}
    \delta\beta_0=q\partial_q\mu.
\end{equation}
This primitive $\beta_0$ establishes the existence of an infinitesimal grading: the correction below removes dependence on this choice up to Hochschild coboundary and ensures compatibility with the closed-string grading through the open--closed map.

Assembling the above, we define geometric infinitesimal gradings.

\begin{defn}\label{def:geometric-inf-grading}
    We define a \emph{geometric infinitesimal grading} on a monotone Fukaya summand $\mathcal{F}(X)_\lambda$ to be of the form
    \begin{equation}
        \eta_\mathrm{geom}=gr+2\beta_\mathrm{geom}
    \end{equation}
    where
    \begin{equation}
        \beta_\mathrm{geom}=CO(C_\lambda)-K_{\Pi,D}-\beta_\Pi\cup CO(D)+K_{1,D}+\lambda\cdot\beta_\Pi-K_{D,1}^\mathrm{inc}-\beta_0\cup CO(1_X)+\beta_0
    \end{equation}
    satisfies $\delta\beta_\mathrm{geom}=q\partial_q\mu$ for the $A_\infty$-structure cocycle $\mu$.
\end{defn}

We verify that the above indeed defines an infinitesimal grading.

\begin{proposition}\label{prop:geometric-grading-is-a-grading}
    The cochain $\eta_\mathrm{geom}$ constructed above is an infinitesimal grading for the monotone Fukaya category $\mathcal{F}(X)_\lambda$.
\end{proposition}
\begin{proof}
    We need to verify that
    \begin{equation}
        \delta\eta_\mathrm{geom}=\kappa_\lambda
    \end{equation}
    at chain level. Writing
    \begin{equation}
        \mu^k=\sum \mu^k_d
    \end{equation}
    where $\mu^k_d$ is the contribution of holomorphic discs with intersection number $d$ with $D$, we compute that
    \begin{equation}
        (\delta gr)^k=\sum (2-k-2d)\mu^k_d=\kappa_\lambda-2q\partial_q\mu.
    \end{equation}
    Therefore it is enough to show that
    \begin{equation}
        \delta\beta_\mathrm{geom}=q\partial_q\mu.
    \end{equation}
    We directly compute
    \begin{align}
        \delta\beta_\mathrm{geom}&=CO(\Pi_\lambda*D)-\lambda CO(\Pi_\lambda)-CO(\Pi_\lambda*D)+CO(\Pi_\lambda)\cup CO(D)\\
        &-(CO(\Pi_\lambda)-CO(1_X))\cup CO(D)+CO(D)-CO(1_X)\cup CO(D)+\lambda(CO(\Pi_\lambda)-CO(1_X))\\
        &-CO(D)+q\partial_q\mu\cup CO(1_X)+\lambda CO(1_X)-q\partial_q\mu\cup CO(1_X)+q\partial_q\mu\\
        &=q\partial_q\mu,
    \end{align}
    which finishes the proof.
\end{proof}

We now discuss the dependence of the geometric infinitesimal grading on choices. Fix the Floer data and the representatives $D,\Pi$, so that the horocycle cochains are fixed. Changing $C_\lambda$ changes $\beta_\mathrm{geom}$ by the closed--open image of a cocycle, which need not be a coboundary. Replacing $\beta_0$ by adding a cocycle $z$ changes $\beta_\mathrm{geom}$ by $z-z\cup CO(1_X)$, which is a Hochschild coboundary since $[CO(1_X)]$ is the unit in Hochschild cohomology. Similarly, replacing $\beta_\Pi$ by $\beta_\Pi+z$, where $\delta z=0$, changes $\beta_\mathrm{geom}$ by
\begin{equation}
    -z\cup CO(D)+\lambda z=\delta(z\cup(K_{D,1}^\mathrm{inc}+\beta_0\cup CO(1_X)))+\lambda(z-z\cup CO(1_X)),
\end{equation}
a coboundary. Therefore with the geometric data fixed, the remaining chain level ambiguity is only by Hochschild coboundaries and the chain level image of the closed--open map.

The following theorem connects the closed-string grading operators studied above with the open string grading operators discussed in Section \ref{sec:algebra}. It shows that the grading operator on Hochschild homology induced from \emph{geometric} infinitesimal gradings can be computed in small quantum cohomology.

\begin{thm}\label{thm:OC-compatibility-reprise}
    Let $(X,\omega)$ be a closed monotone symplectic manifold satisfying $c_1*\pi_\lambda=\lambda\pi_\lambda$, and let $\eta_\mathrm{geom}$ be a choice of geometric infinitesimal grading on $\mathcal{F}(X)\lambda$. Then the open--closed map satisfies
    \begin{equation}
        OC_\lambda\circ\Gamma^{HH_*}_{\eta_\mathrm{geom}}=\pi_\lambda*(\mathrm{gr}_X-n)\circ OC_\lambda
    \end{equation}
    on homology.
\end{thm}

The proof uses horocycle moduli spaces and occupies the following subsection. We first use this comparison to establish the corresponding compatibility for the closed--open map.

\begin{thm}\label{thm:CO-compatibility-reprise}
    Let $(X,\omega)$ and $\lambda$ satisfy the conditions of Theorem \ref{thm:OC-compatibility-reprise}. Then the closed--open map intertwines grading operators on cohomology:
    \begin{equation}
        CO\circ \Gamma_\lambda=\Gamma^{HH^*}_{\eta_\mathrm{geom}}\circ CO.
    \end{equation}
\end{thm}
\begin{proof}  
    For any $x\in QH^*(X)_\lambda$ and $\upsilon\in HH_*(\mathcal{F}(X)_\lambda,\mathcal{F}(X)_\lambda)$, on homology
    \begin{align}
        \pi_\lambda*(\mathrm{gr}_X-n)\circ(x*OC(\upsilon))&=\pi_\lambda*(\mathrm{gr}_X-n)\circ OC(CO(x)\cap\upsilon)\\
        &=OC\circ\Gamma^{HH_*}_\eta(CO(x)\cap\upsilon)\\
        &=OC\circ(\Gamma_\eta^{HH^*}(CO(x))\cap\upsilon)+OC(CO(x)\cap\Gamma_\eta^{HH_*}(\upsilon))\\
        &=OC\circ(\Gamma_\eta^{HH^*}(CO(x))\cap\upsilon)+x*OC(\Gamma_\eta^{HH_*}(\upsilon))\\
        &=OC\circ(\Gamma_\eta^{HH^*}(CO(x))\cap\upsilon)+x*\pi_\lambda*(\mathrm{gr}_X-n)\circ OC_\lambda(\upsilon)
    \end{align}
    where we used the compatibility of infinitesimal gradings and cap products (Proposition \ref{prop:cap-product-matches-grading}). In particular,
    \begin{equation}
        OC\circ(\Gamma_\eta^{HH^*}(CO(x))\cap\upsilon)=\pi_\lambda*\mathrm{gr}_X(x*OC(\upsilon))-x*\mathrm{gr}_X(OC(\upsilon)).
    \end{equation}
    Recall that the monotone Fukaya category carries a weak proper Calabi--Yau structure defined by
    \begin{equation}
        \Phi(\alpha)=\int_X OC_\lambda(\alpha).
    \end{equation}
    Applying this to $\alpha=CO(x)\cap\upsilon$, we obtain
    \begin{align}
        \Phi(\Gamma_\eta^{HH^*}(CO(x))\cap\upsilon)&= \int_X\pi_\lambda*\mathrm{gr}_X(x*OC(\upsilon))-\int_Xx*\mathrm{gr}_X(OC(\upsilon))\\
        &=\int_X\Gamma_\lambda(x)*OC(\upsilon)=\Phi(CO(\Gamma_\lambda(x))\cap\upsilon),
    \end{align}
    where the middle equality follows from
    \begin{equation}
        \langle \mathrm{gr}_X(a),b\rangle_X+\langle a,\mathrm{gr}_X(b)\rangle_X=2n\langle a,b\rangle_X.
    \end{equation}
    To be explicit, note that for $z=OC(\upsilon)$,
    \begin{align}
        \langle\pi_\lambda,\mathrm{gr}_X(x*z)\rangle_X-\langle x,\mathrm{gr}_X(z)\rangle &= 2n\langle \pi_\lambda,x*z\rangle_X-\langle \mathrm{gr}_X(\pi_\lambda),x*z\rangle_X-2n\langle x,z\rangle_X+\langle \mathrm{gr}_X(x),\pi_\lambda*z\rangle_X\\
        &=-\langle \mathrm{gr}_X(\pi_\lambda)*x,z\rangle_X+\langle \mathrm{gr}_X(x)*\pi_\lambda,z\rangle_X=\langle\Gamma_\lambda(x),z\rangle_X
    \end{align}
    where the terms involving $2n$ cancel because $\pi_\lambda*z=z$. By nondegeneracy of $\Phi$, we conclude that
    \begin{equation}
        CO\circ\Gamma_\lambda=\Gamma_\eta^{HH^*}\circ CO. \qedhere
    \end{equation}
\end{proof}

To relate the open--closed comparison to Calabi--Yau dimension, we consider weak proper Calabi-Yau structures that are different from the one used above. These will be modified by even invertible classes in $QH^*(X)_\lambda$. The following proposition constructs such Calabi--Yau structures and gives a criterion for their homogeneity.

\begin{proposition}\label{thm:CY-criterion-reprise}
    Let $\theta\in QH^*(X)_\lambda$ be an even class invertible in $QH^*(X)_\lambda$, satisfying
    \begin{equation}\label{eqn:inf-dim-formula}
        \pi_\lambda*gr_X(\theta)=(n-d)\theta.
    \end{equation}
    Then the functional $\Phi_\theta(\alpha)\coloneqq \int_X \theta*OC_\lambda(\alpha)$ on $HH_*(\mathcal{F}(X)_\lambda,\mathcal{F}(X)_\lambda)$ defines a homogeneous weak proper Calabi--Yau structure of infinitesimal dimension $d$ with respect to the geometric infinitesimal grading.
\end{proposition}
\begin{proof}
    Since $\theta$ is an invertible class, we can find an inverse $\theta^{-1}\in QH^*(X)_\lambda$. Then
    \begin{align}
        \Phi_\theta(CO(\theta^{-1})\cap\alpha)&=\int_X\theta*OC_\lambda(CO(\theta^{-1})\cap\alpha)\\
        &=\int_X\theta*\theta^{-1}*OC_\lambda(\alpha)=\Phi(\alpha)
    \end{align}
    is the standard nondegenerate Calabi--Yau structure, which shows that $\Phi_\theta$ is also nondegenerate.

    By the same argument as in the proof of Theorem \ref{thm:CO-compatibility-reprise}, we obtain
    \begin{align}
        \Phi_\theta(\Gamma^{HH_*}_\eta(\alpha))&=\langle\theta,\pi_\lambda*(\mathrm{gr}_X-n)OC_\lambda(\alpha)\rangle_X\\
        &= \langle n\theta-\pi_\lambda*\mathrm{gr}_X(\theta),OC_\lambda(\alpha)\rangle_X\\
        &=d\Phi_\theta(\alpha)
    \end{align}
    on homology, which shows the required infinitesimal dimension.
\end{proof}

When the base field has characteristic zero, this dimension $d$ is unique.

\begin{corollary}
    Suppose that $\mathrm{char}(k)=0$, and $\theta_1,\theta_2$ are two classes satisfying the conditions of Proposition \ref{thm:CY-criterion-reprise} with infinitesimal dimensions $d_1,d_2$. Then $d_1=d_2$.
\end{corollary}
\begin{proof}
    The derivation formula shows that
    \begin{equation}
        \Gamma_\lambda(\theta_2*\theta_1^{-1})=\Gamma_\lambda(\theta_2)*\theta_1^{-1}+\theta_2*\Gamma_\lambda(\theta_1^{-1})=(d_1-d_2)\theta_2*\theta_1^{-1},
    \end{equation}
    which shows that $\theta_2*\theta_1^{-1}$ is an invertible eigenvector of eigenvalue $d_1-d_2$. Thus its powers give eigenvectors of eigenvalue $n(d_1-d_2)$ for each integer $n$, yet $QH^*(X)_\lambda$ is finite-dimensional. Therefore it follows that when $\mathrm{char}(k)=0$, $d_1-d_2$ must be zero.
\end{proof}

For a fixed infinitesimal grading $\eta$, an $\eta$-graded spherical object determines the infinitesimal dimension of any weak proper Calabi--Yau structure homogeneous with respect to $\eta$. 

\begin{proposition}\label{prop:CY-unique}
    Assume that $\mathrm{char}(k)\neq 2$, and fix an infinitesimal grading $\eta$ on $\mathcal{F}(X)_\lambda$. Suppose that $(V,\xi_V)$ is an $\eta$-graded spherical object of infinitesimal dimension $d$. Then any weak proper Calabi--Yau structure on $\mathcal{F}(X)_\lambda$ homogeneous with respect to $\eta$ must have infinitesimal dimension $d$.
\end{proposition}
\begin{proof}
    Let $\Gamma_V$ be the induced grading operator on $HF^*(V,V)$, and let $\langle~\rangle$ denote the induced pairing of dimension $d'$. Since $\mathrm{char}(k)\neq 2$, it follows that $[V]$ is an eigenvector for $\Gamma_V$, and thus
    \begin{equation}
        \Gamma_V([V])=d[V].
    \end{equation}
    Using the homogeneous pairing and $\Gamma_V(e_V)=0$, we conclude from
    \begin{align}
        \langle \Gamma_V(e_V),[V]\rangle+\langle e_V,\Gamma_V([V])\rangle=d'\langle e_V,[V]\rangle
    \end{align}
    that $d=d'$.
\end{proof}

With this discussion on homogeneous Calabi--Yau structures and the uniqueness of their dimensions, we now state the main theorem of this subsection.

\begin{thm}\label{thm:CY-shift-reprise}
    Assume the conditions of Theorem \ref{thm:OC-compatibility-reprise}, and further assume the existence of $\theta$ as in Proposition \ref{thm:CY-criterion-reprise}. Then
    \begin{equation}
        \theta^{-1}*OC_\lambda\circ\Gamma_\eta^{HH_*}=(\Gamma_\lambda-d)\circ(\theta^{-1}* OC_\lambda)
    \end{equation}
    in homology. Thus, $\theta^{-1}*OC_\lambda$ raises weights by $d$.
\end{thm}
\begin{proof}
    Recall that the class $\theta$ satisfies
    \begin{equation}
        \Gamma_\lambda(\theta)=(n-d)\theta-\mathrm{gr}_X(\pi_\lambda)*\theta.
    \end{equation}
    Therefore by the derivation property, for any $z\in QH^*(X)_\lambda$,
    \begin{align}
        \Gamma_\lambda(\theta^{-1}*z)&=(-(n-d)\theta^{-1}+\mathrm{gr}_X(\pi_\lambda)*\theta^{-1})*z+\theta^{-1}*\Gamma_\lambda(z)\\
        &=\theta^{-1}*(\pi_\lambda*\mathrm{gr}_X(z)-(n-d)z).
    \end{align}
    Now applying this to $z=OC_\lambda(\alpha)$ from the open--closed compatibility
    \begin{equation}
        OC_\lambda\circ\Gamma_\eta^{HH_*}(\alpha)=\pi_\lambda*(\mathrm{gr}_X-n)\circ OC_\lambda(\alpha),
    \end{equation}
    we conclude that
    \begin{align}
        \theta^{-1}*OC_\lambda\circ\Gamma_\eta^{HH_*}=\pi_\lambda*(\theta^{-1}*(\mathrm{gr}_X-n)\circ OC_\lambda)=(\Gamma_\lambda-d)\circ (\theta^{-1}*OC_\lambda).
    \end{align}
\end{proof}

This criterion will be used to determine the infinitesimal dimensions of spherical objects in Section 5.

\subsection{Open-closed map matches gradings} \label{ssec:OC-compatibility}

In this subsection, we prove Theorem \ref{thm:OC-compatibility-reprise}. The key inputs will be the boundary analysis for two and three pointed horocycle moduli spaces, combined with Ganatra--Sheridan's proof \cite{GanatraSheridan} that the cyclic open--closed map preserves connections. We only use the $u^0$ term of the chain level identity for this proof. 

Some conventions: we omit all $\lambda$ and $\mathrm{geom}$ subscripts for brevity throughout. In this subsection only, we shift the closed-string complex by degree $n$ such that the open--closed map has degree zero: this gets rid of the $(-1)^n$ that appears in signs. For brevity, we define for a homogeneous linear map $\Phi:C^*\to D^*$,
\begin{equation}
    \partial_\mathrm{hom}\Phi\coloneqq \partial_D\circ\Phi-(-1)^{\lvert\Phi\rvert}\Phi\circ\partial_C.
\end{equation}

The statement of Theorem \ref{thm:OC-compatibility-reprise} is that
\begin{equation}
    \pi*((\mathrm{gr}_X-n)\circ OC-OC\circ(\mathcal{L}_{\eta_\mathrm{geom}}-\Lambda))
\end{equation}
is zero in homology. From the $\mathbb{Z}$-grading of the relative Fukaya category, it follows that
\begin{equation}
    OC\circ(gr-\Lambda+2q\partial_q)=(\mathrm{gr}_X-n+2q\partial_q)\circ OC,
\end{equation}
since the unshifted open--closed map $OC$ has degree $n$. Thus it is enough to show the following proposition.

\begin{proposition}\label{prop:reduction-0}
    $\pi*(q\partial_q\circ OC-OC\circ(q\partial_q-\mathcal{L}_{\beta_\mathrm{geom}}))$
    is trivial in homology.
\end{proposition}

The proof has four steps. First, we will use that the cyclic open--closed map $OC^-$ preserves connections to reduce to Proposition \ref{prop:reduction-1}. Then we will analyze the boundary of a horocycle moduli space with three interior incidence conditions to reduce to Proposition \ref{prop:reduction-2}. Then Proposition \ref{prop:reduction-3} uses the homotopies $S_{A,A'}$ and $H_C$ to reduce to an expression involving the unit, which Proposition \ref{prop:OC-compatiblity-end-of-proof} shows vanishes on homology.

\begin{figure}[htbp]
\centering
\resizebox{\linewidth}{!}{%
\begin{tikzpicture}[
  line width=0.8pt, line cap=round, line join=round,
  font=\large,
  morse/.style={dashed,->,line width=0.7pt},
  limit/.style={->,line width=0.8pt},
  pics/gsmark/.style={code={%
    \draw (-0.055,-0.055)--(0.055,0.055);
    \draw (-0.055,0.055)--(0.055,-0.055);
  }}
]

\foreach \mode/\xx in {0/0,1/4.4,2/8.8}{%
  \begin{scope}[xshift=\xx cm]
    \draw (0,0) circle (0.87);
    \foreach \a in {35,145,215,325}
      \fill (\a:0.87) circle (0.032);
    \fill (0,-0.87) circle (0.037);
    \node[below=3pt] at (0,-0.87) {$p_\infty$};
    \node at (0,-0.26) {$OC$};

    \draw (0,2.00) circle (0.65);
    \draw[gray!65,line width=0.4pt,densely dashed]
      (-0.65,2.00) arc[start angle=180,end angle=0,
                       x radius=0.65,y radius=0.13];
    \draw[gray!65,line width=0.4pt]
      (-0.65,2.00) arc[start angle=180,end angle=360,
                       x radius=0.65,y radius=0.13];

    \pic at (0.27,2.00) {gsmark};
    \node[font=\normalsize] at (0.32,2.24) {$D$};
    \pic at (-0.20,2.28) {gsmark};
    \draw[morse] (-0.20,2.34)--(-0.20,3.08);
    \node[above=3pt] at (-0.20,3.08) {$\mathrm{out}$};

    \ifnum\mode=0
      \pic at (0,0.13) {gsmark};
      \pic at (0,1.64) {gsmark};
      \draw[morse] (0,0.19)--(0.5,1.25);
      \draw[morse] (0.5,1.25)--(0,1.58);
    \fi
    \ifnum\mode=1
      \pic at (0,0.13) {gsmark};
      \pic at (0,1.64) {gsmark};
      \draw[morse] (0,0.19)--(0,1.58);
      \node[right=4pt] at (0,1.10) {$\ell$};
    \fi
    \ifnum\mode=2
      \draw (0,0.13)--(0,1.64);
      \fill (0,0.13) circle (0.036);
      \fill (0,1.64) circle (0.036);
      \node[right=4pt,font=\normalsize] at (0,1.10) {node};
    \fi
  \end{scope}
}

\begin{scope}[xshift=13.2cm]
  \draw (0,0.65) circle (1.25);
  \foreach \a in {35,145,215,325}
    \fill ({1.25*cos(\a)},{0.65+1.25*sin(\a)}) circle (0.032);

  \draw[line width=0.55pt] (0,0.194) circle (0.794);
  \pic at (-0.65,0.65) {gsmark};
  \pic at (0.65,0.65) {gsmark};
  \node[above=5pt] at (0.65,0.65) {$D$};
  \draw[morse] (-0.65,0.71)--(-0.65,3.08);
  \node[above=3pt] at (-0.65,3.08) {$\mathrm{out}$};
  \fill (0,-0.60) circle (0.037);
\end{scope}

\begin{scope}[xshift=18cm]
  \foreach \xx in {-1.24,0,1.24}
    \draw (\xx,0.65) circle (0.62);
  \fill (-0.62,0.65) circle (0.036);
  \fill (0.62,0.65) circle (0.036);
  \foreach \xx in {-1.24,1.24}{%
    \foreach \a in {235,305}
      \fill ({\xx+0.62*cos(\a)},{0.65+0.62*sin(\a)}) circle (0.030);
  }
  \fill (0,0.03) circle (0.037);
  \node[below=3pt] at (0,0.03) {$p_\infty$};

  \pic at (-1.24,0.82) {gsmark};
  \draw[morse] (-1.24,0.88)--(-1.24,3.08);
  \node[above=3pt] at (-1.24,3.08) {$\mathrm{out}$};
  \node at (-1.24,0.40) {$OC$};
  \node at (0,0.65) {$\mu$};
  \pic at (1.24,0.82) {gsmark};
  \node[right=3pt,font=\normalsize] at (1.14,0.62) {$D$};
\end{scope}

\draw[limit] (3.13,1.10)--(1.25,1.10)
  node[midway,above=6pt] {$\ell\to\infty$};
\draw[limit] (5.65,1.10)--(7.40,1.10)
  node[midway,above=6pt] {$\ell\to0$};
\draw[limit] (11.78,1.10)--(10.18,1.10)
  node[midway,above=6pt] {$t\to0$};
\draw[limit] (14.65,1.10)--(15.90,1.10)
  node[midway,above=6pt] {$t\to1$};

\node at (0,-1.95) {$D*^\mathrm{inc} OC$};
\node at (4.4,-1.95) {$H^1_0$};
\node at (8.8,-1.95) {$H^{12}_0$};
\node at (13.2,-1.95) {$H^2_0$};
\node at (18,-1.95) {$OC\circ\iota_{q\partial_q\mu}$};

\node[font=\normalsize] at (0,-2.40) {$\ell=\infty$};
\node[font=\normalsize] at (4.4,-2.40) {$0<\ell<\infty$};
\node[font=\normalsize] at (8.8,-2.40) {$\ell=0,\quad t=0$};
\node[font=\normalsize] at (13.2,-2.40) {$0<t<1$};
\node[font=\normalsize] at (18,-2.40) {$t=1$};
\node[below=3pt] at (13.2,-0.57) {$p_\infty$};

\end{tikzpicture}%
}
\caption{The two Ganatra--Sheridan families at $k=0$, with their common face $H^{12}_0$ drawn once in the center. At $t=1$, there is an additional $\lambda \cdot OC$ contribution when the $D$-disc has no ordinary inputs and the $\mu$ disc has only the distinguished input $p_\infty$; after forgetting the divisor attachment, the latter reduces to the identity strip. We define $H_\nabla^{-1}=H_0^1+H_0^2$.}
\label{fig:H_nabla^-1}
\end{figure}
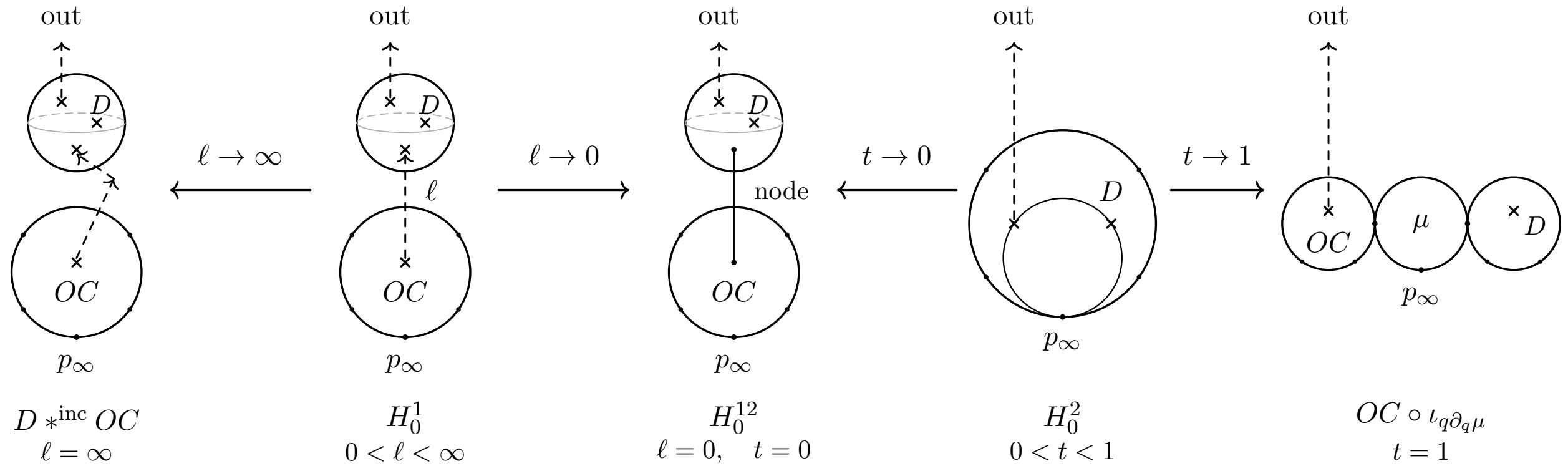

We show that Proposition \ref{prop:reduction-0} follows from the following proposition:

\begin{proposition}\label{prop:reduction-1}
    The proof of Proposition \ref{prop:reduction-0} reduces to showing that
    \begin{equation}\label{eqn:reduction-1}
     (\Pi*H_\nabla^{-1}-\Pi*OC\circ\iota_{\beta_\mathrm{geom}}+(\Pi*m_3^\mathrm{inc}(D,1_X,~)+m_3(\Pi,D,~)+C*)\circ OC)\circ B
    \end{equation}
    is zero in homology, where $H_\nabla^{-1}$ is defined below and $m_3$ is the quantum Massey product, and $m_3^\mathrm{inc}$ is its variant with direct incidence at the $D$-marked point.
\end{proposition}
\begin{proof}
    Recall that Ganatra--Sheridan \cite[Theorem 7.10]{GanatraSheridan} prove that the cyclic open--closed map
    \begin{equation}
        OC^-=OC+u\cdot OC^1+\cdots
    \end{equation}
    preserves $q\partial_q$-connections up to a geometrically defined homotopy. We use the families of domains from \cite[Lemmas 7.3--7.6]{GanatraSheridan} for the uncurved monotone relative Fukaya summand, with the perturbation data chosen above. Throughout this calculation, $\mu^0=0$, and we require the connection preservation only on the $u^{-1}$ and $u^0$ coefficients. The quantum connection in the $q\partial_q$-direction (with sign conventions from \cite[Section 2]{GanatraSheridan}) is
    \begin{equation}
        \nabla_{q\partial_q}^{\mathrm{closed}}=q\partial_q+u^{-1}D*_q^\mathrm{inc},
    \end{equation}
    where $D$ is the chosen chain level representative of $c_1$. This chain-level quantum product $D*_q^\mathrm{inc}$ is not given by the Morse model definition above, but instead counts spheres with an incoming Morse trajectory, an outgoing Morse trajectory, and a marked point directly constrained to $D$. Since $D*^\mathrm{inc}1_X$ represents $c_1(X)$, we choose this cocycle as the Morse representative $D$, so that $D*^\mathrm{inc}1_X=D$ holds at chain level.
    
    The corresponding Getzler--Gauss--Manin connection in the $q\partial_q$-direction is defined as
    \begin{equation}
        \nabla_{q\partial_q}^{\mathrm{GGM}}=q\partial_q-u^{-1}I_{q\partial_q\mu}.
    \end{equation}
    Here $\mu$ is the structure cochain for the uncurved category, retaining $q$-weights, and $I_\phi$ for a Hochschild cochain $\phi$ is defined as 
    \begin{equation}
        I_\phi=(-1)^{\lvert\phi\rvert+1}\iota_\phi+uS_\phi, \quad S_\phi=B^{1,1}(\phi,~)
    \end{equation}
    where $\iota_\phi$ is the cap product by $\phi$ and $S_\phi$ is defined as
    \begin{align}
        S_\phi(z\otimes x_n\otimes\cdots\otimes x_1)&=\sum_{0\le i\le j\le k\le n}(-1)^{\epsilon_{i,k}}\,e^+\otimes x_k\otimes\cdots\otimes x_{j+1}\otimes\phi^{j-i}(x_j,\ldots,x_{i+1})\\&\qquad\qquad\otimes x_i\otimes\cdots\otimes x_1\otimes z\otimes x_n\otimes\cdots\otimes x_{k+1}
    \end{align}
    for the sign
    \begin{equation}
        \epsilon_{i,k}=1+\|\phi\|\bigl(\dagger_{k+1}^{n}+\dagger_i+\|z\|\bigr)+\bigl(|z|+\dagger_{k+1}^{n}\bigr)(1+\dagger_k).
    \end{equation}
    This is the same operator as \cite[(75)]{Hugtenburg}, but in our right-to-left input sign convention. 
    
    Since we work in the uncurved monotone setting, there are additional codimension 1 boundary contributions from Maslov 2 discs. In the second family $H_0^2$ (see Figure \ref{fig:H_nabla^-1}), when the disc carrying $D$ has no ordinary inputs and the middle component reduces to the identity strip, it contributes $\lambda q\cdot OC$. The corresponding configurations with distinguished input $e_L^+$ give the same scalar factor, showing that
    \begin{equation}\label{eqn:connection-homotopy}
      (q\partial_q+u^{-1}(D*_q^\mathrm{inc}-\lambda q\cdot \id) )\circ OC^- =OC^-\circ(q\partial_q-u^{-1}I_{q\partial_q\mu})-\partial_X\circ H_\nabla-H_\nabla\circ(b+uB)    
    \end{equation}
    for a geometrically defined homotopy $H_\nabla$; the sign of $H_\nabla$ is implicitly chosen in the above equation. Now the Cartan homotopy formula \cite[Theorem 2.2]{Getzler} (see also \cite[Proposition 3.15]{Hugtenburg}) states that
    \begin{equation}
        [I_\phi,b+uB]=u\mathcal{L}_\phi+(-1)^{\lvert\phi\rvert}I_{\delta\phi}
    \end{equation}
    for a Hochschild cochain $\phi$, so it follows that
    \begin{equation}
        \nabla_{q\partial_q}^{\mathrm{GGM}}=q\partial_q-\mathcal{L}_{\beta}+u^{-1}[I_\beta,b+uB].
    \end{equation}
    Writing the geometric homotopy
    \begin{equation}
        H_\nabla=u^{-1}H_\nabla^{-1}+H_\nabla^0+\cdots
    \end{equation}
    and expanding \eqref{eqn:connection-homotopy} at the $u^0$-term, we obtain
    \begin{align}
        q\partial_q\circ OC-OC\circ (q\partial_q-\mathcal{L}_\beta)&=-\left(D*_q^\mathrm{inc}-\lambda q\cdot \id\right)\circ OC^1 - (H_\nabla^{-1}-OC\circ\iota_\beta)\circ B \label{eqn:connection-htpy-bdry1} \\
        &+ \partial_X\circ(-H_\nabla^0+OC^1\circ\iota_\beta+OC\circ S_\beta)+(-H_\nabla^0+OC^1\circ\iota_\beta+OC\circ S_\beta)\circ b. \label{eqn:connection-htpy-bdry2}
    \end{align}
    Since the terms in the second line \eqref{eqn:connection-htpy-bdry2} are zero in homology, it is enough to show that the terms in the first line \eqref{eqn:connection-htpy-bdry1} also vanish in homology, after quantum multiplication with the chain level representative $\Pi$ for the idempotent $\pi$.

    We now evaluate the above identity at $q=1$; hence $q\partial_q\mu$ and $q\partial_q\circ OC$ denote the $q$-derivatives, and the remaining operations are specialized at $q=1$.

    Now recall that quantum multiplication is not associative at chain level: the quantum Massey operation $m_3$ measures this failure by
    \begin{equation}
        \Pi*(D*x)-(\Pi*D)*x=\partial_X m_3(\Pi,D,x)+m_3(\Pi,D,\partial_Xx).
    \end{equation}
    Since $\partial_XC=\Pi*D-\lambda \Pi$, the Leibniz rule gives
    \begin{equation}
        (\partial_X C)*=\partial_X\circ(C*)+(C*)\circ\partial_X.
    \end{equation}
    Finally, the variant quantum Massey product $m_3^\mathrm{inc}$ satisfies
    \begin{equation}
        \partial_X\circ m_3^\mathrm{inc}(D,1_X,x)+m_3^\mathrm{inc}(D,1_X,\partial_Xx)=D*^\mathrm{inc}x-D*x.
    \end{equation}
    Combining these identities, we obtain that
    \begin{equation}
        \partial_\mathrm{hom}(\Pi*m_3^\mathrm{inc}(D,1_X,~)+m_3(\Pi,D,~)+C*~)=\Pi*(D*^\mathrm{inc}-\lambda\cdot\id).
    \end{equation}
    Using $\partial_X\circ OC^1=OC^1\circ b+OC\circ B$ we conclude that
    \begin{align}
        &\Pi*(D*^\mathrm{inc}-\lambda\cdot\id)\circ OC^1+\Pi*(H_\nabla^{-1}-OC\circ\iota_\beta)B\\
        &=\bigl( \Pi*(H_\nabla^{-1}-OC\circ\iota_\beta)+(\Pi*m_3^\mathrm{inc}(D,1_X,~)+m_3(\Pi,D,~)+C*~ )\circ OC\bigr)\circ B \label{eqn:4.142}\\
        &+\partial_\mathrm{hom}[\left(\Pi*m_3^\mathrm{inc}(D,1_X,~)+m_3(\Pi,D,~)+C*~\right)\circ OC^1]. \label{eqn:4.143}
    \end{align}
    Since \eqref{eqn:4.143} vanishes in homology, it is enough to show that \eqref{eqn:4.142} is trivial in homology.
\end{proof}

To prove the vanishing of \eqref{eqn:reduction-1} in homology from Proposition \ref{prop:reduction-1}, we recall the geometric construction of the connection homotopy $H_\nabla^{-1}$ from \cite[Lemmas 7.3--7.4]{GanatraSheridan}, for the two $k=0$ families. The first family from \cite[Lemma 7.3]{GanatraSheridan} consists of an open--closed disc joined to a sphere carrying the direct divisor incidence, with an internal Morse edge of length $l\in[0,\infty]$. As $l\to\infty$, the edge breaks and gives $D*^\mathrm{inc} OC$.

The second family is the mixed horocycle family of \cite[Lemma 7.4]{GanatraSheridan}. We choose coherent perturbation data agreeing at the $H_0^{12}$ configuration and restricting to the prescribed data for the operations at the remaining ends. The sum of the algebraic operations defined by the two families is $H_\nabla^{-1}$: see Figure \ref{fig:H_nabla^-1}.

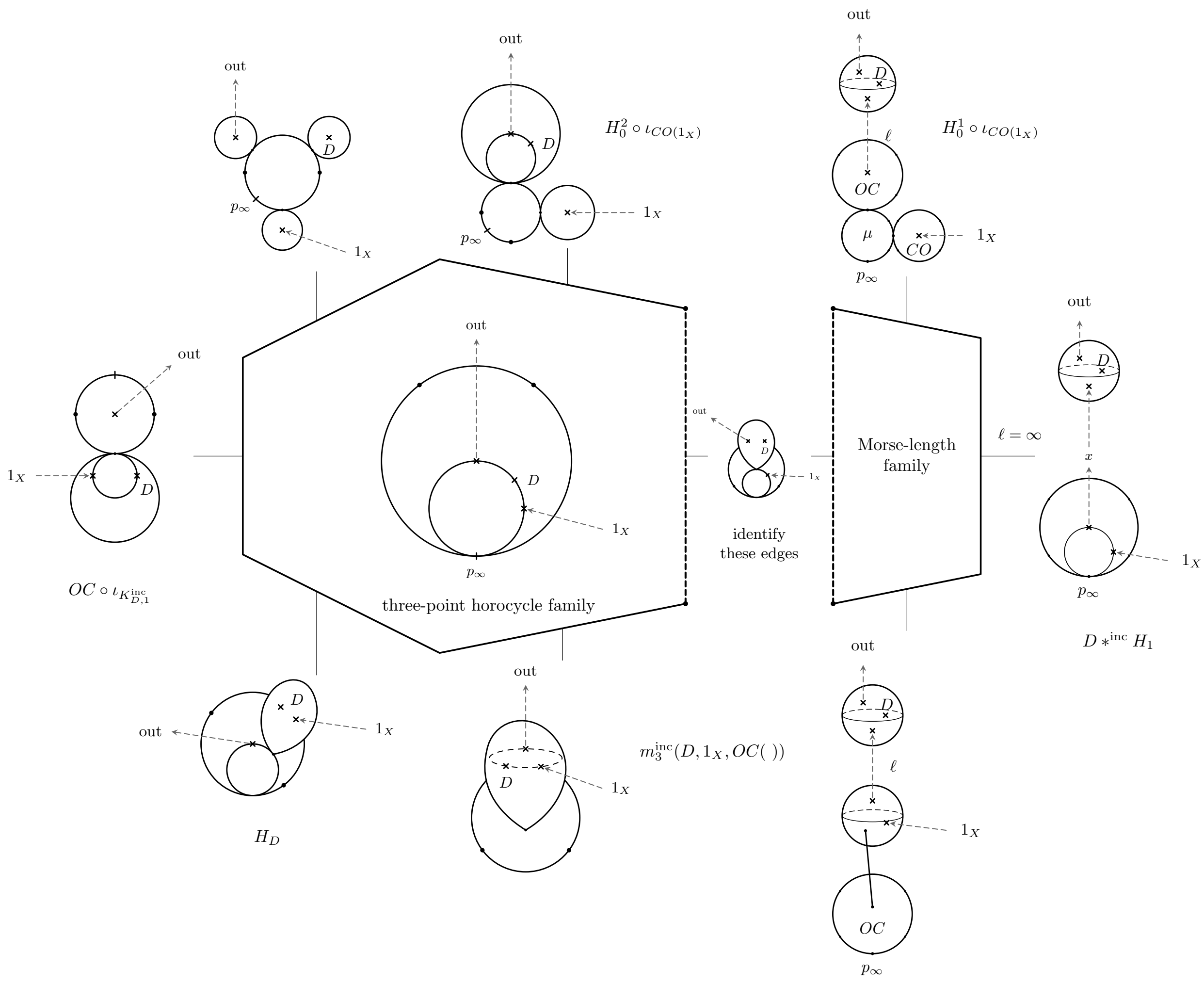
\begin{figure}[!htbp]\centering
\resizebox{\linewidth}{!}{%
\begin{tikzpicture}[
 line width=0.8pt,line cap=round,line join=round,
 morse/.style={draw=black!60,densely dashed,line width=0.55pt,->,>=stealth},
 leglabel/.style={font=\small,inner sep=2pt},
 guide/.style={black,line width=0.35pt},
 seam/.style={black,densely dashed,line width=1.15pt}
]

\draw[line width=1pt] (3,3)--(-2,4)--(-6,2)--(-6,-2)--(-2,-4)--(3,-3);
\draw[line width=1pt] (6,3)--(9,2.4)--(9,-2.4)--(6,-3);
\draw[seam] (3,-3)--(3,3);
\draw[seam] (6,-3)--(6,3);
\fill (3,3) circle (0.05);\fill (6,3) circle (0.05);
\fill (3,-3) circle (0.05);\fill (6,-3) circle (0.05);
\draw[guide] (3,0)--(3.45,0);
\draw[guide] (5.55,0)--(6,0);
\node[font=\small,align=center] at (4.5,-1.8) {identify\\these edges};
\node at (-1.0,-3.05) {three-point horocycle family};
\node[align=center] at (7.5,0) {Morse-length\\family};

\begin{scope}[shift={(-1.2500,-0.1000)},scale=0.84,transform shape]
\begin{scope}[shift={(0.0000,0.0000)},scale=1,transform shape]
\draw (0,0) circle (2.3);
\fill (-1.38,1.84) circle (0.055);
\fill (1.38,1.84) circle (0.055);
\draw (0,-1.15) circle (1.15);
\draw (-0.055,-0.055) -- (0.055,0.055);
\draw (-0.055,0.055) -- (0.055,-0.055);
\draw (0.865,-0.515) -- (0.975,-0.405);
\draw (0.865,-0.405) -- (0.975,-0.515);
\draw (1.095,-1.205) -- (1.205,-1.095);
\draw (1.095,-1.095) -- (1.205,-1.205);
\draw (0,-2.375) -- (0,-2.225);
\node[below=5pt] at (0,-2.3) {$p_\infty$};

\draw[morse,shorten <=2pt] (0,0) -- (0,2.95);
\node[above=3pt] at (0,2.95) {$\mathrm{out}$};
\node[right=5pt] at (0.92,-0.46) {$D$};
\draw[morse,shorten >=2pt] (3.05,-1.65) -- (1.15,-1.15);
\node[right=3pt] at (3.05,-1.65) {$1_X$};
\end{scope}

\end{scope}

\begin{scope}[shift={(4.5000,0.0000)},scale=0.6,transform shape]
\begin{scope}[shift={(6.1500,-2.8000)},scale=1,transform shape]

\draw (-6.25,2.35) circle (0.95);
\fill (-7.01,1.78) circle (0.045);
\fill (-5.49,1.78) circle (0.045);
\draw (-6.25,1.875) circle (0.475);
\filldraw[fill=white,draw=black]
  (-6.25,2.35)
  .. controls (-6.72,2.60) and (-6.98,3.19) .. (-6.84,3.58)
  .. controls (-6.73,4.02) and (-6.13,4.18) .. (-5.85,3.83)
  .. controls (-5.40,3.29) and (-5.72,2.59) .. (-6.25,2.35)
  -- cycle;
\fill (-6.25,2.35) circle (0.028);
\draw (-6.575,3.265) -- (-6.485,3.355);
\draw (-6.575,3.355) -- (-6.485,3.265);
\draw (-6.015,3.265) -- (-5.925,3.355);
\draw (-6.015,3.355) -- (-5.925,3.265);
\draw (-5.915,2.115) -- (-5.825,2.205);
\draw (-5.915,2.205) -- (-5.825,2.115);

\draw[morse,shorten <=2pt] (-6.53,3.31) -- (-7.80,4.10);
\node[leglabel,above left=2pt] at (-7.80,4.10) {$\mathrm{out}$};
\node[leglabel,below=4pt] at (-5.97,3.31) {$D$};
\draw[morse,shorten >=2pt] (-4.65,2.10) -- (-5.87,2.16);
\node[leglabel,right=3pt] at (-4.65,2.10) {$1_X$};

\end{scope}

\end{scope}

\begin{scope}[shift={(0.0000,6.2500)},scale=1.0,transform shape]
\begin{scope}[shift={(-0.2000,-6.0000)},scale=1,transform shape]
\draw (-0.35,6.3) circle (1);
\draw (-0.35,5.8) circle (0.5);
\draw (-0.35,4.7) circle (0.6);
\draw (0.8,4.7) circle (0.55);
\fill (-0.95,4.70) circle (0.045);
\fill (-0.35,4.10) circle (0.045);
\fill (-0.35,5.3) circle (0.028);
\fill (0.25,4.7) circle (0.028);
\draw (-0.395,6.255) -- (-0.305,6.345);
\draw (-0.395,6.345) -- (-0.305,6.255);
\draw (0.005,6.055) -- (0.095,6.145);
\draw (0.005,6.145) -- (0.095,6.055);
\draw (0.755,4.655) -- (0.845,4.745);
\draw (0.755,4.745) -- (0.845,4.655);
\draw (-0.883,4.303) -- (-0.778,4.378);
\node[font=\small,inner sep=1pt,below left=2pt]
  at (-0.83,4.34) {$p_\infty$};

\draw[morse,shorten <=2pt] (-0.35,6.3) -- (-0.35,7.95);
\node[leglabel,above=3pt] at (-0.35,7.95) {$\mathrm{out}$};
\node[leglabel,right=4pt] at (0.05,6.1) {$D$};
\draw[morse,shorten >=2pt] (2.15,4.7) -- (0.8,4.7);
\node[leglabel,right=3pt] at (2.15,4.7) {$1_X$};
\end{scope}

\end{scope}

\begin{scope}[shift={(-5.1500,6.0000)},scale=0.95,transform shape]
\begin{scope}[shift={(-6.3000,-3.6000)},scale=1,transform shape]
\draw (6.25,3.35) circle (0.8);
\draw (5.25,4.1) circle (0.45);
\draw (7.25,4.1) circle (0.45);
\draw (6.25,2.12) circle (0.43);
\fill (5.45,3.35) circle (0.045);
\fill (7.05,3.35) circle (0.045);
\fill (5.61,3.83) circle (0.028);
\fill (6.89,3.83) circle (0.028);
\fill (6.25,2.55) circle (0.028);
\draw (5.205,4.055) -- (5.295,4.145);
\draw (5.205,4.145) -- (5.295,4.055);
\draw (7.205,4.055) -- (7.295,4.145);
\draw (7.205,4.145) -- (7.295,4.055);
\draw (6.205,2.075) -- (6.295,2.165);
\draw (6.205,2.165) -- (6.295,2.075);
\draw (5.738,2.838) -- (5.630,2.730);
\node[font=\small,inner sep=1pt,below left=2pt]
  at (5.684,2.784) {$p_\infty$};

\draw[morse,shorten <=2pt] (5.25,4.1) -- (5.25,5.35);
\node[leglabel,above=3pt] at (5.25,5.35) {$\mathrm{out}$};
\node[leglabel,below=2pt] at (7.25,4.1) {$D$};
\draw[morse,shorten >=2pt] (7.60,1.65) -- (6.25,2.12);
\node[leglabel,right=3pt] at (7.60,1.65) {$1_X$};
\end{scope}

\end{scope}

\begin{scope}[shift={(-8.7000,0.0000)},scale=1.0,transform shape]
\begin{scope}[shift={(-6.0000,2.7000)},scale=1,transform shape]
\draw (6.1,-1.85) circle (0.8);
\draw (6.1,-3.55) circle (0.9);
\fill (5.30,-1.85) circle (0.045);
\fill (6.90,-1.85) circle (0.045);
\draw (6.1,-3.1) circle (0.45);
\fill (6.1,-2.65) circle (0.028);
\draw (6.1,-1.12) -- (6.1,-0.98);
\draw (6.055,-1.895) -- (6.145,-1.805);
\draw (6.055,-1.805) -- (6.145,-1.895);
\draw (6.505,-3.145) -- (6.595,-3.055);
\draw (6.505,-3.055) -- (6.595,-3.145);
\draw (5.605,-3.145) -- (5.695,-3.055);
\draw (5.605,-3.055) -- (5.695,-3.145);

\draw[morse,shorten <=2pt] (6.1,-1.85) -- (7.25,-0.85);
\node[leglabel,above right=2pt] at (7.25,-0.85) {$\mathrm{out}$};
\node[leglabel] at (6.75,-3.38) {$D$};
\draw[morse,shorten >=2pt] (4.50,-3.1) -- (5.65,-3.1);
\node[leglabel,left=3pt] at (4.50,-3.1) {$1_X$};
\end{scope}

\end{scope}

\begin{scope}[shift={(-5.5000,-5.8000)},scale=1.0,transform shape]
\begin{scope}[shift={(-0.3000,5.8000)},scale=1,transform shape]
\draw (0,-5.85) circle (1.05);
\fill (-0.84,-5.22) circle (0.045);
\fill (0.63,-6.69) circle (0.045);
\draw (0,-6.375) circle (0.525);
\filldraw[fill=white,draw=black]
  (0.42,-6.06)
  .. controls (0.12,-5.72) and (0.11,-5.21) .. (0.30,-4.91)
  .. controls (0.52,-4.48) and (1.06,-4.40) .. (1.24,-4.85)
  .. controls (1.48,-5.41) and (1.06,-5.94) .. (0.42,-6.06)
  -- cycle;
\fill (0.42,-6.06) circle (0.028);
\draw (-0.045,-5.895) -- (0.045,-5.805);
\draw (-0.045,-5.805) -- (0.045,-5.895);
\draw (0.525,-5.145) -- (0.615,-5.055);
\draw (0.525,-5.055) -- (0.615,-5.145);
\draw (0.835,-5.395) -- (0.925,-5.305);
\draw (0.835,-5.305) -- (0.925,-5.395);

\draw[morse,shorten <=2pt] (0,-5.85) -- (-1.65,-5.60);
\node[leglabel,left=3pt] at (-1.65,-5.60) {$\mathrm{out}$};
\node[leglabel] at (0.90,-4.94) {$D$};
\draw[morse,shorten >=2pt] (2.30,-5.55) -- (0.88,-5.35);
\node[leglabel,right=3pt] at (2.30,-5.55) {$1_X$};
\end{scope}

\end{scope}

\begin{scope}[shift={(0.2500,-6.5000)},scale=1.0,transform shape]
\begin{scope}[shift={(5.7000,2.4000)},scale=1,transform shape]
\draw (-6.2,-3.25) circle (1.1);
\fill (-7.08,-3.91) circle (0.045);
\fill (-5.32,-3.91) circle (0.045);
\filldraw[fill=white,draw=black]
  (-6.20,-3.50)
  .. controls (-6.83,-3.11) and (-7.15,-2.56) .. (-6.99,-1.88)
  .. controls (-6.88,-1.29) and (-6.21,-1.08) .. (-5.74,-1.49)
  .. controls (-5.11,-2.04) and (-5.35,-3.00) .. (-6.20,-3.50)
  -- cycle;
\fill (-6.20,-3.50) circle (0.028);
\draw[dashed,line width=0.55pt]
  (-6.2,-2.04) ellipse [x radius=0.74,y radius=0.19];
\draw (-6.245,-1.895) -- (-6.155,-1.805);
\draw (-6.245,-1.805) -- (-6.155,-1.895);
\draw (-6.645,-2.245) -- (-6.555,-2.155);
\draw (-6.645,-2.155) -- (-6.555,-2.245);
\draw (-5.935,-2.257) -- (-5.845,-2.167);
\draw (-5.935,-2.167) -- (-5.845,-2.257);

\draw[morse,shorten <=2pt] (-6.2,-1.85) -- (-6.2,-0.55);
\node[leglabel,above=3pt] at (-6.2,-0.55) {$\mathrm{out}$};
\node[leglabel,below=4pt] at (-6.6,-2.2) {$D$};
\draw[morse,shorten >=2pt] (-4.65,-2.65) -- (-5.89,-2.212);
\node[leglabel,right=3pt] at (-4.65,-2.65) {$1_X$};
\end{scope}

\end{scope}

\begin{scope}[shift={(6.7000,4.1000)},scale=0.95,transform shape]
\draw[] (0.0000,1.7000) circle (0.75);
\fill (0.6144,2.1302) circle (0.024);
\fill (-0.6144,2.1302) circle (0.024);
\fill (-0.6144,1.2698) circle (0.024);
\fill (0.6144,1.2698) circle (0.024);
\node[] at (0.0000,1.4000) {$OC$};
\draw[] (0.0000,3.6500) circle (0.6);
\draw[black,line width=0.4pt,densely dashed] (-0.6000,3.6500) arc[start angle=180,end angle=0,x radius=0.6,y radius=0.12];
\draw[black,line width=0.4pt] (-0.6000,3.6500) arc[start angle=180,end angle=360,x radius=0.6,y radius=0.12];
\draw[] (0.2050,3.6050)--(0.2950,3.6950);
\draw[] (0.2050,3.6950)--(0.2950,3.6050);
\node[] at (0.2700,3.8800) {$D$};
\draw[] (-0.2250,3.8550)--(-0.1350,3.9450);
\draw[] (-0.2250,3.9450)--(-0.1350,3.8550);
\draw[morse,shorten <=2pt,shorten >=2pt] (-0.1800,3.9000)--(-0.1800,4.8000);
\node[above=3pt] at (-0.1800,4.8000) {$\mathrm{out}$};
\draw[] (-0.0450,1.7050)--(0.0450,1.7950);
\draw[] (-0.0450,1.7950)--(0.0450,1.7050);
\draw[] (-0.0450,3.2850)--(0.0450,3.3750);
\draw[] (-0.0450,3.3750)--(0.0450,3.2850);
\draw[morse,shorten <=2pt,shorten >=2pt] (0.0000,1.7500)--(0.0000,3.3300);
\node[anchor=west] at (0.2400,2.5400) {$\ell$};
\draw[] (0.0000,0.4000) circle (0.55);
\fill (0.4505,0.7155) circle (0.024);
\fill (-0.4505,0.7155) circle (0.024);
\fill (-0.4505,0.0845) circle (0.024);
\fill (0.4505,0.0845) circle (0.024);
\draw[] (1.1000,0.4000) circle (0.55);
\fill (1.5505,0.7155) circle (0.024);
\fill (0.6495,0.7155) circle (0.024);
\fill (0.6495,0.0845) circle (0.024);
\fill (1.5505,0.0845) circle (0.024);
\fill (0.0000,0.9500) circle (0.032);
\fill (0.5500,0.4000) circle (0.032);
\fill (0.0000,-0.1500) circle (0.032);
\node[] at (0.0000,0.4000) {$\mu$};
\node[below=3pt] at (0.0000,-0.1500) {$p_\infty$};
\draw[] (1.0550,0.3550)--(1.1450,0.4450);
\draw[] (1.0550,0.4450)--(1.1450,0.3550);
\draw[morse,shorten <=2pt,shorten >=2pt] (2.1000,0.4000)--(1.1000,0.4000);
\node[right=3pt] at (2.1000,0.4000) {$1_X$};
\node[] at (1.1000,0.1000) {$CO$};
\end{scope}

\begin{scope}[shift={(11.2000,-2.4000)},scale=0.95,transform shape]
\draw[] (0.0000,4.3500) circle (0.65);
\draw[black,line width=0.4pt,densely dashed] (-0.6500,4.3500) arc[start angle=180,end angle=0,x radius=0.65,y radius=0.13];
\draw[black,line width=0.4pt] (-0.6500,4.3500) arc[start angle=180,end angle=360,x radius=0.65,y radius=0.13];
\draw[] (0.2350,4.3050)--(0.3250,4.3950);
\draw[] (0.2350,4.3950)--(0.3250,4.3050);
\node[] at (0.3000,4.5800) {$D$};
\draw[] (-0.2450,4.5750)--(-0.1550,4.6650);
\draw[] (-0.2450,4.6650)--(-0.1550,4.5750);
\draw[morse,shorten <=2pt,shorten >=2pt] (-0.2000,4.6200)--(-0.2000,5.5000);
\node[above=3pt] at (-0.2000,5.5000) {$\mathrm{out}$};
\draw[] (0.0000,1.0000) circle (1.05);
\fill (0.8601,1.6023) circle (0.024);
\fill (-0.8601,1.6023) circle (0.024);
\fill (-0.8601,0.3977) circle (0.024);
\fill (0.8601,0.3977) circle (0.024);
\fill (0.0000,-0.0500) circle (0.032);
\node[below=3pt] at (0.0000,-0.0500) {$p_\infty$};
\draw[line width=0.5pt] (0.0000,0.4750) circle (0.525);
\draw[] (0.4800,0.4300)--(0.5700,0.5200);
\draw[] (0.4800,0.5200)--(0.5700,0.4300);
\draw[morse,shorten <=2pt,shorten >=2pt] (1.7500,0.2800)--(0.5250,0.4750);
\node[right=3pt] at (1.7500,0.2800) {$1_X$};
\draw[] (-0.0450,0.9550)--(0.0450,1.0450);
\draw[] (-0.0450,1.0450)--(0.0450,0.9550);
\draw[] (-0.0450,3.9750)--(0.0450,4.0650);
\draw[] (-0.0450,4.0650)--(0.0450,3.9750);
\draw[morse,shorten <=2pt,shorten >=2pt] (0.0000,1.0000)--(0.0000,2.3700);
\draw[morse,shorten <=2pt,shorten >=2pt] (0.0000,2.6500)--(0.0000,4.0200);
\node[fill=white,inner sep=1pt,font=\footnotesize] at (0.0000,2.5100) {$x$};
\end{scope}

\begin{scope}[shift={(6.8000,-9.4000)},scale=0.95,transform shape]
\draw[] (0.0000,0.1000) circle (0.85);
\fill (0.6963,0.5875) circle (0.024);
\fill (-0.6963,0.5875) circle (0.024);
\fill (-0.6963,-0.3875) circle (0.024);
\fill (0.6963,-0.3875) circle (0.024);
\fill (0.0000,-0.7500) circle (0.032);
\node[below=3pt] at (0.0000,-0.7500) {$p_\infty$};
\node[] at (0.0000,-0.2000) {$OC$};
\draw[] (0.0000,2.2000) circle (0.65);
\draw[black,line width=0.4pt,densely dashed] (-0.6500,2.2000) arc[start angle=180,end angle=0,x radius=0.65,y radius=0.13];
\draw[black,line width=0.4pt] (-0.6500,2.2000) arc[start angle=180,end angle=360,x radius=0.65,y radius=0.13];
\draw[] (0.2550,2.0050)--(0.3450,2.0950);
\draw[] (0.2550,2.0950)--(0.3450,2.0050);
\draw[morse,shorten <=2pt,shorten >=2pt] (1.6500,1.8800)--(0.3000,2.0500);
\node[right=3pt] at (1.6500,1.8800) {$1_X$};
\draw[] (0.0000,4.3500) circle (0.65);
\draw[black,line width=0.4pt,densely dashed] (-0.6500,4.3500) arc[start angle=180,end angle=0,x radius=0.65,y radius=0.13];
\draw[black,line width=0.4pt] (-0.6500,4.3500) arc[start angle=180,end angle=360,x radius=0.65,y radius=0.13];
\draw[] (0.2350,4.3050)--(0.3250,4.3950);
\draw[] (0.2350,4.3950)--(0.3250,4.3050);
\node[] at (0.3000,4.5800) {$D$};
\draw[] (-0.2450,4.5750)--(-0.1550,4.6650);
\draw[] (-0.2450,4.6650)--(-0.1550,4.5750);
\draw[morse,shorten <=2pt,shorten >=2pt] (-0.2000,4.6200)--(-0.2000,5.5000);
\node[above=3pt] at (-0.2000,5.5000) {$\mathrm{out}$};
\draw[] (0.0000,0.2500)--(-0.1500,1.8800);
\fill (0.0000,0.2500) circle (0.032);
\fill (-0.1500,1.8800) circle (0.032);
\draw[] (-0.0450,2.4750)--(0.0450,2.5650);
\draw[] (-0.0450,2.5650)--(0.0450,2.4750);
\draw[] (-0.0450,3.9750)--(0.0450,4.0650);
\draw[] (-0.0450,4.0650)--(0.0450,3.9750);
\draw[morse,shorten <=2pt,shorten >=2pt] (0.0000,2.5200)--(0.0000,4.0200);
\node[anchor=west] at (0.2400,3.2700) {$\ell$};
\end{scope}

\draw[guide] (0.6,3.48)--(0.6,4.22);
\draw[guide] (-4.5,2.75)--(-4.5,3.75);
\draw[guide] (-6,0)--(-7.0,0);
\draw[guide] (-4.5,-2.75)--(-4.5,-4.45);
\draw[guide] (0.5,-3.5)--(0.5,-4.15);
\draw[guide] (7.5,2.7)--(7.5,3.65);
\draw[guide] (9,0)--(10.1,0);
\draw[guide] (7.5,-2.7)--(7.5,-3.55);

\node[font=\small,fill=white,inner sep=2pt] at (9.8,.45) {$\ell=\infty$};

\node at (2.35,6.65) {$H_0^2\circ\iota_{CO(1_X)}$};
\node at (9.2,6.65) {$H_0^1\circ\iota_{CO(1_X)}$};
\node at (-8.7,-2.8) {$OC\circ\iota_{K^{\mathrm{inc}}_{D,1}}$};
\node at (-5.5,-7.75) {$H_D$};
\node at (3.55,-6.0) {$m_3^\mathrm{inc}(D,1_X,OC(~))$};
\node at (11.8,-3.75) {$D*^\mathrm{inc}H_{1}$};

\end{tikzpicture}

}
\caption{The two families whose sum defines the homotopy $S_{D,1}^\mathrm{inc}$.}
\label{fig:Q_D}
\end{figure}

We compare $H_\nabla^{-1}$ with the Morse-Floer module homotopy $H_D$ using a three-point horocycle family with mixed incidence conditions, namely one direct divisor incidence and one Morse unit input. Define the operation $S_{D,1}^\mathrm{inc}$ as described above: the domain is drawn in Figure \ref{fig:Q_D}, with its six boundary degenerations. The homotopy $H_\nabla^{-1}$ used is from \eqref{eqn:connection-homotopy}, with orientations chosen such that its boundary is
\begin{equation}
    \partial_{\mathrm{hom}}H_\nabla^{-1}=OC\circ\iota_{q\partial_q\mu}-(D*^{\mathrm{inc}}-\lambda\,\mathrm{id})\circ OC.
\end{equation}
The term $\lambda\cdot OC$ is the Maslov-two divisor-bubble contribution described above, specialized at $q=1$.

We claim that the following equation holds:
\begin{equation}\label{eqn:S_D^inc}
    \partial_\mathrm{hom}(S_{D,1}^\mathrm{inc}\circ B)=(H_D-OC\circ\iota_{K_{D,1}^\mathrm{inc}}-H_\nabla^{-1}\circ\iota_{CO(1_X)}-m_3^\mathrm{inc}(D,1_X,OC(~))-D*^\mathrm{inc}H_1)\circ B.
\end{equation}

For Figure \ref{fig:Q_D}, the leftmost face represents $OC\circ\iota_{K_{D,1}^\mathrm{inc}}$, and the top-left face consisting of four discs contributes zero after composition with the Connes $B$ operator. Since $D*^\mathrm{inc}1_X=D$, the bottom-left face in Figure \ref{fig:Q_D} can be identified with $H_D$. The two families on the top-right can be identified with the two components that define $H_\nabla^{-1}$ precomposed with $\iota_{CO(1_X)}$, while the two families on the bottom right add up to the variant quantum Massey product $m_3^\mathrm{inc}(D,1_X,OC(~))$. On the four-disc face after composition with $B$, we forget the $e^+$-marking before contracting any input-free divisor disc, so the strict unit with three attaching nodes gives zero from the forgetful data. Finally, breaking the internal Morse edge in the attached Morse-length family gives the composition $D*^\mathrm{inc}\circ H_1$, contributing the term $-D*^\mathrm{inc}\circ H_1\circ B$. 

\begin{proposition}\label{prop:reduction-2}
    The proof of Proposition \ref{prop:reduction-1} reduces to showing that
    \begin{equation}\label{eqn:reduction-2-main-conclusion}
        (\Pi*H_D+(m_3(\Pi,D,~)+C*)\circ OC-\lambda\Pi*H_1
        -\Pi*OC\circ\iota_{CO(C)-K_{\Pi,D}-\beta_\Pi\cup CO(D)+K_{1,D}+\lambda\beta_\Pi})\circ B
    \end{equation}
    is zero in homology.
\end{proposition}
\begin{proof}
    Recall that we aim to show
    \begin{equation}\label{eqn:reduction-11}
     (\Pi*H_\nabla^{-1}-\Pi*OC\circ\iota_{\beta_\mathrm{geom}}+(\Pi*m_3^\mathrm{inc}(D,1_X,~)+m_3(\Pi,D,~)+C*)\circ OC)\circ B
    \end{equation}
    is zero in homology. We first verify that the operation before $B$ is closed:
    \begin{align}
        &\partial_\mathrm{\hom}(\Pi*H_\nabla^{-1}-\Pi*OC\circ\iota_{\beta_\mathrm{geom}}+(\Pi*m_3^\mathrm{inc}(D,1_X,~)+m_3(\Pi,D,~)+C*)\circ OC)\\
        &=-\Pi*(D*^\mathrm{inc}-\lambda\cdot\id)\circ OC+\Pi*(D*^\mathrm{inc}-\lambda\cdot\id)\circ OC=0.
    \end{align}
    Since $\iota_{CO(1_X)}$ induces the identity on Hochschild homology, inserting it before $B$ does not change the induced map on homology. Thus it is enough to show that the precomposition
    \begin{equation}\label{eqn:reduction-12}
     (\Pi*H_\nabla^{-1}-\Pi*OC\circ\iota_{\beta_\mathrm{geom}}+(\Pi*m_3^\mathrm{inc}(D,1_X,~)+m_3(\Pi,D,~)+C*)\circ OC)\iota_{CO(1_X)}\circ B
    \end{equation}
    is trivial in homology. We first note that by the Leibniz rule,
    \begin{align}
        &\partial_\mathrm{hom}\bigl[ (\Pi*m_3^\mathrm{inc}(D,1_X,~)+m_3(\Pi,D,~)+C*)\circ H_1\circ B \bigr]\\
        &= \Pi*(D*^\mathrm{inc}-\lambda)\circ H_1\circ B-(\Pi*m_3^\mathrm{inc}(D,1_X,~)+m_3(\Pi,D,~)+C*)\circ (OC\circ\iota_{CO(1_X)}-OC)\circ B.
    \end{align}
    Thus \eqref{eqn:reduction-12} vanishes on homology if and only if the following expression does:
    \begin{equation}\label{eqn:reduction-13}
     (\Pi*H_\nabla^{-1}\circ\iota_{CO(1_X)}-\Pi*OC\circ\iota_{\beta_\mathrm{geom}\cup CO(1_X)}+\Pi*(D*^\mathrm{inc}-\lambda)\circ H_1+(\Pi*m_3^\mathrm{inc}(D,1_X,~)+m_3(\Pi,D,~)+C*)\circ OC)\circ B.
    \end{equation}
    Now from \eqref{eqn:S_D^inc}, we find that it is enough to show
    \begin{equation}
        (\Pi*H_D+(m_3(\Pi,D,~)+C*)\circ OC-\lambda\Pi*H_1-\Pi*OC\circ\iota_{K_{D,1}^\mathrm{inc}+\beta_{geom}\cup CO(1_X)})\circ B
    \end{equation}
    is trivial in homology. We recall from the definition of $\beta_\mathrm{geom}$ that
    \begin{align}
        &K_{D,1}^\mathrm{inc}+\beta_\mathrm{geom}\cup CO(1_X)-(CO(C)-K_{\Pi,D}-\beta_\Pi\cup CO(D)+K_{1,D}+\lambda\beta_\Pi)\\
        =&(\beta_\mathrm{geom}-\beta_0)\cup CO(1_X)-(\beta_\mathrm{geom}-\beta_0),
    \end{align}
    and this is a Hochschild coboundary because $\beta_\mathrm{geom}-\beta_0$ is closed ($\delta\beta_\mathrm{geom}=\delta\beta_0=q\partial_q\mu$) and the operation $\cup CO(1_X)$ is homotopic to the identity. This completes our reduction to \eqref{eqn:reduction-2-main-conclusion}.
\end{proof}

\begin{proposition}\label{prop:reduction-3}
    The proof of Proposition \ref{prop:reduction-2} reduces to showing that
    \begin{equation}\label{eqn:reduction-3-main-conclusion}
        \lambda\bigl( (H_\Pi-H_1-OC\circ\iota_{\beta_\Pi})\circ(\id-\iota_{CO(1_X)})+(\id-\Pi*)\circ H_1 \bigr)\circ B
    \end{equation}
    is zero in homology.
\end{proposition}
\begin{proof}
    Recall that we have to show
    \begin{equation}\label{eqn:reduction-3-main-cancellation}
        (\Pi*H_D+(m_3(\Pi,D,~)+C*)\circ OC-\lambda\Pi*H_1
        -\Pi*OC\circ\iota_{CO(C)-K_{\Pi,D}-\beta_\Pi\cup CO(D)+K_{1,D}+\lambda\beta_\Pi})\circ B
    \end{equation}
    acts trivially on homology. From the module homotopy equations,
    \begin{align}
        \partial_\mathrm{hom}(S_{\Pi,D}\circ B)&=\left(H_\Pi\circ\iota_{CO(D)}+\Pi*H_D-H_{\Pi*D}+OC\circ\iota_{K_{\Pi,D}}+m_3(\Pi,D,~)\circ OC\right)\circ B\\
        \partial_\mathrm{hom}(-H_C\circ B)&=\left( H_{\Pi*D}-\lambda H_\Pi+C*OC-OC\circ\iota_{CO(C)} \right)\circ B\\
        \partial_\mathrm{hom}(-S_{1,D}\circ B)&=-\left(H_1\circ\iota_{CO(D)}+1_X*H_D-H_{1_X*D}+OC\circ\iota_{K_{1,D}}+m_3(1_X,D,~)\circ OC\right)\circ B,
    \end{align}
    and from our choice of perturbation data we see that
    \begin{equation}
        1_X*H_D=H_D,\quad 1_X*D=D,\quad m_3(1_X,D,~)=0.
    \end{equation}
    Therefore we add the above equations to find
    \begin{align}
        \partial_\mathrm{hom}((S_{\Pi,D}-S_{1,D}-H_C)\circ B)= \Bigl[ &\Bigl( \Pi*H_D+(m_3(\Pi,D,~)+C*)\circ OC-\lambda H_\Pi \label{eqn:reduction-31} \\
        &+(H_\Pi-H_1)\circ\iota_{CO(D)}+OC\circ\iota_{K_{\Pi,D}-K_{1,D}-CO(C)}\Bigr) \Bigr]\circ B. \label{eqn:reduction-32}
    \end{align}
    Comparing this to \eqref{eqn:reduction-3-main-cancellation} we see that there is a missing $\Pi*$ term for the cap product with $K_{\Pi,D}$ and similar terms. The homotopies below will insert the missing $\Pi*$ projection into the cap-product terms:
    \begin{align}
         \partial_\mathrm{hom}(H_\Pi-H_1-OC\circ\iota_{\beta_\Pi})&=(\id-\Pi*)\circ OC,\\
        \delta\bigl(CO(C)-K_{\Pi,D}-\beta_\Pi\cup CO(D)+K_{1,D}+\lambda\beta_\Pi\bigr)&=CO(D)-\lambda CO(1_X).
    \end{align}
    This leads to 
    \begin{align}
        &\partial_\mathrm{hom}\Bigl[(H_\Pi-H_1-OC\circ\iota_{\beta_\Pi})\circ\iota_{CO(C)-K_{\Pi,D}-\beta_\Pi\cup CO(D)+K_{1,D}+\lambda\beta_\Pi}\circ B\Bigr] \label{eqn:reduction-33} \\
        &=(\id-\Pi*)\circ OC\circ\iota_{CO(C)-K_{\Pi,D}-\beta_\Pi\cup CO(D)+K_{1,D}+\lambda\beta_\Pi}\circ B \label{eqn:reduction-34} \\
        &\quad-(H_\Pi-H_1-OC\circ\iota_{\beta_\Pi})\circ\iota_{CO(D)-\lambda CO(1_X)}\circ B. \label{eqn:reduction-35}
    \end{align}
    Subtracting the sum of the boundaries \eqref{eqn:reduction-31}-\eqref{eqn:reduction-32} and \eqref{eqn:reduction-33}--\eqref{eqn:reduction-35} from \eqref{eqn:reduction-3-main-cancellation} gives \eqref{eqn:reduction-3-main-conclusion}.
\end{proof}

The two-unit homotopy $S_{1,1}$ provides the final remaining cancellation.

\begin{proposition}\label{prop:OC-compatiblity-end-of-proof}
    The previous expression
    \begin{equation}\label{eqn:cancellation}
        \lambda\bigl( (H_\Pi-H_1-OC\circ\iota_{\beta_\Pi})\circ(\id-\iota_{CO(1_X)})+(\id-\Pi*)\circ H_1 \bigr)\circ B
    \end{equation}
    vanishes in homology.
\end{proposition}
\begin{proof}
    We first check that the operator before $B$ in \eqref{eqn:cancellation} is closed. Indeed,
    \begin{equation}
        \partial_\mathrm{hom}(H_\Pi-H_1-OC\circ\iota_{\beta_\Pi})=(\id-\Pi*)\circ OC,
    \end{equation}
    and therefore
    \begin{align}
        &\partial_\mathrm{hom}\left ((H_\Pi-H_1-OC\circ\iota_{\beta_\Pi})\circ(\id-\iota_{CO(1_X)})+(\id-\Pi*)\circ H_1\right)\\
        &=(\id-\Pi*)\circ OC\circ(\id-\iota_{CO(1_X)})+(\id-\Pi*)\circ(OC\circ\iota_{CO(1_X)}-OC)=0.
    \end{align}
    Again since $\iota_{CO(1_X)}$ induces the identity on Hochschild homology, it is enough to show that
    \begin{equation}\label{eqn:cancellation2}
        \lambda\bigl( (H_\Pi-H_1-OC\circ\iota_{\beta_\Pi})\circ(\id-\iota_{CO(1_X)})+(\id-\Pi*)\circ H_1 \bigr)\circ \iota_{CO(1_X)} \circ B
    \end{equation}
    vanishes in homology. Recall that
    \begin{align}
        \partial_\mathrm{hom}(S_{1,1}\circ B)&=(H_1\circ\iota_{CO(1_X)}+OC\circ\iota_{K_{1,1}})\circ B,\\
        \partial_\mathrm{hom}(\iota_{K_{1,1}}\circ B)
        &=(\iota_{CO(1_X)}-\iota_{CO(1_X)}\circ\iota_{CO(1_X)})\circ B.
    \end{align}
    Therefore by the Leibniz rule, we conclude that
    \begin{align}
        &\partial_\mathrm{hom}\Bigl[ (\id-\Pi*)\circ S_{1,1}\circ B-(H_\Pi-H_1-OC\circ\iota_{\beta_\Pi})\circ\iota_{K_{1,1}}\circ B \Bigr]\\
        &=(\id-\Pi*)\circ (H_1\circ\iota_{CO(1_X)}+OC\circ\iota_{K_{1,1}})\circ B-(\id-\Pi*)\circ OC\circ\iota_{K_{1,1}}\circ B\\
        &+(H_\Pi-H_1-OC\circ\iota_{\beta_\Pi})\circ(\iota_{CO(1_X)}-\iota_{CO(1_X)\cup CO(1_X)})\circ B\\
        &=\Bigl((H_\Pi-H_1-OC\circ\iota_{\beta_\Pi})\circ(\id-\iota_{CO(1_X)})+(\id-\Pi*)\circ H_1\Bigr)\circ\iota_{CO(1_X)}\circ B.
    \end{align}
    Thus the expression after insertion of $\iota_{CO(1_X)}$ is homotopic to zero, and \eqref{eqn:cancellation} vanishes in homology.
\end{proof}

This concludes the proof of Theorem \ref{thm:OC-compatibility-reprise}.

\begin{remark}
    Although the proof above is written for Lagrangians that admit a lift to the relative Fukaya category, over $\mathbb{C}$ this hypothesis can be removed by using the capping construction as in the proof of Theorem \ref{thm:monotone-fukaya-gradeability-reprise} and an area-weighted deformation.

    Weight each holomorphic curve by its area $q^{\int u^*\omega}$, and let $G$ be defined as
    \begin{equation}
        G^1(x)=(\deg(x,w)-2\langle \omega,w\rangle)x
    \end{equation}
    for each generator $x$ and capping $w$. Then
    \begin{align}
        \delta G&=\kappa-2q\partial_q\mu,\\
        OC\circ(\mathcal L_G-\Lambda)
        &=(\mathrm{gr}_X-n)\circ OC+2q\partial_qOC.
    \end{align}
    Moreover, replacing the direct divisor incidence with $D$ by insertion of the closed two-form $\omega$, we can show
    \begin{equation}
        \delta K_{D,1}^{\mathrm{inc}}=CO(D)-(q\partial_q\mu)\cup CO(1_X)-\lambda CO(1_X)
    \end{equation}
    holds again, so defining $\beta_\mathrm{geom}$ by the same formula shows
    \begin{equation}
        \delta\beta_{\mathrm{geom}}=q\partial_q\mu,
    \end{equation}
    and thus $G+2\beta_\mathrm{geom}$ is an infinitesimal grading. The proof of the open--closed compatibility goes through for this grading as well, after taking the insertion by the form $\omega$ instead of direct incidence with $D$.
\end{remark}

\section{Examples} \label{sec:examples}

In this section, we establish the existence of infinitesimal gradings on monotone Fukaya categories of certain Fano varieties, and explore their topological consequences. We will consider these Fano varieties as complex Kähler manifolds, equipped with a monotone symplectic form $\omega$ normalized such that $[\omega]=c_1$. For the quantum products, we will denote by $\alpha*_q\beta$ the small quantum product valued in cohomology with Novikov coefficients, and $\alpha*\beta$ the small quantum product with the Novikov variable $q$ set to $1$. Throughout this section, the base field is $\mathbb{C}$.

Each example will follow the same roadmap. We first identify the relevant $c_1*$-eigensummand and its ring structure. The key step in each example is to compute the closed-string infinitesimal degrees of a convenient basis. Finally, using the closed-open and open-closed string maps we determine the infinitesimal dimensions of graded Lagrangian spheres and apply the results from previous sections.

\subsection{MM 2-25} \label{ssec:MM2-25}

Our first example will be the Fano threefold MM 2-25, as described in \cite{Fano3}. This is the blowup of $\mathbb{P}^3$ along an elliptic curve, given by the intersection of two quadric hypersurfaces. Alternatively, this can be described as a bidegree $(1,2)$ hypersurface in $\mathbb{P}^1\times\mathbb{P}^3$. We note that \cite[Example 1.12]{Seidel23} provides an interesting observation about the arithmetic properties of the Fukaya category of this Fano variety.

Our main result for the MM 2-25 threefold is as follows.

\begin{thm} \label{thm:MM-2-25-reprise}
    Let $X$ be a smooth bidegree $(1,2)$ hypersurface in $\mathbb{P}^1\times\mathbb{P}^3$, equipped with its monotone Kähler form. With its closed-string infinitesimal grading, the generalized $(-1)$-eigensummand $QH^*(X)_{-1}$ of quantum cohomology admits an isomorphism of graded $\mathbb{C}$-algebras
    \begin{equation}
        QH^*(X)_{-1}\cong H^*(E;\mathbb{C})
    \end{equation}
    for an elliptic curve $E$ with its standard cohomological grading.

    There is a choice of geometric infinitesimal grading on $\mathcal{F}(X)_{-1}$, for which homogeneous weak proper Calabi--Yau structures exist precisely in infinitesimal dimension $1$. Moreover, every Lagrangian sphere $L\subset X$ with $[L]\neq0$ in $H_3(X;\mathbb{C})$ which admits a graded lift with respect to the chosen geometric infinitesimal grading has infinitesimal dimension $1$.
\end{thm}

To prove the above theorem, we will first understand the infinitesimal gradings associated to the $(-1)$-eigensummand of the quantum cohomology of $X$. The proposition below follows from Proposition \ref{prop:qh-MM225}:

\begin{proposition} \label{prop:QH^*_-1-MM225}
    The $c_1*$ action on the quantum cohomology of $X$ has the eigenvalues
    \begin{equation}
        \lambda=-1,0,8,-4\pm4i,
    \end{equation}
    and is diagonalizable. The $(-1)$-eigensummand has rank $4$, and the other eigensummands have rank $1$. A convenient basis for the $(-1)$-eigensummand is
    \begin{align}
        &\pi=\frac{1}{225}(88-17h_1-30h_2+68h_1h_2-76h_2^2+32[\mathrm{pt}]),\\
        &\nu=-16-16h_1+4h_1h_2-8h_2^2+16[\mathrm{pt}],
    \end{align}
    and the two odd-degree classes $\psi_1$, $\psi_2$. Their quantum products are as follows:
    \begin{equation}
    \begin{array}{c|cccc}
    * 
    & \pi 
    & \psi_1 
    & \psi_2
    & \nu 
    \\
    \hline
    \pi 
    & \pi 
    & \psi_1 
    & \psi_2
    & \nu
    \\[2mm]

    \psi_1
    & \psi_1
    & 0 
    & -\frac{\nu}{16} 
    & 0 
    \\[2mm]

    \psi_2
    & \psi_2
    & \frac{\nu}{16} 
    & 0 
    & 0    
    \\[2mm]

    \nu
    & \nu
    & 0
    & 0
    & 0
    \\[2mm]

    \end{array}
    \end{equation}
    In other words, $\pi$ is the idempotent for the $(-1)$-eigensummand, and $c_1(X)*\pi_{-1}=-\pi_{-1}$: therefore the monotone Fukaya category $\mathcal{F}(X)_{-1}$ admits geometric infinitesimal gradings as in Definition \ref{def:geometric-inf-grading}.
\end{proposition}

Proposition \ref{prop:QH^*_-1-MM225} now allows us to compute the grading operator.

\begin{proposition}
    The grading operator on the $-1$-summand satisfies
    \begin{align}
        &\Gamma_{-1}(\pi)=0,\\
        &\Gamma_{-1}(\psi_i)=\psi_i,\\
        &\Gamma_{-1}(\nu)=2\nu.
    \end{align}
\end{proposition}
    
\begin{proof}
    Recall from Definition \ref{def:closed-string-grading-operator} and \eqref{eqn:GW_4-qdq} that the grading operator can be computed as 
    \begin{equation} \label{eqn:closed-string-grading-operator-MM225}
        \Gamma_{-1}(x)=\pi*(gr(x)-2q\frac{d}{dq}(\pi*_qx)\vert_{q=1}).
    \end{equation}
    The idempotent has closed-string infinitesimal grading zero from Proposition \ref{prop:idempotent-zero}, so it remains to compute the other three classes. We start with the odd classes: from the computation of the small quantum products of the odd classes in Proposition \ref{prop:qh-MM225},
    \begin{align}
        \pi*_q\psi_i&= \frac{1}{225}(88-17h_1-30h_2+68h_1h_2-76h_2^2+32[\mathrm{pt}])*_q\psi_i\\
        &=\frac{1}{225}(88+17q+0+0+152q^2-32q^3)\psi_i,
    \end{align}
    and therefore the closed-string infinitesimal gradings of $\psi_i$ can be computed as
    \begin{align}
        \Gamma_{-1}(\psi_i)&=\pi*(3\psi_i-\frac{2}{225}(17+304-96)\psi_i)\\
        &=\psi_i.
    \end{align}

    Similarly, for $\nu$ we first compute that
    \begin{align}
        \pi*_q\nu&=\frac{1}{225}(88-17h_1-30h_2+68h_1h_2-76h_2^2+32[\mathrm{pt}])*_q(-16-16h_1+4h_1h_2-8h_2^2+16[\mathrm{pt}])\\
        &=\frac{16}{225}(192q^6-432q^5+24q^4+79q^3-88)\cdot 1+\frac{16}{225}(96q^5-184q^4-79q^3+30q^2-17q-71)\cdot h_1\\
        &+\frac{32}{225}(32q^5+30q^4-109q^3+15q^2+17q+15)\cdot h_2-\frac{4}{225}(64q^3-64q^2-289q+64)h_1h_2\\
        &+\frac{8}{225}(64q^4-64q^3-289q+64)h_2^2+16[pt],
    \end{align}
    and therefore
    \begin{align}
        q\frac{d}{dq}(\pi*_q\nu)&=\frac{16}{225}(1152q^6-2160q^5+96q^4+237q^3)\cdot 1+\frac{16}{225}(480q^5-736q^4-237q^3+60q^2-17q)\cdot h_1\\
        &+\frac{32}{225}(160q^5+120q^4-327q^3+30q^2+17q)\cdot h_2+\frac{4}{225}(-192q^3+128q^2+289q)h_1h_2\\
        &+\frac{8}{225}(256q^4-192q^3-289q)h_2^2.
    \end{align}
    Therefore the closed-string infinitesimal grading of $\nu$ is
    \begin{align}
        \Gamma_{-1}(\nu)&=\pi*(-32h_1+16h_1h_2-32h_2^2+96[pt]-2(-48-32h_1+4h_1h_2-8h_2^2))\\
        &=\frac{1}{225}(88-17h_1-30h_2+68h_1h_2-76h_2^2+32[\mathrm{pt}])*(96+32h_1+8h_1h_2-16h_2^2+96[pt])\\
        &=\frac{88}{225}(96+32h_1+8h_1h_2-16h_2^2+96[pt])-\frac{17}{225}(160+160h_1+256h_2-8h_1h_2+16h_2^2-32[pt])\\
        &-\frac{30}{225}(64+64h_1+256h_2+128h_1h_2)-\frac{76}{225}(832+192h_1+64h_2+48h_1h_2+288h_2^2+64[pt])\\
        &+\frac{68}{225}(512+128h_2+128h_1h_2+256h_2^2)+\frac{32}{225}(544+288h_1+256h_2-8h_1h_2+80h_2^2+96[pt])\\
        &=-32-32h_1+8h_1h_2-16h_2^2+32[pt]=2\nu.
    \end{align} 
    Finally, since $\pi*gr(\pi)=2\pi+2/225\nu$ and $\pi*gr(\nu)=4\nu$, the class
    \begin{equation}
        \theta=\pi-\frac{1}{225}\nu
    \end{equation}
    satisfies $\pi*gr(\theta)=2\theta$, providing the required class for Theorem \ref{thm:intro-cy-criterion}. Therefore we conclude that this class provides a weak proper Calabi-Yau structure, homogeneous of infinitesimal dimension $d=3-2=1$.
\end{proof}

\begin{corollary}
    With the closed-string infinitesimal grading given by (\ref{eqn:closed-string-grading-operator-MM225}), the (-1)-eigensummand of quantum cohomology is isomorphic as a graded algebra to the cohomology ring of an elliptic curve.
\end{corollary}

The algebra assertion follows from the derivation property (Proposition \ref{prop:derivation}). We now prove Theorem \ref{thm:MM-2-25-reprise}.
\begin{proof}
    We first construct a Lagrangian sphere in the Fano threefold $X$ under discussion. By Picard-Lefschetz theory, the Poincaré duals of vanishing cycles of nodal degenerations of bidegree $(1,2)$ hypersurfaces span $H^3(X;\mathbb{Q})$, since $H^3(\mathbb{P}^1\times\mathbb{P}^3;\mathbb{Q})=0$. As $H^3(X;\mathbb{C})$ is nonzero, there exists a Lagrangian vanishing sphere $L$ with $[L]\neq0$.
    
    Since $OC([1_L])=PD([L])\neq0$, this $L$ is a nonzero object in the monotone Fukaya category. Moreover, a class in $H^3(X)$ pairing nontrivially with $[L]$ must map under the closed-open map $QH^*(X)\to HF^*(L,L)$ to a nonzero multiple of the point class, for degree reasons. Together with unitality, this shows closed-open surjectivity. 

    In particular, it follows by parity of degree that
    \begin{equation}
        QH^{odd}_\lambda(X)\neq 0,
    \end{equation}
    forcing $m_0(L)=-1$.

    It remains to show that $L$ is gradable with respect to a choice of geometric infinitesimal grading. Recall that there is a choice of a cochain $C\in CM^*(X)$ and chain-level representative $\Pi$ of $\pi_{-1}$ such that
    \begin{equation}
        \Pi*D+\Pi=\partial_X C
    \end{equation}
    holds at the chain level. In particular, the choice of a geometric infinitesimal grading may be changed by the closed--open image of a cocycle (see the discussion after Proposition \ref{prop:geometric-grading-is-a-grading}). Since we have established that the closed--open map is surjective for $L$ in cohomology, we may choose the geometric infinitesimal grading $\eta_\mathrm{geom}$ for $X$ such that for this specific Lagrangian $L$,
    \begin{equation}
        [\eta_\mathrm{geom}^0\vert_L]=0
    \end{equation}
    in $HF^*(L,L)$, and thus $L$ is gradable for that choice of geometric infinitesimal grading $\eta_\mathrm{geom}$.

    Pick $\alpha\in QH^*(X)_{-1}$ such that $CO^0(\alpha)=[pt_L]$. By Lemma \ref{lem:smith-lemma}, the composition
    \begin{equation}
        QH^*(X)\to HF^*(L,L)\to QH^*(X)
    \end{equation}
    is equal to the small quantum product with $PD[L]$. Since this class lies entirely in $H^3(X)$, $PD[L]$ is a linear combination of $\psi_1$ and $\psi_2$. In particular,
    \begin{equation}
        OC^0([pt_L])=\alpha*PD[L],
    \end{equation}
    but we have seen from the quantum multiplication table of the $(-1)$-summand that the only even degree classes in $QH^*(X)_{-1}$ that are multiples of linear combinations of $\psi_1$ and $\psi_2$ are scalar multiples of $\nu$. Therefore we conclude
    \begin{equation}
        OC^0([pt_L])=c\cdot\nu
    \end{equation}
    for some constant $c\in\mathbb{C}$.

    Now Ritter--Smith \cite[Lemma 12.17]{RitterSmith} show that the leading order term of the image $OC^0([pt_L])$ must be $[pt_X]$. Thus $c=1/16$, and
    \begin{equation}\label{eqn:MM2-25-CY}
        OC^0([pt_L])=\frac{1}{16}\nu.
    \end{equation}
    Since $\theta^{-1}*\nu=\nu$, Theorem \ref{thm:CY-shift-reprise} and \eqref{eqn:MM2-25-CY} show that $L$ has infinitesimal dimension
    \begin{equation}
        \mathrm{weight}_{\Gamma_{-1}}(\theta^{-1}*\frac{1}{16}\nu)-d=2-1=1.
    \end{equation}

    Thus $L$ is a Lagrangian sphere with infinitesimal dimension $1$, even though its geometric dimension is $3$. Restricting any homogeneous weak proper Calabi--Yau structure to this sphere forces its infinitesimal dimension to be $1$ by nondegeneracy (Proposition \ref{prop:CY-unique}).
\end{proof}

We give a criterion for surjectivity of the closed-open map of a Lagrangian sphere.

\begin{lemma} \label{lem:CO^0-A2}
    Let $L_0$ and $L_1$ be Lagrangian spheres in a monotone symplectic manifold $X$ that intersect transversely at a single point. Then the closed-open map
    \begin{equation}
        CO^0:QH^*(X)\to HF^*(L_0,L_0)
    \end{equation}
    is surjective.
\end{lemma}
\begin{proof}
    Since the closed-open map is unital, it hits the unit of $HF^*(L_0,L_0)$:
    \begin{equation}
        CO^0_{L_0}(1_X)=1_{L_0}.
    \end{equation}
    For the inclusion $\iota:L_0\xhookrightarrow{}X$, the classical contribution to the closed-open image $CO(PD[L_1])$ is
    \begin{equation}
        \iota^*PD[L_1]=\pm[pt_{L_0}],
    \end{equation}
    since $L_0$ and $L_1$ intersect at one point. Therefore there exists some $c\in\mathbb{C}$ such that
    \begin{equation}
        CO(PD[L_1])=\pm[pt_{L_0}]+c\cdot 1_{L_0},
    \end{equation}
    which shows that the images of $1_X$ and $PD[L_1]$ under the closed-open map span $HF^*(L_0,L_0)$.
\end{proof}

\subsection{Gushel-Mukai fourfolds} \label{ssec:Gushel-Mukai}

Our second example will be Gushel-Mukai fourfolds. Gushel-Mukai fourfolds have attracted recent interest because of their connections with K3 surfaces \cite{DebarreKuznetsov}, \cite{KuznetsovPerry}. In particular, Kuznetsov--Perry \cite[Proposition 2.6]{KuznetsovPerry} show that the Kuznetsov component of a Gushel-Mukai fourfold can be regarded as a ``noncommutative K3 surface'', in the sense that its associated Serre functor is a shift by degree $2$. In terms of their Hodge structures, \cite{Fano4K3} note that the middle-dimensional Hodge structure of Gushel-Mukai fourfolds and Verra fourfolds resembles that of a K3 surface.

Our main result is the following.

\begin{thm} \label{thm:Gushel-Mukai-reprise}
    Let $(X,\omega)$ be a smooth Gushel-Mukai fourfold equipped with a monotone Kähler form. The geometric infinitesimal grading on $\mathcal{F}(X)_0$ admits homogeneous weak proper Calabi--Yau structures precisely in infinitesimal dimension $2$.

    Every Lagrangian sphere in $X$ has infinitesimal dimension $2$, and those arising as vanishing cycles of nodal degenerations have infinite order Dehn twist in $\pi_0\mathrm{Symp}(X,\omega)$.

    Moreover, $X$ contains an $A_4$-configuration of Lagrangian spheres, each of infinitesimal dimension $2$. The first three Dehn twists induce a faithful action of $B_4$ on $D^\pi\mathcal{F}(X)_0$, even after forgetting the infinitesimal grading. Hence they induce an injective homomorphism $B_4\to \pi_0\mathrm{Symp}(X,\omega)$.
\end{thm}

First, we recall the topology of Gushel-Mukai fourfolds.

Ordinary Gushel-Mukai fourfolds are obtained by intersecting $\mathrm{Gr}(2,5)\subset\mathbb{P}^9$ with a hyperplane and a quadric hypersurface. Additively, their cohomology can be computed by the Lefschetz hyperplane theorem: following the notation of \cite{BenedettiManivelPerrin}, we may write the additive decomposition as
\begin{equation}
    H^*(X)=\mathbb{C}\langle \sigma_0,\sigma_1,\sigma_2,\sigma_{11},\sigma_3,\sigma_{31}\rangle\oplus H^4_{prim}(X),
\end{equation}
where the $\sigma$ classes are the ambient Schubert classes coming from the Grassmannian of cohomological degree $0,2,4,4,6,8$, and $H^4_{prim}(X)$ is the middle-dimensional primitive cohomology of rank $22$. For clarification, we note that $\sigma_0=1$ is the unit, $\sigma_1=H$ is the hyperplane class, and $\sigma_{31}=2[pt_X]$.

The first Chern class can be computed as twice the hyperplane class $c_1(X)=2H$, indeed matching that $X$ is a Fano variety. It follows that any Gushel-Mukai fourfold also admits a smooth anticanonical divisor: since $-K_X=2H$ one can take a smooth divisor in the linear system $\lvert 2H\rvert$.

We recall the ambient small quantum cohomology of a Gushel-Mukai fourfold from \cite[Theorem 6]{BenedettiManivelPerrin}: since we require the Novikov variable $q$ to have degree $2$, the powers of $q$ in the formula are doubled.

\begin{proposition} \label{prop:qh-GM4}
    The small ambient quantum cohomology ring of a Gushel-Mukai fourfold is
    \begin{equation}
        QH^*_{amb}(X)=\frac{\mathbb{C}[H,\sigma_{11}]}{(5H\sigma_{11}-2H^3+14q^2H,5\sigma_{11}^2+20q^2\sigma_{11}-H^4+12q^2H^2+20q^4,H^5-44q^2H^3-16q^4H)}.
    \end{equation}
\end{proposition}

We also include a small quantum product table of the ambient classes for the reader.

\begin{figure}[htbp]
\centering
\renewcommand{\arraystretch}{1.35}
\scriptsize
\[
\begin{array}{c|cccccc}
\star 
& \sigma_0 
& \sigma_1 
& \sigma_2 
& \sigma_{11} 
& \sigma_3 
& \sigma_{31} \\
\hline
\sigma_0
& \sigma_0
& \sigma_1
& \sigma_2
& \sigma_{11}
& \sigma_3
& \sigma_{31}
\\
\sigma_1
& \sigma_1
& 6\sigma_0+\sigma_2+\sigma_{11}
& 10\sigma_1+3\sigma_3
& 6\sigma_1+2\sigma_3
& 24\sigma_0+4\sigma_2+2\sigma_{11}+\sigma_{31}
& 24\sigma_1+6\sigma_3
\\
\sigma_2
& \sigma_2
& 10\sigma_1+3\sigma_3
& 80\sigma_0+8\sigma_2+12\sigma_{11}+2\sigma_{31}
& 52\sigma_0+8\sigma_2+4\sigma_{11}+\sigma_{31}
& 60\sigma_1+10\sigma_3
& 176\sigma_0+28\sigma_2+24\sigma_{11}
\\
\sigma_{11}
& \sigma_{11}
& 6\sigma_1+2\sigma_3
& 52\sigma_0+8\sigma_2+4\sigma_{11}+\sigma_{31}
& 32\sigma_0+6\sigma_2+\sigma_{31}
& 40\sigma_1+6\sigma_3
& 112\sigma_0+20\sigma_2+12\sigma_{11}
\\
\sigma_3
& \sigma_3
& 24\sigma_0+4\sigma_2+2\sigma_{11}+\sigma_{31}
& 60\sigma_1+10\sigma_3
& 40\sigma_1+6\sigma_3
& 120\sigma_0+20\sigma_2+20\sigma_{11}
& 120\sigma_1+24\sigma_3
\\
\sigma_{31}
& \sigma_{31}
& 24\sigma_1+6\sigma_3
& 176\sigma_0+28\sigma_2+24\sigma_{11}
& 112\sigma_0+20\sigma_2+12\sigma_{11}
& 120\sigma_1+24\sigma_3
& 368\sigma_0+64\sigma_2+48\sigma_{11}
\end{array}
\]
\caption{The small quantum products between the ambient classes of a Gushel-Mukai fourfold. The Novikov variable $q$ is set to 1 for readability, but can be easily recovered from the cohomological degrees.}
\label{fig:qh-GM4}
\end{figure}

The kernel of quantum multiplication with $c_1$ is determined in \cite[Lemma 7]{BenedettiManivelPerrin} to be
\begin{equation}
    \ker(c_1*)=\langle \alpha,\beta\rangle\oplus H^4_{prim}(X)
\end{equation}
for the two classes
\begin{align}
    &\alpha=-2q^2+2\sigma_2-3\sigma_{11},\\
    &\beta=-4q^4-2q^2\sigma_2+2[pt].
\end{align}
Moreover, $c_1*$ acts trivially on the $0$-generalized eigenspace of quantum cohomology. The other eigensummands of $c_1*$ can be determined to be $1$-dimensional. In particular, the infinitesimal gradings on those summands will be zero (Proposition \ref{prop:idempotent-zero}), so we will restrict our attention to the $0$-summand.

To determine the idempotent of the $0$-generalized eigensummand of quantum cohomology, we will determine the small quantum products involving primitive classes. First, since the small quantum products of ambient classes did not involve any primitive classes, it follows from Frobenius symmetry
\begin{equation}
    GW_{0,3,A}(\alpha,\beta,\gamma)=GW_{0,3,A}(\alpha,\gamma,\beta)
\end{equation}
that the small quantum product of an ambient class with a primitive class must be a primitive cohomology class. Now because the Gushel-Mukai fourfold is defined as a hyperplane section of a projective variety, the monodromy action on its primitive cohomology is irreducible \cite[Theorem 3.27]{Voisin}. Therefore, for any ambient class $\alpha$ and primitive class $\gamma$, the map
\begin{equation}
    \Phi:H^4_{prim}(X)\to H^4_{prim}(X),~\Phi(\gamma)\coloneqq \alpha*\gamma
\end{equation}
is a map of irreducible representations over $\mathbb{C}$, and thus is multiplication by a constant (Schur's lemma).

Since the Fano index of a Gushel-Mukai fourfold is $2$, we immediately conclude that for a primitive class $\gamma$,
\begin{equation}
    \sigma_1*\gamma=0,~\sigma_3*\gamma=0.
\end{equation}
Denote the scalars corresponding to quantum multiplication by the ambient classes $\sigma_{11},\sigma_2,\sigma_{31}$ by $\delta_{11},\delta_2,\delta_{31}$.

To determine the remaining small quantum products, recall that
\begin{align}
    \sigma_1*\sigma_1&=\sigma_2+\sigma_{11}+6,\\
    \sigma_1*\sigma_3&=\sigma_{31}+4\sigma_2+2\sigma_{11}+24,\\
    \sigma_{11}*\sigma_{11}&= \sigma_{31}+6\sigma_2+32.
\end{align}
Therefore by associativity of the quantum product, we have for any $\gamma\in H^*_{prim}(X)$,
\begin{equation}
    (\delta_2+\delta_{11}+6)\gamma=\sigma_1*(\sigma_1*\gamma)=0,
\end{equation}
and similarly
\begin{align}
    &\delta_{31}+4\delta_2+2\delta_{11}+24=0,\\
    &\delta_{11}^2=\delta_{31}+6\delta_2+32.
\end{align}
The above three equations uniquely determine the values of $\delta_2,\delta_{11},\delta_{31}$ to be
\begin{equation}
    \delta_2=-4,~\delta_{11}=-2,~\delta_{31}=-4.
\end{equation}
Therefore this fully determines the quantum product between the ambient classes and the primitive classes, and by Frobenius symmetry this determines the ambient parts of the small quantum products between the primitive classes. 

It remains to rule out a primitive component in the quantum product of two primitive classes. Equivalently, we need to show that nonzero monodromy-invariant pairings
\begin{equation}
    H^4_{prim}(X)\otimes H^4_{prim}(X)\to H^4_{prim}(X)
\end{equation}
do not exist. For this, we recall that the monodromy action on $H^4_{\mathrm{prim}}(X)$ is irreducible, and its invariant intersection form is indefinite. The monodromy group is therefore infinite. The  big-monodromy argument as in \cite{Voisin} then implies that the Picard–Lefschetz transformations generate a Zariski-dense subgroup of the orthogonal group. Since this group contains $-\id\in\mathbf{O}(n)$, which acts as $-\id$ on odd-degree tensors, we conclude that there is no nonzero monodromy-equivariant map
\begin{equation}
    H^4_{\mathrm{prim}}(X)\otimes H^4_{\mathrm{prim}}(X)\to H^4_{\mathrm{prim}}(X).
\end{equation}

In conclusion, for two primitive classes $\gamma_i,\gamma_j\in H^4_{prim}(X)$ whose cup product is $\gamma_i\cup\gamma_j=c_{ij}[pt_X]$, we have
\begin{equation}
    \gamma_i*\gamma_j=c_{ij}(-2q^4-q^2\sigma_2+[pt_X]).
\end{equation}
Note that the cohomology class in parentheses is half of $\beta$.

From this table, the idempotent for the generalized 0-eigensummand can be computed as
\begin{equation}
    \pi_0=-\frac{1}{4}\alpha-\frac{5}{16}\beta=\frac{7}{4}\sigma_0+\frac{1}{8}\sigma_2+\frac{3}{4}\sigma_{11}-\frac{5}{16}\sigma_{31},
\end{equation}
after setting $q=1$. We summarize the previous computations into the following proposition:

\begin{proposition} \label{prop:qh_0-GM4}
    The $0$-generalized eigensummand of the quantum cohomology of a Gushel-Mukai fourfold has the following basis
    \begin{equation}
        \pi_0,\gamma_1,\cdots,\gamma_{22},\nu,
    \end{equation}
    where $\pi_0$ is the idempotent, $\gamma_i$ are a basis of primitive classes of $H^4_{prim}(X)$, and $\nu=\beta\vert_{q=1}$. Moreover, $c_1*$ acts trivially on this generalized eigensummand. The quantum products between these classes satisfy
    \begin{align}
        &\pi_0*\pi_0=\pi_0,~\pi_0*\gamma_i=\gamma_i,~\pi_0*\nu=\nu,\\
        &\gamma_i*\gamma_j=\frac{c_{ij}}{2}\nu,\\
        &\gamma_i*\nu=0,~\nu*\nu=0.
    \end{align}
\end{proposition}

With this computation in hand, we are finally ready to compute the infinitesimal gradings of the $0$-generalized eigensummand.

\begin{proposition}
    The grading operator on the $0$-summand satisfies
    \begin{align}
        &\Gamma_0(\pi_0)=0,\\
        &\Gamma_0(\gamma_i)=2\gamma_i,\\
        &\Gamma_0(\nu)=4\nu.
    \end{align}
\end{proposition}
    
\begin{proof}
    Recall that the grading operator on the $0$-summand is defined as 
    \begin{equation}
        \Gamma_0(x)=\pi_0*(gr(x)-2q\frac{d}{dq}(\pi_0*_qx)\vert_{q=1}).
    \end{equation}
    We have seen that the infinitesimal grading of the idempotent is always zero. For the $\nu$ class, the small quantum product table shows
    \begin{align}
        \pi_0*_q\nu=&(\frac{7}{4}\sigma_0+\frac{1}{8}\sigma_2+\frac{3}{4}\sigma_{11}-\frac{5}{16}\sigma_{31})*_q(-4\sigma_0-2\sigma_2+\sigma_{31})\\
        =&(-115q^8+216q^6-98q^4-7)\sigma_0+(-20q^6+36q^4-14q^2-4)\sigma_2\\
        &+(-15q^6+27q^4-9q^2-3)\sigma_{11}+\sigma_{31}.
    \end{align}
    Therefore it follows that
    \begin{equation}
        q\frac{d}{dq}(\pi_0*_q\nu)=(-920q^8+1296q^6-392q^4)\sigma_0+(-120q^6+144q^4-28q^2)\sigma_2+(-90q^6+108q^4-18q^2)\sigma_{11}.
    \end{equation}
    Thus we compute
    \begin{align}
        \Gamma_0(\nu)&=\pi_0*(gr(\nu)-2q\frac{d}{dq}(\pi_0*_q\nu)\vert_{q=1})\\
        &=\pi_0*(-8\sigma_2+8\sigma_{31}-(-32\sigma_0-8\sigma_2))\\
        &=\pi_0*(32\sigma_0+8\sigma_{31})=32\pi_0+8\pi_0*\sigma_{31}\\
        &=32(\frac{7}{4}\sigma_0+\frac{1}{8}\sigma_2+\frac{3}{4}\sigma_{11}-\frac{5}{16}\sigma_{31})+8(-9\sigma_0-\frac{3}{2}\sigma_2-3\sigma_{11}+\frac{7}{4}\sigma_{31})\\
        &=-16\sigma_0-8\sigma_2+4\sigma_{31}=4\nu.
    \end{align}
    Similarly for the $\gamma_i$ classes,
    \begin{align}
        \pi_0*_q\gamma_i&=(\frac{7}{4}\sigma_0+\frac{1}{8}\sigma_2+\frac{3}{4}\sigma_{11}-\frac{5}{16}\sigma_{31})*_q\gamma_i\\
        &=\frac{7}{4}\gamma_i-\frac{1}{2}q^2\gamma_i-\frac{3}{2}q^2\gamma_i+\frac{5}{4}q^4\gamma_i,
    \end{align}
    and therefore
    \begin{align}
        \Gamma_0(\gamma_i)&=\pi_0*(4\gamma_i-2\gamma_i)=2\gamma_i.
    \end{align}
    In particular, for the class $\theta=\pi_0+5\nu/32$, we compute
    \begin{equation}
        \pi_0*gr(\theta)=2\theta,
    \end{equation}
    and thus Theorem \ref{thm:intro-cy-criterion} provides a weak proper Calabi--Yau structure of infinitesimal dimension $4-2=2$.
\end{proof}

\begin{remark}
    Kuznetsov--Perry \cite[Proposition 2.6]{KuznetsovPerry} show that the Kuznetsov component of the derived category of a Gushel-Mukai fourfold can be regarded as a ``noncommutative $K3$ surface''. Our computations suggest that the infinitesimally graded $0$-eigensummand of a Gushel-Mukai fourfold is isomorphic, as a graded $\mathbb{C}$-algebra, to the graded cohomology ring of a manifold of dimension $4$ with rank 22 middle-dimensional cohomology, which provides an interesting comparison.
\end{remark}

Now that we have computed the closed-string infinitesimal grading associated to the $0$-eigensummand of the quantum cohomology of a Gushel-Mukai fourfold $X$, we turn to its applications to the $0$-summand of the monotone Fukaya category $\mathcal{F}_0(X)$.

We first provide some lemmas about Lagrangian spheres in Fano varieties, which we will use repeatedly for other examples as well.

\begin{lemma} \label{lem:lagrangian-vanishing-spheres}
    Let $X$ be a Fano variety of Fano index $\rho\geq2$, and let $L$ be a Lagrangian sphere of dimension $n\geq2$. Then $L$ belongs to the $0$-summand of the monotone Fukaya category of $X$.
\end{lemma}
\begin{proof}
    We show that there are no Maslov index $2$ discs with boundary on $L$. Since $n\geq2$, $L$ is simply connected, and thus $\pi_2(X)\to\pi_2(X,L)$ is surjective. Now the Maslov index formula shows
    \begin{equation}
        \mu(u_\beta)=2\langle c_1(X),\beta\rangle.
    \end{equation}
    But since $X$ has Fano index $\rho\geq2$, the Chern number of any lift of $\beta\in\pi_2(X,L)$ to $\pi_2(X)$ must be divisible by $\rho$. Also by monotonicity, any such disc with positive area must have positive Chern number. Therefore, any disc with boundary on $L$ and positive area must have Maslov index at least $2\rho\geq4$, and thus there are no Maslov 2 discs with boundary on $L$. Thus we conclude that $L$ is a monotone Lagrangian in the $0$-summand of the monotone Fukaya category.
\end{proof}

From this, we conclude that the vanishing cycles constructed from the degenerating family of Gushel-Mukai fourfolds (constructed below) are all monotone Lagrangian spheres in the $0$-summand of the monotone Fukaya category. To determine their infinitesimal dimensions, we show the following lemma about the images of the closed-open maps.

\begin{lemma} \label{lem:image-of-CO}
    Let $X$ be a monotone symplectic manifold of even complex dimension, and let $L\subset X$ be a monotone Lagrangian sphere. Then the closed-open map $CO^0:QH^*(X)\to HF^*(L,L)$ is surjective.
\end{lemma}
\begin{proof}
    Since the closed-open map is unital, the unit $e_X\in QH^*(X)$ is sent to $e_L\in HF^*(L,L)$. Thus it is enough to show that the point class $[pt_L]\in HF^*(L,L)$ is in the image. We claim that
    \begin{equation}
        CO^0(PD[L])=\iota^*(PD[L])+\mathrm{lower~order~terms}
    \end{equation}
    for the map $\iota:L\xhookrightarrow{}X$. Indeed, this is the contribution of the constant discs to the map $CO^0$. Since $L$ is Lagrangian, an almost complex structure identifies the normal bundle of $L$ with $TL$, and therefore
    \begin{equation}
        \iota^*PD[L]=e(N_L)=\pm e(TL)=\pm\chi(L)[pt_L].
    \end{equation}
    But since $L$ is an even-dimensional Lagrangian sphere, $\chi(L)$ is nonzero, and therefore we conclude that $[pt_L]$ is also in the image of $CO^0$.
\end{proof}

\begin{lemma} \label{lem:2-dim-qh-OC}
    Let $X$ be a monotone symplectic manifold, and suppose that the $\lambda$-eigensummand of the quantum cohomology $QH^*(X)$ with respect to $c_1*$ multiplication has the algebra structure
    \begin{equation}
        QH^*(X)_\lambda=\mathbb{C}\langle\pi\rangle\oplus\mathfrak{m},\quad \mathfrak{m}^2=\mathbb{C}\langle\nu\rangle
    \end{equation}
    where $\pi$ is the $\lambda$-idempotent and $\mathfrak{m}$ is the unique maximal ideal. Then for any Lagrangian $L$ in the $\lambda$-summand of the monotone Fukaya category $\mathcal{F}_\lambda(X)$ for which $PD[L]\in\mathfrak{m}$ and the closed-open map $CO^0_L:QH^*(X)\to HF^*(L,L)$ hits the point class $[pt_L]$, the open-closed map $OC^0_L:HF^*(L,L)\to QH^*(X)$ satisfies
    \begin{equation}
        OC^0_L([pt_L])\in\mathbb{C}\langle PD([L])\rangle+\mathbb{C}\langle\nu\rangle.
    \end{equation}
    In particular, if $[pt_L]\in HF^*(L,L)$ is nilpotent, we conclude that
    \begin{equation}
        OC^0_L([pt_L])=c\cdot\nu
    \end{equation}
    for some $c\in\mathbb{C}$.
\end{lemma}
\begin{proof}
    Since $CO^0_L$ hits $[pt_L]$, there exists a class $\alpha\in QH^*(X)_\lambda$ such that $CO^0(\alpha)=[pt_L]$. Write
    \begin{equation}
        \alpha=b\cdot\pi+r,\quad r\in\mathfrak{m}
    \end{equation}
    for some $b\in\mathbb{C}$. Recall from Lemma \ref{lem:smith-lemma} that the composition
    \begin{equation}
        QH^*(X)\to HF^*(L,L)\to QH^*(X)
    \end{equation}
    is the quantum product with $PD([L])$. Specializing to $\alpha$, we conclude that
    \begin{equation}
        OC^0([pt_L])=\alpha*PD([L])=b\cdot PD([L])+r*PD([L]),
    \end{equation}
    and the assumption that $\mathfrak{m}^2=\mathbb{C}\cdot\nu$ shows that the second term lies in $\mathbb{C}\langle \nu\rangle$.

    When $[pt_L]$ is a nilpotent element in $HF^*(L,L)$, $\alpha$ cannot be an invertible class since $CO$ preserves products. Thus $b=0$. 
\end{proof}

With the above lemma, we can determine the open-closed image of the point class of a Lagrangian sphere constructed as a vanishing cycle. We remark that this lemma crucially uses the quantum product structure of the $0$-eigensummand of the quantum cohomology of Gushel-Mukai fourfolds, computed in Proposition \ref{prop:qh_0-GM4}.

\begin{lemma} \label{lem:OC^-GM4}
    Let $X$ be a Gushel-Mukai fourfold, and let $L\subset X$ be a Lagrangian sphere. Then the open-closed map $OC^0_L:HF^*(L,L)\to QH^*(X)$ satisfies
    \begin{equation}\label{eqn:GM4-CY}
        OC^0_L([pt_L])=\frac{1}{2}\nu
    \end{equation}
    for the nilpotent generator $[pt_L]$ of $HF^*(L,L)$. In particular, any Lagrangian sphere has infinitesimal dimension $2$ induced from the open-string infinitesimal grading operator.
\end{lemma}
\begin{proof}
    We first understand the cohomology class $PD([L])\in QH^4(X)$. The ambient middle cohomology is spanned by $h^2$ and $c_2(TX)$, both of which restrict trivially to $L$; the former because $H^2(L)=0$ and the latter because $TX\vert_L\cong TL\otimes_{\mathbb{R}}\mathbb{C}$, and $TS^4$ is stably trivial. Therefore $PD([L])$ must belong to the primitive part of the cohomology of $X$.

    Now from Lemma \ref{lem:image-of-CO}, the closed-open map $CO^0_L$ is surjective. Also from Proposition \ref{prop:qh_0-GM4}, the 0-eigensummand of the quantum cohomology of a Gushel-Mukai fourfold satisfies the assumptions of Lemma \ref{lem:2-dim-qh-OC}, and the point class of a spherical object is nilpotent. Therefore, we conclude that
    \begin{equation}
        OC^0([pt_L])=c\cdot\nu
    \end{equation}
    for some scalar $c\in\mathbb{C}$. From \cite[Lemma 12.17]{RitterSmith}, the leading order term of $OC^0_L([pt_L])$ is the point class $[pt_X]$. Therefore we can determine the value of the scalar $c$ to be exactly $1/2$.

    Note that $HF^\mathrm{odd}(L,L)=0$, so $L$ automatically admits a graded lift. Since $\theta^{-1}*\nu=\nu$, Theorem \ref{thm:CY-shift-reprise} and \eqref{eqn:GM4-CY} show that $L$ has infinitesimal dimension
    \begin{equation}
        \mathrm{weight}_{\Gamma_0}(\theta^{-1}*\frac{1}{2}\nu)-d=4-2=2.
    \end{equation}
     Moreover, restricting any homogeneous weak proper Calabi--Yau structure to these spheres forces its infinitesimal dimension to be $2$ by nondegeneracy.
\end{proof}

To conclude the proof of Theorem \ref{thm:Gushel-Mukai-reprise}, we construct a family of Gushel-Mukai fourfolds with $A_n$-singularity degenerations. The idea is to modify the equation of the quadric that appears in the definition of Gushel-Mukai fourfolds to introduce such a singularity.

To provide an explicit construction, consider the point $p\in Gr(2,5)$ spanned by the column vectors $e_1,e_2$. Then locally near $p$, $Gr(2,5)$ can be parametrized by the 2-planes given by
\begin{equation}
\begin{pmatrix}
    1 &0 &a &b &c\\
    0 &1 &d &e &f
\end{pmatrix}    
\end{equation}
for $a,b,c,d,e,f\in\mathbb{C}$.

Then the local Plücker coordinates of the image of the Grassmannian under the Plücker embedding are the $2\times2$ minors of this matrix. Now choose the hyperplane that defines the Gushel-Mukai variety to be $b=f$, which cuts down the number of free variables to 5. We now have the freedom to choose the quadric hypersurface that defines the Gushel-Mukai fourfold.

Choose the equation of the quadric to be $Q=(a+b^2-ce)^2+c^2+d^2+e^2+a(ab-cd)$. Solving the equations for $\partial Q/\partial x$ for the four variables $x=a,c,d,e$ gives
\begin{equation}
    \{a=-\frac{b^2}{1+b},c=d=e=0\}
\end{equation}
and in these coordinates the equation $Q$ takes the form
\begin{equation}
    Q(-\frac{b^2}{1+b},b,0,0,0)=\frac{b^5}{1+b},
\end{equation}
and therefore there exists a holomorphic change of coordinates near the origin such that
\begin{equation}
    Q(z_1,z_2,z_3,z_4,z_5)=z_1^2+z_2^2+z_3^2+z_4^2+z_5^5,
\end{equation}
an $A_4$ singularity. We expect this method to be able to construct families of Gushel-Mukai fourfolds that degenerate into a central fiber with higher isolated $A_n$-singularities for $n>4$, but we do not determine the largest possible $n$ here.

From this $A_4$-degeneration, standard arguments for Lefschetz fibrations (\cite[Section 6]{KhovanovSeidel}, \cite[Section 20]{Seidel08}) show that there exist vanishing cycles in a regular fiber of the family that form an $A_4$-configuration of Lagrangian spheres.

For each Lagrangian sphere $L_i$ that appears in the $A_4$ configuration constructed from a vanishing cycle, we have demonstrated that its infinitesimal dimension is $2$. Therefore by Theorem \ref{thm:ADE-faithfulness-reprise}, we conclude that the Dehn twist along any such Lagrangian sphere has infinite order, and moreover the first three Dehn twists coming from the $A_4$ configuration of Lagrangian spheres generate a faithful braid group action of $B_{4}$ on the derived monotone Fukaya category summand $D^\pi\mathcal{F}_0(X)$ of the Gushel-Mukai fourfold. Hence they induce an injective group homomorphism to the symplectic mapping class group $B_4\to\pi_0\mathrm{Symp}(X,\omega)$.

Moreover, since the open locus of quadric sections of the fixed smooth hyperplane section of $\mathrm{Gr}(2,5)$ having exactly one ordinary node is connected, any two Lagrangian vanishing cycles in this family are carried to one another by a symplectomorphism of the smooth fiber. Thus their Dehn twists are conjugate to those in the configuration above and therefore have infinite order.

\begin{remark}
    For a $K3$ surface, Seidel \cite{Seidel08b} shows that for a choice of symplectic form, one can construct a braid group action of $B_{15}$ on the Fukaya category. Since Gushel-Mukai fourfolds and the Verra fourfold (which we discuss next) admit Kuznetsov components of ``K3 type'', it would be interesting to see if a similar $B_{15}$ action can be constructed on the corresponding $0$-summand of the monotone Fukaya category.
\end{remark}

\subsection{Verra fourfolds} \label{ssec:Verra}

We next consider Verra fourfolds, another class of Fano fourfolds of index $2$ and $K3$ type \cite{Fano4K3}.

Our main result for Verra fourfolds is the following.

\begin{thm} \label{thm:Verra-reprise}
    Let $(X,\omega)$ be a smooth Verra fourfold, equipped with a monotone Kähler form. The geometric infinitesimal grading on $\mathcal{F}(X)_0$ admits homogeneous weak proper Calabi--Yau structures precisely in infinitesimal dimension $2$.

    Every Lagrangian sphere has infinitesimal dimension $2$, and moreover vanishing cycles arising from nodal degenerations have infinite order Dehn twists in the symplectic mapping class group $\pi_0\mathrm{Symp}(X,\omega)$.

    Moreover, $X$ contains an $A_5$-configuration of Lagrangian spheres, each of infinitesimal dimension $2$. The Dehn twists along the first four spheres define an injective group homomorphism
    \begin{equation}
        B_5\to\pi_0\mathrm{Symp}(X,\omega).
    \end{equation}
\end{thm}

A Verra fourfold is defined to be a double cover of $\mathbb{P}^2\times\mathbb{P}^2$ branched along a smooth divisor of bidegree $(2,2)$. Define the pullbacks along the projections of the two hyperplane classes of $\mathbb{P}^2$ to be $h_1,h_2$. Then the cohomology can be additively computed as
\begin{equation}
    H^*(X;\mathbb{C})=\mathbb{C}\langle 1,h_1,h_2,h_1^2,h_1h_2,h_2^2,h_1^2h_2,h_1h_2^2,h_1^2h_2^2\rangle\oplus H^4_{prim}(X),
\end{equation}
where the primitive part $H^4_{prim}(X)$ has rank 21, and $h_1^2h_2^2=2[pt_X]$ (see \cite[Lemma 3]{BenedettiGuereManivelPerrin}).

Any Verra fourfold also admits a smooth anticanonical divisor by Bertini's theorem, since the anticanonical divisor $-K_X=\pi^*\mathcal{O}(2,2)$ is globally generated. We also note that $c_1(X)=2h_1+2h_2$.

The ambient small quantum cohomology of a Verra fourfold was computed in \cite[Section 3]{BenedettiGuereManivelPerrin} and \cite[Proposition 6.1]{Fay}, using different methods which yield the same result. We provide the result in the following table (Figure \ref{fig:qh-Verra4}).

\begin{figure}[htbp]
\centering
\renewcommand{\arraystretch}{1.4}
\scriptsize
\[
\begin{array}{c|cccc} \star & h_1 & h_1^2 & h_1h_2 & h_1^2h_2 \\ 
\hline 
h_1 & h_1^2+2 & 2h_1+4h_2 & 2h_2+h_1^2h_2 & 8+2h_1h_2+4h_2^2 \\
h_2 & h_1h_2 & h_1^2h_2 & 2h_1+h_1h_2^2 & 8+2h_1^2+h_1^2h_2^2 \\
h_1^2 & 2h_1+4h_2 & 4+4h_1h_2 & 8+2h_1h_2+4h_2^2 & 8h_1+4h_2+4h_1h_2^2 \\ 
h_1h_2 & 2h_2+h_1^2h_2 & 8+2h_1h_2+4h_2^2 & 12+2h_1^2+2h_2^2+h_1^2h_2^2 & 20h_1+16h_2+2h_1h_2^2 \\ 
h_2^2 & h_1h_2^2 & 8+h_1^2h_2^2 & 8+4h_1^2+2h_1h_2 & 8h_1+16h_2+2h_1^2h_2 \\ 
h_1^2h_2 & 8+2h_1h_2+4h_2^2 & 8h_1+4h_2+4h_1h_2^2 & 20h_1+16h_2+2h_1h_2^2 & 40+16h_1^2+16h_1h_2+4h_2^2 \\ 
h_1h_2^2 & 8+2h_2^2+h_1^2h_2^2 & 16h_1+8h_2+2h_1h_2^2 & 16h_1+20h_2+2h_1^2h_2 & 32+8h_1^2+20h_1h_2+8h_2^2 \\ 
h_1^2h_2^2 & 8h_1+8h_2+2h_1h_2^2 & 32+8h_1^2+8h_1h_2+4h_2^2 & 32+8h_1^2+12h_1h_2+8h_2^2 & 32h_1+40h_2+8h_1^2h_2+8h_1h_2^2 
\end{array}
\]

\[
\begin{array}{c|cc} \star & h_1h_2^2 & h_1^2h_2^2 \\ 
\hline 
h_1 & 8+2h_2^2+h_1^2h_2^2 & 8h_1+8h_2+2h_1h_2^2 \\ 
h_2 & 8+4h_1^2+2h_1h_2 & 8h_1+8h_2+2h_1^2h_2 \\ 
h_1^2 & 16h_1+8h_2+2h_1h_2^2 & 32+8h_1^2+8h_1h_2+4h_2^2 \\ 
h_1h_2 & 16h_1+20h_2+2h_1^2h_2 & 32+8h_1^2+12h_1h_2+8h_2^2 \\ 
h_2^2 & 4h_1+8h_2+4h_1^2h_2 & 32+4h_1^2+8h_1h_2+8h_2^2 \\ 
h_1^2h_2 & 32+8h_1^2+20h_1h_2+8h_2^2 & 32h_1+40h_2+8h_1^2h_2+8h_1h_2^2 \\ 
h_1h_2^2 & 40+4h_1^2+16h_1h_2+16h_2^2 & 40h_1+32h_2+8h_1^2h_2+8h_1h_2^2 \\ 
h_1^2h_2^2 & 40h_1+32h_2+8h_1^2h_2+8h_1h_2^2 & 144+32h_1^2+32h_1h_2+32h_2^2 \end{array} 
\]

\caption{The small quantum products between the ambient classes of a Verra fourfold. One can recover the omitted small quantum products from symmetry between $h_1$ and $h_2$. The Novikov variables $q_1$ and $q_2$ are both set to 1 for readability.}
\label{fig:qh-Verra4}
\end{figure}

Again, we need to extend this computation to the primitive components. The argument will be similar to the case of Gushel-Mukai fourfolds.

First, since the small quantum product classes between ambient classes are ambient, it follows from Frobenius symmetry that the small quantum product between an ambient class and a primitive class must be primitive. Now recall that the monodromy action on the primitive cohomology of Verra fourfolds is shown to be irreducible in \cite[Theorem 1.3]{LeePirola}. Therefore, we again obtain that for any ambient class $\alpha$, the action
\begin{equation}
    \gamma\in H^4_{prim}(X) \mapsto \alpha*\gamma\in H^4_{prim}(X)
\end{equation}
is multiplication by a scalar $\delta_\alpha$. Since the Fano index of a Verra fourfold is $2$, the same arguments as in the Gushel-Mukai fourfold case show that the only nontrivial quantum products between ambient and primitive classes are given by 
\begin{equation}
    h_1^2*\gamma,\quad h_1h_2*\gamma,\quad h_2^2*\gamma,\quad h_1^2h_2^2*\gamma.
\end{equation}
Multiplication by these ambient classes acts by scalar multiples of the identity; we denote the corresponding scalars by $\delta_{11},\delta_{12},\delta_{22},\delta_4$.

From the quantum product ring, we have
\begin{align}
    h_1*h_1&=h_1^2+2,\\
    h_1^2h_2*h_1&=2h_1h_2+4h_2^2+8,\\
    h_1h_2^2*h_1&=h_1^2h_2^2+2h_2^2+8.
\end{align}
The above relations uniquely determine the values of the $\delta$ to be
\begin{equation}
    \delta_{11}=\delta_{22}=-2,\delta_{12}=0,\delta_4=-4.
\end{equation}
Again, we have for two primitive classes $\gamma_i,\gamma_j$ whose cup product is $\gamma_i\cup\gamma_j=c_{ij}[pt_X]$,
\begin{equation}
    \gamma_i*\gamma_j=c_{ij}(-2-h_1^2-h_2^2+\frac{1}{2}h_1^2h_2^2).
\end{equation}
Primitive components are excluded by the same big-monodromy argument as in the Gushel-Mukai case.

From this, the kernel of $c_1*$ can be computed as
\begin{equation}
    \ker(c_1*)=\mathbb{C}\langle 4-h_1^2-h_2^2+h_1h_2,-4-2h_1^2-2h_2^2+h_1^2h_2^2,2h_1-2h_2+h_1^2h_2-h_1h_2^2\rangle\oplus H^4_{prim},
\end{equation}
and moreover one can show that the $0$-generalized eigensummand of the $c_1*$-action agrees with the kernel (see \cite[Proposition 14]{BenedettiGuereManivelPerrin}), and the other generalized eigensummands are at most $2$-dimensional. Therefore, we again restrict our attention to the infinitesimal gradings of the $0$-eigensummand.

From the above computations, the idempotent for the $0$-generalized eigensummand can be computed as
\begin{equation}
    \pi_0=\frac{11}{16}-\frac{1}{32}(h_1^2+h_2^2)+\frac{1}{8}h_1h_2-\frac{3}{64}h_1^2h_2^2.
\end{equation}
Let
\begin{align}
    \alpha&=2h_1-2h_2+h_1^2h_2-h_1h_2^2,\\
    \nu&=-2-h_1^2-h_2^2+\frac{1}{2}h_1^2h_2^2
\end{align}
be the complementary basis elements of the 0-generalized eigenspace. Then the previous computations of the small quantum product can be summarized as follows:

\begin{proposition} \label{prop:qh_0-Verra4}
    The $0$-generalized eigensummand of the quantum cohomology of a Verra fourfold has the following basis
    \begin{equation}
        \pi_0,\alpha, \gamma_1,\cdots,\gamma_{21},\nu,
    \end{equation}
    where $\pi_0$ is the idempotent, $\gamma_i$ are a basis of primitive classes of $H^4_{prim}(X)$, and $\alpha,\nu$ are defined above. Moreover, $c_1*$ acts trivially on this generalized eigensummand. The quantum products between these classes satisfy
    \begin{align}
        &\alpha*\alpha=-16\nu,~\alpha*\gamma_i=0,\\
        &\gamma_i*\gamma_j=c_{ij}\nu,\\
        &\alpha*\nu=0,~\gamma_i*\nu=0,~\nu*\nu=0.
    \end{align}
\end{proposition}

\begin{proposition}
    The grading operator on the 0-generalized eigenspace of quantum cohomology satisfies
    \begin{align}
        \Gamma_0(\pi_0)&=0,\\
        \Gamma_0(\gamma_i)&=2\gamma_i,\\
        \Gamma_0(\alpha)&=2\alpha,\\
        \Gamma_0(\nu)&=4\nu.
    \end{align}
\end{proposition}
\begin{proof}
    This time, there are two Novikov variables, one for each $\mathbb{P}^2$ factor of the base of the double covering. So the grading operator becomes
    \begin{equation}
        \Gamma_0(x)=\pi_0*(gr(x)-2(q_1\frac{d}{dq_1}+q_2\frac{d}{dq_2})(\pi_0*_qx)\vert_{q_1=q_2=1}).
    \end{equation}
    Before we begin, we list the small quantum products that we will use, this time with the Novikov variables. There is a symmetry between $h_1$ and $h_2$, so we exclude the products that can be deduced from symmetry.
    \begin{align}
        h_1^2*_qh_1^2&=4q_1^4+4q_1^2h_1h_2\\
        h_1^2*_qh_1h_2&=8q_1^2q_2^2+2q_1^2h_1h_2+4q_1^2h_2^2\\
        h_1^2*_qh_2^2&=8q_1^2q_2^2+h_1^2h_2^2\\
        h_1h_2*_qh_1h_2&=12q_1^2q_2^2+2q_2^2h_1^2+2q_1^2h_2^2+h_1^2h_2^2\\
        h_1^2*_qh_1^2h_2^2&=32q_1^4q_2^2+8q_1^2q_2^2h_1^2+8q_1^2q_2^2h_1h_2+4q_1^4h_2^2\\
        h_1h_2*_qh_1^2h_2^2&=16q_1^4q_2^2+16q_1^2q_2^4+8q_1^2q_2^2h_1^2+12q_1^2q_2^2h_1h_2+8q_1^2q_2^2h_2^2\\
        h_1^2h_2^2*_qh_1^2h_2^2&=144q_1^4q_2^4+32q_1^2q_2^4h_1^2+16(q_1^4q_2^2+q_1^2q_2^4)h_1h_2+32q_1^4q_2^2h_2^2.
    \end{align}
    
    We start with $\nu$. The small quantum multiplication table shows
    \begin{align}
        &\pi_0*_q\nu=(\frac{11}{16}-\frac{1}{32}(h_1^2+h_2^2)+\frac{1}{8}h_1h_2-\frac{3}{64}h_1^2h_2^2)*_q(-2-h_1^2-h_2^2+\frac{1}{2}h_1^2h_2^2)\\
        &=-\frac{27q_1^4q_2^4-16q_1^4q_2^2-q_1^4-16q_1^2q_2^4+12q_1^2q_2^2-q_2^4+11}{8}-\frac{6q_1^2q_2^4-6q_1^2q_2^2-q_2^4+4q_2^2+5}{8}h_1^2\\
        &-\frac{3q_1^4q_2^2+3q_1^2q_2^4-10q_1^2q_2^2+q_1^2+q_2^2+2}{8}h_1h_2-\frac{6q_1^4q_2^2-q_1^4-6q_1^2q_2^2+4q_1^2+5}{8}h_2^2+\frac{1}{2}h_1^2h_2^2
    \end{align}
    so applying $2(q_1\frac{d}{dq_1}+q_2\frac{d}{dq_2})$, we get
    \begin{align}
        \Gamma_0(\nu)&=\pi_0*(-4h_1^2-4h_2^2+4h_1^2h_2^2+16+4(h_1^2+h_2^2))\\
        &=\pi_0*(16+4h_1^2h_2^2)\\
        &=16\pi_0+4(\frac{11}{16}-\frac{1}{32}(h_1^2+h_2^2)+\frac{1}{8}h_1h_2-\frac{3}{64}h_1^2h_2^2)*h_1^2h_2^2.
    \end{align}
    To finish, we plug in the small quantum products from the beginning:
    \begin{align}
        \Gamma_0(\nu)&=16\pi_0+\frac{11}{4}h_1^2h_2^2-\frac{1}{8}(h_1^2+h_2^2)*h_1^2h_2^2+\frac{1}{2}h_1h_2*h_1^2h_2^2-\frac{3}{16}h_1^2h_2^2*h_1^2h_2^2\\
        &=16(\frac{11}{16}-\frac{1}{32}(h_1^2+h_2^2)+\frac{1}{8}h_1h_2-\frac{3}{64}h_1^2h_2^2)+\frac{11}{4}h_1^2h_2^2-\frac{1}{8}(64+12h_1^2+16h_1h_2+12h_2^2)\\
        &+\frac{1}{2}(32+8h_1^2+12h_1h_2+8h_2^2)-\frac{3}{16}(144+32h_1^2+32h_1h_2+32h_2^2)\\
        &=-8-4(h_1^2+h_2^2)+2h_1^2h_2^2=4\nu.
    \end{align}
    For the primitive classes, we have that
    \begin{equation}
        \pi_0*_q\gamma_i=\frac{11}{16}\gamma_i+\frac{1}{16}(q_1^2+q_2^2)\gamma_i+\frac{3}{16}q_1^2q_2^2\gamma_i,
    \end{equation}
    and therefore
    \begin{align}
        \Gamma_0(\gamma_i)&=\pi_0*(4\gamma_i-2(\frac{1}{4}\gamma_i+\frac{3}{4}\gamma_i))\\
        &=\pi_0*2\gamma_i=2\gamma_i.
    \end{align}
    Finally, for $\alpha$ we have
    \begin{align}
        \pi_0*_q\alpha&=(\frac{11}{16}-\frac{1}{32}(h_1^2+h_2^2)+\frac{1}{8}h_1h_2-\frac{3}{64}h_1^2h_2^2)*_q(2h_1-2h_2+h_1^2h_2-h_1h_2^2)\\
        &=\frac{-6q_1^4q_2^2+9q_1^2q_2^4+4q_1^2q_2^2-q_1^2+q_2^4-2q_2^2+11}{8}h_1+h_1^2h_2\\
        &-\frac{9q_1^4q_2^2+q_1^4-6q_1^2q_2^4+4q_1^2q_2^2-2q_1^2-q_2^2+11}{8}h_2-h_1h_2^2.
    \end{align}
    Thus we compute
    \begin{align}
        \Gamma_0(\alpha)&=\pi_0*(4h_1-4h_2+6h_1^2h_2-6h_1h_2^2-8h_1+8h_2)\\
        &=(\frac{11}{16}-\frac{1}{32}(h_1^2+h_2^2)+\frac{1}{8}h_1h_2-\frac{3}{64}h_1^2h_2^2)*(-4h_1+4h_2+6h_1^2h_2-6h_1h_2^2)\\
        &=\frac{11}{16}(-4h_1+4h_2+6h_1^2h_2-6h_1h_2^2)-\frac{1}{32}(-16h_1+16h_2-8h_1^2h_2+8h_1h_2^2)\\
        &+\frac{1}{8}(32h_1-32h_2-16h_1^2h_2+16h_1h_2^2)-\frac{3}{64}(-48h_1+48h_2+8h_1^2h_2-8h_1h_2^2)\\
        &=4h_1-4h_2+2h_1^2h_2-2h_1h_2^2=2\alpha.
    \end{align}
    In particular, we compute that the class $\theta=\pi_0+3\nu/64$ satisfies
    \begin{equation}
        \pi_0*gr(\theta)=2\theta,
    \end{equation}
    and thus Theorem \ref{thm:intro-cy-criterion} provides the desired Calabi--Yau structure of infinitesimal dimension $2$.
\end{proof}

\begin{remark}
    Our computations show that the $0$-eigensummand of the quantum cohomology of a Verra fourfold is (noncanonically) isomorphic, as a graded $\mathbb{C}$-algebra, to the graded cohomology ring of a K3 surface: since we can rescale the basis element $\alpha$ corresponding to the hyperplane class by any complex scalar, our computations do not determine the degree of the K3 surface. To compare the result to Kuznetsov components, \cite[Example 3.11]{PerryPertusiZhao} explains that the Kuznetsov component of the derived category of a Verra fourfold is equivalent to the derived category of a degree 2 K3 surface with a Brauer class that is 2-torsion. Kuznetsov--Shinder \cite[Section 1.3]{KuznetsovShinder} explain that there is a derived equivalence between the derived category of a degree 8 K3 surface and the twisted derived category of a degree 2 K3 surface with a 2-torsion Brauer class.
\end{remark}

Now again, Lemma \ref{lem:lagrangian-vanishing-spheres} shows that any Lagrangian vanishing cycle belongs to the $0$-summand of the monotone Fukaya category, and Lemma \ref{lem:image-of-CO} shows that the associated closed-open map 
\begin{equation}
    CO^0:QH^*(X)\to HF^*(L,L)
\end{equation}
is surjective. Now we prove a similar lemma about the image of the open-closed map for Verra fourfolds.

\begin{lemma} \label{lem:OC^0-Verra4}
    Let $X$ be a Verra fourfold, and let $L\subset X$ be a Lagrangian sphere. Then the open-closed map $OC^0_L:HF^*(L,L)\to QH^*(X)$ maps the point class of $L$ to
    \begin{equation}
        OC^0_L([pt_L])=\nu
    \end{equation}
    for the nilpotent generator $[pt_L]$ of $HF^*(L,L)$. In particular, any Lagrangian sphere in a Verra fourfold has infinitesimal dimension $2$ induced from the open-string infinitesimal grading operator.
\end{lemma}
\begin{proof}
    We claim that for any Lagrangian sphere $L\subset X$, $PD([L])$ is an element in the primitive part of the cohomology of $X$. The ambient middle cohomology of $X$ is generated by products of degree $2$ classes, which restrict trivially to $L$ since $H^2(L)$=0. Thus $PD([L])$ lies in the primitive cohomology of $X$. Therefore we may apply Lemma \ref{lem:2-dim-qh-OC} and conclude that
    \begin{equation}
        OC^0([pt_L])=c\cdot\nu
    \end{equation}
    for some scalar $c\in\mathbb{C}$. From \cite[Lemma 12.17]{RitterSmith}, the leading order term of $OC^0_L([pt_L])$ is the point class $[pt_X]$. Therefore we can determine the constant $c$ to be exactly $1$. 
    
    Similarly to Lemma \ref{lem:OC^-GM4}, we conclude that $L$ has infinitesimal dimension $2$, and thus any homogeneous weak proper Calabi--Yau structure must also have infinitesimal dimension $2$.
\end{proof}

To finish the proof of Theorem \ref{thm:Verra-reprise}, we exhibit an explicit $A_5$-degeneration of Verra fourfolds. This degeneration is constructed by choosing an explicit equation for the branch locus of the double cover in $\mathbb{P}^2\times\mathbb{P}^2$.

To begin the construction, pick local coordinates
\begin{equation}
    ([1;x_1;x_2],[1;y_1;y_2])\in\mathbb{P}^2\times\mathbb{P}^2,
\end{equation}
and define the local equation for the defining polynomial of the branch locus to be
\begin{equation}
    p(x_1,x_2,y_1,y_2)=x_1x_2+y_1^2+x_1y_2^2+x_2y_1y_2.
\end{equation}
Then the kernel of the Hessian is in the $\frac{\partial}{\partial y_2}$-direction, and in that direction we can reparametrize the equation by solving
\begin{equation}
    \frac{\partial p}{\partial x_1}=\frac{\partial p}{\partial x_2}=\frac{\partial p}{\partial y_1}=0,
\end{equation}
which gives
\begin{align}
    x_2=-y_2^2,\\
    x_1+y_1y_2=0,\\
    2y_1+x_2y_2=0,
\end{align}
equivalently
\begin{equation}
    x_2=-y_2^2,\quad y_1=\frac{1}{2}y_2^3,\quad x_1=-\frac{1}{2}y_2^4.
\end{equation}
Thus, the polynomial $p$ grows in the $\frac{\partial}{\partial y_2}$ direction as
\begin{equation}
    p(-\frac{1}{2}y_2^4,-y_2^2,\frac{1}{2}y_2^3,y_2)=-\frac{1}{4}y_2^6,
\end{equation}
and thus the branch locus acquires an isolated $A_5$ singularity. Then the local defining equation for the Verra fourfold, which is a double branched cover, would be
\begin{equation}
    z^2=p(x_1,x_2,y_1,y_2),
\end{equation}
which also acquires an $A_5$-singularity.

Therefore we can construct an $A_5$-configuration of Lagrangian vanishing cycles in a Verra fourfold, and each Lagrangian $L_i$ in this configuration has infinitesimal dimension 2. Therefore it follows that the first four Dehn twists along this $A_5$ configuration of Lagrangian spheres generate a faithful braid group action of $B_5$ on the derived monotone Fukaya category $D^\pi\mathcal{F}_0(X)$ of the Verra fourfold, and hence an injective map to the symplectic mapping class group $B_5\to\pi_0\mathrm{Symp}(X,\omega)$.

Similarly to Gushel-Mukai fourfolds, vanishing cycles of nodal degenerations of Verra fourfolds lie in a single symplectic monodromy orbit, so their Dehn twists are conjugate to those in the configuration above and therefore have infinite order.

\subsection{Symplectic Grassmannians} \label{ssec:IG(2,2n)}

Let $(V,\omega)$ be a complex vector space of dimension $2n$ with a symplectic form. The symplectic Grassmannian $X=\operatorname{IG}(2,2n)$ parametrizes isotropic two-planes in $V$ and has complex dimension $4n-5$. 

Our main theorem for this subsection is that the family of symplectic Grassmannians exhibit fractional infinitesimal dimensions.

\begin{thm} \label{thm:IG-reprise}
    Let $(X,\omega)$ be a symplectic Grassmannian $\operatorname{IG}(2,2n)$ equipped with a monotone Kähler form. For $n\geq3$, one can choose a geometric infinitesimal grading on $\mathcal{F}(X)_0$ for which the standard weak proper Calabi--Yau structure is homogeneous of infinitesimal dimension $\frac{n-2}{n}$, and every homogeneous weak proper Calabi--Yau structure for this grading has this dimension.
\end{thm}

To do so, we first describe its small quantum cohomology. Take the tautological rank two bundle $\mathcal{E}\to\operatorname{IG}(2,2n)$, and define
\begin{equation}
    a_1=c_1(\mathcal{E}),\quad a_2=c_2(\mathcal{E}),\quad b_j=c_{2j}(\mathcal{E}^\perp/\mathcal{E}).
\end{equation}
The first Chern class of symplectic Grassmannians is $c_1(TX)=-(2n-1)a_1$.

The quantum cohomology of $\operatorname{IG}(2,2n)$ was first computed in \cite{BKT}; we use the description given in \cite[Corollary 4.2]{MMPS}.

\begin{proposition}\label{prop:QH-IG(2,2n)}
    The small quantum cohomology of the symplectic Grassmannian is given by
    \begin{equation}\label{eqn:QH-IG(2,2n)}
        QH^*(\operatorname{IG}(2,2n))\cong \frac{\mathbb{C}[a_1,a_2]}{(P_{n-1}(a_1,a_2),P_n(a_1,a_2)-a_1)}
    \end{equation}
    for the polynomial $P_n(a_1,a_2)$ defined recursively as
    \begin{equation}\label{eqn:P_n-recursion}
        P_k=(a_1^2-2a_2)P_{k-1}-a_2^2P_{k-2},\quad P_0=1,\quad P_1=a_1^2-2a_2.
    \end{equation}
    The ring structure used throughout this subsection is quantum multiplication, not ordinary cup product.
\end{proposition}

The generalized eigensummands of $c_1*$ were computed in \cite[Proposition 4.3]{MMPS}: there are $(2n-1)(n-1)$ summands of rank $1$ and the zero generalized eigensummand has rank $n-1$. Moreover, writing the idempotent of the zero summand as $\pi_0$, the proof of \cite[Proposition 4.3]{MMPS} shows that there is an isomorphism for the zero generalized eigensummand
\begin{equation}
    QH^*(\operatorname{IG}(2,2n))_0\cong\mathbb{C}[\varepsilon]/(\varepsilon^{n-1}),
\end{equation}
where $\varepsilon=\pi_0*a_2$. In particular, the $c_1*$ action on the zero generalized eigensummand is trivial.

The following proposition provides an explicit formula for the idempotent.

\begin{proposition}
    The idempotent of the zero generalized eigensummand $QH^*(\operatorname{IG}(2,2n))_0$ is given by
    \begin{equation}
        \pi_0=1+\frac{1}{4n}(4a_2-a_1^2)\frac{\partial P_{n-1}}{\partial a_1},
    \end{equation}
    and satisfies $a_1*\pi_0=0$.
\end{proposition}
\begin{proof}
    Set $P_{-1}=0$. For each $k\geq0$, define the polynomial $Q_k(a_1,a_2)$ as
    \begin{equation}
        Q_{k}(a_1,a_2)=\frac{4a_2-a_1^2}{4}\frac{\partial P_{k}}{\partial a_1}.
    \end{equation}
    We claim the following general identity
    \begin{equation}\label{eqn:Q_n}
        a_1Q_{n}=(n+1)a_2^2P_{n-1}-\frac{n}{2}(a_1^2-2a_2)P_n.
    \end{equation}
    Once this holds, it follows that in quantum cohomology
    \begin{equation}
        a_1Q_{n-1}=na_2^2P_{n-2}-\frac{n-1}{2}(a_1^2-2a_2)P_{n-1}=-na_1.
    \end{equation}
    Since the polynomial $P_n$ is even in the variable $a_1$, its derivative is divisible by $a_1$, and therefore $\pi_0\equiv 1~(\mathrm{mod}~a_1)$. Thus multiplication with $\pi_0$ is identity on $QH^*(X)_0$ and is zero on the other summands, which shows that $\pi_0$ is the idempotent.
    
    We prove \eqref{eqn:Q_n} by induction: at $n=0$ both sides are zero, and at $n=1$ using $P_1=a_1^2-2a_2$ we compute
    \begin{equation}
        a_1Q_1=a_1\cdot\frac{4a_2-a_1^2}{4}\cdot 2a_1=\frac{a_1^2(4a_2-a_1^2)}{2}=2a_2^2-\frac{1}{2}(a_1^2-2a_2)^2.
    \end{equation}
    Now for the induction step differentiate the recursive relation \eqref{eqn:P_n-recursion} in $a_1$:
    \begin{equation}
        Q_k=2a_1\cdot\frac{4a_2-a_1^2}{4}P_{k-1}+(a_1^2-2a_2)Q_{k-1}-a_2^2Q_{k-2},
    \end{equation}
    so by inductive hypotheses
    \begin{align}
        a_1Q_k&=(2a_2^2-\frac{1}{2}(a_1^2-2a_2)^2)P_{k-1}+(a_1^2-2a_2)(ka_2^2P_{k-2}\\
        &-\frac{k-1}{2}(a_1^2-2a_2)P_{k-1})-a_2^2((k-1)a_2^2P_{k-3}-\frac{k-2}{2}(a_1^2-2a_2)P_{k-2})\\
        &=(2a_2^2-\frac{k}{2}(a_1^2-2a_2)^2)P_{k-1}+\frac{3k-2}{2}(a_1^2-2a_2)a_2^2P_{k-2}-(k-1)a_2^4P_{k-3}\\
        &=(k+1)a_2^2P_{k-1}-\frac{k}{2}(a_1^2-2a_2)P_k,
    \end{align}
    where the last equality uses the recursive relation for $P_k$ again.
\end{proof}

With this, we now compute the grading operator on the generalized zero summand.

\begin{proposition}
    The closed-string grading operator on the zero summand can be computed as
    \begin{equation}
        \Gamma_0(\varepsilon^k)=\frac{2k}{n}\varepsilon^k,\qquad 0\leq k\leq n-2.
    \end{equation}
    Moreover, the criterion \eqref{eqn:inf-dim-formula} for the infinitesimal dimension is satisfied by the idempotent:
    \begin{equation}
        \pi_0*\mathrm{gr}_X(\pi_0)=\frac{(4n-2)(n-1)}{n}\pi_0
    \end{equation}
    and this determines the infinitesimal dimension to be $\frac{n-2}{n}$.
\end{proposition}
\begin{proof}
    Since the grading operator is a derivation, it is enough to check that
    \begin{equation}
        \Gamma_0(\varepsilon)=\frac{2}{n}\varepsilon.
    \end{equation}
    We recall from the quantum Pieri rule (\cite[Theorem 1.4]{BKT}) that
    \begin{align}
        (4a_2-a_1^2)*\frac{\partial P_{n-1}}{\partial a_1}&=(4a_2-a_1^2)\cup\frac{\partial P_{n-1}}{\partial a_1}-4(n-1)\cdot 1_X,\label{eqn:quantum-relation-IG}\\
        a_2*\left((4a_2-a_1^2)\cup\frac{\partial P_{n-1}}{\partial a_1}\right)&=a_2\cup(4a_2-a_1^2)\cup\frac{\partial P_{n-1}}{\partial a_1}-4(a_1^2-a_2).
    \end{align}
    The first identity gives a homogeneous degree decomposition of the idempotent
    \begin{equation}
        \pi_0=\frac{1}{n}\cdot1_X+\frac{1}{4n}(4a_2-a_1^2)\cup\frac{\partial P_{n-1}}{\partial a_1},
    \end{equation}
    and thus the grading operator $\mathrm{gr}_X$ applied to the idempotent is
    \begin{equation}\label{eqn:IG-grpi}
        \mathrm{gr}_X(\pi_0)=(4n-2)(\pi_0-\frac{1}{n}\cdot 1_X),
    \end{equation}
    which shows that
    \begin{equation}\label{eqn:IG-grading-1}
        \mathrm{gr}_X(\pi_0)*\varepsilon=\frac{(4n-2)(n-1)}{n}\varepsilon.
    \end{equation}
    Similarly, the second identity provides a homogeneous degree decomposition of $\varepsilon$:
    \begin{equation}
        \varepsilon=\frac{1}{n}(-a_1^2+2a_2)+\frac{1}{4n}a_2\cup (4a_2-a_1^2)\cup\frac{\partial P_{n-1}}{\partial a_1},
    \end{equation}
    and each term has cohomological degrees $4$ and $4n+2$. Since $\pi_0*a_1=0$ and $\pi_0*a_2=\varepsilon$, multiplying both sides with $\pi_0$ shows that
    \begin{equation}
         \pi_0*\left(a_2\cup(4a_2-a_1^2)\cup\frac{\partial P_{n-1}}{\partial a_1}\right)=4(n-2)\varepsilon.
    \end{equation}
    Thus we compute
    \begin{equation}\label{eqn:IG-grading-2}
    \begin{aligned}
        \pi_0*\mathrm{gr}_X(\varepsilon)&=\pi_0*(\frac{4}{n}(-a_1^2+2a_2)+\frac{4n+2}{4n}a_2\cup(4a_2-a_1^2)\cup\frac{\partial P_{n-1}}{\partial a_1})\\
        &=\frac{8}{n}\varepsilon+\frac{(4n+2)(n-2)}{n}\varepsilon.
    \end{aligned}    
    \end{equation}
    Subtracting \eqref{eqn:IG-grading-1} from \eqref{eqn:IG-grading-2}  shows that
    \begin{equation}
        \Gamma_0(\varepsilon)=\pi_0*\mathrm{gr}_X(\varepsilon)-\mathrm{gr}_X(\pi_0)*\varepsilon=(\frac{8+(4n+2)(n-2)-(4n-2)(n-1)}{n})\varepsilon=\frac{2}{n}\varepsilon.
    \end{equation}
    Finally from \eqref{eqn:IG-grpi} we compute that
    \begin{equation}
        \pi_0*\mathrm{gr}_X(\pi_0)=\pi_0*(4n-2)(\pi_0-\frac{1}{n}\cdot 1_X)=\frac{(4n-2)(n-1)}{n}\pi_0,
    \end{equation}
    and thus the infinitesimal dimension is computed as
    \begin{equation}
        d=(4n-5)-\frac{(4n-2)(n-1)}{n}=\frac{n-2}{n}.
    \end{equation}
\end{proof}

Therefore \ref{thm:CY-criterion-reprise} shows that the canonical homogeneous weak proper Calabi--Yau structure on the monotone Fukaya summand $\mathcal{F}(\operatorname{IG}(2,2n))_0$ has infinitesimal dimension $\frac{n-2}{n}$ with respect to a geometric infinitesimal grading. To show that it is unique for a chosen grading, we invoke the uniqueness criterion (Proposition \ref{prop:CY-dim-unique}), and exhibit a nonzero gradeable object in the monotone Fukaya summand.

Since $\operatorname{IG}(2,2n)$ is defined as a hyperplane section of the Grassmannian $\operatorname{Gr}(2,2n)$, by considering a deformation of the equation we can construct a family of symplectic Grassmannians with nodal degenerations. This family constructs a Lagrangian vanishing sphere $L$ in $\operatorname{IG}(2,2n)$, and since $c_1(\operatorname{IG}(2,2n))=-(2n-1)a_1$, by Lemma \ref{lem:lagrangian-vanishing-spheres} this is a monotone Lagrangian sphere in the zero summand of the monotone Fukaya category since there are no Maslov two discs. Moreover, since the minimal Maslov number satisfies
\begin{equation}
    N_L=4n-2>4n-5+1,
\end{equation}
the Oh spectral sequence shows that
\begin{equation}
    HF^*(L,L)\cong H^*(S^{4n-5}),
\end{equation}
as graded vector spaces and thus $L$ is a nonzero object.

It remains to show that $L$ is gradeable with respect to the geometric infinitesimal grading $\eta$. To do so, recall that $HF^*(L,L)$ also admits a $\mathbb{Z}/(4n-2)$-grading from its Maslov cover. If one shows that
\begin{equation}
    [\eta^0_L]\in HF^{1}(L,L)=0
\end{equation}
with respect to the $\mathbb{Z}/(4n-2)$-grading, it follows that $[\eta^0_L]=0$ since the cohomology $HF^*(L,L)$ vanishes in that degree.

Recall that we may choose the components of the geometric infinitesimal grading to each lie in degree $1$ for $C$, $-1$ for $\beta_\Pi$, and $1$ for $\beta_0$ and the homotopies $K_{\Pi,D},K_{1,D},K_{D,1}^\mathrm{inc}$. Fixing these choices, the corresponding geometric infinitesimal grading lies in degree $1$ with respect to the $\mathbb{Z}/(4n-2)$-grading, and in particular so does $\eta^0_L$. It follows that the vanishing cycle $L$ constructed above is a nonzero gradeable object in the infinitesimally graded category $(\mathcal{F}(X)_0,\eta_\mathrm{geom})$, and thus its infinitesimal dimension is uniquely determined as $\frac{n-2}{n}$, matching the Serre dimension of the $A_{n-1}$ Kuznetsov component.

\begin{appendices}

\section{Quantum cohomology of MM 2-25} \label{app:MM2-25}

In this appendix, we compute the small quantum cohomology ring of the Fano threefold MM 2-25.

Recall that $X$ can be realized as the blowup of $\mathbb{P}^3$ along an elliptic curve, an intersection of two quadrics $Q_0$ and $Q_\infty$ which form a pencil of quadrics $\{Q_t\}$ on $\mathbb{P}^3$. The cohomology of $X$ can be computed additively as
\begin{equation}
    H^*(X)\cong H^*(\mathbb{P}^3)\oplus H^{*-2}(E),
\end{equation}
where $E$ is the center of the blowup. Alternatively, $X$ can also be realized as a bidegree $(1,2)$ hypersurface in $\mathbb{P}^1\times\mathbb{P}^3$, from which we obtain the following description of the generators:
\begin{equation}
    H^*(X)\cong\mathbb{C}\langle1,h_1,h_2,\psi_1,\psi_2,h_1h_2,h_2^2,[pt]\rangle,
\end{equation}
where $h_1$ and $h_2$ are the pullbacks of the hyperplane classes from $\mathbb{P}^1$ and $\mathbb{P}^3$, and $\psi_1$, $\psi_2$ are degree $3$ classes from the $H^1$ classes of the blowup locus, which is an elliptic curve. Also, one computes the classical relations
\begin{equation}
    h_1^2=0,\quad h_1h_2^2=2[pt],\quad h_2^3=[pt],\quad \psi_1\psi_2=-[pt].
\end{equation}

We first state our computation of the quantum cohomology of $X$.

\begin{proposition} \label{prop:qh-MM225}
    The small quantum products between even degree cohomology classes of $X$ are as follows:
    \begin{equation*}
    \resizebox{\textwidth}{!}{$
    \begin{array}{c|cccccc}
        * 
        & 1 
        & h_1 
        & h_2 
        & h_1h_2 
        & h_2^2 
        & [\mathrm{pt}]
        \\
        \hline
        1 
        & 1 
        & h_1 
        & h_2 
        & h_1h_2 
        & h_2^2 
        & [\mathrm{pt}]
        \\[2mm]

        h_1    
        & h_1 
        & 2h_2-h_1 
        & h_1h_2 
        & 4+2h_2^2-h_1h_2 
        & 4+2[\mathrm{pt}] 
        & 2+h_1+2h_2
        \\[2mm]

        h_2 
        & h_2 
        & h_1h_2 
        & h_2^2+2 
        & 4+2h_1+2[\mathrm{pt}] 
        & 4+3h_1+2h_2+[\mathrm{pt}] 
        & 1+h_1+2h_2+h_1h_2
        \\[2mm]

        h_1h_2 
        & h_1h_2 
        & 4+2h_2^2-h_1h_2 
        & 4+2h_1+2[\mathrm{pt}] 
        & 4+4h_1+8h_2 
        & 2+2h_1+8h_2+2h_1h_2 
        & 4+2h_2+h_1h_2+2h_2^2
        \\[2mm]

        h_2^2 
        & h_2^2 
        & 4+2[\mathrm{pt}] 
        & 4+3h_1+2h_2+[\mathrm{pt}] 
        & 2+2h_1+8h_2+2h_1h_2 
        & 5+h_1+6h_2+4h_1h_2 
        & 8+2h_1+h_2+h_1h_2+2h_2^2
        \\[2mm]

        [\mathrm{pt}] 
        & [\mathrm{pt}] 
        & 2+h_1+2h_2 
        & 1+h_1+2h_2+h_1h_2 
        & 4+2h_2+h_1h_2+2h_2^2 
        & 8+2h_1+h_2+h_1h_2+2h_2^2 
        & 6+3h_1+2h_2+h_2^2
        \end{array}
    $}    
    \end{equation*}
    The small quantum products between the odd and even degree classes are as follows: for any $\psi\in H^3(X)$,
    \begin{align}
        h_1*\psi=-\psi,\\
        h_2*\psi=0,\\
        h_1h_2*\psi=0,\\
        h_2^2*\psi=-2\psi,\\
        [pt]*\psi=-\psi,
    \end{align}
    and from this, one can compute the small quantum product between the odd classes to be
    \begin{equation}
        \psi_1*\psi_2=1+h_1+\frac{1}{2}h_2^2-\frac{1}{4}h_1h_2-[pt].
    \end{equation}
\end{proposition}

The proof will occupy the rest of this section. First, we enumerate all possible irreducible rational curve images for the effective curves of anticanonical degree at most $6$: these are the only curves that can contribute to the small quantum product. 

Suppose that $C$ is a rational curve of degree $d$ in $\mathbb{P}^3$ that meets the blowup locus $E$ in scheme length $l$. Then its strict transform has bidegree $(2d-l,d)$ in $\mathbb{P}^1\times\mathbb{P}^3$. Since the first Chern class is $h_1+2h_2$, it follows that the effective curve classes of anticanonical degree at most $6$ must satisfy $4d-l\leq6$. From this, we can enumerate the effective curve classes as (strict transforms of) the following curve classes:
\[
\begin{array}{c@{\quad:\quad}l}
(1,0) & \text{exceptional fibers} \\
(0,1) & \text{bisecant lines to } E \\
(1,1) & \text{lines in } \mathbb{P}^3 \text{ meeting } E \text{ once} \\
(2,1) & \text{lines in } \mathbb{P}^3 \text{ disjoint from } E \\
(0,2) & \text{4-secant conics} \\
(1,2) & \text{3-secant conics} \\
(2,2) & \text{2-secant conics} \\
(0,3) & \text{rational cubics in a quadric fiber}
\end{array}
\]

Now we examine each effective curve class case by case.

\subsection{(1,0)-curves}

These are the exceptional fibers in the blowup; denote their class by $f$. First, the divisor axiom shows that for any cohomology classes $x$, $y$,
\begin{equation}
    GW_{0,3,f}(x,y,h_2)=(h_2\cdot f) GW_{0,2,f}(x,y)=0.
\end{equation}
We then claim that the only nonzero Gromov-Witten invariant for the class $f$ among products with even degree insertions is
\begin{equation}
    GW_{0,3,f}(h_1,h_1,h_1h_2)=4.
\end{equation}
Indeed, the Poincaré dual of $h_1$ can be chosen to be a generic fiber $Q_\lambda$ of the pencil, and the Poincaré dual of $h_1h_2$ can be chosen to be an intersection of a fiber $Q_\nu$ with a generic hyperplane $H\subset\mathbb{P}^3$. Let the fiber class $f$ be represented by a generic exceptional fiber $f_e$ supported over a point $e\in E$ in the blowup locus. The incidence conditions for the Gromov-Witten invariant force the point $e$ to be in the intersection of the blowup locus $E$ with the chosen hyperplane $H$: since $E$ has degree $4$, there are 4 such intersection points, and each exceptional fiber $f_e$ over those points contributes with multiplicity $1$ to the Gromov-Witten invariant.

The other Gromov-Witten invariant
\begin{equation}
    GW_{0,3,f}(h_1,h_1,h_2^2)=0
\end{equation}
follows because the Poincaré dual of $h_2^2$ is represented by a generic line in $\mathbb{P}^3$, so the incidence condition with a generic fiber $f_e$ cannot be satisfied. 

Now suppose that the small quantum product $h_1*_fh_1$ is given by
\begin{equation}
    h_1*_fh_1=ah_1+bh_2.
\end{equation}
Then from the above computations, we find that
\begin{align}
    4&=\langle h_1*_fh_1, h_1h_2\rangle=a\int_Xh_1^2h_2+b\int_Xh_1h_2^2=2b,\\
    0&=\langle h_1*_fh_1,h_2^2\rangle=a\int_Xh_1h_2^2+b\int_Xh_2^3=2a+b,
\end{align}
which determines the small quantum product
\begin{equation}
    h_1*_{(1,0)}h_1=2h_2-h_1.
\end{equation}
Symmetry of the Gromov-Witten invariants determines the other small quantum product
\begin{align}
    h_1*_{(1,0)}h_1h_2=2h_2^2-h_1h_2,
\end{align}
and these are the only nonzero contributions to the small quantum products of the class $f$ among products with even degree insertions.

\subsection{(0,1)-curves}

These are the contributions of the bisecant lines $b$ in $\mathbb{P}^3$ to the blowup locus $E$, whose strict transforms are ruling lines in the quadric $Q_\lambda$ of the pencil. The divisor axiom shows
\begin{equation}
     GW_{0,3,b}(h_1,x,y)=0,
\end{equation}
for cohomology classes $x,y$, and it follows that the small quantum product $x*_{(0,1)}y$ has no $h_2^2$ term.

To compute the Gromov-Witten invariants, we use \cite[Proposition 2.1]{ABT}, which states that the bisecant lines to an elliptic quartic in $\mathbb{P}^3$ form a smooth congruence of bidegree $(2,6)$ in the Grassmannian $G(1,3)$. Therefore, through a general point $p\in\mathbb{P}^3$, there exist two bisecant lines to $E$. This determines one Gromov-Witten invariant from the divisor axiom:
\begin{equation}
    GW_{0,3,b}(h_2,h_2,[pt])=(b\cdot h_2)^2\cdot GW_{0,1,b}([pt])=2.
\end{equation}
We now determine the remaining three Gromov-Witten invariants
\begin{equation}
     GW_{0,3,b}(h_2,h_1h_2,h_1h_2),\quad  GW_{0,3,b}(h_2,h_1h_2,h_2^2),\quad  GW_{0,3,b}(h_2,h_2^2,h_2^2).
\end{equation}
For this we prefer to understand the curve class $b$ as ruling lines in the quadric fiber $Q_t$. Recall that the Poincaré dual of $h_1h_2$ is represented by a quadric $Q_\lambda$ intersected with a generic hyperplane $H$. Since the representatives for two different $h_1h_2$ classes can be generically chosen to be supported in two different quadric fibers $Q_\lambda,Q_{\lambda'}$, and ruling lines are completely contained in quadric fibers, we have
\begin{equation}
    GW_{0,3,b}(h_2,h_1h_2,h_1h_2)=0.
\end{equation}
Now the Poincaré dual of $h_2^2$ can be represented by a generic line $l\subset\mathbb{P}^3$, and this meets a generic quadric $Q_\lambda$ at two points. Through each of these two points, there exist two ruling lines of the quadric, which determines
\begin{equation}
    GW_{0,3,b}(h_2,h_1h_2,h_2^2)=4.
\end{equation}
The final undetermined Gromov-Witten invariant can be determined via Schubert calculus from the Grassmannian description: we count the number of lines in the congruence of bidegree $(2,6)$ that meet two generic lines in $\mathbb{P}^3$, so
\begin{align}
    GW_{0,3,b}(h_2,h_2^2,h_2^2)&=\int_{G(1,3)}\sigma_1\cdot\sigma_1\cdot(2\sigma_2+6\sigma_{1,1})\\
    &=\int_{G(1,3)}(\sigma_2+\sigma_{1,1})\cdot(2\sigma_2+6\sigma_{1,1})\\
    &=\int_{G(1,3)}2\sigma_2\cdot\sigma_2+8\sigma_2\cdot\sigma_{1,1}+6\sigma_{1,1}\cdot\sigma_{1,1}\\
    &=2+0+6=8.
\end{align}
Now a computation similar to that in the $(1,0)$ case determines the small quantum products
\begin{align}
    h_2*_{(0,1)}h_2&=2,\\
    h_2*_{(0,1)}h_1h_2&=2h_1,\\
    h_2*_{(0,1)}h_2^2&=3h_1+2h_2,\\
    h_2*_{(0,1)}[pt]&=h_1h_2,\\
    h_1h_2*_{(0,1)}h_2^2&=2h_1h_2,\\
    h_2^2*_{(0,1)}h_2^2&=4h_1h_2.
\end{align}

\subsection{(1,1)-curves}

These are the strict transforms of lines meeting the blowup locus $E$ at 1 point. 

We first determine the Gromov-Witten invariants involving the point class: the result is
\begin{align}
    GW_{0,3,(1,1)}(h_1,h_1h_2,[pt])&=4,\\
    GW_{0,3,(1,1)}(h_1,h_2^2,[pt])&=4,\\
    GW_{0,3,(1,1)}(h_2,h_1h_2,[pt])&=4,\\
    GW_{0,3,(1,1)}(h_2,h_2^2,[pt])&=4.
\end{align}

By the divisor axiom, the above reduces to showing that
\begin{equation}
    GW_{0,2,(1,1)}(h_1h_2,[pt])=GW_{0,2,(1,1)}(h_2^2,[pt])=4.
\end{equation}

The Poincaré dual of $h_1h_2$ is a quadric $Q_\lambda$ intersected with a generic hyperplane $H$. Fix a generic point $q\not\in Q_\lambda$, and a point $e$ in the blowup locus $E$. The line $l_e$ joining $q$ and $e$ meets $Q_\lambda$ at exactly two points: since $E\subset Q_\lambda$, one of them is $e$. Denote the other intersection point by $r(e)$. As $e$ varies in $E$, $r(e)$ sweeps out a quartic curve $R$. A generic hyperplane $H$ meets $R$ at four points away from $E$, which we name $\{p_1,p_2,p_3,p_4\}$. The corresponding four lines $p_iq$ are precisely the one-secant lines through $q$ whose strict transforms meet the strict transform of $Q_\lambda\cap H$. This determines
\begin{equation}
    GW_{0,2,(1,1)}(h_1h_2,[pt])=4.
\end{equation}

The Poincaré dual of $h_2^2$ is a generic line $l$ in $\mathbb{P}^3$, and the dual of $[pt]$ determines a generic point $q\in\mathbb{P}^3$. There is a unique plane $P$ passing through $q$ and $l$, and this again meets the blowup locus $E$ at 4 points $q_1,q_2,q_3,q_4$. Thus the 4 lines contributing to the Gromov-Witten invariants are determined as the strict transforms of $\overline{qq_i}$. We conclude
\begin{equation}
    GW_{0,2,(1,1)}(h_2^2,[pt])=4.
\end{equation}
These determine the small quantum products
\begin{align}
    h_1*_{(1,1)}h_1h_2&=4,\\
    h_2*_{(1,1)}h_1h_2&=4,\\
    h_1*_{(1,1)}h_2^2&=4,\\
    h_2*_{(1,1)}{h_2^2}&=4,\\
    h_1*_{(1,1)}[pt]&=h_1+2h_2,\\
    h_2*_{(1,1)}[pt]&=h_1+2h_2,\\
    h_1h_2*_{(1,1)}[pt]&=h_1h_2+2h_2^2,\\
    h_2^2*_{(1,1)}[pt]&=h_1h_2+2h_2^2.
\end{align}

The remaining Gromov-Witten invariants can be determined as follows. We claim that
\begin{equation}
    GW_{0,3,(1,1)}(h_2^2,h_2^2,h_2^2)=8.
\end{equation}
Indeed, this is the number of lines in $\mathbb{P}^3$ meeting the degree 4 blowup locus $E$ and three generic lines, so Schubert calculus tells us that
\begin{equation}
    GW_{0,3,(1,1)}(h_2^2,h_2^2,h_2^2)=\int_{G(1,3)}(4\sigma_1)\cdot\sigma_1\cdot\sigma_1\cdot\sigma_1=8.
\end{equation}
For the Gromov-Witten invariants involving $h_1h_2$, it is more convenient to compute the Gromov-Witten invariants involving the exceptional fiber class $f$. Since $h_1h_2$ is the Poincaré dual of the strict transform of a conic
\begin{equation}
    C=Q_\lambda\cap H\subset\mathbb{P}^3
\end{equation}
which meets the blowup locus $E$ at 4 points, we conclude that
\begin{equation}
    h_1h_2=2h_2^2-4f,
\end{equation}
where $f$ is the Poincaré dual of a fiber class. We now claim that
\begin{equation}
    GW_{0,3,(1,1)}(h_2^2,h_2^2,f)=1.
\end{equation}
Take the two lines $l_1,l_2$ in $\mathbb{P}^3$ that represent the Poincaré duals of the two $h_2^2$ classes, and suppose that the fiber class $f$ is supported over a point $e\in E$. Then the pairs $(e,l_1),(e,l_2)$ each determine a unique plane $P_1,P_2$, and therefore the 1-secant line to $E$ is uniquely determined as the intersection of these two planes. Finally, we claim
\begin{equation}
    GW_{0,3,(1,1)}(h_2^2,f,f)=GW_{0,3,(1,1)}(f,f,f)=0,
\end{equation}
since each fiber class is supported over a single point in the base locus, and two generic fiber incidences cannot be satisfied by strict transforms of 1-secant lines in $\mathbb{P}^3$. Using $h_1h_2=2h_2^2-4f$, we conclude
\begin{align}
     GW_{0,3,(1,1)}(h_2^2,h_2^2,h_1h_2)&=12,\\
      GW_{0,3,(1,1)}(h_2^2,h_1h_2,h_1h_2)&=16,\\
       GW_{0,3,(1,1)}(h_1h_2,h_1h_2,h_1h_2)&=16.
\end{align}
These determine the remaining small quantum products
\begin{align}
    h_1h_2*_{(1,1)}h_1h_2&=4h_1+8h_2,\\
    h_1h_2*_{(1,1)}h_2^2&=2h_1+8h_2,\\
    h_2^2*_{(1,1)}h_2^2&=h_1+6h_2.
\end{align}

\subsection{(2,1)-curves}

These are lines in $\mathbb{P}^3$ that are disjoint from the blowup locus $E$. We first show
\begin{equation}
    GW_{0,3,(2,1)}(h_1,[pt],[pt])=2,\quad GW_{0,3,(2,1)}(h_2,[pt],[pt])=1,
\end{equation}
which both follow from the divisor axiom applied to
\begin{equation}
    GW_{0,2,(2,1)}([pt],[pt])=1.
\end{equation}
Indeed, through 2 general points in $\mathbb{P}^3$, there exists a unique line passing through these points. This determines the small quantum products
\begin{align}
     h_1*_{(2,1)}[pt]&=2,\\
     h_2*_{(2,1)}[pt]&=1,\\
     [pt]*_{(2,1)}[pt]&=h_2^2.
\end{align}

We now claim that
\begin{align}
    GW_{0,3,(2,1)}(h_1h_2,h_1h_2,[pt])&=4,\\
    GW_{0,3,(2,1)}(h_1h_2,h_2^2,[pt])&=2,\\
    GW_{0,3,(2,1)}(h_2^2,h_2^2,[pt])&=1.
\end{align}
These all follow from the fact that the lines in $\mathbb{P}^3$ through a generic point $p\in\mathbb{P}^3$ meeting a degree $d$ curve form a degree $d$ curve in $\mathbb{P}^2$. Therefore by Bezout's theorem, the number of lines through a generic point $p\in\mathbb{P}^3$ meeting a degree $d_1$ curve and a degree $d_2$ curve is $d_1d_2$.

Applying this to our situation, recall that the  Poincaré dual of the class $h_1h_2$ is a conic, and the class $h_2^2$ is a line. Therefore the $(2,1)$-contribution to the quantum products between such classes counts the number of lines that meet the prescribed conics or lines, which justifies the above computation.

This leads to the computation of the small quantum products
\begin{align}
    h_1h_2*_{(2,1)}h_1h_2&=4,\\
    h_1h_2*_{(2,1)}h_2^2&=2,\\
    h_2^2*_{(2,1)}h_2^2&=1,\\
    h_1h_2*_{(2,1)}[pt]&=2h_2,\\
    h_2^2*_{(2,1)}[pt]&=h_2.
\end{align}

\subsection{(0,2)-curves}

These are $4$-secant conics to $E$, and thus their strict transforms must lie in a fiber of the quadric pencil $X\to\mathbb{P}^1$. Thus, these are in one-to-one correspondence with (1,1)-conic curves in a quadric fiber $Q_t\cong\mathbb{P}^1\times\mathbb{P}^1$. We first claim that
\begin{equation}
    GW_{0,3,(0,2)}(h_2^2,h_2^2,[pt])=4.
\end{equation}
Indeed, the  Poincaré dual of $h_2^2$ is a general line $l\subset\mathbb{P}^3$, which meets $Q_t$ at two points. Thus for two general lines $l_1,l_2$ which meet $Q_t$ at the points $\{a_1,a_2\}$, $\{b_1,b_2\}$, and a generic point $p$, there exists a unique conic in $Q_t$ through the three points $a_i,b_j,p$, because $Q_t$ is isomorphic to $\mathbb{P}^1\times\mathbb{P}^1$, and $\lvert\mathcal{O}_{Q_t}(1,1)\rvert\cong\mathbb{P}^3$.

We then claim
\begin{equation}
    GW_{0,3,(0,2)}(h_1,[pt],[pt])=0,\quad GW_{0,3,(0,2)}(h_2,[pt],[pt])=0,
\end{equation}
because the representatives of the Poincaré duals of the two point classes project to two general points in $\mathbb{P}^3$, which lie in two different quadrics $Q_\lambda,Q_{\lambda'}$ of the pencil. Since each $(0,2)$-curve class is completely contained in one of the quadrics, it cannot satisfy the two incidence conditions given by the point classes.

The remaining Gromov-Witten invariants are
\begin{equation}
    GW_{0,3,(0,2)}(h_1h_2,h_1h_2,[pt])=0,\quad GW_{0,3,(0,2)}(h_1h_2,h_2^2,[pt])=0,
\end{equation}
since again the representative of the Poincaré dual of $h_1h_2$ also induces an incidence condition with a general quadric $Q_\lambda$, which is incompatible with the other incidence condition coming from the point class.

Therefore the small quantum products can be computed as
\begin{align}
    h_2^2*_{(0,2)}h_2^2&=4,\\
    h_2^2*_{(0,2)}[pt]&=2h_1.
\end{align}

\subsection{(1,2)-curves}

These are the strict transforms of 3-secant conics to the blowup locus $E$. For degree reasons, the Gromov-Witten invariants relevant to the small quantum product are
\begin{equation}
    GW_{0,3,(1,2)}(h_1h_2,[pt],[pt]),\quad GW_{0,3,(1,2)}(h_2^2,[pt],[pt]).
\end{equation}
We claim that
\begin{align}
    GW_{0,3,(1,2)}(h_1h_2,[pt],[pt])&=4, \label{eqn:GW-invt-MM225} \\
    GW_{0,3,(1,2)}(h_2^2,[pt],[pt])&=8. 
\end{align}
We first pick two generic points $p,q\in\mathbb{P}^3$, given as the Poincaré duals of the point classes. Then there exists a unique line $l$ passing through $p$ and $q$. The planes $P$ passing through $l$ form a $\mathbb{P}^1$-family, and each plane $P$ meets the blowup locus $E$ at 4 points $r_1,r_2,r_3,r_4$. Choosing any three of these points $r_i,r_j,r_k$ uniquely determines a 3-secant conic in the plane $P$ passing through the five points $p,q,r_i,r_j,r_k$.

The class that the one-dimensional family of conics sweeps out can be determined as
\begin{equation}
    S=8H-3D
\end{equation}
for the hyperplane class $H$ and the exceptional divisor $D$. Indeed, a generic point in $E$ lies on three conics in the one-dimensional family, which determines the coefficient of the $D$ term. To determine the coefficient of $H$, take a plane $P'$ in $\mathbb{P}^3$ that contains the line $l$. Then the intersection of $P'$ with $S$ is the union of the four conics determined above, which determines the $8H$ term.

Therefore it follows that
\begin{align}
    GW_{0,3,(1,2)}(f,[pt],[pt])&=3,\quad GW_{0,3,(1,2)}(h_2^2,[pt],[pt])=8,
\end{align}
and we use the relation $h_1h_2=2h_2^2-4f$ to conclude (\ref{eqn:GW-invt-MM225}).

From the Gromov-Witten invariants, the small quantum products are determined as
\begin{align}
    h_1h_2*_{(1,2)}[pt]&=4,\\
    h_2^2*_{(1,2)}[pt]&=8,\\
    [pt]*_{(1,2)}[pt]&=3h_1+2h_2.
\end{align}

\subsection{(2,2)-curves}

These are the strict transforms of $2$-secant conics to $E$. For degree reasons, the only Gromov-Witten invariant relevant to the small quantum product that these classes can contribute to is
\begin{equation}
    GW_{0,3,(2,2)}([pt],[pt],[pt]).
\end{equation}
Since the projection of a representative of the Poincaré dual of the point class is a generic point in $\mathbb{P}^3$, the above invariant counts the number of $2$-secant conics through three general points $a,b,c$ in $\mathbb{P}^3$.

The three general points determine a plane $P\subset\mathbb{P}^3$, and this plane meets the blowup locus $E$ at four points $\{x_1,x_2,x_3,x_4\}$. Now for any choice of two distinct points $x_i$, $x_j$ among these four points, there exists a unique conic in the plane $P$ passing through the five points $a,b,c,x_i,x_j$. Therefore we conclude that
\begin{equation}
    GW_{0,3,(2,2)}([pt],[pt],[pt])=\begin{pmatrix}4\\2
    \end{pmatrix}=6,
\end{equation}
and therefore
\begin{equation}
    [pt]*_{(2,2)}[pt]=6\cdot 1_X.
\end{equation}

\subsection{(0,3)-curves}

These are the rational cubic curves in a single fiber $Q_t\cong\mathbb{P}^1\times\mathbb{P}^1$ of the quadric pencil. Again, for degree reasons, the only Gromov-Witten invariant relevant is
\begin{equation}
    GW_{0,3,(0,3)}([pt],[pt],[pt])=0.
\end{equation}
Indeed, the three point classes induce incidence relations with three general points which lie in different $Q_\lambda$s, but any $(0,3)$-curve class lies entirely in a single quadric $Q_\lambda$, and thus cannot satisfy the incidence conditions.

The remaining task is to exclude the contribution from nonsimple stable maps. For a multiple cover that could contribute, the underlying simple curve class $\beta$ has anticanonical degree at most $3$. Such curves are exceptional fibers or strict transforms of bisecant and once-secant lines. The corresponding families of curve images have dimensions $1,2,3$, equal to the anticanonical degrees. Therefore the images of their three-point evaluation maps can have dimension at most $c_1(X)\cdot\beta+3$, and thus cannot meet the general incidence conditions defining a zero-dimensional moduli space for a multiple cover. A similar dimension argument also excludes the contribution of stable maps with multiple irreducible components in their images. Therefore we have determined that the only contributions to the small quantum product come from embedded curves.

Before we consider the odd classes, we briefly justify that the moduli spaces we considered for the embedded curves are regular: it suffices to show that each embedded image $C$ satisfies
\begin{equation}
    H^1(C,N_{C/X})=0.
\end{equation}
For exceptional fibers and strict transforms of lines this is standard. For the strict transforms of smooth plane conics meeting a generic elliptic quartic $E$ at $r$ points with the plane transverse to $E$, the normal bundle fits into the short exact sequence
\begin{equation}
    0\longrightarrow \mathcal{O}_{\mathbb{P}^1}(4-r)\longrightarrow N_{C/X} \longrightarrow \mathcal{O}_{\mathbb{P}^1}(2)\longrightarrow 0,
\end{equation}
and thus for the cases $r=2,3,4$ contributing to the small quantum product, the desired vanishing of cohomology holds.

We now move on to the small quantum products between the odd classes and the even classes. Since all the odd cohomology classes are supported in degree $3$, only the curves of anticanonical degree $\leq3$ can contribute to these small quantum products.

Take the odd cohomology classes $\psi_1$, $\psi_2$ where $\int_X\psi_1\cup\psi_2=-1$. We first compute
\begin{equation}
    GW_{0,2,(1,0)}(\psi_1,\psi_1)=0,\quad  GW_{0,2,(1,0)}(\psi_1,\psi_2)=1.
\end{equation}
To do so, define the maps $j:D\to X$, $p:D\to E$ for the exceptional divisor $D$. Then the classes $\psi_1=j_*p^*\alpha$, $\psi_2=j_*p^*\beta$ are obtained by the pull-push construction from the $H^1$ classes $\alpha,\beta$ of the blowup locus, which is an elliptic curve. Since the $(1,0)$ class is an exceptional fiber, there exists exactly one exceptional fiber supported over the intersection point of the $\alpha$ and $\beta$ curves.

Thus by the divisor axiom, we conclude
\begin{equation}
     GW_{0,3,(1,0)}(h_1,\psi_1,\psi_2)=1,\quad  GW_{0,3,(1,0)}(h_2,\psi_1,\psi_2)=0.
\end{equation}

Moving on to the anticanonical degree $2$ classes, we claim
\begin{equation}
    GW_{0,3,(0,1)}(h_1h_2,\psi_1,\psi_2)=0,\quad GW_{0,3,(0,1)}(h_2^2,\psi_1,\psi_2)=2.
\end{equation}
Recall that the $(0,1)$ classes are ruling lines entirely contained in a fiber $Q_t$ of the quadric pencil. The incidence condition with the Poincaré dual of the $h_1h_2$ class fixes this fiber, and imposes an additional incidence condition with a generic hyperplane $H\subset\mathbb{P}^3$ which is satisfied at exactly one point for ruling lines. Therefore the moduli space of such ruling lines has two connected components (one for each family of rulings), each isomorphic to $\mathbb{P}^1$. Thus, the incidence condition with $\psi_1$ after the pull-push becomes a degree $1$ cohomology class of $\mathbb{P}^1$, which is necessarily trivial. Therefore the $h_1h_2$-Gromov-Witten invariant is zero.

For the fiber class $f$, we will show that
\begin{equation}
    GW_{0,3,(0,1)}(f,\psi_1,\psi_2)=1.
\end{equation}
The Poincaré dual of the fiber class is supported over a single point $e\in E$, and the $(0,1)$-curve (which is a bisecant line to $E$ in $\mathbb{P}^3$) passes through $e$. Now the moduli space of such bisecant lines is parametrized by the second secant point $e'\in E$, and this, together with the $\psi_1,\psi_2$ incidence conditions, determines the value to be $1$. This determines the Gromov-Witten invariant involving $h_2^2$ from the relation $h_1h_2=2h_2^2-4f$.

Finally, we claim that
\begin{equation}
     GW_{0,3,(1,1)}([pt],\psi_1,\psi_2)=1.
\end{equation}
We recall that the Poincaré dual of the point class projects to a generic point $p\in\mathbb{P}^3$, and the incidence condition of $\psi_1$ and $\psi_2$ determines a unique point $q$ on the blowup locus $E$. Therefore the unique 1-secant line satisfying these incidence conditions is $\overline{pq}$.

This concludes our computation of the small quantum products involving the odd classes.   \hfill $\Box$

\section{Noncommutative formulae} \label{app:HH-calculus}

In this appendix, we verify some of the formulae that appear in the proofs of Sections 2 and 3.

The first identity to check is (\ref{eqn:App-1}): we first recall the relevant definition.

\begin{defn}
    For $\upsilon=z\otimes x_k\otimes\cdots\otimes x_1$, we define the operation
    \begin{align}
        \iota_\phi^\beta(\upsilon)=&\sum(-1)^{\vartheta_\phi}\beta^{i+k-t-s+2}(x_i,\cdots,\phi^s(x_{j+s},\cdots,x_{j+1}),\cdots,z,\cdots,x_{t+1})\otimes\cdots\otimes x_{i+1}
    \end{align}
    for the sign
    \begin{align}\label{eqn:Appendix-sign}
        \vartheta_\phi=\norm{\beta}+\norm{\phi}\dagger_j+(\dagger_i+\norm{\phi})(\norm{z}+\sum_{r=i+1}^k\norm{x_r})+\norm{\beta}\sum_{r=i+1}^t\norm{x_r}.
    \end{align}
\end{defn}

\begin{lemma} \label{lem:App-1}
    For the operation $\iota$ defined above, the following equation holds:
    \begin{equation}
        b\circ\iota^\eta_\phi-(-1)^{\norm{\phi}+\norm{\eta}}\iota^\eta_\phi\circ b=\iota_\phi^{\delta\eta}-(-1)^{\norm{\eta}}\iota^\eta_{\delta\phi}-\mathcal{L}_{\eta}(\phi\cap~)+(-1)^{\norm{\eta}\cdot\lvert\phi\rvert}\phi\cap\mathcal{L}_\eta+(-1)^{\norm{\eta}}[\eta,\phi]\cap~.
    \end{equation}
\end{lemma}
\begin{proof}
    This is the standard Cartan homotopy computation carried out in \cite[Lemma 2.3]{Getzler}. Our operation $\iota^\eta_\phi$ corresponds to Getzler's operation $\rho\{\eta,\phi\}$. Both operations insert $\eta$ across the distinguished Hochschild-chain input, with $\phi$ appearing as an inner input of $\eta$, and with the corresponding Koszul signs from the suspended bar complex. Our convention $\norm{x}=\lvert x\rvert-1$ records the degree of a suspended input, and we order the bar inputs from right to left. The differentials we use are $\delta=-[\mu,~]$, and $b=\mathcal{L}_\mu$.

    The computation with signs for the Cartan homotopy formula is carried out in \cite[Lemma 2.3]{Getzler}. Another reference closer to our sign conventions is \cite[(4.12)--(4.13)]{GetzlerFormulae}, where the homotopy is denoted as $T(A,B)$. The authors explicitly note that their homotopy $T(A,B)$ agrees with $\rho$ from Getzler. We remark that their equation $(4.12)$ with reduced degrees for inputs reproduces our lemma.
\end{proof}

The identities to check now are (\ref{eqn:App-3}), (\ref{eqn:App-4}), and ($\ref{eqn:App-5}$).

Recall the setup: $\mathcal{F},\mathcal{G}:\mathcal{C}\to\mathcal{D}$ are two $A_\infty$-functors between $\mathbb{Z}/2$-graded $A_\infty$-categories, and $T\in CC^*(\mathcal{C},\prescript{~}{\mathcal{F}}{\Delta}_{\mathcal{G}})$ is a pre-natural transformation. Then $\iota_T:CC^*(\mathcal{D},\mathcal{D})\to CC^*(\mathcal{C},\prescript{~}{\mathcal{F}}{\Delta}_{\mathcal{G}})$ is defined as
    \begin{equation}
        (\iota_T\phi)^k(x_k,\cdots,x_1)=\sum(-1)^{\norm{T}\cdot\dagger}\phi^{s+t+1}(\mathcal{G}^{j_s}(x_k,\cdots,x_{k-j_s+1}),\cdots,T^l(x_q,\cdots,x_{i_1+\cdots+i_t+1}),\cdots,\mathcal{F}^{i_1}(x_{i_1},\cdots,x_1))
    \end{equation}
    for the sign
    \begin{equation}
        \dagger=\norm{x_{i_1+\cdots+i_t}}+\cdots+\norm{x_1}.
    \end{equation}

We also recall the sign conventions for the $A_\infty$-functor category:
\begin{align}
    (\mu^1_{fun}T)^d(x_d,\cdots,x_1)&=\sum(-1)^{\norm{T}\cdot\dagger_i}\mu^{r}(\mathcal{G}^{s_r}(x_d,\cdots),\cdots,T^{s_j}(\cdots),\mathcal{F}^{s_{j-1}}(x_i,\cdots),\cdots,\mathcal{F}^{s_1}(\cdots,x_1))\\
    &-\sum(-1)^{\norm{T}+\dagger_i} T^{d-k+1}(x_d,\cdots,\mu^k(x_{i+k},\cdots,x_{i+1}),x_i,\cdots,x_1)
\end{align}
and for $T_1:\mathcal{F}_0\to\mathcal{F}_1$, $T_2:\mathcal{F}_1\to\mathcal{F}_2$, the composition is defined as
\begin{align}
    \mu^2_{fun}(T_2,T_1)^d(x_d,\cdots,x_1)&=\sum(-1)^{\norm{T_2}\cdot\dagger_j+\norm{T_1}\cdot\dagger_i}\mu^r(\mathcal{F}_2(~),\cdots,T_2(~),\mathcal{F}_1(x_j,\cdots),\cdots,T_1(~),\mathcal{F}_0(x_i,\cdots),\cdots,\mathcal{F}_0(~)).
\end{align}

\begin{lemma} \label{lem:App-3}
    In the above setup, for a Hochschild cochain $\phi\in CC^*(\mathcal{D},\mathcal{D})$, the following holds:
    \begin{equation}
        \iota_{\mu^1T}(\phi)-(-1)^{\norm{\phi}}\mu^1_{fun}(\iota_T\phi)=(-1)^{\norm{\phi}}\iota_T(\delta\phi)+(-1)^{\norm{\phi}\cdot\lvert T\rvert}\mu^2_{fun}(\mathcal{G}^*\phi,T)+(-1)^{\norm{\phi}}\mu^2_{fun}(T,\mathcal{F}^*\phi).
    \end{equation}
    In particular, specializing to $\phi=\eta_\mathcal{D}$ recovers (\ref{eqn:App-3}).
\end{lemma}
\begin{proof}
    We suppress variables that don't affect the sign to make the proof clearer. For the first term,
    \begin{align}
        \iota_{\mu^1_{fun}T}(\phi)^k(x_k,\cdots,x_1)&=\sum(-1)^{\lvert T\rvert\cdot\dagger_i}\phi(\mathcal{G}(~),\cdots,\mu^1_{fun}T(~),\mathcal{F}(x_i,\cdots),\cdots,\mathcal{F}(~))\\
        &=\sum(-1)^{A_1}\phi(\mathcal{G}(~),\cdots,\mu(\mathcal{G}(~),\cdots,T(~),\mathcal{F}(x_j,\cdots,x_{i+1})),\cdots,\mathcal{F}())\\
        &+\sum(-1)^{A_2}\phi(\mathcal{G}(~),\cdots,T(\cdots,\mu(\cdots),x_r,\cdots,x_{i+1}),\cdots,\mathcal{F}(~))
    \end{align}
    for the signs
    \begin{align}
        A_1&=\dagger_i+\norm{T}\cdot\dagger_j\\
        A_2&=\lvert T\rvert\cdot\dagger_i+\lvert T\rvert+\dagger_{i+1}^{r}
    \end{align}
    where $\dagger_{j+1}^i=\norm{x_i}+\cdots+\norm{x_{j+1}}$. The second term on the left-hand side is
    \begin{align}
        -(-1)^{\norm{\phi}}\mu^1_{fun}(\iota_T\phi)&=\sum(-1)^{\lvert\phi\rvert+\norm{T}\cdot\dagger_i}\mu^1_{fun}(\phi(\mathcal{G}(~),\cdots,T(~),\mathcal{F}(x_i,\cdots),\cdots,\mathcal{F}(~)))\\
        &=\sum(-1)^{A_3}\mu(\mathcal{G}(~),\cdots,\phi(\mathcal{G}(~),\cdots,T(\cdots,x_{i+1}),\cdots,\mathcal{F}(\cdots,x_{j+1})),\cdots,\mathcal{F}(~))\\
        &+\sum(-1)^{A_4}\phi(\mathcal{G}(~),\cdots,\mathcal{G}(\cdots,\mu(\cdots,x_{r+1}),\cdots),\cdots,T(\cdots,x_{i+1}),\cdots,\mathcal{F}(~))\\
        &+\sum(-1)^{B_2}\phi(\mathcal{G}(~),\cdots,T(\cdots,\mu(~),x_r,\cdots,x_{i+1}),\cdots,\mathcal{F}(~))\\
        &+\sum(-1)^{A_5}\phi(\mathcal{G}(~),\cdots,T(\cdots,x_{i+1}),\cdots,\mathcal{F}(\cdots,\mu(\cdots,x_{j+1}),\cdots),\cdots,\mathcal{F}(~))
    \end{align}
    for the signs
    \begin{align}
        A_3&=\lvert\phi\rvert+\norm{T}\cdot\dagger_i+\norm{\phi}\cdot \dagger_j\\
        A_4&=\norm{T}\cdot\dagger_i+\dagger_r+\norm{T}\\
        B_2&=\norm{T}\cdot \dagger_i+\dagger_r+\norm{T}\\
        A_5&=\norm{T}\cdot \dagger_i+\dagger_j.
    \end{align}
    Here the cancellation of $\lvert\phi\rvert$ in $A_4,B_2,A_5$ occurs because the degree of $\iota_T\phi$ is $\lvert\phi\rvert+\lvert T\rvert-1$, which feeds into the formula for $\mu^1_{fun}$. In particular, $B_2=A_2-1$, so these two terms cancel.
    
    For the first term on the right hand side,
    \begin{align}
        (-1)^{\norm{\phi}}\iota_T(\delta\phi)&=\sum(-1)^{\norm{\phi}+\norm{T}\cdot\dagger_i}\delta\phi\bigl(\mathcal{G}(~),\cdots,T(\cdots,x_{i+1}),\cdots,\mathcal{F}(~)\bigr)\\
        &=\sum(-1)^{B_3}\mu(\mathcal{G}(~),\cdots,\phi(\mathcal{G}(~),\cdots,T(\cdots,x_{i+1}),\cdots,\mathcal{F}(\cdots,x_{j+1})),\cdots,\mathcal{F}(~))\\
        &+\sum(-1)^{C_G}\mu(\mathcal{G}(~),\cdots,\phi(\cdots,\mathcal{G}(\cdots,x_{r+1})),\cdots,T(\cdots,x_{i+1}),\cdots,\mathcal{F}(~)\bigr)\\
        &+\sum(-1)^{C_F}\mu(\mathcal{G}(~),\cdots,T(\cdots,x_{i+1}),\cdots,\phi(\cdots,\mathcal{F}(\cdots,x_{j+1})),\cdots,\mathcal{F}(~))\\
        &+\sum(-1)^{B_1}\phi(\mathcal{G}(~),\cdots,\mu(\mathcal{G}(~),\cdots,T(\cdots,x_{i+1}),\cdots,\mathcal{F}(\cdots,x_{j+1})),\cdots,\mathcal{F}(~))\\
        &+\sum(-1)^{B_4}\phi(\mathcal{G}(~),\cdots,\mu(\cdots,\mathcal{G}(\cdots,x_{r+1})),\cdots,T(\cdots,x_{i+1}),\cdots,\mathcal{F}(~))\\
        &+\sum(-1)^{B_5}\phi(\mathcal{G}(~),\cdots,T(\cdots,x_{i+1}),\cdots,\mu(\cdots,\mathcal{F}(\cdots,x_{j+1})),\cdots,\mathcal{F}(~)).
    \end{align}
    The signs are given by
    \begin{align}
        B_3&=\norm{\phi}+\norm{T}\cdot\dagger_i+1+\norm{\phi}\cdot\dagger_j,\\
        C_G&=\norm{\phi}+\norm{T}\cdot\dagger_i+1+\norm{\phi}\cdot(\dagger_r+\norm{T}),\\
        C_F&=\norm{\phi}+\norm{T}\cdot\dagger_i+1+\norm{\phi}\cdot\dagger_j,\\
        B_1&=\norm{T}\cdot\dagger_j+\dagger_i,\\
        B_4&=\norm{T}\cdot\dagger_i+\dagger_r+\norm{T},\\
        B_5&=\norm{T}\cdot\dagger_i+\dagger_j.
    \end{align}
    In particular, the $B_1,B_3,B_4,B_5$ agree with their counterparts $A_1,A_3,A_4,A_5$, so these terms cancel. Finally, the $C_G$ and $C_F$ terms are exactly the signs that cancel out the two remaining $\mu^2_{fun}$ terms.
\end{proof}

We now recall the operation $\rho_T:CC^*(\mathcal{C},\mathcal{C})\to CC^*(\mathcal{C},\prescript{~}{\mathcal{F}}{\Delta}_{\mathcal{G}})$ defined as
\begin{equation}
        (\rho_T\phi)^k(x_k,\cdots,x_1)=\sum(-1)^{\norm{\phi}\cdot\dagger_i}T^{k-j+1}(x_k,\cdots,\phi^j(x_{i+j},\cdots,x_{i+1}),x_i,\cdots,x_1).
    \end{equation}
\begin{lemma} \label{lem:App-4}
    For a Hochschild cochain $\psi\in CC^*(\mathcal{C},\mathcal{C})$, the following identity holds:
    \begin{equation}
        \rho_{\mu^1_{fun}T}(\psi)-\mu^1_{fun}(\rho_T\psi)=(-1)^{\norm{T}}\rho_T(\delta\psi)+\mu^2_{fun}(\mathcal{G}_*\psi,T)+(-1)^{\norm{T}\cdot\norm{\psi}}\mu^2_{fun}(T,\mathcal{F}_*\psi).
    \end{equation}
\end{lemma}
\begin{proof}
    The proof is identical to that of the previous lemma. Expanding the above formula, one sees the same families of terms as in the proof for $\iota_T$, with the Hochschild cochain now inserted into an input of $T$, instead of an output. The terms involving the source \(A_\infty\)-operations combine to give $(-1)^{\norm{T}}\rho_T(\delta\psi)$, while the two boundary terms give $\mu^2_{fun}(\mathcal{G}_*\psi,T)$ and $(-1)^{\norm{T}\norm{\psi}}\mu^2_{fun}(T,\mathcal{F}_*\psi)$.
\end{proof}

Finally, we verify the equation (\ref{eqn:App-5}).
\begin{lemma} \label{lem:App-5}
    \begin{equation}
        \Lambda(\mu^1_{fun}T)-\mu^1_{fun}(\Lambda T)=\mu^1_{fun}T-\iota_T(\kappa_\mathcal{D})+(-1)^{\norm{T}}\rho_T(\kappa_\mathcal{C})-\mu^2_{fun}(H_\mathcal{G},T)-\mu^2_{fun}(T,H_\mathcal{F}).
    \end{equation}
\end{lemma}
\begin{proof}
    We start with the suppressed notation
    \begin{align}
        \mu^1_{fun}T&=\sum(-1)^{\norm{T}\cdot\dagger_i}\mu(\mathcal{G}(~),\cdots,T(~),\mathcal{F}(x_i,\cdots),\cdots)\\
    &-\sum(-1)^{\norm{T}+\dagger_i} T(\cdots,\mu(~),x_i,\cdots).
    \end{align}
    Then
    \begin{align}
        \Lambda(\mu^1_{fun}T)^d&=\sum(-1)^{\norm{T}\cdot\dagger_i}d\cdot\mu(\mathcal{G}(~),\cdots,T(~),\mathcal{F}(x_i,\cdots),\cdots)\\
    &-\sum(-1)^{\norm{T}+\dagger_i} d\cdot T(\cdots,\mu^k(~),x_i,\cdots),
    \end{align}
    and
    \begin{align}
        \mu^1_{fun}(\Lambda T)^d&=\sum(-1)^{\norm{T}\cdot\dagger_i}r\cdot \mu(\mathcal{G}(~),\cdots,T^r(~),\mathcal{F}(x_i,\cdots),\cdots)\\
    &-\sum(-1)^{\norm{T}+\dagger_i} (d-k+1)T(\cdots,\mu^k(~),x_i,\cdots).
    \end{align}

    On the right hand side, we compute
    \begin{align}
        -\iota_T(\kappa_\mathcal{D})&=-\sum(-1)^{\norm{T}\cdot\dagger_i}(2-r)\mu^r(\mathcal{G}^{j_r}(~),\cdots,T^{j_s}(\cdots,x_{i+1}),\cdots,\mathcal{F}^{j_1}(~)),\\
        (-1)^{\norm{T}}\rho_T(\kappa_\mathcal{C})&=\sum(-1)^{\norm{T}+\dagger_i}(2-k)T(\cdots,\mu^k(\cdots,x_{i+1}),\cdots),
    \end{align}
    and
    \begin{align}
        -\mu^2_{fun}(H_\mathcal{G},T)&=-\sum(-1)^{\norm{T}\cdot\dagger_i}(1-s_c)\mu(\mathcal{G}^{s_b}(~),\cdots,\mathcal{G}^{s_1}(~),T(~),\mathcal{F}^{t_a}(x_i,\dots,),\cdots,\mathcal{F}^{t_1}()),\\
        -\mu^2_{fun}(T,H_\mathcal{F})&=-\sum(-1)^{\norm{T}\cdot\dagger_i}(1-t_c)\mu(\mathcal{G}^{s_b}(~),\cdots,\mathcal{G}^{s_1}(~),T(~),\mathcal{F}^{t_a}(x_i,\dots,),\cdots,\mathcal{F}^{t_1}()).
    \end{align}
    Now we confirm that the sums match up. First, for the particular term
    \begin{equation}
        \mu^{1+a+b}(\mathcal{G}^{s_b}(x_d,\cdots),\cdots,T^r(\cdots),\cdots,\mathcal{F}^{t_1}(\cdots,x_1)),
    \end{equation}
    the left-hand side contribution is
    \begin{equation}
        d-r,
    \end{equation}
    while the right-hand side contribution is
    \begin{equation}
        1-(2-1-a-b)-\sum(1-s_c)-\sum(1-t_c),
    \end{equation}
    which agree since $s_b+\cdots+s_1+r+t_a+\cdots+t_1=d$. The Koszul signs of the two contributions agree.

    Now for the particular source-insertion term
    \begin{equation}
        T^{d-k+1}(x_d,\cdots,\mu^k(x_{i+k},\cdots,x_{i+1}),x_i,\cdots,x_1),
    \end{equation}
    factor out the common signed multiplier $-(-1)^{\|T\|+\dagger_i}$. The coefficient on the left is
    \begin{equation}
        d-(d-k+1)=k-1,
    \end{equation}
    while the coefficient on the right is
    \begin{equation}
        1-(2-k)=k-1.
    \end{equation}
    Thus the two contributions agree.
\end{proof}

\end{appendices}


\begin{thebibliography}{40}

\bibitem[AboSmi16]{AbouzaidSmith}
M.~Abouzaid, \& I.~Smith.
\newblock The symplectic arc algebra is formal.
\newblock {\em Duke Math. J.}, 165(6):985--1060, 2016.

\bibitem[ABT01]{ABT}
E.~Arrondo, M.~Bertolini, \& C.~Turrini.
A focus on focal surfaces.
\newblock {\em Asian J. Math.}, 5(3):535--560, 2001.

\bibitem[BGMP26]{BenedettiGuereManivelPerrin}
V.~Benedetti,  J.~Guéré, L.~Manivel, \& N.~Perrin.
\newblock Quantum cohomology and birational geometry of Verra fourfolds.
\newblock Preprint arXiv:2605.30450.

\bibitem[BMP26]{BenedettiManivelPerrin}
V.~Benedetti, L.~Manivel, \& N.~Perrin.
\newblock Quantum cohomology and irrationality of Gushel-Mukai fourfolds.
\newblock Preprint arXiv:2603.17487.

\bibitem[BMFT26]{Fano4K3}
M.~Bernardara, E.~Fatighenti, L.~Manivel, \& F.~Tanturri.
\newblock Fano fourfolds of K3 type.
\newblock {\em Ann. Inst. Fourier (Grenoble)}, Online first, 121 p., 2026.

\bibitem[BT11]{BravThomas}
C.~Brav, \& H.~Thomas.
\newblock Braid groups and Kleinian singularities.
\newblock {\em Math. Ann.}, 351(4):1005--1017, 2011.

\bibitem[BKT09]{BKT}
A.~S.~Buch, A.~Kresch, \& H.~Tamvakis.
\newblock Quantum Pieri rules for isotropic Grassmannians.
\newblock {\em Invent. Math.}, 178(2):345--405, 2009.

\bibitem[CLT21]{CaldararuLiTu}
A.~Căldăraru, S.~Li, \& J.~Tu.
\newblock Categorical primitive forms and Gromov-Witten invariants of $A_n$ singularities.
\newblock {\em Int. Math. Res. Not.}, (24):18489--18519, 2021.

\bibitem[CP25]{CP}
I.~Cheltsov, \& V.~Przyjalkowski.
\newblock Katzarkov-Kontsevich-Pantev conjecture for Fano threefolds.
\newblock {\em Proc. Steklov Inst. Math.}, 328(1):1--156, 2025.

\bibitem[CCGK16]{Fano3}
T.~Coates, A.~Corti, S.~Galkin, \& A.~Kasprzyk.
\newblock Quantum periods for 3-dimensional Fano manifolds.
\newblock {\em Geom. Topol.}, 20(1):103--256, 2016.

\bibitem[Deb20]{Debarre20}
O.~Debarre.
\newblock Gushel-Mukai varieties.
\newblock Preprint arXiv:2001.03485.

\bibitem[DK18]{DebarreKuznetsov}
O.~Debarre, \& A. Kuznetsov.
\newblock Gushel-Mukai varieties: classification and birationalities.
\newblock {\em Algebr. Geom.}, 5(1):15--76, 2018.

\bibitem[EL17]{EtguLekili}
T.~Etgü, \& Y.~Lekili.
\newblock Koszul duality patterns in Floer theory.
\newblock {\em Geom. Topol.}, 21(6):3313--3389, 2017.

\bibitem[Fay26]{Fay}
A.~Fay.
\newblock Equivariant irrationality of very general symmetric Verra fourfolds.
\newblock Preprint arXiv:2605.30439.

\bibitem[Gan13]{GanatraThesis}
S.~Ganatra.
\newblock Symplectic cohomology and duality for the wrapped Fukaya category.
\newblock Ph.D. thesis, arXiv:1304.7312.

\bibitem[Gan23]{Ganatra23}
S.~Ganatra.
\newblock Cyclic homology, $S^1$-equivariant Floer cohomology and Calabi-Yau structures.
\newblock {\em Geom. Topol.}, 27(9):3461--3584, 2023.

\bibitem[GS25]{GanatraSheridan}
S.~Ganatra, \& N.~Sheridan.
\newblock The cyclic open--closed map and variations of Hodge structures.
\newblock Preprint arXiv:2511.04498.

\bibitem[Get93]{Getzler}
E.~Getzler.
\newblock Cartan homotopy formulas and the Gauss-Manin connection in cyclic homology.
\newblock In {\em Quantum deformations of algebras and their representations}, pages 65--78. Israel Math. Conf. Proc. 7, 1993.

\bibitem[Hug24]{Hugtenburg}
K.~Hugtenburg.
\newblock The cyclic open-closed map, u-connections and R-matrices.
\newblock {\em Selecta Math. (N.S.)}, 30(2), Paper No. 29, 2024.

\bibitem[Iri23]{Iritani}
H.~Iritani.
\newblock Quantum cohomology of blowups.
\newblock Preprint arXiv:2307.13555.

\bibitem[KKPY25]{KKPY}
L.~Katzarkov, M.~Kontsevich, T.~Pantev, \& T.~Y.~Yu.
\newblock Birational Invariants from Hodge Structures and Quantum Multiplication.
\newblock Preprint arXiv:2508.05105.

\bibitem[Kea14]{Keating14}
A.~Keating.
\newblock Dehn twists and free subgroups of symplectic mapping class groups.
\newblock {\em J. Topol.}, 7(2):436--474, 2014.

\bibitem[KRW23]{KeatingRandalWilliams}
A.~Keating, \& O.~Randal-Williams.
\newblock On the order of Dehn twists.
\newblock {\em New York J. Math.}, 29:203--212, 2023.

\bibitem[KS02]{KhovanovSeidel}
M.~Khovanov, \& P.~Seidel.
\newblock Quivers, Floer cohomology, and braid group actions.
\newblock {\em J. Amer. Math. Soc.}, 15(1):203--271, 2002.

\bibitem[Kim26]{Kim}
Y.~Kim.
\newblock Simple homotopy theory for Fukaya categories.
\newblock {\em J. Topol.}, 19(3), Paper No. e70099, 2026.

\bibitem[KP18]{KuznetsovPerry}
A. Kuznetsov, \& A.~Perry.
\newblock Derived categories of Gushel-Mukai varieties.
\newblock {\em Compos. Math.}, 154(7):1362--1406, 2018.

\bibitem[KS18]{KuznetsovShinder}
A. Kuznetsov, \& E.~Shinder.
\newblock Grothendieck ring of varieties, D- and L-equivalence, and families of quadrics.
\newblock {\em Selecta Math. (N.S.)}, 24(4):3475--3500, 2018.

\bibitem[LP21]{LeePirola}
Y.~Lee, \& G.~P.~Pirola.
\newblock Vanishing cohomology on a double cover.
\newblock {\em Bull. Lond. Math. Soc.}, 53(2):370--379, 2021.

\bibitem[LU21]{LekiliUeda}
Y.~Lekili, \& K.~Ueda.
\newblock Homological mirror symmetry for Milnor fibers of simple singularities.
\newblock {\em Algebr. Geom.}, 8(5):562--586, 2021.

\bibitem[LW25]{LiuWang}
J.~Liu, \& Z.~Wang.
\newblock $A_\infty$-deformations of zigzag algebras via Ginzburg dg algebras.
\newblock {\em J. Algebra}, 677:360--371, 2025.

\bibitem[MMPS19]{MMPS}
J.~A.~Cruz Morales, A.~Mellit, N.~Perrin, \& M.~Smirnov.
\newblock On quantum cohomology of Grassmannians of isotropic lines, unfoldings of $A_n$-singularities, and Lefschetz exceptional collections.
\newblock {\em Ann. Inst. Fourier (Grenoble)}, 69(3):955--991, 2019.

\bibitem[NV19]{NordskovaVolkov}
A.~Nordskova, \& Y.~Volkov.
\newblock Faithful actions of braid groups by twists along ADE-configurations of spherical objects.
\newblock Preprint arXiv:1910.02401.

\bibitem[Orl93]{Orlov}
D.~Orlov.
\newblock Projective bundles, monoidal transformations, and derived categories of coherent sheaves.
\newblock {\em  Russian Acad. Sci. Izv. Math.}, 41(1):133--141, 1993.

\bibitem[Par04]{Paris}
L.~Paris.
\newblock Artin groups of spherical type up to isomorphism.
\newblock {\em J. Algebra}, 281(2):666--678, 2004.

\bibitem[PPZ26]{PerryPertusiZhao}
A.~Perry, L.~Pertusi, \& X.~Zhao.
\newblock Moduli spaces of stable objects in Enriques categories.
\newblock {\em Selecta Math. (N.S.)}, 32, Paper No. 47, 2026.

\bibitem[PS23]{PerutzSheridan}
T.~Perutz, \& N.~Sheridan.
\newblock Constructing the relative Fukaya category.
\newblock {\em J. Symplectic Geom.}, 21(5):997--1076, 2023.

\bibitem[RitSmi17]{RitterSmith}
A.~Ritter, \& I.~Smith.
\newblock The monotone wrapped Fukaya category and the open-closed string map.
\newblock {\em Selecta Math.}, 23(1):533--642, 2017.

\bibitem[Sei00]{Seidel00}
P.~Seidel.
\newblock Graded Lagrangian submanifolds.
\newblock {\em Bull. Soc. Math. France}, 128(1):103--149, 2000.

\bibitem[Sei03]{Seidel03}
P.~Seidel.
\newblock A long exact sequence for symplectic Floer homology.
\newblock {\em Topology}, 42(5):1003--1063, 2003.

\bibitem[Sei08a]{Seidel08}
P.~Seidel.
\newblock {\em Fukaya categories and Picard-Lefschetz theory}.
\newblock European Math. Soc., 2008.

\bibitem[Sei08b]{Seidel08b}
P.~Seidel.
\newblock Lectures on four-dimensional Dehn twists.
\newblock In {\em Symplectic 4-manifolds and algebraic surfaces}, pages 231--267. Lecture Notes in Math., 1938, Springer-Verlag, Berlin, 2008.

\bibitem[Sei08c]{Seidel08c}
P.~Seidel.
\newblock $A_\infty$-subalgebras and natural transformations.
\newblock {\em Homology Homotopy Appl.}, 10(2):83--114, 2008.

\bibitem[Sei12]{Seidel12}
P.~Seidel.
\newblock Lagrangian homology spheres in $A_m$ Milnor fibres via $\C^*$-equivariant $A_\infty$-modules.
\newblock {\em Geom. Topol.}, 16(4):2343--2389, 2012.

\bibitem[Sei13]{Seidel13}
P.~Seidel.
\newblock {\em Lectures on categorical dynamics and symplectic topology}.
\newblock Notes, available at \url{https://math.mit.edu/~seidel}.

\bibitem[Sei14]{Seidel14}
P.~Seidel.
\newblock Disjoinable Lagrangian spheres and dilations.
\newblock {\em Invent. Math.}, 197(2):299--359, 2014.

\bibitem[Sei15a]{Seidel15}
P.~Seidel.
\newblock Picard-Lefschetz theory and dilating $\C^*$-actions.
\newblock {\em J. Topol.}, 8(4):1167--1201, 2015.

\bibitem[Sei23]{Seidel23}
P.~Seidel.
\newblock Formal Groups and Quantum Cohomology.
\newblock {\em Geom. Topol.}, 27(8):2937--3060, 2023.

\bibitem[Seidel]{SeidelErratum}
P.~Seidel.
\newblock {\em Erratum to: Fukaya categories and Picard-Lefschetz theory, Lemma 2.1}.
\newblock Available at \url{https://math.mit.edu/~seidel/errata/strictly-unital.pdf}.

\bibitem[SS12]{SeidelSolomon}
P.~Seidel, \& J.~Solomon.
\newblock Symplectic cohomology and q-intersection numbers.
\newblock {\em Geom. Funct. Anal.}, 22:443--477, 2012.

\bibitem[ST01]{SeidelThomas}
P.~Seidel, \& R.~Thomas.
\newblock Braid group actions on derived categories of coherent sheaves.
\newblock {\em Duke Math. J.}, 108(1):37--108, 2001.

\bibitem[SS22]{GetzlerFormulae}
A.~Sharapov, \& E.~Skvortsov.
\newblock Homotopy Cartan calculus and inner deformations of $A_\infty$-algebras.
\newblock Preprint arXiv:2202.10925.

\bibitem[She16]{Sheridan}
N.~Sheridan.
\newblock On the Fukaya category of a Fano hypersurface in projective space.
\newblock {\em Publ. Math. Inst. Hautes Études Sci.}, 124:165--317, 2016.

\bibitem[She20]{SheridanFormulae}
N.~Sheridan.
\newblock Formulae in noncommutative Hodge theory.
\newblock {\em J. Homotopy Relat. Struct.}, 15(1):249--299, 2020.

\bibitem[Smi12]{Smith}
I.~Smith.
\newblock Floer cohomology and pencils of quadrics.
\newblock {\em Invent. Math.}, 189(1):149--250, 2012.

\bibitem[ST24]{SolomonTukachinsky}
J.~Solomon, \& S.~Tukachinsky.
\newblock Relative quantum cohomology.
\newblock {\em J. Eur. Math. Soc.}, 26(9):3497--3573, 2024.

\bibitem[VWX26]{VenWoodXu}
S.~Venugopalan, C.~Woodward, \& G.~Xu.
\newblock Fukaya categories of blowups.
\newblock {\em J. Inst. Math. Jussieu}, 25(1):85--214, 2026.

\bibitem[Voi03]{Voisin}
C.~Voisin.
\newblock Hodge theory and complex algebraic geometry. II.
\newblock {\em Cambridge Stud. Adv. Math.}, Vol. 77. Cambridge University Press, Cambridge, 2003.

\end{thebibliography}
\end{document}